\documentclass[11pt]{article}
\usepackage{amsthm,amsfonts,amsmath,amssymb}
\usepackage{multirow}
\newtheorem{thm}{Theorem}[section]
\newtheorem{prop}[thm]{Proposition}
\newtheorem{lem}[thm]{Lemma}
\newtheorem{cor}[thm]{Corollary}

\newtheorem{rem}[thm]{Remark}

\def\di{\bigm|} \def\lg{\langle} \def\rg{\rangle}
\def\nd{\mathrel{\bigm|\kern-.7em/}}

\def\f{\noindent}
\def\maxsgp{<\!\!\!\cdot\ }

\def\mod{\hbox{\rm mod }}

\def\demo{\f{\bf Proof}\hskip10pt}
\def\cha{\hbox{\rm\ char }}

\def\qed{\hfill $\Box$}

\def\Z{{\rm Z}}

\def\lg{\langle}
\def\rg{\rangle}

\def\rr#1{\item[{\rm (#1)}]}

\def\A{\mathcal{A}$$}

\begin{document}
\title{Finite $p$-groups all
of whose subgroups of index $p^3$ are abelian
\thanks{This work was supported by NSFC (no.11371232) and
by NSF of Shanxi Province (no.2012011001). }}
\author{Qinhai Zhang\thanks{Corresponding author. e-mail: zhangqh@dns.sxnu.edu.cn\ \
zhangqhcn@hotmail.com},\ \ Libo Zhao,\ Miaomiao Li and Yiqun Shen\\
Department of Mathematics, Shanxi Normal University, Linfen, Shanxi 041004\\
People's Republic of China}

\maketitle

\begin{abstract}
Suppose that $G$ is a finite $p$-group. If all subgroups of index
$p^t$ of $G$ are abelian and at least one subgroup of index
$p^{t-1}$ of $G$ is not abelian, then $G$ is called an
$\mathcal{A}_t$-group. In this paper, some information about
$\mathcal{A}_t$-groups are obtained and $\mathcal{A}_3$-groups are
completely classified. This solves an {\it old problem} proposed by
Berkovich and Janko in their book \cite{Ber2}. Abundant information about $\mathcal{A}_3$-groups are given.


\medskip
\noindent{\bf 2010MSC} 20D15
\end{abstract}

\baselineskip=16pt

\section{Introduction}

\hskip0.2in Finite $p$-groups are an important class of finite
groups. After the classification of finite simple groups was finally
completed, the study of finite $p$-groups becomes more and more
active. Many leading group theorists, for example, Glauberman and
Janko have turned their attentions to the study of finite
$p$-groups. Although finite simple groups are classified, it is
impossible to classify finite $p$-groups in the classical sense. The
reason is that a finite $p$-group has `` too many'' normal subgroups
and consequently there is an extremely large number of
non-isomorphic $p$-groups of a given fixed order. In fact,
Higman\cite{Higman1, Higman2}  gave a formula for the number
$f(n,p)$ of non-isomorphic $p$-groups of order $p^n$:
\begin{center}when $n\to\infty$,\ \ $f(n,p)=p^{n^3({2/27+O(n^{-1/3})})}.$
\end{center}
It is easy to see that when $n$ becomes large, the number of
non-isomorphic $p$-groups of order $p^n$ becomes large in exponent
speed. For example, up to now, the known  results about the
classification of $2$-groups \cite{BEO} are:

 {\normalsize\begin{center}
\begin{tabular}{|c||c|c|c|c|c|c|c|c|c|c|} \hline
$n$ &$2$ &$2^2$ &$2^3$ &$2^4$ &$2^5$ &$2^6$ &$2^7$  &$2^8$ &$2^9$
&$2^{10}$\\ \hline $f(n,2)$ &$1$ &$2$ &$5$ &$14$ &$51$ &$267$
&$2328$ &$56092$ &$10494213$ &$49487365422$\\ \hline
\end{tabular}
\end{center}}

For $p>2$, $p$-groups of order $p^7$ are classified by \cite{OL}.
The result is

$f(7,3)=9310$, $f(7,5)=34297$. For $p >5$,

$f(7,p)=3p^5 +12p^4 +44p^3 +170p^2 +707p + 2455 + (4p^2 +44p
+291)\gcd(p-1, 3) +(p^2 + 19p + 135)\gcd(p -1, 4)+ (3p +31)
\gcd(p-1, 5)+4 \gcd(p-1, 7) +5 \gcd(p- 1, 8)+ \gcd(p-1, 9)$.

\medskip
Because of the difficult of the classification of finite $p$-groups
in the classical sense,  Janko and Berkovich sponsored and led an
research project that aims at classifying certain classes of finite
$p$-groups defined by their subgroup structure. As Janko mentioned
in the Foreword of \rm \cite{Ber1}, to study $p$-groups with
``large'' abelian subgroups is another approach to finite
$p$-groups. We know that non-abelian $p$-groups with ``largest''
abelian subgroups are minimal non-abelian groups. A non-abelian
group $G$ is said to be {\it minimal non-abelian} if every proper
subgroup of $G$ is abelian. Minimal non-abelian groups were
classified in \cite{MIller}, and in more detail for finite
$p$-groups in \cite{R}. Berkovich and Janko in  \cite{BJ} introduced
a new concept, $\mathcal{A}_t$-groups, which is a more general
concept than that of minimal non-abelian $p$-groups. For a positive
integer $t$, a finite $p$-group $G$ is called an {\it
$\mathcal{A}_t$-group} if all subgroups of index $p^t$ of $G$ are
abelian, and at least one subgroup of index $p^{t-1}$ of $G$ is not
abelian. In this paper, an $\mathcal{A}_0$-group is an abelian
group. Obviously, $\mathcal{A}_1$-groups are exactly the minimal
non-abelian $p$-groups. For small $t$, $\mathcal{A}_t$-groups can be
considered as
groups having ``large'' abelian subgroups. Many scholars studied and
classified $\mathcal{A}_2$-groups, see
\cite{BJ,Ber2,Dra,Kaz,She,ZSAX}. Obviously, classifying $\mathcal
{A}_3$-groups is a fascinating problem. In fact, this is a problem
proposed by Berkovich and Janko in their joint book \cite{Ber2}.

\medskip
{\bf Problem 1278.} {\it  {\rm (}Old problem{\rm )} Classify
${\mathcal{A}}_3$-groups.}
\medskip

In this paper, we completely classify $\mathcal{A}_3$-groups in
classical sense. The groups described in the title are the totality
of abelian groups, ${\mathcal{A}}_1$-groups,
${\mathcal{A}}_2$-groups and ${\mathcal{A}}_3$-groups. Thus the
groups described in the title are completely classified.


\medskip


Related to ${\mathcal{A}}_3$-groups, Berkovich and Janko proposed
the following

\medskip
{\bf Problem 893(\cite{Ber2}).} {\it Describe the set
$\{\alpha_1(G)\di G \  is \ an\  {\mathcal{A}}_3$-group \}.}

\medskip
{\bf Problem 1595(\cite{Ber3}).} {\it Classify the
${\mathcal{A}}_3$-groups $G$ such that $\alpha_1(G)<p^2+p+1$, where
$\alpha_1(G)$ denotes the number of ${\mathcal{A}}_1$-subgroups in a
$p$-group $G$.}

\medskip
{\bf Problem 2829(\cite{Ber4}).} {\it Find ${\rm
max}\{\alpha_1(G)\}$, where a $p$-group $G$ runs over all
${\mathcal{A}}_3$-groups.}

As a corollary of the classification of ${\mathcal{A}}_3$-groups, we
get all $\alpha_1(G)$ for ${\mathcal{A}}_3$-groups $G$. Hence these
problems mentioned above are also solved. Moreover, we give the
triple $(\mu_0,\mu_1,\mu_2)$, where $\mu_i$ denotes the number of
$\mathcal{A}_i$-subgroups of index $p$ in $\mathcal{A}_3$-groups.

\smallskip

We classify $\mathcal{A}_3$-groups $G$ in two parts: $G$ has an
$\mathcal{A}_1$-subgroup of index $p$, and $G$ has no
$\mathcal{A}_1$-subgroup of index $p$. The sketch of the
classification of ${\mathcal{A}}_3$--groups are as follows.

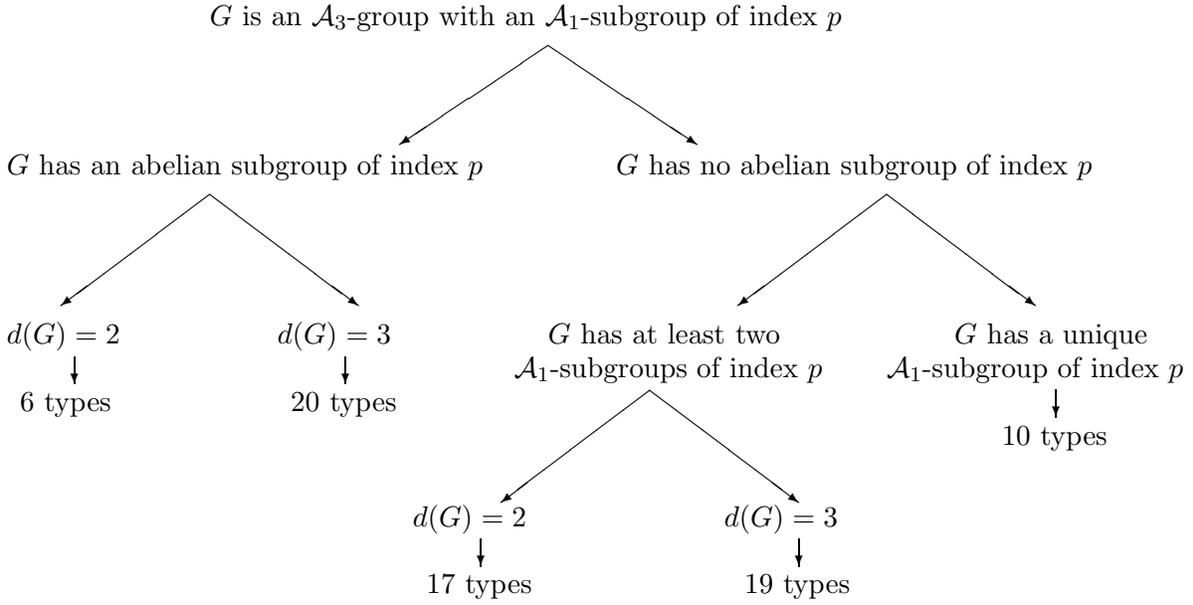
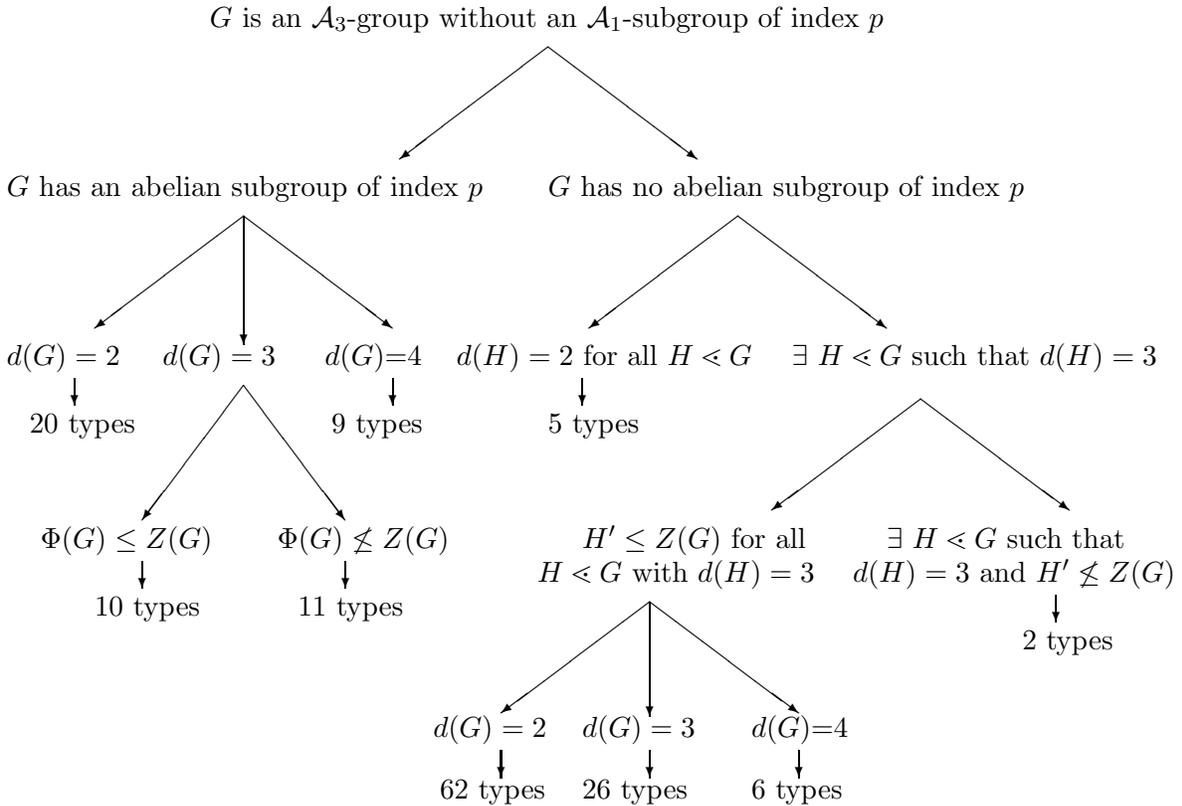
\begin{center}
\vspace{1cm} \setlength{\unitlength}{.90mm}
\begin{picture}(140, 260)(-10, 0)

\put(10, 247) {$G$ is an $\mathcal{A}_3$-group with an
$\mathcal{A}_1$-subgroup of index $p$}
    \put(60, 244){\vector(-3, -2){22}} \put(60, 244){\vector(3, -2){22}}
    \put(-20, 225){$G$ has an  abelian subgroup of index $p$} \put(70, 225){$G$ has no abelian subgroup of index $p$}
        \put(10, 222){\vector(-4, -3){22}} \put(10, 222){\vector(4, -3){22}}\put(110, 222){\vector(-4, -3){22}} \put(110, 222){\vector(4, -3){22}}
        \put(-20, 200){$d(G)=2$} \put(20, 200){$d(G)=3$}\put(60, 200){$G$ has at least two } \put(120, 200){$G$ has a unique }
          \put(55, 195){$\mathcal{A}_1$-subgroups of index $p$} \put(110, 195){$\mathcal{A}_1$-subgroup of index $p$}

       \put(75, 193){\vector(4, -3){22}}\put(75, 193){\vector(-4, -3){22}}
           \put(40, 173){$d(G)=2$} \put(86, 173){$d(G)=3$}

           \put(-10, 198){\vector(0, -3){4}}\put(30, 198){\vector(0, -3){4}}
            \put(-18, 190){6 types }\put(22, 190){20 types }

            \put(135, 193){\vector(0, -3){4}}
            \put(127, 185){10 types }

            \put(50, 171){\vector(0, -3){4}}\put(97, 171){\vector(0, -3){4}}
            \put(42, 163){17 types }\put(89, 163){19 types }

    \put(10, 135){$G$ is an $\mathcal{A}_3$-group without an $\mathcal{A}_1$-subgroup of index $p$}
        \put(60, 132){\vector(-4, -3){22}}
        \put(-20, 110){$G$ has an abelian subgroup of index $p$}
        \put(60, 132){\vector(4, -3){22}}
        \put(60, 110){$G$ has no abelian subgroup of index $p$}

            \put(15, 107){\vector(4, -3){22}}\put(15, 107){\vector(-4, -3){22}}\put(15, 107){\vector(0, -2){19}}
            \put(-20, 85){$d(G)=2$}\put(3, 85){$d(G)=3$}\put(27, 85){$d(G)$=4}

            \put(88, 107){\vector(4, -3){22}}\put(88, 107){\vector(-4, -3){22}}
           \put(45, 85){ $d(H)=2$ for all $H\lessdot G$} \put(96, 85){$\exists\ H\lessdot G$ such that $d(H)=3$}

            \put(-10, 83){\vector(0, -3){4}}\put(37, 83){\vector(0, -3){4}}
            \put(-18, 75){ 20 types }\put(28, 75){9 types}

             \put(65, 83){\vector(0, -3){4}}
            \put(60, 75){5 types }

            \put(15, 82){\vector(3, -4){15}}\put(15, 82){\vector(-3, -4){15}}
            \put(-15,58){$\Phi(G)\le Z(G)$}         \put(20, 58){$\Phi(G)\not\le Z(G)$}

             \put(0, 56){\vector(0, -3){4}}\put(30, 56){\vector(0, -3){4}}
                   \put(-7, 48){10 types }\put(23, 48){11 types}

               \put(115, 80){\vector(4, -3){22}}\put(115, 80){\vector(-4, -3){22}}
               \put(65, 58){$H'\le Z(G)$ for all} \put(110, 58){$\exists\ H\lessdot G$ such that}
               \put(57,53){ $H\lessdot G$ with $d(H)=3$ }         \put(105, 53){$d(H)=3$ and $H'\not\le Z(G)$}
                                                        \put(135, 51){\vector(0, -3){4}}
                                                      \put(130, 43){2 types}

                   \put(75, 50){\vector(4, -3){22}}\put(75, 50){\vector(-4, -3){22}}\put(75, 50){\vector(0, -2){17}}
                   \put(43, 30){$d(G)=2$}\put(65, 30){$d(G)=3$}\put(90, 30){$d(G)$=4}
                       \put(53, 28){\vector(0, -3){4}} \put(75, 28){\vector(0, -3){4}}\put(97, 28){\vector(0, -3){4}}
                       \put(44, 21){62 types}\put(65, 21){26 types}\put(90, 21){6 types}

\end{picture}
\end{center}

\section{Preliminaries}

In this paper, $p$ is always a prime. We use $F_p$ to denote the
finite field containing $p$ elements. $F_p^*$ is the multiplicative
group of $F_p$. $(F_p^*)^2=\{ a^2\di a\in F_p^*\}$ is a subgroup of
$F_p^*$.

Let $G$ be a finite group. We use $c(G)$, $\exp(G)$ and $d(G)$ to
denote the nilpotency class, the exponent and the minimal number of
generators of $G$ respectively.  We use $C_{p^m}$, $C_{p^m}^{n}$ and
$H*K$ to denote the cyclic group of order $p^m$, the direct product
of $n$ cyclic groups of order $p^m$, and a central product of $H$
and $K$ respectively. We use $M\lessdot G$ to denote $M$ is a
maximal subgroup of $G$ and
$$G>G'=G_2>G_3>\dots>G_{c+1}=1$$
denote the lower central series of $G$, where $c=c(G)$.

Let $G$ be a finite $p$-group. We use $G\in \mathcal{A}_t$ to denote
$G$ is an $\mathcal{A}_t$-group. For any positive integer $s$, we
define
$$\Omega_s(G)=\lg a\in G\di a^{p^s}=1\rg\ {\rm and}\ \mho_s(G)=\lg
a^{p^s}\di a\in G\rg.$$

We use $M_p(n,m)$ to denote the $p$-groups
$$\langle a, b \di a^{p^n}=b^{p^m}=1, a^b=a^{1+p^{n-1}}\rangle,
{\rm where}\ n\geq 2.$$

We use $M_p(n,m,1)$ to denote the $p$-groups
$$\langle a, b; c \di a^{p^n}=b^{p^m}=c^p=1, [a, b]=c, [c, a]=[c,
b]=1 \rangle,$$ {\rm where}$$\ n\geq m,\ {\rm and\ if}\ p=2,\ {\rm
then}\ m+n\geq 3.$$

For other notation and terminology the reader is referred to
\cite{Hup}.


\medskip
The following results are used in this paper, we gather them
together.

\begin{lem}\label{phi(G')G_3}{\rm (\cite[Theorem 2.3]{Bla})}
\label{metacyclic} Suppose that $G$ is a finite $p$-group. Then $G$
is metacyclic if and only if $G/\Phi(G')G_3$ is metacyclic.
\end{lem}

\begin{lem}{\rm (\cite[Lemma 2.2]{Alj})}
\label{minimal non-abelian equivalent conditions} Suppose that $G$
is a finite non-abelian $p$-group. Then the following conditions are
equivalent:

{\rm (1)} $G$ is minimal non-abelian;

{\rm (2)} $d(G)=2$ and $|G'|=p$;

{\rm (3)} $d(G)=2$ and $\Phi(G)=Z(G)$.
\end{lem}

\begin{lem}\label{thm=Redei}{\rm (\cite {R})}
Suppose that $G$ is an $\mathcal{A}_1$-group. Then $G$ is one of the
following groups: $Q_8$, $M_p(n,m)$ or $M_p(n,m,1)$.
\end{lem}

\begin{lem}{\rm (\cite[\S 9, Exercise 10]{Ber1})}\label{jidalei3group}
Let $G$ be a $3$-group of maximal class. Then the fundamental
subgroup $G_1$ of $G$ is either abelian or minimal nonabelian.
\end{lem}

\begin{lem}\label{A_2}{\rm(\cite{ZSAX})} A finite $p$-group $G$ is an $\mathcal{A}_2$-group if and only if  $G$ is
one of the following pairwise non-isomorphic groups:
\begin{enumerate}
\rr{I} $d(G)=2$ and $G$ has an abelian subgroup of index $p$. In this case, $\alpha_1(G)=p$.
\begin{enumerate}
\rr{1} $\lg a,b\di a^{8}=b^{2^m}=1, [a,b]=a^{-2}\rg$, where $m\ge
1$;

\rr{2} $\lg a,b\di a^{8}=b^{2^m}=1, [a,b]=a^2\rg$, where $m\ge 1$;

\rr{3} $\lg a,b\di a^{8}=1, b^{2^m}=a^{4}, [a,b]=a^{-2}\rg$, where
$m\ge 1$;

\rr{4} $\lg a_1,b;a_2,a_3\di
a_1^p=a_2^p=a_3^p=b^{p^m}=1,[a_1,b]=a_{2},[a_2,b]=a_3,[a_3,b]=1,[a_i,a_j]=1\rg,$
where $p\ge 5$ for $m=1$, $p\ge 3$ and $1\leq i,j\leq 3$;

\rr{5} $\lg a_1,b;a_2\di
a_1^p=a_2^p=b^{p^{m+1}}=1,[a_1,b]=a_{2},[a_2,b]=b^{p^m},[a_1,a_2]=1\rg$,
where $p\ge 3$;

\rr{6} $\lg a_1,b;a_2\di
a_1^{p^2}=a_2^p=b^{p^m}=1,[a_1,b]=a_{2},[a_2,b]=a_1^{\nu
p},[a_1,a_2]=1\rg,$ where $p\ge 3$ and $\nu=1$ or a fixed quadratic
non-residue modulo $p$.

\rr{7} $\lg a_1, b; a_2\di
a_1^9=a_2^3=1,b^3=a_1^3,[a_1,b]=a_2,[a_2,b]=a_1^{-3},[a_2,a_1]=1\rg$.

\end{enumerate}

\rr{II} $d(G)=3$, $|G'|=p$ and $G$ has an abelian subgroup  of index
$p$. In this case, $\alpha_1(G)=p^2$.

\begin{enumerate}
\rr{8} $\langle {a},{b},{x} \di {a}^4={x}^{2}=1,
{b}^2={a}^2=[{a},{b}], [{x},{a}]=[{x},{b}]=1 \rangle \cong Q_8
\times C_{2}$;

\rr{9}  $\langle {a}, {b},{x} \di {a}^{p^{n+1}}={b}^{p^m}={x}^{p}=1,
[{a},{b}]={a}^{p^{n}}, [{x},{a}]=[{x},{b}]=1\rangle\cong
M_p(n+1,m)\times C_{p}$;

\rr{10}  $\langle {a}, {b}, {x}; {c} \di
{a}^{p^n}={b}^{p^m}={c}^{p}={x}^p=1, [{a},{b}]={c},
[{c},{a}]=[{c},{b}]=[{x},{a}]=[{x},{b}]=1\rangle\cong
M_p(n,m,1)\times C_{p}$, where $n\ge m$, and $n\ge 2$ if $p=2$;

\rr{11}  $\langle {a}, {b}, {x} \di {a}^4=1,
{b}^2={x}^{2}={a}^2=[{a},{b}], [{x},{a}]=[{x},{b}]=1\rangle \cong
Q_8 \ast C_{4}$;

\rr{12}  $\langle {a}, {b}, {x} \di
{a}^{p^n}={b}^{p^m}={x}^{p^{2}}=1, [{a},{b}]={x}^{p},
[{x},{a}]=[{x},{b}]=1 \rangle\cong M_p(n,m,1)\ast C_{p^{2}}$, where
$n\ge 2$ if $p=2$ and $n\ge m$.

\end{enumerate}

\rr{III} $d(G)=3$, $|G'|=p^2$ and $G$ has an abelian subgroup  of
index $p$. In this case, $\alpha_1(G)=p^2+p$.

\begin{enumerate}

\rr{13} $\lg {a},{b},{c}\di {a}^{4}={b}^{4}=1,
{c}^2={a}^2{b}^2,[{a},{b}]={b}^2,[{c},{a}]={a}^2,[{c},{b}]=1\rg$;

\rr{14} $\lg {a},{b},{d}\di {a}^{p^{m}}={b}^{p^2}={d}^{p}=1,
[{a},{b}]={a}^{p^{m-1}},[{d},{a}]={b}^p,[{d},{b}]=1\rg$, where
$m\geq 3$ if $p=2$;

\rr{15} $\lg {a},{b},{d}\di {a}^{p^m}={b}^{p^2}={d}^{p^2}=1,
[{a},{b}]={d}^p,[{d},{a}]={b}^{jp},[{d},{b}]=1\rg$, where $(j,p)=1$,
$p>2$, $j$ is a fixed quadratic non-residue modulo $p$, and $-4j$ is
a quadratic non-residue modulo $p$;

\rr{16} $\lg {a},{b},{d}\di
 {a}^{p^m}={b}^{p^2}={d}^{p^2}=1,[{a},{b}]={d}^p,[{d},{a}]={b}^{jp}{d}^p,[{d},{b}]=1\rg$,
where if $p$ is odd, then $4j =1-\rho^{2r+1}$ with $1\le
r\le\frac{p-1}{2}$ and $\rho$ the smallest positive integer which is
a primitive root $(\mod p)$; if $p = 2$, then $j = 1$.

\end{enumerate}

\rr{IV} $d(G)=2$ and $G$ has no abelian subgroup  of index $p$. In this case, $\alpha_1(G)=1+p$.

\begin{enumerate}

\rr{17} $\lg a, b\di a^{p^{r+2}}=1,b^{p^{r+s+t}}=a^{p^{r+s}},[a,
b]=a^{p^r}\rg$, where $r\ge 2$ for $p=2$, $r\ge 1$ for $p\ge 3$,
$t\ge 0$, $0\le s\le 2$ and $r+s\ge 2$;

\rr{18} $\lg a,b;c\di a^{p^2}=b^{p^2}=c^p=1,[a,b]=c,[c,a]=b^{\nu p},
[c,b]=a^{p}\rg$, where $p\ge 5$, $\nu$ is a fixed quadratic
non-residue modulo $p$;

\rr{19} $\lg a,b;c\di a^{p^2}=b^{p^2}=c^p=1,[a,
b]=c,[c,a]=a^{-p}b^{-lp},[c, b]=a^{-p}\rg$, where $p\ge 5$,
$4l=\rho^{2r+1}-1$, $r=1, 2, \dots, \frac{1}{2}(p-1)$, $\rho$ is the
smallest positive integer which is a primitive root modulo $p$;

\rr{20} $\lg a, b;c\di a^9=b^9=c^3=1,[a,b]=c,[c,a]=b^{-3}, [c,
b]=a^3\rg$;

\rr{21} $\lg a, b;c\di a^9=b^9=c^3=1, [a,b]=c, [c,a]=b^{-3}, [c,
b]=a^{-3}\rg$.

\end{enumerate}

\rr{V} $d(G)=3$ and $G$ has no abelian subgroup  of index $p$. In this case, $\alpha_1(G)=1+p+p^2$.

\begin{enumerate}

\rr{22} $\lg a,b,d\di a^4=b^4=d^4=1, [a,b]=d^2, [d,a]=b^2d^2,
[d,b]=a^2b^2,[a^2,b]=[b^2,a]=1\rg$.

\end{enumerate}
\end{enumerate}

\end{lem}

By analyzing the groups in Lemma {\rm\ref{A_2}}, we have following
lemma.

\begin{lem}\label{A_2-property} Suppose that $G$ is an $\mathcal{A}_2$-group of order
$p^n$. Then

\begin{enumerate}
\rr{1} $d(G)\le 3$ and $c(G)\le 3$;

\rr{2} if $d(G)=3$, then $c(G)=2$ and $G'\le C_p^3$. Moreover, if
$G'=C_p^3$, then $|G|=2^6$, $G'=\Omega_1(G)=Z(G)\cong C_2^3$ and $G$
has no abelian subgroup  of index $p$;

\rr{3} if $d(G)=3$ and $|G'|=p$, then $|G:Z(G)|=p^2$ and the type of
$G/G'$ is $(p^u,p^v,p)$ where $u+v=n-2$;

\rr{4} if $d(G)=3$ and $|G'|=p^2$, then $c(G)=2$,
$\Phi(G)=\mho_1(G)=Z(G)$, $G$ has a unique abelian subgroup $A$  of
index $p$ and $\exp(A)>p$, $G/G'$ has the type $(p^{n-4},p,p)$ and
$A/G'$ has the type $(p^{n-5},p,p)$.

\rr{5} if $G$ has an abelian subgroup  of index $p$, then $G$ is
metacyclic if and only if $p=2$;

\rr{6} if $p=2$, then $G$ is metacyclic if and only if $d(G)=2$;

\rr{7} if $|G'|=p$, then $\alpha_1(G)=p^2$;

\rr{8} if $d(G)=2$ and $|G'|=p^3$, then $\alpha_1(G)=1+p$;

\rr{9} if $d(G)=2$ and $G'\cong C_p^2$, then $\alpha_1(G)=p$.
\end{enumerate}
\end{lem}

\begin{lem}{\rm (\cite[$p_{27}$, Exercise 6]{Ber1})}
\label{number of abelian maximal subgroups} Let $G$ be a non-abelian
 $p$-group. Then the number of abelian subgroups of index $p$ in $G$ is $0,1$ or
$p+1$.

\end{lem}

\begin{lem}{\rm \cite[ $p_{
259}$, Aufgabe 2]{Hup})} \label{Ber} Suppose that a finite
non-abelian $p$-group $G$ has an abelian normal subgroup $A$, and
$G/A=\lg bA\rg$ is cyclic. Then the map $a\mapsto [a,b]$, $a\in A$,
is an epimorphism from $A$ to $G'$, and $G'\cong A/A\cap Z(G)$. In
particular, if a non-abelian $p$-group $G$ has an abelian subgroup
of index $p$, then $|G|=p|G'||Z(G)|$.
\end{lem}

\begin{lem}
\label{equivalent equation 1} Let $p$ be an odd prime. Then the
equation $x^2+ry^2-u=0$ about $x,y$ over $F_p$ has a solution, and
the following conclusions hold:

$(1)$ If $-r\in (F_p^*)^2$ and $u\in F_p^*$, then the
equation $x^2+ry^2-u=0$ has exactly $p-1$ solutions;

$(2)$ If $-r\not\in (F_p^*)^2$ and $u\in F_p^*$, then the
equation $x^2+ry^2-u=0$ has exactly $p+1$ solutions.
\end{lem}
\demo Let $A=\{ a^2\di a\in F_p\}$ and $B=\{u-rb^2\di b\in F_p\}$.
Then $|A|=|B|=\frac{p+1}{2}$. It follows that $A\cap B\ne
\emptyset$. Hence there exist $a,b\in F_p$ such that $a^2=u-rb^2$.

(1) Let $-r=\alpha^2$. Then the
equation $x^2+ry^2-u=0$ is $x^2-\alpha^2y^2=u$, and hence is equivalent to
$$\{\begin{array}{c}
x+\alpha y=k\\
x-\alpha y=k^{-1}u
\end{array}
$$
where $k=1,2,\dots,p-1$. Hence the
equation $x^2+ry^2-u=0$ has exactly $p-1$ solutions.

(2) Let $n(u)$ be the number of solutions of the equation
$x^2+ry^2-u=0$. Since $-r\not\in (F_p^*)^2$, $n(0)=1$. It is obvious
that $n(u_1)=n(u_2)$ for $u_1u_2^{-1}\in (F_p^*)^2$. Let $\nu\not\in
(F_p^*)^2$. Since $F_p\times F_p$ has a partition
$$\{(x,y)\in F_p\times F_p\mid x^2+ry^2\in (F_p^*)^2\}\cup \{(x,y)\in F_p\times F_p\mid x^2+ry^2\not\in (F_p^*)^2\},$$
$p^2=n(0)+\frac{p-1}{2}n(1)+\frac{p-1}{2}n(\nu)$. It follows that $n(1)+n(\nu)=2p+2$.
If we prove that $n(1)=p+1$, then $n(u)=p+1$ for all $u\in F_p^*$.

Now we calculate $n(1)$. If $y=0$, then the
equation $x^2+ry^2-1=0$ has two solutions $(1,0)$ and $(-1,0)$. If $y\neq 0$, then
the
equation $x^2+ry^2-1=0$ is equivalent to $(xy^{-1})^2-y^{-2}+r=0$. By (1),
the later has $p-1$ solutions. Hence $n(1)=p+1$.\qed

\begin{lem}
\label{equivalent equation} Let $p$ be an odd prime. If
$s^2-4r\not\equiv 0(\mod p)$, then the equation
$x^2+sxy+ry^2+wx+vy+u=0$ about $x,y$ over $F_p$ has a solution.
\end{lem}
\demo Let $x_1=x+2^{-1}sy$, $y_1=y$, $r_1=r-4^{-1}s^2$ and
$v_1=v-2^{-1}ws$. Then the equation is turned to
$x_1^2+r_1y_1^2+wx_1+v_1y_1+u=0$. Let $x_2=x_1+2^{-1}w$,
$y_2=y_1+2^{-1}r_1^{-1}v_1$ and
$u_2=u-4^{-1}w^2-4^{-1}r_1^{-1}v_1^2$. Then the equation is turned
to $x_2^2+r_1y_2^2+u_2=0$. By Lemma \ref{equivalent equation 1}, the
last equation has a solution $(a,b)$. Thus
$(a-2^{-1}w-2^{-1}sb+4^{-1}sr_1^{-1}v_1,b-2^{-1}r_1^{-1}v_1)$ is a
solution of the first equation.\qed

\begin{lem}{\rm (\cite[Lemma 3.1]{Alj})}
\label{alj1} Let $G$ be a non-abelian two-generator $p$-group having
an abelian subgroup $A$ of index $p$. Assume that $|G/G'|=p^{m+1}$
and $c(G)=c$. Then $m\ge 1$, $c\ge 2$ and
\begin{enumerate}

\rr1 $G$ has the lower central complexion
$(m+1,\underbrace{1,\dots,1}_{c-1})$ and hence $|G'|=p^{c-1}$,
$|G|=p^{m+c}$;

\rr2 $|\Z(G)|=p^m$ and $G/\Z(G)$ is of maximal class;

\rr3 $\Z(G)\leq \Phi(G)$, $\Phi(G)=G'\Z(G)$ and $G'\cap\Z(G)=G_c$;

\rr4 Let $M$ be a non-abelian subgroup  of index $p$ of $G$. Then
$\Z(M)=\Z(G)$ and $$M'=G_3,M_3=G_4,\cdots,M_{c-1}=G_c.$$

\end{enumerate}

\end{lem}

\begin{lem}{\rm (\cite[Corollary 3.8]{Alj})}\label{alj2}
Suppose that $G$ is a finite $p$-group having an abelian subgroup of
index $p$, and all non-abelian subgroups of $G$ are generated by two
elements. Then
\begin{enumerate}
\rr1 $G\in\mathcal{A}_2$ if and only if $c(G)=3$;

\rr2 If $p=2$, then $G$ is metacyclic;

\rr3 If $c(G)\le p$, then $d(G')=c(G)-1$. If $c(G)\ge p+1$, then
$d(G')=p-1$.
\end{enumerate}
\end{lem}

\begin{lem}{\rm (\cite[Theorem 4.1]{Alj})}
\label{Dp3=>A_2} Suppose that $G$ is a three-generator non-abelian
$p$-group, and all non-abelian proper subgroups of $G$ are generated
by two elements. If $G$ has an abelian subgroup of index $p$, then
$G\in \A_2$.
\end{lem}

\begin{lem}
\label{G' is abelian} {\rm (\cite[{\rm Theorem 4}]{Bla2})} Let $G$
be a $p$-group. If both $G$ and $G'$ can be generated by two
elements, then $G'$ is abelian.
\end{lem}

\begin{lem}\label{G'}
Suppose that $G$ is a finite $p$-group. $M_1$ and $M_2$ are two
distinct maximal subgroups of $G$. Then $|G'|\leq p|M_1'M_2'|$.
\end{lem}
\demo Let $\bar{G}=G/{M_1'M_2'}$. Then $\bar{G}$ has abelian
subgroups $\bar{M}_1$ and $\bar{M}_2$  of index $p$. Hence
$Z(\bar{G})\ge \bar{M}_1\cap\bar{M}_2$. By Lemma \ref{Ber},
$|\bar{G}'|\le p$. Thus $|G'|\leq p|M_1'M_2'|$.\qed

\begin{lem}\label{d(K) leq k+1 3}{\rm (\cite[{\rm Theorem 5.6}]{ADZ})}
If $G$ is a finite non-abelian $p$-group and $|G'|=p^k$, then $G$
has a subgroup $K$ such that $d(K)\le k+1$ and $K_n=G_n$ for all
$2\le n\le c(G)$.
\end{lem}

\begin{prop}{\rm (\cite{X})} \label{commutator fomula} Let $G$ be a metabelian group and
$a,b\in G$. For any positive integers $i$ and $j$, let
$$[ia,jb]=[a,b,\underbrace{a,\dots ,a}_{i-1},
\underbrace{b,\dots ,b}_{j-1}].$$ Then, for any positive integers
$m$ and $n$,
$$[a^m,b^n]=\prod_{i=1}^m\prod_{j=1}^n [ia,jb]^{{m\choose i} {n\choose j}}.$$

\end{prop}

The following lemma is equivalent to \cite[Theorem 3.4]{LQZ1}

\begin{lem}\label{class 2}
Assume $G$ is a finite $p$-group. If $c(H)\le 2$ for all $H<G$, then
$c(G)\le 3$.
\end{lem}

Let $G$ be a group of order $p^m$, $|G:\Phi(G)|=p^d$,
$$\Gamma_i=\{H<G\mid \Phi(G)\le H,\ |G:H|=p^i\}$$ and $\mathfrak{M}$ be
a set of proper subgroups of $G$. For $K\le G$, we denote by
$\alpha(K)$ the number of members of the set $\mathfrak{M}$ that are
subgroups of $K$. Obviously, $\alpha(G)=|\mathfrak{M}|$.

\begin{thm}\label{Hall}
$($Hall's enumeration principle$)$ In the above notation,
$$\alpha(G)=\sum_{i=1}^{d}\sum_{H\in \Gamma_i}^{}(-1)^{i-1}p^{i \choose 2}\alpha(H).$$

\end{thm}

\section{Some general properties of $\mathcal{A}_t$-groups with $t\geq 3$}

Y. Berkovich was the first where $\mathcal{A}_t$-groups were
investigated. He obtained many significant results. The following
are some related to $\mathcal{A}_3$-groups.

\begin{lem}\label{metacyclic An}{\rm(\cite[Lemma J(i)]{Yak0})}
Let $G$ be a metacyclic $p$-group. Then $G$ is an
$\mathcal{A}_t$-group if and only if $|G'|=p^t$.
\end{lem}

\begin{lem}\label{A3-groups property1}{\rm(\cite[Proposition 72.2]{Ber2})}
If a $p$-group $G$ is an $\mathcal{A}_3$-group, then $|G'|\le p^4$.
For $p=5$, there exists  an $\mathcal{A}_3$-group $G$ with
$|G'|=5^4$.
\end{lem}

\begin{lem}\label{A3-groups property2}{\rm(\cite[Proposition 72.3]{Ber2})}
Suppose that a $p$-group $G$ is an $\mathcal{A}_3$-group with
$|G'|=p^4$. Then

{\rm (a)}\ $G'$ is abelian.

{\rm (b)}\ If $\exp(G')=p$, then $p>2$, $G/G'$ is abelian of type
$(p^n, p)(n\geq 1)$ and $G'\cong C_p^m$. \mbox{\hskip0.45in}If, in
addition, $n>1$, then $\Omega_1(G)=G'$.

{\rm (c)}\ $G'$ is not metacyclic.
\end{lem}

\begin{thm}\label{A3-groups property3}{\rm(\cite[Theorem 72.6]{Ber2})}
If a $p$-group $G$ is an $\mathcal{A}_4$-group, then $|G'|\le p^6$.
\end{thm}

In following we give more information about  $\mathcal
{A}_t$-groups.

\begin{lem}{\rm (\cite[{\rm Corollary 2.4}]{ALQZ})}
\label{A_t+k} {\rm (1)} Let $M$ be an $\mathcal {A}_t$-group, and
$A$ be an abelian group of order $p^k$. Then $G=M\times A$ is an
$\mathcal {A}_{t+k}$-group.

{\rm (2)} Let $M$ be an $\mathcal {A}_t$-group with $|M'|=p$,
$G=M\ast A$, where $A$ is an abelian group of order $p^{k+1}$ and
$M\cap A=G'$. Then $G$ is an $\mathcal {A}_{t+k}$-group.
\end{lem}

\begin{lem}{\rm (\cite[{\rm Theorem 5.4}]{QZGA})}\label{54} Assume $G$ is a finite $2$-group
having a unique $\mathcal{A}_1$-subgroup $A$ of index $p$. If
$|G|=2^t\ge 2^9$ and $G/A'$ has a unique abelian subgroup of index
$p$, then one of the following is true:

\medskip

$(${\rm 1}$)$ $G$ is metacyclic;

$(${\rm 2}$)$ $G$ is an $\mathcal{A}_{t-2}$-group or an
$\mathcal{A}_{t-3}$-group.

\end{lem}

\begin{lem}{\rm (\cite[{\rm Corollary 5.5}]{QZGA})}\label{55}
Assume that $G$ is a finite $2$-group
having a unique $\mathcal{A}_1$-subgroup $A$  of index $p$. If $G$
is an $\mathcal{A}_3$-group and $|G|\ge 2^9$, then $G$ is
metacyclic.
\end{lem}

\begin{lem}\label{alj3}
Suppose that $G$ is a finite $p$-group having an abelian subgroup of
index $p$. If all non-abelian subgroups of $G$ are generated by two
elements, then $G$ is an $\mathcal{A}_t$-group if and only if
$c(G)=t+1$.
\end{lem}

\demo Suppose that $c(G)=t+1$ and $M$ is a non-abelian subgroup of
index $p$. By Lemma \ref{alj1}(4), $c(M)=t$. By hypothesis,
$d(M)=2$. By induction on $t$, $M$ is an $\mathcal{A}_{t-1}$-group.
Hence $G$ is an $\mathcal{A}_t$-group. \qed

\begin{lem}\label{alj7}
If $G$ is a $p$-group of maximal class of order $p^n$ with $n\geq
3$, then $G$ is an $\mathcal{A}_{n-2}$-group.
\end{lem}

\demo By \cite[Theorem 9.6(f)]{Ber1}, there exists a subgroup $A$ of
order $p^2$ of $G$ such that $C_G(A)=A$. Let $B$ be a subgroup of
$G$ with $B\geq A$ and $|B|=p^3$. Then $B$ is not abelian. It
follows that $G\in {\mathcal{A}}_{n-2}$. \qed

\begin{lem}\label{alj11}
If a $p$-group $G$ is an $\mathcal{A}_t$-group and $G'$ is regular,
then $\exp(G')\leq p^t$.
\end{lem}

\demo By induction on $t$ we have $\exp(M')\le p^{t-1}$ for all
$M\lessdot G$. Let $N=\prod\limits_{M\lessdot G}M'$. Since $G'$ is
regular, $\exp(N)\le p^{t-1}$ by \cite[Theorem 7.2(b)]{Ber1}. It is
easy to see that $G/N$ is abelian or an $\mathcal{A}_1$-group. It
follows by Lemma \ref{minimal non-abelian equivalent conditions}(2)
that $|G'N/N|=|(G/N)'|\le p$. Thus $\exp(G')\le p\cdot p^{t-1}=p^t$.
\qed

\begin{cor}
Assume $G$ is a finite $p$-group.

{\rm (1)}\ If $G$  is an $\mathcal{A}_3$-group, then $\exp(G')\leq
p^3$.

{\rm (2)}\ If $p\ge 7$ and  $G$  is an $\mathcal{A}_4$-group, then
$\exp(G')\leq p^4$.
\end{cor}

\demo By Lemma \ref{A3-groups property2}(a), Theorem \ref{A3-groups
property3} and \cite[Theorem 7.1(b)]{Ber1} we get $G'$ is regular.
Thus the conclusion follows by Lemma \ref{alj11}. \qed

\begin{thm}\label{alj6}
If a $p$-group $G$ is an $\mathcal{A}_3$-group, then $c(G)\leq 4$.
\end{thm}

\demo Let $N=\prod\limits_{M\lessdot G}M_3$.  Since $M_3$ is
characteristic in $M$, $M_3\lhd G$. Thus $N\lhd G$. By Lemma
\ref{class 2}, $c(G/N)\le 3$. If $N\le Z(G)$, then $c(G)\le 4$.
Assume $N\nleq Z(G)$. Then there exists $M\lessdot G$ such that
$M_3\nleq Z(G)$. Thus $|M_3|\ge p^2$. Notice that $M$ is an
$\mathcal{A}_2$-group. Hence $|M|=p^5$ by checking the list of
groups in Lemma \ref{A_2}. Thus $|G|=p^6$. If $c(G)=5$, then $G$ is
of maximal class. It follows by Lemma \ref{alj7} that $G$ an
$\mathcal{A}_4$-group. This is a contradiction. \qed

\begin{rem}
There exist $\mathcal{A}_3$-groups of class $4$. In fact, those
groups listed in Theorem \ref{d=3-1} are all $\mathcal{A}_3$-groups
of class 4.
\end{rem}

\begin{lem}\label{alj5}
Assume a $p$-group $G$ is an $\mathcal{A}_3$-group.

$(1)$\ If $\Phi(G)$ is non-abelian, then $\Phi(G)$ is metacyclic.

$(2)$\ If $p>2$, then $\Phi(G)$ is non-abelian if and only if $G$ is
metacyclic.
\end{lem}

\demo (1) Since $G$ is an $\mathcal{A}_3$-group, $d(G)=2$ and
$\Phi(G)\in \mathcal{A}_1$. By Lemma \ref{minimal non-abelian
equivalent conditions}, $d(\Phi(G))=2$. It follows from
\cite[Theorem 44.12]{Ber1} that $\Phi(G)$ is metacyclic.

(2)\ $\Longleftarrow$ Let $G=\langle a, b\rangle$.  Then
$o([a,b])=p^3$ by Lemma \ref{metacyclic An}. Since $G$ is regular,
$[a^p,b^p]=1$ if and only if $[a,b]^{p^2}=1$ by \cite[Theorem
7.2(e)]{Ber1}. Since $o([a,b])=p^3$, $[a^p,b^p]\neq 1$. It follows
that $\Phi(G)$ is non-abelian.

$\Longrightarrow$ Assume $G$ is not metacyclic. It follows by Lemma
\ref{metacyclic An} and Lemma \ref{A3-groups property2}(a) and (c)
that $|G': \Omega_1(G')|\le p$. Thus $G_3\le \Omega_1(G')$. It
follows that $\exp (G_3)\le p$.  on the other hand, $c(G)\le 4$ by
Lemma \ref{alj6}. Since $G'$ is abelian, $G$ is metabelian. By lemma
\ref{commutator fomula} we have
$$[x^p,y^p]=\prod\limits_{i=1}^{p}\prod\limits_{j=1}^{p}[ix,jy]^{{p\choose
i} {p\choose j}},\ \ \ \
[x^p,z]=\prod\limits_{i=1}^{p}[ix,z]^{{p\choose i}}.$$ where  $x,
y\in G, z\in G'$.

It follows from $\exp(G')\le p^2$, $\exp (G_3)\le p$, $c(G)\le 4$
and $p>2$ that $[x^p,y^p]=1$ and $[x^p,z]=1$. Thus $\Phi(G)$ is
abelian. This is a contradiction. \qed


\begin{rem}
In Lemma $\ref{alj5}$, if $p=2$, then $(2)$  is not true. For
example, $D_{32}$ is a metacyclic $\mathcal{A}_3$-group, but
$\Phi(D_{32})\cong C_8$.
\end{rem}

\begin{lem}\label{lem a0,b0,M' notin Z(G)}
Suppose that $G$ is an $\mathcal{A}_3$-group. If $G$ has a maximal
subgroup $M$ such that $d(M)=3$ and $M'\not\le Z(G)$, then $G$ has a
three-generator maximal subgroup $L$ such that $d(L)=3$, $|L'|\le
p^2$ and $L'\not\le Z(G)$.
\end{lem}

\demo Otherwise, for every three-generator maximal subgroup $L$ with
$d(L)=3$ and $|L'|\le p^2$, we have $L'\le Z(G)$. Hence $|M'|\ge
p^3$. It follows from Lemma \ref{A_2-property}(2) that $$p=2,\
|M|=2^6\ {\rm and}\ M'=\Omega_1(M)=Z(M)\cong C_2^3.$$ Let $d(G)=r$
where $2\le r\le 4$. Then $G$ has $s=1+2+\dots+2^{r-1}$ maximal
subgroups. Let $H_1,H_2,\cdots,H_{s-1}$ and $M$ be all maximal
subgroups of $G$. Since $M'\le H_i$, $H_i$ is not metacyclic. It
follows from Lemma \ref{A_2-property}(6) that $d(H_i)=3$. If
$|H'|=8$, then we also have $H'=\Omega_1(H)=Z(H)$. Since $M'\cong
C_2^3\le H$, $M'=H'$. It follows that $[M',G]=[M',HM]=1$ and hence
$M'\le Z(G)$, a contradiction. Hence $|H_i'|\le 4$ and $H_i'\le
Z(G)$. Let $N=H_1'H_2'\cdots H_{s-1}'$. Then $G/N$ has $s-1$ abelian
subgroups of index $p$. Since $s-1\ge 2$, we have, by Lemma
\ref{number of abelian maximal subgroups}, the number of abelian
subgroups of index $p$ of $G/N$ is $3$. Hence $s-1\le 3$. It follows
that $d(G)=2$ and $M/N$ is also abelian. Thus $M'\le N\le Z(G)$, a
contradiction. \qed

\begin{lem}\label{d=4}Suppose that $G$ is a four-generator
$\mathcal{A}_3$-group. Then $c(G)=2$, $\Phi(G)\le Z(G)$, $G'\le
C_p^3$ and  all $\mathcal{A}_1$-subgroups of $G$ contain $\Phi(G)$.
\end{lem}

\demo Let $M$ be an $\mathcal{A}_1$-subgroup of $G$ and $L$ be a
maximal subgroup containing $M$. Since $d(G)=4$, we have $d(L)\ge
3$. Hence $L\in\mathcal{A}_2$. By Lemma \ref{A_2-property}(2), we
get $d(L)=3$, $c(L)=2$ and $L'\le C_p^3$. By comparing the indexes
of $\Phi(M),\Phi(L)$ and $\Phi(G)$, we get
$$Z(M)=\Phi(M)=\Phi(L)=\Phi(G).$$ By the arbitrariness of $M$ we get
$\Phi(G)\le Z(G)$. Hence $c(G)=2$. Let $N=\prod\limits_{L\maxsgp G}
L'$. Then $N\le Z(G)$ and $\exp(N)=p$. Since all subgroups of $G/N$
are abelian, we have $G/N$ is also abelian. It follows that $G'=N$.
Hence $\exp(G')=p$. Since $G'\le
\Omega_1(\Phi(G))=\Omega_1(\Phi(M))$, $G'\le C_p^3$. \qed

\begin{lem}\label{alpha of d=3}
Let $G$ be an $\mathcal{A}_3$-group such that $d(G)=3$ and $\Phi(G)\le Z(G)$.
Then $\alpha_1(G)=\mu_1+p^2\mu_2$.
\end{lem}
\demo Groups satisfying $d(G)=3$ and $\Phi(G)\le Z(G)$ were classified in \cite{ALQZ,QXA}.
The derived groups of non-abelian maximal subgroups of $G$ are of order $p$.

Let $H\in \Gamma_2$. Then, by $\Phi(G)\le Z(G)$, $H$ is abelian and
hence $\alpha_1(H)=0$. Let $M\in \Gamma_1$. If $M\in\mathcal{A}_2$,
then, by Lemma \ref{A_2-property}(7), $\alpha_1(M)=p^2$. By Hall's
enumeration principle, $ \alpha_1(G)=\sum_{H\in \Gamma_1}
\alpha_1(H)=\mu_1+p^2\mu_2.$ \qed

\section{$\mathcal{A}_3$-groups with an $\mathcal{A}_1$-subgroup of index $p$}

In this section, we classify $\mathcal{A}_3$-groups with an
$\mathcal{A}_1$-subgroup of index $p$. Since finite $p$-groups with
an $\mathcal{A}_1$-subgroup of index $p$ are classified by
\cite{ALQZ,AHZ,QYXA,QXA,QZGA}, it is enough to check whether those
groups in \cite{ALQZ,AHZ,QYXA,QXA,QZGA} are $\mathcal{A}_3$-groups
or not. In this case, $\mathcal{A}_3$-groups have 72 non-isomorphic
types.

For convenience, in this section, assume $G$ is an
$\mathcal{A}_3$-group with an $\mathcal{A}_1$-subgroup of index $p$
in Theorem \ref{d=2-1}, \ref{d=2-2}, \ref{d=2-3}, \ref{d=2-4} and
\ref{d=2-5}. It is easy to see that $d(G)\le 3$.

\begin{thm}\label{d=2-1}
$G$ has an abelian subgroup of index $p$ and $d(G)=2$ if and only if
$G$ is isomorphic to one of the following pairwise non-isomorphic
groups:
\begin{enumerate}
\rr{Ai} $c(G)=3$ and $G'\cong C_4$.
\begin{enumerate}
\rr{A1}  $\langle a, b; c \di a^{4}=b^2=c^4=1,
[a,b]=c,[c,a]=[c,b]=c^2 \rangle$. Moreover $|G|=2^5$, $\Phi(G)=\lg
a^2,c\rg\cong C_4 \times C_2$, $G'=\lg c\rg$ and $Z(G)=\lg
a^2,c^2\rg\cong C_2^2$.

\rr{A2} $\langle a, b; c \di a^{4}=b^{4}=1, c^2=a^{2},
[a,b]=c,[c,a]=[c,b]=c^2 \rangle$. Moreover, $|G|=2^5$, $\Phi(G)=\lg
b^2,c\rg \cong C_4 \times C_2$,
$G'=\lg c\rg$ and $Z(G)=\lg b^2,c^2\rg\cong C_2^2$.

\rr{A3}  $\langle a, b; c \di a^{8}=b^2=1,c^{2}=a^4,
[a,b]=c,[c,a]=[c,b]=c^2 \rangle$. Moreover, $|G|=2^5$, $\Phi(G)=\lg
a^2,c\rg\cong C_4 \times C_2$,
$G'=\lg c\rg$ and $Z(G)=\lg a^2\rg\cong C_4$.
\end{enumerate}
\rr{Aii}$c(G)=3$ and $G'\cong C_p^2$ where $p>2$.

\begin{enumerate}

\rr{A4}  $\langle a, b; c \di a^{p^{3}}=b^{p}=c^p=1,
[a,b]=c,[c,a]=1,[c,b]=a^{\nu p^2} \rangle$, where $\nu=1$ or a fixed
quadratic non-residue modulo $p$. Moreover, $|G|=p^5$, $\Phi(G)=\lg
a^p,c\rg\cong C_{p^2}\times C_p$, $G'=\lg a^{p^2},c\rg$ and
$Z(G)=\lg a^p\rg\cong C_{p^2}$.

\rr{A5} $\langle a, b; c \di a^{p^{2}}=b^{p^2}=c^p=1, 
[a,b]=c,[c,a]=1,[c,b]=b^{p} \rangle$. Moreover, $|G|=p^5$, $\Phi(G)=\lg
a^p,b^p,c\rg\cong C_p^3$, $G'=\lg b^{p},c\rg$ and $Z(G)=\lg
a^p,b^p\rg\cong C_p^2$.

\rr{A6}  $\langle a, b; c,d \di a^{p^{2}}=b^{p}=c^p=d^p=1,
[a,b]=c,[c,a]=1,[c,b]=d, [d,a]=[d,b]=1 \rangle$.  Moreover,
$|G|=p^5$, $\Phi(G)=\lg a^p,c,d\rg\cong C_p^3$,
$G'=\lg c,d\rg$ and $Z(G)=\lg a^p,d\rg\cong C_p^2$.
\end{enumerate}
\end{enumerate}
Moreover, $(\mu_0,\mu_1,\mu_2)=(1,p-1,1)$ and
 $\alpha_1(G)=p^2+p-1$.
\end{thm}

\demo  Let $A$ be an abelian subgroup of index $p$ and $B$ be an
$\mathcal{A}_1$-subgroup of index $p$. By Lemma \ref{alj1},
$G_3=B'$, $c(G)=3$ and $|G'|=p^2$. Let $\bar{G}=G/G_3$. Then
$|\bar{G}'|=p$.  Since $d(\bar{G})=d(G)=2$, $\bar{G}\in
\mathcal{A}_1$ by Lemma \ref{minimal non-abelian equivalent
conditions}. Obviously, $G_3=\Phi(G')G_3$. If $\bar{G}$ is
metacyclic, then $G$ is also metacyclic by Lemma \ref{metacyclic}.
Hence $G$ is an $\mathcal{A}_2$-group by Lemma \ref{metacyclic An}.
This contradicts that $G$ is an $\mathcal{A}_3$-group. Thus
$\bar{G}$ is a non-metacyclic $\mathcal{A}_1$-group. By Lemma
\ref{thm=Redei}, $\bar{G}\cong M_p(n,m,1)$. Now we have
$\Phi(G')G_3\le C_p^2$, $\Phi(G')G_3\le Z(G)$ and
$G/\Phi(G')G_3\cong M_p(n,m,1)$. Thus $G$ is one of the groups
classified in \cite{AHZ}. Since $|G'|=p^2$, $G$ is either one of the
groups listed in \cite[Theorem 3.5]{AHZ} or \cite[Theorem 4.6]{AHZ}.
By hypothesis, we only need pick out those groups satisfying the
following three conditions:

(1) the minimal index of $\mathcal{A}_1$-subgroups is $1$;

(2) the maximal index of $\mathcal{A}_1$-subgroups is $2$;

(3) $G$ has an abelian subgroup of index $p$.

\medskip

If $G'\cong C_{p^2}$, then $G$ is one of the groups listed in
\cite[Theorem 3.5]{AHZ}. \cite[Theorem 3.1 \& 3.6]{AHZ} tell us
those groups in \cite[Theorem 3.5]{AHZ} satisfy the conditions (1)
and (2). Moreover, it is easy to check which one satisfies the
condition (3) in these groups. Thus we get the groups (A1)--(A3). If
$G'\cong C_p^2$, then $G$ is one of the groups listed in
\cite[Theorem 4.6]{AHZ}. \cite[Theorem 4.1 \& 4.7]{AHZ} tell us
those groups in \cite[Theorem 4.6]{AHZ}  satisfy the conditions (1)
and (2). Moreover, in these groups we check which one satisfies the
condition (3). Thus we get the groups (A4)--(A6). For convenience,
in Table \ref{table1} we give the relationship between the groups in
Theorem \ref{d=2-1} and those groups in paper \cite{AHZ}.
\begin{table}[h]
  \centering
{ \tiny
\begin{tabular}{c|c||c|c}
\hline
Groups in Theorem \ref{d=2-1}& Groups in \cite[Theorem 3.5]{AHZ}& Groups in Theorem \ref{d=2-1}& Groups in \cite[Theorem 4.6]{AHZ} \\
\hline
(A1)&(A1)& (A4)&(J1) where $n=2$ and $m=1$\\
(A2)&(A2)& (A5)&(J3) where $n=2$ and $m=1$\\
(A3)&(A3)& (A6)&(J5) where $n=2$ and $m=1$\\
\hline
\end{tabular}
\caption{The correspondence from Theorem \ref{d=2-1} to
\cite[Theorem 3.5 \& 4.6]{AHZ}} \label{table1}}
\end{table}

In following, we calculate the $(\mu_0,\mu_1,\mu_2)$ and
$\alpha_1(G)$.

Since $d(G)=2$, $G$ has $1+p$ maximal subgroups. For any
$\mathcal{A}_3$-group, we always have $\mu_2\ge 1$. By the
hypothesis, $\mu_0\ge 1$ and $\mu_1\ge 1$.

If $G$ is one of the groups {\rm(A1)--(A3)}, then $p=2$. Hence $\mu_0+\mu_1+\mu_2=1+p=3$.
It follows that $(\mu_0,\mu_1,\mu_2)=(1,1,1)$.

If $G$ is one of the groups {\rm(A4)--(A6)}, then it is easy to see that $\lg c,a,\Phi(G)\rg$
and $\lg c,b,\Phi(G)\rg$ are the unique abelian subgroup of index $p$ and
the $\mathcal{A}_2$-subgroup of index $p$ respectively. Hence $(\mu_0,\mu_1,\mu_2)=(1,p-1,1)$.

Now we calculate $\alpha_1(G)$. It is obvious that $\alpha_1(\Phi(G))=0$. Let $M\in\Gamma_1$. If
$M$ is the $\mathcal{A}_2$-subgroup of index $p$, then $|M'|=p$. By Lemma \ref{A_2-property} (7), $\alpha_1(M)=p^2$.

By Hall's enumeration principle,
$$
\alpha_1(G)=\sum_{H\in \Gamma_1}\alpha_1(H) -p\alpha_1(\Phi(G))=\sum_{H\in
\Gamma_1} \alpha_1(H)=\mu_1+p^2\mu_2=p^2+p-1.
$$
\qed

\begin{thm}\label{d=2-2}
$G$ has an abelian subgroup of index $p$ and $d(G)=3$ if and only if
$G$ is isomorphic to one of the following pairwise non-isomorphic
groups:

\begin{enumerate}
 \rr{Bi} $c(G)=2$ and $G'\cong C_p$.
\begin{enumerate}
\rr{B1}  $\langle a, b, c \di a^4=c^{4}=1,\ b^2=a^2=[a,b],\
[c,a]=[c,b]=1 \rangle  \cong  Q_8 \times C_{4}$; where  $|G|=2^5$,
$\Phi(G)=\lg a^2,c^2\rg\cong C_2 \times C_2$, $G'=\lg a^2\rg$ and
$Z(G)=\lg a^2,c\rg\cong C_4 \times C_2$.

\rr{B2}  $\langle a, b, c \di a^{p^{n+1}}=b^{p^m}=c^{p^2}=1,\
[a,b]=a^{p^{n}},\ [c,a]=[c,b]=1\rangle\cong M_p(n+1,m)\times
C_{p^2}$, where $\min\{n,m \}=1$; where $|G|=p^{m+n+3}$,
$\Phi(G)=\lg a^p, b^p,c^p\rg \cong C_{p^n} \times C_{p^{m-1}}\times
C_{p}$ if $m>1$, $\Phi(G) \cong C_{p^n} \times C_{p}$ if $m=1$,
$G'=\lg a^{p^n}\rg$, $Z(G)=\lg a^p, b^p,c\rg \cong C_{p^n} \times
C_{p^{m-1}}\times C_{p^2}$ if $m>1$ and $Z(G) \cong C_{p^n} \times
C_{p^2}$ if $m=1$.

\rr{B3}  $\langle a, b, c;d \di a^{p^n}=b^{p}=c^{p^2}=d^p=1,\
[a,b]=d,\ [c,a]=[c,b]=1\rangle\cong M_p(n,1,1)\times C_{p^2}$, where
$n\ge 2$ for $p=2$; where  $|G|=p^{n+4}$, $\Phi(G)=\lg a^p,c^p,d\rg
\cong C_{p^{n-1}} \times C_{p}\times C_{p}$ if $n>1$, $\Phi(G)\cong
C_{p}\times C_{p}$ if $n=1$, $G'=\lg d\rg$, $Z(G)=\lg a^p,c,d\rg
\cong C_{p^{n-1}} \times C_{p^{2}}\times C_{p}$ if $n>1$, $Z(G)
\cong C_{p^{2}}\times C_{p}$ if $n=1$.

\rr{B4}  $\langle a, b, c \di a^4=1,\ b^2=c^{4}=a^2=[a,b],\
[c,a]=[c,b]=1\rangle  \cong  Q_8 \ast C_{8}$; where  $|G|=2^{5}$,
$\Phi(G)=\lg c^2\rg \cong C_4$, $G'=\lg c^4\rg$ and $Z(G)=\lg c\rg
\cong C_8$.

\rr{B5}  $\langle a, b, c \di a^{p^n}=b^{p}=c^{p^{3}}=1,\
[a,b]=c^{p^2},\ [c,a]=[c,b]=1 \rangle\cong M_p(n,1,1)\ast
C_{p^{3}}$, where $n\ge 2$ for $p=2$. Moreover,  $|G|=p^{n+4}$,
$\Phi(G)=\lg a^p,c^p\rg \cong C_{p^{n-1}} \times C_{p^{2}}$ if
$n>1$, $\Phi(G)\cong C_{p^{2}}$ if $n=1$, $G'=\lg c^{p^2}\rg$,
$Z(G)=\lg a^p,c\rg \cong C_{p^{n-1}} \times C_{p^{3}}$ if $n>1$,
$Z(G)\cong C_{p^{3}}$ if $n=1$.
\end{enumerate}
\rr{Bii} $c(G)=2$ and $G'\cong C_p^2$.
\begin{enumerate}
\rr{B6} $\langle a, b, c \di a^{p}=b^{p^{2}}=c^{p^{2}}=1,
[b,c]=1,[c,a]=c^{p},[a,b]=b^{-p}\rangle$, where $p$ is odd; and
$|G|=p^{5}$, $\Phi(G)=G'=Z(G)=\lg b^{p},c^p\rg \cong C_p^2$.

\rr{B7} $\langle a, b, c \di a^{p^{l}}=b^{p^{2}}=c^{p^{2}}=1,
[b,c]=1,[c,a]=b^{p}c^{p},[a,b]=b^{-p}\rangle$, where $p$ is odd; and
$|G|=p^{l+4}$, $\Phi(G)=Z(G)=\lg a^p, b^p,c^p\rg \cong C_{p^{l-1}}
\times C_{p}\times C_{p}$ if $l>1$, $\Phi(G)=Z(G) \cong C_{p}\times
C_{p}$ if $l=1$, $G'=\lg b^{p},c^p\rg$.

\rr{B8} $\langle a, b, c \di a^{p^{l}}=b^{p^{2}}=c^{p^{2}}=1,
[b,c]=1,[c,a]=b^{p}c^{tp},[a,b]=b^{-tp}c^{\nu p}\rangle$, where $p$
is odd, $\nu=1$ or a fixed quadratic non-residue modulo $p$,
$-\nu\in (F_p^*)^2$, and $t\in\{0,1,\dots,\frac{p-1}{2}\}$ such that
$t^2\neq -\nu$; and $|G|=p^{l+4}$, $\Phi(G)=Z(G)=\lg a^p, b^p,c^p\rg
\cong C_{p^{l-1}} \times C_{p}\times C_{p}$ if $l>1$, $\Phi(G)=Z(G)
\cong C_{p}\times C_{p}$ if $l=1$, $G'=\lg b^{p},c^p\rg$.

\rr{B9} $\langle a, b, c \di a^{p}=b^{p^3}=c^{p^3}=1,
[b,c]=1,[c,a]=b^{p^{2}}c^{tp^{2}},[a,b]=b^{-tp^2}c^{\nu
p^2}\rangle$, where $p$ is odd, $\nu=1$ or a fixed quadratic
non-residue modulo $p$, $-\nu\not\in (F_p^*)^2$, and
$t\in\{0,1,\dots,\frac{p-1}{2}\}$; moreover, $|G|=p^{7}$,
$\Phi(G)=Z(G)=\lg b^p,c^p\rg\cong C_{p^2}\times  C_{p^2}$, $G'=\lg
b^{p^2},c^{p^2}\rg$.

\rr{B10} $\langle a, b, c \di a^{2^{l}}=b^{4}=c^{4}=1,
[b,c]=1,[c,a]=b^{2},[a,b]=c^{2}\rangle$; where $|G|=2^{l+4}$,
$\Phi(G)=Z(G)=\lg a^2, b^2,c^2\rg \cong C_{2^{l-1}} \times
C_{2}\times C_{2}$ if $l>1$, $\Phi(G)=Z(G) \cong C_{2}\times C_{2}$
if $l=1$, $G'=\lg b^{2},c^2\rg$.

\rr{B11}  $\langle a, b, c \di a^{2}=b^{8}=c^{8}=1,
[b,c]=1,[c,a]=b^{4},[a,b]=b^{4}c^{4}\rangle$; where  $|G|=2^{7}$,
$\Phi(G)=Z(G)=\lg  b^2,c^2\rg \cong C_4\times C_4$, $G'=\lg
b^{4},c^4\rg$.

 \rr{B12}  $\langle
a, b, c \di a^{p^{l}}=b^{p^{3}}=c^{p^{2}}=1,[b,c]=1, [a,b]=c^{\nu
p},[c,a]=b^{p^2} \rangle$, where $\nu=1$ or a fixed quadratic
non-residue modulo $p$;  and $|G|=p^{l+5}$, $\Phi(G)=Z(G)=\lg a^p,
b^p,c^p\rg \cong C_{p^{l-1}} \times C_{p^2}\times C_{p}$ if $l>1$,
$\Phi(G)=Z(G) \cong C_{p^2}\times C_{p}$ if $l=1$, $G'=\lg
b^{p^2},c^p\rg$.

\rr{B13}  $\langle a, b, c \di
a^{p^{l+1}}=b^{p^{m}}=c^{p^{2}}=1,[b,c]=1,[c,a]=c^{p} ,
[a,b]=a^{p^l}\rangle$, where $m\le 2$; moreover, $|G|=p^{l+m+3}$,
$\Phi(G)=Z(G)=\lg a^p, b^p,c^p\rg \cong C_{p^{l}}\times C_{p}^2$ if $m=2$, $\Phi(G)=Z(G) \cong C_{p^l}\times
C_{p}$ if $m=1$, $G'=\lg a^{p^l},c^p\rg$.

\rr{B14}  $\langle a, b, c \di
a^{p^{l+1}}=b^{p^{m+1}}=c^{p^{n}}=1,[b,c]=1,[c,a]=b^{p^{m}},
[a,b]=a^{p^l} \rangle$, where $\min\{l,m,n\}=1$ and $\max\{m,n\}=2$;
moreover, $|G|=p^{l+m+n+2}$, $\Phi(G)=Z(G)=\lg a^p, b^p,c^p\rg \cong
C_{p^{l}} \times C_{p^{m}}\times C_{p^{n-1}}$ if $n>1$,
$\Phi(G)=Z(G)\cong C_{p^{l}} \times C_{p^{m}}$ if $n=1$, $G'=\lg
a^{p^l},b^{p^m}\rg$.

\rr{B15}  $\langle a, b, c \di a^{4}=b^{4}=c^{4}=1,[b,c]=1,
[c,a]=a^2c^2,[a,b]=c^2 \rangle$; where $|G|=2^{6}$,
$\Phi(G)=Z(G)=\lg a^2, b^2,c^2\rg \cong C_2^3$, $G'=\lg
a^{2},c^2\rg$.

\rr{B16}  $\langle a, b, c \di a^{4}=b^{4}=c^4=1,[b,c]=1,
[c,a]=a^2=c^2, [a,b]=b^{2} \rangle$; where $|G|=2^{5}$,
$\Phi(G)=Z(G)=G'=\lg a^2, b^2\rg \cong C_2^2$.

\rr{B17}  $\langle a, b, c; x\di a^{p^{l}}=b^{p}=c^{p^{2}}=x^p=1,
[a,b]=x,[a,c]=c^{p},[b,c]=[x,a]=[x,b]=[x,c]=1 \rangle$, where $l\ge
2$ if $p=2$; moreover, $|G|=p^{l+4}$, $\Phi(G)=Z(G)=\lg a^p, c^p,
x\rg \cong C_{p^{l-1}} \times C_{p}\times C_{p}$ if $l>1$,
$\Phi(G)=Z(G) \cong C_{p}\times C_{p}$ if $l=1$, $G'=\lg
c^{p},x\rg$.

\rr{B18}  $\langle a, b, c; x\di
a^{p^{l}}=b^{p^{m}}=c^{p^{2}}=x^p=1,
[a,b]=c^{p},[a,c]=x,[b,c]=[x,a]=[x,b]=[x,c]=1 \rangle$, where $l\ge
2$ if $p=2$, $m\le 2$ and $\min\{l,m\}=1$; moreover,
$|G|=p^{l+m+3}$, $\Phi(G)=Z(G)=\lg a^p, b^p,c^p,x\rg\cong
C_{p^{l-1}} \times C_{p}\times C_{p}$ if $m=1$, $\Phi(G)=Z(G)\cong
C_{p^{m-1}}\times C_{p}\times C_{p}$ if $l=1$, $G'=\lg c^{p},x\rg$.

\rr{B19}  $\langle a, b, c; x\di a^{p^{l+1}}=b^{p^{m}}=c^{p}=x^p=1,
[a,b]=a^{p^{l}},[a,c]=x,[b,c]=[x,a]=[x,b]=[x,c]=1 \rangle$, where
$m\le 2$; moreover, $|G|=p^{l+m+3}$, $\Phi(G)=Z(G)=\lg a^p,
b^p,x\rg\cong C_{p^{l}} \times C_{p^{m-1}}\times C_{p}$ if $m>1$,
$\Phi(G)=Z(G)\cong C_{p^{l}} \times C_{p}$ if $m=1$, $G'=\lg
a^{p^l},x\rg$.

\rr{B20}  $\langle a, b,c; x \di a^{4}=b^{2}=c^{4}=x^2=1,
[a,b]=x,[a,c]=a^2=c^2,[b,c]=[x,a]=[x,b]=[x,c]=1 \rangle$; where
$|G|=2^{5}$, $\Phi(G)=Z(G)=G'=\lg a^2, x\rg\cong C_2^2$.
\end{enumerate}
\end{enumerate}

Moreover, Table $\ref{table2}$ gives $(\mu_0,\mu_1,\mu_2)$ and
 $\alpha_1(G)$ for the groups {\rm (B1)--(B20)}.

\begin{table}[h]
\centering
 {
\tiny
\begin{tabular}[t]{c|c||c}
\hline
 $(\mu_0,\ \mu_1,\ \mu_2)$&$\alpha_1(G)$      & types of  $\mathcal{A}$$_3$ groups \\

\hline $(p+1,\ p^2-1,\ 1)$  & $2p^2-1$&  (B1); (B2) where $m=n= 1$; (B3) where $n=1$; (B4); (B5)  where $n=1$\\

\hline $(p+1,\ p^2-p,\ p)$  & $p^3+p^2-p$& (B2) where $m>1=n$ or $n>1=m$; (B3) where $n>1$; (B5)  where $n>1$ \\

\hline \multirow{2}*{$(1,\ p^2-1,\ p+1 )$} & \multirow{2}*{$p^3+2p^2-1$} &  (B6); (B9); (B11); (B13) where $p=2$ and $m=l=1$\\
& &(B17) where $l=1$; (B19) where $p=2$ and $m=l=1$\\

\hline \multirow{3}*{$(1,\ p^2,\ p )$ } & \multirow{3}*{$p^3+p^2$}&  (B7) where $l>1$; (B10); (B12) where $l>1$; (B13) where $m=1$ and $l>1$\\
& &(B14) where $n=1$ or $m=1$; (B16); (B18) where $m=1$\\
 &                             &(B19) where $p=2$ and $l>1=m$; (B19) where $p>2$ and $m=1$\\

\hline $(1,\ p^2+p-1,\ 1)$  & {$2p^2+p-1$}& (B7) where $l=1$; (B12) where $l=1$; (B13) where $p>2$ and $m=l=1$; (B20)\\

\hline $(1,\ p^2+p-2,\ 2)$  & $3p^2+p-2$& (B8) where $l=1$\\

\hline \multirow{2}*{$(1,\ p^2-p,\ 2p)$}  &\multirow{2}*{$2p^3+p^2-p$} & (B8) where $l>1$; (B13) where $m=2$; (B14) where $n=m=2$\\
&&(B15); (B17) where $l>1$; (B18) where $m=2$; (B19) where $m=2$\\
\hline
\end{tabular}}
\caption{The enumeration of (B1)--(B20)\label{table2}}
\end{table}
\end{thm}

\demo  Let $A$ be an abelian subgroup of index $p$ and $B$ be an
$\mathcal{A}_1$-subgroup of index $p$. By Lemma \ref{G'}, $|G'|\le
p|A'||B'|=p^2$. It is easy to see that $\Phi(B)=\Phi(G)$. By Lemma
\ref{minimal non-abelian equivalent conditions}, $\Phi(B)=Z(B)$.
Since $[\Phi(G),A]=[\Phi(G),B]=1$ and $G=AB$, we have $\Phi(G)\le
Z(G)$. Moreover, it is easy to prove that $G'\le C_p^2$. If
$|G'|=p$, then $G$ is one of the groups listed in \cite[Theorem
3.1]{ALQZ}. Since $G$ has an $\mathcal{A}_1$-subgroup of index $p$,
they are the groups (B1)--(B5). If $|G'|=p^2$, then $G$ is one of
the groups listed in \cite[Theorem 4.7]{ALQZ}. By checking
\cite[Table 4]{ALQZ} we get the groups (B6)--(B20). The
correspondence shows as Table \ref{table3}.

\begin{table}[h]
  \centering
{
\tiny
\begin{tabular}{c|c||c|c}
\hline Groups  & Groups in \cite[Theorem 4.8]{ALQZ}& Groups &
Groups in \cite[Theorem
4.8]{ALQZ} \\
\hline
(B6)&(A1) where $l=m=1$ & (B14)&(A10) where $\min\{l,m,n\}=1$ and $\max\{m,n\}=2$\\
(B7)&(A2) where $m=1$ & (B15)&(A11) where $n=2$\\
(B8)&(A3) where $m=1$ and $-\nu\in (F_p^*)^2$ & (B16)&(A12) where $h=1$\\
(B9)&(A3) where $l=1$, $m=2$ and $-\nu\not\in (F_p^*)^2$& (B17)&(B1) where $m=n=2=1$\\
(B10)&(A4) where $m=1$ & (B18)&(B2) where $n=1$, $m\le 2$ and $\min\{l,m\}=1$\\
(B11)&(A6) where $m=2$ & (B19)&(B3) where $n=1$ and $m\le 2$\\
(B12)&(A8) where $m=2$ and $n=1$& (B20)&(B4) where $h=1$\\
(B13)&(A9) where $n=1$ and $m\le 2$& &\\
\hline
\end{tabular}
\caption{The correspondence from Theorem \ref{d=2-2} to
\cite[Theorem 4.7]{ALQZ}}\label{table3}}
\end{table}
We calculate the $(\mu_0,\mu_1,\mu_2)$ and $\alpha_1(G)$ of those
groups in Theorem \ref{d=2-2} as follows.

\medskip
By Lemma \ref{alpha of d=3}, $\alpha_1(G)=\mu_1+p^2\mu_2$. Hence we
only need to calculate $(\mu_0,\mu_1,\mu_2)$.  Since $d(G)=3$, $G$
has $1+p+p^2$ maximal subgroups. They are respectively

$M=\lg b, c, \Phi(G)\rg$;

         $M_i=\lg ab^i, c, \Phi(G)\rg$, where $0\leq i\leq p-1$;

          $M_{ij}=\lg ac^i, bc^j, \Phi(G)\rg$, where $0\leq i,j\leq
          p-1$.

\medskip

{\bf Case 1.} $G'\cong C_p$. That is, $G$ is one of the groups (B1)--(B5).

In this case, $M$ and $M_i$ are abelian. Notice that $|M_{ij}'|=p$.
If $d(M_{ij})=2$, then, by Lemma \ref{minimal non-abelian equivalent
conditions}, $M_{ij}\in\mathcal{A}_1$. Since $\Phi(M_{ij})\le
\Phi(G)$, we get that $M_{ij}\in\mathcal{A}_1$ if and only if
$\Phi(G)=\Phi(M_{ij})$.

By calculations, $\Phi(G)=\lg a^p,b^p,c^p,G'\rg$ and
$\Phi(M_{ij})=\lg a^pc^{ip}, b^pc^{jp},G'\rg$. Since $\Phi(G)/G'=\lg
\bar{a}^p\rg\times \lg \bar{b}^p\rg\times \lg \bar{c}^p\rg$ where
$\bar{a}^p=1$ or $\bar{b}^p=1$, and $\lg \bar{c}^p\rg\cong C_{p}$,
the following conclusions hold:

(i) If $\bar{a}^p=1$ and $\bar{b}^p\neq 1$, then $\Phi(G)=\Phi(M_{ij})$ if and only if $(i,p)=1$;

(ii) If $\bar{a}^p\neq 1$ and $\bar{b}^p=1$, then $\Phi(G)=\Phi(M_{ij})$ if and only if $(j,p)=1$;

(iii) If $\bar{a}^p=\bar{b}^p=1$, then $\Phi(G)=\Phi(M_{ij})$ if and only if $(i,p)=1$ or $(j,p)=1$.

\medskip

If $|G|=p^5$, then case (iii) happens. In this case, $M_{00}$ is the unique $\mathcal{A}_2$-subgroup of $G$.
Hence $(\mu_0,\mu_1,\mu_2)=(p+1,p^2-1,1)$.

If $|G|>p^5$, then case (i) or (ii) happens. If (i) happens, then $M_{ij}\in\mathcal{A}_2$
if and only if $i=0$. If (ii) happens, then $M_{ij}\in \mathcal{A}_2$ if and only if $j=0$.
Hence $(\mu_0,\mu_1,\mu_2)=(p+1,p^2-p,p)$.

\medskip

{\bf Case 2.} $G'\cong C_p^2$. That is, $G$ is one of the groups (B6)--(B20).

In this case, $M$ is the unique abelian subgroup of index $p$.
Notice that $|M_{i}'|=p$. If $d(M_{i})=2$, then, by Lemma
\ref{minimal non-abelian equivalent conditions},
$M_{i}\in\mathcal{A}_1$. Since $\Phi(M_{i})\le \Phi(G)$, we get that
$M_{i}\in\mathcal{A}_1$ if and only if $\Phi(G)=\Phi(M_{ij})$.
Similar argument gives that $M_{ij}\in\mathcal{A}_1$ if and only if
$\Phi(G)=\Phi(M_{ij})$.

\medskip
Subcase 1. $G$ is the group (B6).

In this case, $M_0=\lg a,c\rg\times \lg b^p\rg\in\mathcal{A}_2$. For
$i\ne 0$, $M_i=\lg ab^i,c\rg\in \mathcal{A}_1$. By calculations,
$\Phi(G)=\lg b^p,c^p\rg$ and $\Phi(M_{ij})=\lg c^{ip},
b^pc^{jp}\rg$. Hence $\Phi(G)=\Phi(M_{ij})$ if and only if
$(i,p)=1$. In Summary, $G$ has $p+1$ $\mathcal{A}_2$-subgroups. They
are $M_0$ and $M_{0j}$. Hence $(\mu_0,\mu_1,\mu_2)=(1, p^2-1, p+1)$.

\medskip
Subcase 2. $G$ is the group (B7).

In this case, $M_i=\lg ab^i,c\rg\in \mathcal{A}_1$.
By calculations, $\Phi(G)=\lg a^p,b^p,c^p\rg$ and
$\Phi(M_{ij})=\lg a^pc^{ip}, b^pc^{jp}, b^{-p}b^{-jp}c^{-jp}\rg=\lg a^pc^{ip}, b^pc^{jp}, b^{-jp}\rg$.

If $l=1$, then $a^p=1$. In this case, $\Phi(G)=\Phi(M_{ij})$ if and
only if $(i,p)=1$ or $(j,p)=1$. That is, $G$ has a unique
$\mathcal{A}_2$-subgroup $M_{00}$. Hence $(\mu_0,\mu_1,\mu_2)=(1,
p^2+p-1, 1)$.

If $l>1$, then $a^p\neq 1$. In this case, $\Phi(G)=\Phi(M_{ij})$ if
and only if $(j,p)=1$. That is, $G$ has $p$
$\mathcal{A}_2$-subgroups $M_{i0}$. Hence $(\mu_0,\mu_1,\mu_2)=(1,
p^2, p)$.

\medskip
Subcase 3. $G$ is the group (B8).

In this case, $M_i=\lg ab^i,c\rg\in \mathcal{A}_1$.
By calculations, $\Phi(G)=\lg a^p,b^p,c^p\rg$ and
$\Phi(M_{ij})=\lg a^pc^{ip}, b^pc^{jp}, b^{-tp}c^{\nu p}b^{-jp}c^{-tjp}\rg=\lg a^pc^{ip}, b^pc^{jp}, c^{\nu p}b^{-jp}\rg$.

If $l=1$, then $a^p=1$. In this case, $\Phi(G)=\Phi(M_{ij})$ if and only if $(i,p)=1$ or $\left|
\begin{array}{cc}
       1 &  j \\
       -j & \nu\\
    \end{array}
\right|\neq 0.$ That is, $G$ has $2$ $\mathcal{A}_2$-subgroups $M_{0j}$ where $j^2=-\nu$. Hence
 $(\mu_0,\mu_1,\mu_2)=(1, p^2+p-2, 2)$.

If $l>1$, then $a^p\neq 1$. In this case, $\Phi(G)=\Phi(M_{ij})$ if and only if $\left|
\begin{array}{ccc}
       1 &  j \\
       -j & \nu\\
    \end{array}
\right|\neq 0.$ That is, $G$ has $2p$ $\mathcal{A}_2$-subgroups $M_{ij}$ where $j^2=-\nu$. Hence
 $(\mu_0,\mu_1,\mu_2)=(1, p^2-p, 2p)$.

\medskip
Subcase 4. $G$ is the group (B9).

In this case, $M_0=\lg a,c\rg\ast \lg b^pc^{tp}\rg\in\mathcal{A}_2$.
For $i\neq 0$, $M_i=\lg ab^i,c\rg\in \mathcal{A}_1$. By
calculations, $\Phi(G)=\lg b^p,c^p\rg$ and $\Phi(M_{ij})=\lg c^{ip},
b^pc^{jp}, G'\rg$. Hence $\Phi(G)=\Phi(M_{ij})$ if and only if
$(i,p)=1$. In summary, $G$ has $p$ $\mathcal{A}_2$-subgroups. They
are $M_0$ and $M_{0j}$. Hence $(\mu_0,\mu_1,\mu_2)=(1, p^2-1, p+1)$.

\medskip
Subcase 5. $G$ is the group (B10).

In this case, $M_i=\lg ab^i,c\rg\in \mathcal{A}_1$. By calculations,
$\Phi(G)=\lg a^2,b^2,c^2\rg$ and $\Phi(M_{ij})=\lg a^2c^{2i}b^{2i},
b^2c^{2j}, c^2b^{2j}\rg$. Hence $\Phi(G)=\Phi(M_{ij})$ if and only
if $j=0$. In summary, $G$ has $p$ $\mathcal{A}_2$-subgroups. They
are $M_{i1}$. Hence $(\mu_0,\mu_1,\mu_2)=(1, p^2, p)$.

\medskip
Subcase 6. $G$ is the group (B11).

In this case,  $M_0=\lg a,c\rg\ast \lg b^2\rg\in\mathcal{A}_2$ and $M_1=\lg ab,c\rg\in \mathcal{A}_1$.
By calculations, $\Phi(G)=\lg b^2,c^2\rg$ and $\Phi(M_{ij})=\lg c^{2i}, b^2c^{2j}, G'\rg$.
Hence $\Phi(G)=\Phi(M_{ij})$ if and only if $i=1$. In summary, $G$ has $p+1$ $\mathcal{A}_2$-subgroups.
They are $M_0$ and $M_{0j}$. Hence
 $(\mu_0,\mu_1,\mu_2)=(1, p^2-1, p+1)$.

\medskip
Subcase 7. $G$ is the group (B12).

By calculations, $\Phi(M_{ij})=\lg a^pc^{ip}, b^pc^{jp}, c^{\nu
p}b^{-jp^2}\rg=\lg a^p, b^p, c^{p}\rg=\Phi(G)$. Hence $M_{ij}\in
\mathcal{A}_1$. Notice that $\Phi(M_i)=\lg a^pb^{ip},c^p,
b^{p^2}\rg$.

If $l=1$, then $a^p=1$. In this case, $\Phi(G)=\Phi(M_{i})$ if and
only if $(i,p)=1$. That is, $G$ has a unique
$\mathcal{A}_2$-subgroup $M_{0}$. Hence $(\mu_0,\mu_1,\mu_2)=(1,
p^2+p-1, 1)$.

If $l>1$, then $a^p\neq 1$ and $\Phi(G)>\Phi(M_{i})$. Hence
$M_{i}\in \mathcal{A}_2$ and $(\mu_0,\mu_1,\mu_2)=(1, p^2, p)$.

\medskip
Subcase 8. $G$ is the group (B13).

By calculations, $\Phi(G)=\lg a^p, b^p,c^{p}\rg$. If $p>2$, then
$\Phi(M_i)=\lg a^pb^{ip}, c^p\rg$. Hence $M_{i}\in \mathcal{A}_2$
for $m=2$ and $M_{i}\in \mathcal{A}_1$ for $m=1$. If $p=2$, then
$\Phi(M_i)=\lg a^2b^{2i}a^{i2^l}, c^2\rg$. Hence $M_{i}\in
\mathcal{A}_2$ for $m=2$ and $M_{i}\in \mathcal{A}_1$ for $m=1\le
l$. If $l=m=1$ and $p=2$, then $M_0\in\mathcal{A}_1$ and
$M_1\in\mathcal{A}_2$.

If $p>2$, then $\Phi(M_{ij})=\lg a^pc^{ip},b^pc^{jp},
a^{p^l}c^{-jp}\rg$. We discuss the value of $m$ and $l$.

If $m=2$ and $l>1$, then $\Phi(M_{ij})=\lg a^pc^{ip},b^pc^{jp},
c^{-jp}\rg$. In this case, $\Phi(G)=\Phi(M_{ij})$ if and only if
$(j,p)=1$ and $G$ has $2p$ $\mathcal{A}_2$-subgroups $M_{i}$ and
$M_{i0}$. Hence $(\mu_0,\mu_1,\mu_2)=(1, p^2-p,2p)$.

If $m=2$ and $l=1$, then $\Phi(M_{ij})=\lg a^pc^{ip},b^pc^{jp},
c^{(i+j)p}\rg$. In this case, $\Phi(G)=\Phi(M_{ij})$ if and only if
$(i+j,p)=1$ and $G$ has $2p$ $\mathcal{A}_2$-subgroups $M_{i}$ and
$M_{ij}$ where $i+j=p$. Hence $(\mu_0,\mu_1,\mu_2)=(1, p^2-p,2p)$.

If $m=1$ and $l>1$, then $\Phi(M_{ij})=\lg a^pc^{ip},c^{jp}\rg$. In
this case, $\Phi(G)=\Phi(M_{ij})$ if and only if $(j,p)=1$ and $G$
has $p$ $\mathcal{A}_2$-subgroups $M_{i0}$. Hence
$(\mu_0,\mu_1,\mu_2)=(1, p^2,p)$.

If $l=m=1$,then $\Phi(M_{ij})=\lg a^pc^{ip},c^{jp},
a^pc^{-jp}\rg=\lg a^p,c^{ip},c^{jp}\rg$. In this case,
$\Phi(G)=\Phi(M_{ij})$ if and only if $(i,p)=1$ or $(j,p)=1$, and
$G$ has a unique $\mathcal{A}_2$-subgroup $M_{00}$. Hence
$(\mu_0,\mu_1,\mu_2)=(1, p^2+p-1,1)$.

If $p=2$, then $\Phi(M_{ij})=\lg a^2,b^2c^{2j},
a^{2^l}c^{-2j}\rg=\lg a^2,b^2, c^{-2j}\rg$. In this case,
$\Phi(G)=\Phi(M_{ij})$ if and only if $(j,2)=1$ and $M_{i0}\in
\mathcal{A}_2$. Hence $(\mu_0,\mu_1,\mu_2)=(1, p^2-p,2p)$ for $m=2$,
$(\mu_0,\mu_1,\mu_2)=(1, p^2,p)$ for $m=1<l$, and
$(\mu_0,\mu_1,\mu_2)=(1, p^2-1,p+1)$ for $l=m=1$.

\medskip
Subcase 9. $G$ is the group (B14).

By calculations, $\Phi(G)=\lg a^p, b^p,c^{p}\rg$. If $p>2$,
then $\Phi(M_i)=\lg a^pb^{ip}, c^p, b^{p^m}\rg$.
Hence $M_{i}\in \mathcal{A}_2$ for $m=2$ and $M_{i}\in \mathcal{A}_1$ for $m=1$.
If $p=2$, then $\Phi(M_i)=\lg a^2b^{2i}a^{i2^l}, c^2, b^{2^m}\rg$.
Hence $M_{i}\in \mathcal{A}_2$ for $m=2$ and $M_{i}\in \mathcal{A}_1$
for $m=1\le l$. If $l=m=1$ and $p=2$, then $M_0\in\mathcal{A}_1$ and $M_1\in\mathcal{A}_2$.

If $p>2$, then $\Phi(M_{ij})=\lg a^pc^{ip},b^pc^{jp},
a^{p^l}b^{-jp^m}\rg$. We discuss the value of $m$ and $l$.

If $n=1$, then $\Phi(M_{ij})=\lg a^p,b^p\rg=\Phi(G)$ and hence
$\Phi(M_{ij})\in \mathcal{A}_1$. Since $\max\{ m,n\}=2$, $m=2$.
Hence $(\mu_0,\mu_1,\mu_2)=(1, p^2,p)$.

If $m=n=2$, then $l=1$ since $\min\{ l,m,n\}=1$. Notice that
$\Phi(M_{ij})=\lg c^{ip},b^pc^{jp}, a^{p}\rg$. In this case,
$\Phi(G)=\Phi(M_{ij})$ if and only if $(i,p)=1$ and $G$ has $2p$
$\mathcal{A}_2$-subgroups $M_{i}$ and $M_{0j}$. Hence
$(\mu_0,\mu_1,\mu_2)=(1, p^2-p,2p)$.

If $n=2$ and $l=m=1$, then $\Phi(M_{ij})=\lg a^pc^{ip},b^pc^{jp},
a^pb^{-jp}\rg$. In this case, $\Phi(G)=\Phi(M_{ij})$ if and only if
$j^2-i\neq 0$ and $G$ has $p$ $\mathcal{A}_2$-subgroups $M_{ij}$
where $i=j^2$. Hence $(\mu_0,\mu_1,\mu_2)=(1, p^2,p)$.

If $n=2$ and $m=1<l$, then $\Phi(M_{ij})=\lg a^pc^{ip},b^pc^{jp},
b^{-jp}\rg$. In this case, $\Phi(G)=\Phi(M_{ij})$ if and only if
$(j,p)=1$ and $G$ has $p$ $\mathcal{A}_2$-subgroups $M_{i0}$. Hence
$(\mu_0,\mu_1,\mu_2)=(1, p^2,p)$.

If $p=2$, then $\Phi(M_{ij})=\lg a^2c^{2i}b^{i2^m},b^2c^{2j},
a^{2^l}b^{-j2^m}\rg$. We discuss the value of $m$ and $l$.

If $n=1$, then $\Phi(M_{ij})=\lg a^2,b^2\rg=\Phi(G)$ and hence
$\Phi(M_{ij})\in \mathcal{A}_1$. Since $\max\{ m,n\}=2$, $m=2$.
Hence $(\mu_0,\mu_1,\mu_2)=(1, p^2,p)$.

If $m=n=2$, then $l=1$ since $\min\{ l,m,n\}=1$. Notice that
$\Phi(M_{ij})=\lg c^{2i},b^2c^{2j}, a^{2}\rg$. In this case,
$\Phi(G)=\Phi(M_{ij})$ if and only if $(i,2)=1$ and $G$ has $2p$
$\mathcal{A}_2$-subgroups $M_{i}$ and $M_{0j}$. Hence
$(\mu_0,\mu_1,\mu_2)=(1, p^2-p,2p)$.

If $n=2$ and $l=m=1$, then $\Phi(M_{ij})=\lg
a^2c^{2i}b^{2i},b^2c^{2j}, a^2b^{2j}\rg$. In this case,
$\Phi(G)=\Phi(M_{ij})$ if and only if $i=1$ or $j=1$, and $G$ has
$2$ $\mathcal{A}_2$-subgroups $M_1$ and $M_{00}$. Hence
$(\mu_0,\mu_1,\mu_2)=(1, p^2,p)$.

If $n=2$ and $m=1<l$, then $\Phi(M_{ij})=\lg
a^2c^{2i}b^{2i},b^2c^{2j}, b^{-2j}\rg$. In this case,
$\Phi(G)=\Phi(M_{ij})$ if and only if $(j,2)=1$ and $G$ has $2$
$\mathcal{A}_2$-subgroups $M_{i0}$. Hence $(\mu_0,\mu_1,\mu_2)=(1,
p^2,p)$.

\medskip
Subcase 10. $G$ is the group (B15).

By calculations, $\Phi(G)=\lg a^2, b^2,c^{2}\rg$, $M_0=\lg
c,a\rg\times\lg b^2\rg\in\mathcal{A}_2$ and $M_1=\lg
c,ab\rg\in\mathcal{A}_1$. Since $\Phi(M_{ij})=\lg
a^{2(1+i)},b^2c^{2j}, c^{2(1+j)}a^{2j}\rg$, $\Phi(G)=\Phi(M_{ij})$
if and only if $i=j=0$. Hence $G$ has $4$ $\mathcal{A}_2$-subgroups
$M_{0}$ and $M_{ij}$ where $i\ne 0$ or $j\neq 0$, and
$(\mu_0,\mu_1,\mu_2)=(1, 2,4)$.

\medskip
Subcase 11. $G$ is the group (B16).

By calculations, $\Phi(G)=\lg a^2, b^2\rg$, $M_i=\lg c,ab^i\rg\times\lg b^2\rg\in\mathcal{A}_2$.
Since $\Phi(M_{ij})=\lg a^{2}=c^2,b^2c^{2j}, b^{2}a^{2j}\rg=\Phi(G)$, $M_{ij}\in\mathcal{A}_1$.
Hence $G$ has $p$ $\mathcal{A}_2$-subgroups $M_{i}$ and $(\mu_0,\mu_1,\mu_2)=(1,p^2,p)$.

\medskip
Subcase 12. $G$ is the group (B17).

By calculations, $\Phi(G)=\lg a^p,c^{p},x\rg$ and $M_{i}=\lg ab^i,c\rg\times \lg x\rg\in \mathcal{A}_2$.

If $p>2$, then $\Phi(M_{ij})=\lg a^pc^{ip},c^{jp}, xc^{-jp}\rg=\lg
a^pc^{ip},c^{jp},x\rg$. In this case, $\Phi(G)=\Phi(M_{ij})$ if and
only if $(j,p)=1$ and $G$ has $2p$ $\mathcal{A}_2$-subgroups $M_{i}$
and $M_{i0}$. Hence $(\mu_0,\mu_1,\mu_2)=(1, p^2-p,2p)$.

If $p=2$, then $\Phi(M_{ij})=\lg a^2,c^{2j}, xc^{2j}\rg=\lg
a^2,c^{2j},x\rg$. In this case, $\Phi(G)=\Phi(M_{ij})$ if and only
if $(j,2)=1$ and $G$ has $4$ $\mathcal{A}_2$-subgroups $M_{i}$ and
$M_{i0}$. Hence $(\mu_0,\mu_1,\mu_2)=(1, p^2-p,2p)$.

\medskip
Subcase 13. $G$ is the group (B18).

If $l=1$ and $m=2$, then $p>2$ and $\Phi(G)=\lg b^{p},c^p,x\rg$. In
addition, $\Phi(M_i)=\lg b^{ip}, c^p, x\rg$ and $\Phi(M_{ij})=\lg
c^{ip},b^pc^{jp}, c^px^{-j}\rg$. In this case, $\Phi(M_i)=\Phi(G)$
if and only if $(i,p)=1$, and $\Phi(M_{ij})=\Phi(G)$ if and only if
$(i,p)=1$ and $(j,p)=1$. Hence $G$ has $2p$
$\mathcal{A}_2$-subgroups $M_{0}$, $M_{0j}$ and $M_{i0}$, and
$(\mu_0,\mu_1,\mu_2)=(1, p^2-p,2p)$.

If $l=m=1$, then $p>2$ and $\Phi(G)=\lg c^p,x\rg$. Since
$\Phi(M_i)=\lg c^p, x\rg=\Phi(G)$, $M_i\in\mathcal{A}_1$. Since
$\Phi(M_{ij})=\lg c^{ip},c^{jp}, c^px^{-j}\rg$,
$\Phi(M_{ij})=\Phi(G)$ if and only if $(j,p)=1$. Hence $G$ has $p$
$\mathcal{A}_2$-subgroups $M_{i0}$, and $(\mu_0,\mu_1,\mu_2)=(1,
p^2,p)$.

If $l>1$, then $m=1$ since $\min\{ l,m\}=1$. In this case,
$\Phi(G)=\lg a^{p},c^p,x\rg$. In addition, if $p>2$, then
$\Phi(M_i)=\lg a^{p}, c^p, x\rg$ and $\Phi(M_{ij})=\lg
a^pc^{ip},c^{jp}, c^px^{-j}\rg$. In this case, $M_i\in\mathcal{A}_1$
and $\Phi(M_{ij})=\Phi(G)$ if and only if $(j,p)=1$. Hence $G$ has
$p$ $\mathcal{A}_2$-subgroups $M_{i0}$ and $(\mu_0,\mu_1,\mu_2)=(1,
p^2,p)$; if $p=2$, then $\Phi(M_i)=\lg a^2c^{2i}, c^2, x\rg$ and
$\Phi(M_{ij})=\lg a^2c^{2i}x^i,c^{2j}, c^2x^{j}\rg$. In this case,
$M_i\in\mathcal{A}_1$ and $\Phi(M_{ij})=\Phi(G)$ if and only if
$(j,2)=1$. Hence $G$ has $2$ $\mathcal{A}_2$-subgroups $M_{i0}$ and
$(\mu_0,\mu_1,\mu_2)=(1, p^2,p)$.

\medskip
Subcase 14. $G$ is the group (B19).

If $m=2$, then $\Phi(G)=\lg a^p, b^{p},x\rg$. In addition, if $p>2$,
then $\Phi(M_i)=\lg a^pb^{ip}, x\rg$ and $\Phi(M_{ij})=\lg a^p,b^p,
a^{p^l}x^{-j}\rg$. In this case, $\Phi(M_i)\in\mathcal{A}_2$ and
$\Phi(M_{ij})=\Phi(G)$ if and only $(j,p)=1$. Hence $G$ has $2p$
$\mathcal{A}_2$-subgroups $M_{i}$ and $M_{i0}$, and
$(\mu_0,\mu_1,\mu_2)=(1, p^2-p,2p)$; if $p=2$, then $\Phi(M_i)=\lg
a^2b^{2i}a^{i2^l}, x\rg$ and $\Phi(M_{ij})=\lg a^2x^i,b^2,
a^{2^l}x^{-j}\rg$. In this case, $\Phi(M_i)\in\mathcal{A}_2$.
$\Phi(M_{ij})=\Phi(G)$ if and only $(j,2)=1$ for $l>1$. Hence $G$
has $2p$ $\mathcal{A}_2$-subgroups $M_{i}$ and $M_{i0}$ for $l>1$.
$\Phi(M_{ij})=\Phi(G)$ if and only $i\neq j$ for $l=1$. Hence $G$
has $2p$ $\mathcal{A}_2$-subgroups $M_{i}$ and $M_{ij}$ where $i\neq
j$ for $l>1$. In summary, $(\mu_0,\mu_1,\mu_2)=(1, p^2-p,2p)$.

If $l=m=1$, then $\Phi(G)=\lg a^p,x\rg$. In addition, if $p>2$, then
$\Phi(M_i)=\lg a^p, x\rg$ and $\Phi(M_{ij})=\lg a^p,
a^{p}x^{-j}\rg$. In this case, $\Phi(M_i)\in\mathcal{A}_1$ and
$\Phi(M_{ij})=\Phi(G)$ if and only $(j,p)=1$. Hence $G$ has $p$
$\mathcal{A}_2$-subgroups $M_{i0}$, and $(\mu_0,\mu_1,\mu_2)=(1,
p^2,p)$; if $p=2$, then $\Phi(M_i)=\lg a^2a^{i2}, x\rg$ and
$\Phi(M_{ij})=\lg a^2x^i, a^{2}x^{-j}\rg$. In this case,
$\Phi(M_0)\in\mathcal{A}_1$, $\Phi(M_1)\in\mathcal{A}_2$ and
$\Phi(M_{ij})=\Phi(G)$ if and only $i\neq j$. Hence $G$ has $p+1$
$\mathcal{A}_2$-subgroups $M_{1}$ and $M_{ij}$ where $i\neq j$, and
$(\mu_0,\mu_1,\mu_2)=(1, p^2-1,p+1)$.

If $l>1=m$, then $\Phi(G)=\lg a^p,x\rg$. In addition, if $p>2$, then
$\Phi(M_i)=\lg a^p, x\rg$ and $\Phi(M_{ij})=\lg a^p,
a^{p^l}x^{-j}\rg$. In this case, $\Phi(M_i)\in\mathcal{A}_1$ and
$\Phi(M_{ij})=\Phi(G)$ if and only $(j,p)=1$. Hence $G$ has $p$
$\mathcal{A}_2$-subgroups $M_{i0}$, and $(\mu_0,\mu_1,\mu_2)=(1,
p^2,p)$; if $p=2$, then $\Phi(M_i)=\lg a^2a^{i2^l}, x\rg$ and
$\Phi(M_{ij})=\lg a^2x^i, a^{2^l}x^{-j}\rg$. In this case,
$\Phi(M_i)\in\mathcal{A}_1$ and $\Phi(M_{ij})=\Phi(G)$ if and only
$(j,2)=1$. Hence $G$ has $p$ $\mathcal{A}_2$-subgroups $M_{i0}$, and
$(\mu_0,\mu_1,\mu_2)=(1, p^2,p)$.

\medskip
Subcase 15. $G$ is the group (B20).

By calculations, $\Phi(G)=\lg a^2=c^{2},x\rg$, $M_{0}=\lg
a,c\rg\times\lg x\rg\in \mathcal{A}_2$, and $M_{1}=\lg ab,c\rg\in
\mathcal{A}_1$. Since $\Phi(M_{ij})=\lg a^2,c^{2j}, xc^{2j}\rg=\lg
a^2,x\rg=\Phi(G)$, $M_{ij}\in\mathcal{A}_1$. Hence
$(\mu_0,\mu_1,\mu_2)=(1, p^2+p-1,1)$.

To sum up, we get the Table \ref{table3}. \qed

\begin{lem}\label{c7-c15}
Suppose that $G$ is a finite $p$-group such that $\Phi(G')=1$,
$C_p^2\cong G_3\le Z(G)$ and $G/G_3\cong M_p(1,1,1)$, where $p\ge
3$. Let $ G=\langle a, b, c;x,y\rangle$ such that
$$a^{p}=x^{w_{11}}y^{w_{12}},b^{p}=x^{w_{21}}y^{w_{22}},\,\,
c^{p}=x^p=y^p=1, [a,b]=c, [c,a]=x,[c,b]=y,$$ $N=\lg c,b,x,y\rg$ and
$M_i=\lg ab^{i},c,x,y\rg$. Then $N$ and $M_i$ are all maximal
subgroups of $G$, and the following conclusions hold:

$(1)$ $N\in\mathcal{A}_2$ if and only if $w_{21}=0$;

$(2)$ If $p>3$, then $M_i\in\mathcal{A}_2$ if and only if $w_{12}+i(w_{22}-w_{11})-i^2w_{21}=0$;

$(3)$ If $p=3$, then $M_i\in\mathcal{A}_2$ if and only if $w_{12}+i(w_{22}-w_{11})-i^2w_{21}-i^2=0$;

$(4)$ $\alpha_1(G)=\mu_1+\mu_2 p^2$.
\end{lem}
\demo If $w_{21}=0$, then $N=\lg c,b\rg\times\lg x\rg\in\mathcal{A}_2$. If $w_{21}\neq 0$, then
$N=\lg c,b\rg\in\mathcal{A}_1$. Hence (1) holds.

It is easy to see that $M_i\in\mathcal{A}_2$ if and only if
$\lg(ab^i)^p\rg=\lg [c,ab^i]\rg=\lg xy^i\rg$.

If $p>3$, then $(ab^i)^p=a^pb^{ip}=x^{w_{11}+iw_{21}}y^{w_{12}+iw_{22}}$. Hence
$M_i\in\mathcal{A}_2$ if and only if ${w_{12}+iw_{22}}=i(w_{11}+iw_{21})$. The later is $w_{12}+i(w_{22}-w_{11})-i^2w_{21}=0$. Hence (2) holds.

If $p=3$, then $(ab^i)^3=a^3b^{3i}[a,b^{-i},a][a,b^{-i},b^{-i}]=x^{w_{11}+iw_{21}-i}y^{w_{12}+iw_{22}+i^2}$. Hence
$M_i\in\mathcal{A}_2$ if and only if ${w_{12}+iw_{22}+i^2}=i(w_{11}+iw_{21}-i)$. The later is $w_{12}+i(w_{22}-w_{11})-i^2w_{21}-i^2=0$. Hence (3) holds.

Let $M\in \Gamma_1$. Then, by above argument, $|M'|=p$. If $M\in\mathcal{A}_2$, then, by Lemma \ref{A_2-property} (7), $\alpha_1(M)=p^2$.
By Hall's enumeration principle, $
\alpha_1(G)=\sum_{H\in
\Gamma_1} \alpha_1(H)=\mu_1+p^2\mu_2.
$ Hence (4) holds.
\qed

\begin{thm}\label{d=2-3}
$G$ has no abelian subgroup of index $p$, $G$ has  at least two
distinct $\mathcal{A}_1$-subgroups of index $p$ and $d(G)=2$ if and
only if $G$ is isomorphic to one of the following pairwise
non-isomorphic groups:
\begin{enumerate}
\rr{Ci} $c(G)=3$ and $G'\cong C_2^2$. In this case, $(\mu_0,\mu_1,\mu_2)=(0,2,1)$ and
 $\alpha_1(G)=6$.

\begin{enumerate}
\rr{C1}  $\langle a, b; c,d \di a^{4}=b^{2}=c^2=d^2=1,
[a,b]=c,[c,a]=d,[c,b]=[d,a]=[d,b]=1 \rangle$; where  $|G|=2^{5}$,
$\Phi(G)=\lg a^2,c,d\rg\cong C_2^3$, $G'=\lg c,d\rg$ and $Z(G)=\lg d\rg\cong C_2$.

\rr{C2} $\langle a, b; c \di a^{8}=b^{2}=c^2=1,
[a,b]=c,[c,a]=a^4,[c,b]=1 \rangle$; where $|G|=2^{5}$, $\Phi(G)=\lg
a^2,
c\rg\cong C_4\times C_2$, $G'=\lg c,a^4\rg$ and $Z(G)=\lg a^4\rg\cong C_2$.

\rr{C3}  $\langle a, b; c \di a^{8}=c^{2}=1,b^2=a^4,
[a,b]=c,[c,a]=[a^2,b]=b^2,[c,b]=1 \rangle$; where $|G|=2^{5}$,
$\Phi(G)=\lg a^2,
c\rg\cong C_4\times C_2$, $G'=\lg c,a^4\rg$ and $Z(G)=\lg a^4\rg\cong C_2$.

\rr{C4}  $\langle a, b; c \di a^{2^{n+1}}=b^{2}=c^2=1,
[a,b]=c,[c,a]=a^{2^n},[c,b]=1 \rangle$, where $n\ge 3$; moreover,
$|G|=2^{n+3}$, $\Phi(G)=\lg a^2,
c\rg\cong C_{2^n}\times C_2$, $G'=\lg c,a^{2^n}\rg$ and $Z(G)=\lg a^4\rg\cong C_{2^{n-1}}$.

\rr{C5}  $\langle a, b; c \di a^{2^{n}}=b^{4}=c^2=1,
[a,b]=c,[c,a]=b^{2},[c,b]=1 \rangle$, where $n\ge 3$; moreover,
$|G|=2^{n+3}$, $\Phi(G)=\lg a^2,b^2,
c\rg\cong C_{2^{n-1}}\times C_2\times C_2$, $G'=\lg c,b^2\rg$ and $Z(G)=\lg a^4,b^2\rg \cong C_{2^{n-2}}\times C_2$.

\rr{C6}  $\langle a, b; c,d \di a^{2^{n}}=b^{2}=c^2=d^2=1,
[a,b]=c,[c,a]=d,[c,b]=1,[d,a]=[d,b]=1 \rangle$, where $n\ge 3$;
moreover, $|G|=2^{n+3}$, $\Phi(G)=\lg a^2,
c,d\rg\cong C_{2^{n-1}}\times C_2\times C_2$, $G'=\lg c,d\rg$ and  $Z(G)=\lg a^4,d\rg\cong C_{2^{n-2}}\times C_2$.
\end{enumerate}

\rr{Cii} $\Phi(G')\le G_3\cong C_p^2$.

\begin{enumerate}
\rr{C7}  $\langle a, b; c,d,e \di a^{3}=b^{3}=c^3=d^3=e^3=1,
[a,b]=c,[c,a]=d,[c,b]=e,[d,a]=[d,b]=[e,a]=[e,b]=1 \rangle$, where
$|G|=3^{5}$, $\Phi(G)=G'=\lg
c,d,e\rg$ and $Z(G)=\lg d,e\rg\cong C_3^2$.

\rr{C8}  $\langle a, b; c,d\di a^{9}=c^{3}=d^3=1, b^3=a^3,
[a,b]=c,[c,a]=d,[c,b]=a^{3},[d,a]=[d,b]=1 \rangle$; where
$|G|=3^{5}$, $\Phi(G)=G'=\lg a^3,
c,d\rg$, $Z(G)=\lg a^3,d\rg\cong C_3^2$.

\rr{C9}  $\langle a, b; c,d \di a^{9}=b^{3}=c^3=d^3=1,
[a,b]=c,[c,a]=d,[c,b]=a^{-3},[d,a]=[d,b]=1 \rangle$; where
$|G|=3^{5}$, $\Phi(G)=G'=\lg a^3,
c,d\rg$ and $Z(G)=\lg a^3,d\rg\cong C_3^2$.

\rr{C10}  $\langle a, b; c \di a^{9}=b^9=c^3=1,
[a,b]=c,[c,a]=a^3,[c,b]=b^{3} \rangle$; where $|G|=3^{5}$,
$\Phi(G)=G'=\lg a^3,b^3,
c\rg$ and $Z(G)=\lg a^3,b^3\rg\cong C_3^2$.

\rr{C11} $\langle a, b, c \di a^{p^{2}}=b^{p^{2}}=c^{p}=1,
[a,b]=c,[c,a]=a^{p}b^{\nu p},[c,b]=b^{p}\rangle$,
where $p>3$, $\nu=1$ or a
fixed quadratic non-residue modular $p$; moreover,
$|G|=p^{5}$, $\Phi(G)=G'=\lg a^p,
b^p,c\rg$ and $Z(G)=\lg a^p,b^p\rg \cong C_p^2$.

\rr{C12} $\langle a, b, c \di a^{p^{2}}=b^{p^{2}}=c^{p}=1,
[a,b]=c,[c,a]=b^{\nu p},[c,b]=a^{-p},[a^p,b]=1\rangle$, where
$p>3$, $\nu=1$ or a
fixed quadratic non-residue modular $p$ such that $-\nu \in F_p^2$; moreover,
$|G|=p^{5}$, $\Phi(G)=G'=\lg a^p,
b^p,c\rg$ and $Z(G)=\lg a^p,b^p\rg \cong C_p^2$.

\rr{C13} $\langle a, b, c \di a^{p^{2}}=b^{p^{2}}=c^{p}=1,
[a,b]=c,[c,a]^{1+r}=a^{p}b^{p},[c,b]^{1+r}=a^{-r p}b^{
p},[a^p,b]=1\rangle$, where $p>3$, $r\neq 0,-1$ and $-r\in (F_p^*)^2$;
moreover,
$|G|=p^{5}$, $\Phi(G)=G'=\lg a^p,
b^p,c\rg$ and $Z(G)=\lg a^p,b^p\rg \cong C_p^2$.

\rr{C14} $\langle a, b; c,d \di a^{p^{2}}=b^{p}=c^{p}=d^p=1,
[a,b]=c,[c,a]=a^{p},[c,b]=d\rangle$, where $p>3$; moreover,
$|G|=p^{5}$, $\Phi(G)=G'=\lg a^p,
c,d\rg$ and $Z(G)=\lg a^p,d\rg \cong C_p^2$.

\rr{C15} $\langle a, b; c,d \di a^{p}=b^{p^2}=c^{p}=d^p=1,
[a,b]=c,[c,a]=b^{\nu p},[c,b]=d,[d,a]=[d,b]=1\rangle$, where $p>3$,
$\nu=1$ or a fixed quadratic non-residue modulo $p$; moreover,
$|G|=p^{5}$, $\Phi(G)=G'=\lg b^p,
c,d\rg$ and $Z(G)=\lg b^p,d\rg\cong C_p^2$.

\rr{C16} $\langle a, b; c \di a^8=c^4=1,b^2=a^4,
[a,b]=c,[c,a]=a^4,[c,b]=c^2\rangle$; where $|G|=2^{6}$, $\Phi(G)=\lg
a^2,
c\rg\cong C_4^2$, $G'=\lg c,a^4\rg$ and $Z(G)=\lg a^4,c^2\rg\cong C_2^2$.

\rr{C17} $\langle a, b; c \di a^{p^{3}}=b^{p^{2}}=c^{p}=1,
[a,b]=c,[c,a]=b^{\nu_1 p},[c,b]=a^{-\nu_2 p^2}\rangle$, where $p\ge
3$, and $\nu_1,\nu_2=1$ or a fixed quadratic non-residue modulo $p$;
moreover, $|G|=p^{6}$, $\Phi(G)=\lg a^p,b^p,
c\rg \cong C_{p^2}\times C_p\times C_p$, $G'=\lg c,a^{p^2},b^p\rg$ and $Z(G)=\lg a^p,b^p\rg\cong C_{p^2}\times C_p$.

\end{enumerate}

Moreover, Table $\ref{table4}$ gives $(\mu_0,\mu_1,\mu_2)$ and
$\alpha_1(G)$ for the groups {\rm (C7)--(C17)}.
\end{enumerate}

\begin{table}[h]
\centering
{
\begin{tabular}[t]{c|c||c}
\hline
 $(\mu_0,\ \mu_1,\ \mu_2)$&$\alpha_1(G)$      & types of  $\mathcal{A}$$_3$ groups \\

\hline $(0,\ p-1,\ 2)$  & {$2p^2+p-1$}& (C7); (C10); (C12)--(C14)\\

\hline $(0,\ p,\ 1)$  & $p^2+p$& (C8)--(C9); (C11); (C15);(C17) \\

\hline $(0,\ p,\ 1)$  & $p^2+2p$& (C16) \\

\hline

\end{tabular}}
\caption{The enumeration of (C7)--(C17)\label{table4}}
\end{table}
\end{thm}
\demo Let $A$ and $B$ be two distinct $\mathcal{A}_1$-subgroups of
index $p$. Then $|A'|=|B'|=p$ and hence $A'B'\le Z(G)$. Let
$\bar{G}=G/A'B'$. Then $\bar{G}$ has two abelian subgroups $\bar{A}$
and $\bar{B}$  of index $p$. By Lemma \ref{number of abelian maximal
subgroups}, $\bar{G}$ has $1+p$ abelian subgroups. Since
$d(\bar{G})=2$, $\bar{G}$ is an $\mathcal{A}_1$-group. It is easy to
see that $A'B'=\Phi(G')G_3$. If $\bar{G}$ is metacyclic, then $G$ is
also metacyclic by Lemma \ref{phi(G')G_3} and hence $A'=B'$. It
follows that $|G'|=p^2$. by Lemma \ref{metacyclic An}, $G$ is an
$\mathcal{A}_2$-group. This contradicts that $G$ is an
$\mathcal{A}_3$-group. Thus $\bar{G}$ is a non-metacyclic
$\mathcal{A}_1$-group.

If $G'\cong C_{p^2}$, then $G_3\leq \Phi(G')\cong C_p$. Thus $G$ is
one of the groups listed in \cite[Theorem 3.5]{AHZ}. By
\cite[Theorem 3.6]{AHZ} we get no group.

If $G'\cong C_p^2$, then $G$ is one of the groups listed in
\cite[Theorem 4.6]{AHZ}. By \cite[Theorem 4.7]{AHZ} we get groups
(C1)--(C6).

If $G_3\cong C_p$ and $\Phi(G')G_3\cong C_p^2$, then, by
\cite[Theorem 5.1, 5.2, 5.5 \& 5.6]{AHZ}, we get no group.

If $\Phi(G')\le G_3\cong C_{p}^2$, then $G$ is one of the groups
listed in \cite[Theorem 6.5]{AHZ}. By \cite[Theorem 6.1]{AHZ}, we
get the groups (C7)--(C17).

Table \ref{table5} gives the correspondence from Theorem \ref{d=2-3}
to \cite[Theorem 4.6 \& 6.5]{AHZ}.

\begin{table}[h]
  \centering
{
\small
\begin{tabular}{c|c||c|c}
\hline Groups &Groups in \cite[Theorem 4.6]{AHZ} & Groups &Groups in
\cite[Theorem
6.5]{AHZ} \\
\hline
(C1)&(H1) &          (C7)&(O1) \\
(C2)&(H2) &          (C8)&(O3)  \\
(C3)&(H3) &          (C9)&(O4)  \\
(C4)&(I1) &           (C10)&(O5)\\
(C5)&(I2) &           (C11)&(P2) where $n=1$\\
(C6)&(I3) &           (C12)&(P3) where $n=1$ and $-\nu\in F_p^2$\\
&         &            (C13)&(P4) where $n=1$ and $-r\in (F_p^*)^2$\\
&         &           (C14)&(P8) where $n=1$\\
 &         &          (C15)&(P9) where $n=1$\\
 &         &          (C16)&(R4)\\
  &         &            (C17) &(S2) where $n=2$\\

\hline
\end{tabular}
\caption{The correspondence from Theorem \ref{d=2-3} to
\cite[Theorem 4.6 \& 6.5]{AHZ}} \label{table5}}
\end{table}

\medskip
We calculate the $(\mu_0,\mu_1,\mu_2)$ and $\alpha_1(G)$ of those
groups in Theorem \ref{d=2-3} as follows.

Since $d(G)=2$, $G$ has $1+p$ maximal subgroups. For any
$\mathcal{A}_3$-group, we always have $\mu_2\ge 1$. By the
hypothesis of Theorem \ref{d=2-3}, $\mu_0=0$ and $\mu_1\ge 2$.

\medskip

{\bf Case 1.} $c(G)=3$ and $G'\cong C_2^2$. That is , $G$ is one of the groups (C1)--(C6).

Obviously, $(\mu_0,\mu_1,\mu_2)=(0,2,1)$. Let $M\in\Gamma_1$. If $M$
is the unique $\mathcal{A}_2$-subgroup, then $|M'|=p$. By Lemma
\ref{A_2-property}, $\alpha_1(M)=p^2=4$. By Hall's enumeration
principle, $ \alpha_1(G)=\sum_{H\in \Gamma_1}
\alpha_1(H)=\mu_1+4\mu_2=6.$

\medskip

{\bf Case 2.} $G$ is one of the groups (C7)--(C15).

In this case, $p\ge 3$, $\Phi(G')=1$, $G_3\cong C_p^2$, $G_3\le Z(G)$ and $G/G_3\cong M_p(1,1,1)$. $
G=\langle a, b, c;x,y\rangle$ such that $$a^{p}=x^{w_{11}}y^{w_{12}},b^{p}=x^{w_{21}}y^{w_{22}},\,\, c^{p}=x^p=y^p=1,
[a,b]=c, [c,a]=x,[c,b]=y.$$
Let $N=\lg c,b,x,y\rg$ and $M_i=\lg ab^{i},c,x,y\rg$. Then $N$ and $M_i$ are all maximal subgroups of $G$.
Let $w(G)=\left(
 \begin{array}{cc}
 w_{11}& w_{12}\\
 w_{21}& w_{22}
 \end{array}
 \right)$.

If $G$ is the group (C7), then $p=3$ and $w(G)=\left(
 \begin{array}{cc}
 0& 0\\
0&0
 \end{array}
 \right)$. By Lemma \ref{c7-c15}, $N\in \mathcal{A}_2$ and $M_i\in \mathcal{A}_2$ if and only if
 $i=0$. Hence $(\mu_0,\mu_1,\mu_2)=(0,p-1,2)$ and $\alpha_1(G)=2p^2+p-1$.

If $G$ is the group (C8), then $p=3$ and $w(G)=\left(
 \begin{array}{cc}
 0& 1\\
0&1
 \end{array}
 \right)$. By Lemma \ref{c7-c15}, $N\in \mathcal{A}_2$ and $M_i\in \mathcal{A}_1$. Hence $(\mu_0,\mu_1,\mu_2)=(0,p,1)$ and $\alpha_1(G)=p^2+p$.

If $G$ is the group (C9), then $p=3$ and $w(G)=\left(
 \begin{array}{cc}
 0& -1\\
0&0
 \end{array}
 \right)$. By Lemma \ref{c7-c15}, $N\in \mathcal{A}_2$ and $M_i\in \mathcal{A}_1$.
 Hence $(\mu_0,\mu_1,\mu_2)=(0,p,1)$ and $\alpha_1(G)=p^2+p$.

If $G$ is the group (C10), then $p=3$ and $w(G)=\left(
 \begin{array}{cc}
 1& 0\\
0&1
 \end{array}
 \right)$. By Lemma \ref{c7-c15}, $N\in \mathcal{A}_2$ and $M_i\in \mathcal{A}_2$ if and only if $i=0$.
 Hence $(\mu_0,\mu_1,\mu_2)=(0,p-1,2)$ and $\alpha_1(G)=2p^2+p-1$.

If $G$ is the group (C11), then $p>3$ and $w(G)=\left(
 \begin{array}{cc}
 1& -\nu\\
0&1
 \end{array}
 \right)$. By Lemma \ref{c7-c15}, $N\in \mathcal{A}_2$ and $M_i\in \mathcal{A}_1$. Hence $(\mu_0,\mu_1,\mu_2)=(0,p,1)$ and $\alpha_1(G)=p^2+p$.

If $G$ is the group (C12), then $p>3$ and $w(G)=\left(
 \begin{array}{cc}
 0& -1\\
\nu^{-1}&0
 \end{array}
 \right)$. By Lemma \ref{c7-c15}, $N\in \mathcal{A}_1$ and $M_i\in \mathcal{A}_2$ if and only if $i^2=-\nu$.
 Hence $(\mu_0,\mu_1,\mu_2)=(0,p-1,2)$ and $\alpha_1(G)=2p^2+p-1$.

If $G$ is the group (C13), then $p>3$ and $w(G)=\left(
 \begin{array}{cc}
 1& -1\\
r&1
 \end{array}
 \right)$. By Lemma \ref{c7-c15}, $N\in \mathcal{A}_1$ and $M_i\in \mathcal{A}_2$ if and only if
 $i^2=-r^{-1}$. Hence $(\mu_0,\mu_1,\mu_2)=(0,p-1,2)$ and $\alpha_1(G)=2p^2+p-1$.

If $G$ is the group (C14), then $p>3$ and $w(G)=\left(
 \begin{array}{cc}
 1& 0\\
0&0
 \end{array}
 \right)$. By Lemma \ref{c7-c15}, $N\in \mathcal{A}_2$ and $M_i\in \mathcal{A}_2$ if and only if
 $i=0$. Hence $(\mu_0,\mu_1,\mu_2)=(0,p-1,2)$ and $\alpha_1(G)=2p^2+p-1$.

If $G$ is the group (C15), then $p>3$ and $w(G)=\left(
 \begin{array}{cc}
 0& 0\\
\nu^{-1}&0
 \end{array}
 \right)$. By Lemma \ref{c7-c15}, $N\in \mathcal{A}_1$ and $M_i\in \mathcal{A}_2$ if and only
 if $i=0$. Hence $(\mu_0,\mu_1,\mu_2)=(0,p,1)$ and $\alpha_1(G)=p^2+p$.

\medskip

{\bf Case 3.} $G$ is the group (C16).

It is obvious that $(\mu_0,\mu_1,\mu_2)=(0,2,1)=(0,p,1)$. The unique
$\mathcal{A}_2$-subgroup of index $p$ is $M=\lg c,b,a^2\rg$, where
$d(M)=3$ and $|M'|=4$. Notice that $\alpha_1(M)=p^2+p=6$. By Hall's
enumeration principle, $$ \alpha_1(G)=\sum_{H\in \Gamma_1}
\alpha_1(H)=\mu_1+(p^2+p)\mu_2=p^2+2p=8.$$

\medskip

{\bf Case 4.} $G$ is the group (C17).

Let $N=\lg c,b,a^p,b^p\rg$ and $M_i=\lg c,ab^{i},a^p,b^p\rg$. Then
$N$ and $M_i$ are all maximal subgroups of $G$. Since $N=\lg
c,b\rg\ast\lg a^p\rg$, $N\in\mathcal{A}_2$ with $d(N)=3$ and
$|N'|=p$. By Lemma \ref{A_2-property} (7), $\alpha_1(N)=p^2$. Since
$M_i=\lg c,ab^i\rg$ and $|M_i'|=p$, $M_i\in\mathcal{A}_1$. Hence
$(\mu_0,\mu_1,\mu_2)=(0,p,1)$ and, by Hall's enumeration principle,
$$ \alpha_1(G)=\sum_{H\in \Gamma_1}
\alpha_1(H)=\mu_1+p^2\mu_2=p^2+p.$$

To sum up, we get the Table \ref{table4}.\qed

\begin{lem}\label{mu of d} Suppose that $G$ is an $\mathcal{A}_3$-group with $d(G)=3$, $\Phi(G)=
Z(G)$ and $G'\cong C_p^3$. Then we may assume that $G=\lg a_1,a_2,a_3;x,y,z\rg$ with
$$x=[a_2,a_3],y=[a_3,a_1],z=[a_1,a_2],x^p=y^p=z^p=1,a_i^{p^{m_i}}=x^{w_{i1}}y^{w_{i2}}z^{w_{i3}}, i=1,2,3.$$
Let $M=\langle a_2, a_3, \Phi(G)\rangle$, $M_i=\langle a_1a_2^i,
a_3, \Phi(G)\rangle$ and $M_{ij}=\langle a_1a_3^i, a_2a_3^j,
\Phi(G)\rangle$, where $0\leq i,j\leq p-1$. Then $M$, $M_i$ and
$M_{ij}$ are all maximal subgroups of $G$, and the following
conclusions hold:

$(1)$ If $m_1=2$, $m_2=m_3=1$ and $a_1^{p^2}=x$, then $M\in\mathcal{A}_2$, $M_i\in \mathcal{A}_2$ if and only if $w_{33}=0$,
and $M_{ij}\in\mathcal{A}_2$ if and only if $w_{22}+j(w_{23}+w_{32})+j^2w_{33}=0$;

$(2)$ If $m_1=m_2=2$, $m_3=1$, $a_3^{p}=z$ and $w_{13}=w_{23}=0$, then $M\in\mathcal{A}_2$,
$M_i\in \mathcal{A}_2$, and $M_{ij}\in\mathcal{A}_2$ if and only if $w_{11}w_{22}-w_{12}w_{21}=0$;

$(3)$ If $m_1=m_2=m_3=1$ and $p>2$, then $M\in\mathcal{A}_2$ if and only if $\left|
 \begin{array}{cc}
 w_{22}& w_{23}\\
 w_{32}& w_{33}
 \end{array}
 \right|=0$, $M_i\in \mathcal{A}_2$ if and only if ${\scriptsize\left|
 \begin{array}{ccc}
 i& -1 &0\\
 w_{11}+iw_{21}&w_{12}+iw_{22}&w_{13}+iw_{23}\\
 w_{31}&w_{32}&w_{33}
 \end{array}
 \right|}=0$, and $M_{ij}\in\mathcal{A}_2$ if and only if ${\scriptsize\left|
 \begin{array}{ccc}
 -i& -j &1\\
 w_{11}+iw_{31}&w_{12}+iw_{32}&w_{13}+iw_{33}\\
  w_{21}+jw_{31}&w_{22}+jw_{32}&w_{23}+jw_{33}
 \end{array}
 \right|}=0$;

$(4)$ If $m_1=m_2=m_3=1$ and $p=2$, then $M\in\mathcal{A}_2$ if and only if $\left|
 \begin{array}{cc}
 w_{22}& w_{23}\\
 w_{32}& w_{33}
 \end{array}
 \right|=0$, $M_i\in \mathcal{A}_2$ if and only if ${\scriptsize\left|
 \begin{array}{ccc}
 i& 1 &0\\
 w_{11}+iw_{21}&w_{12}+iw_{22}&w_{13}+iw_{23}+i\\
 w_{31}&w_{32}&w_{33}
 \end{array}
 \right|}=0$, and $M_{ij}\in\mathcal{A}_2$ if and only if ${\scriptsize\left|
 \begin{array}{ccc}
 i& j &1\\
 w_{11}+iw_{31}&w_{12}+iw_{32}+i&w_{13}+iw_{33}\\
  w_{21}+jw_{31}+j&w_{22}+jw_{32}&w_{23}+jw_{33}
 \end{array}
 \right|}=0$.
\end{lem}
\medskip
\demo Since $G\in \mathcal{A}_3$, any maximal subgroup is either $\mathcal{A}_1$-group or
$\mathcal{A}_2$-group. Let $H\in\Gamma_1(G)$. Then $|H'|=p$. If $d(H)=2$, then,
by Lemma \ref{minimal non-abelian equivalent conditions}, $H\in\mathcal{A}_1$.
Since $\Phi(H)\le \Phi(G)$, we get that $H\in\mathcal{A}_1$ if and only if $\Phi(G)=\Phi(H)$.

(1) If $m_1=2$, $m_2=m_3=1$ and $a_1^{p^2}=x$, then it is obvious
that $M=\lg a_2,a_3\rg\ast \lg a_1\rg\in\mathcal{A}_2$. Since
$x=(a_1a_2^i)^{p^2}\in \Phi(M_i)$ and $y=[a_3,a_1a_2^i]x^i\in
\Phi(M_i)$, $\Phi(M_i)=\Phi(G)$ if and only if $a_3^p\not\in\lg
x,y\rg$. Hence $M_i\in \mathcal{A}_2$ if and only if $w_{33}=0$.
Similar reason gives that $M_{ij}\in\mathcal{A}_1$ if and only if
$(a_2a_3^j)^p\not\in \lg x,[a_1a_3^i,a_2a_3^j]\rg=\lg x,y^{-j}z\rg$.
Hence $M_{ij}\in \mathcal{A}_2$ if and only if ${\scriptsize\left|
 \begin{array}{ccc}
 1& 0 &0\\
 0&-j&1\\
 *&w_{22}+jw_{32}&w_{23}+jw_{33}
 \end{array}
 \right|}=0$. The later is, $w_{22}+j(w_{23}+w_{32})+j^2w_{33}=0$.

\medskip
(2) If $m_1=m_2=2$, $m_3=1$, $a_3^{p}=z$ and $w_{13}=w_{23}=0$, then
it is obvious that $a_1^p\not\in \lg a_2,a_3\rg$ and $\lg
a_1^p,a_2^p\rg\not\le \lg a_1a_2^i,a_3\rg$. Hence $M\in
\mathcal{A}_2$ and $M_i\in\mathcal{A}_2$. Since
$(a_1a_3^i)^{p^2}=a_1^{p^2}=x^{w_{11}}y^{w_{12}}$ and
$(a_2a_3^j)^{p^2}=a_2^{p^2}=x^{w_{21}}y^{w_{22}}$ and
$[a_1a_2^i,a_2a_3^j]\not\in\lg x,y\rg$ $M_{ij}\in\mathcal{A}_1$ if
and only if $\left|
\begin{array}{cc}
w_{11} &w_{12}\\
w_{21}&w_{22}
\end{array}
\right|=0$. The later is, $w_{11}w_{22}-w_{12}w_{21}=0$.

\medskip
(3) If $m_1=m_2=m_3=1$ and $p>2$, then $M\in\mathcal{A}_1$ if and
only if $\lg a_2^p, a_3^p,[a_2,a_3]\rg=G'$. The later is
${\scriptsize\left|
 \begin{array}{ccc}
 w_{21}& w_{22}&w_{23}\\
 w_{31}&w_{32}&w_{33}\\
 1&0&0
 \end{array}
 \right|}\neq 0$. Hence
 $M\in\mathcal{A}_2$ if and only if $\left|
 \begin{array}{cc}
 w_{22}& w_{23}\\
 w_{32}& w_{33}
 \end{array}
 \right|=0$. Similar arguments give that $M_i\in \mathcal{A}_2$ if and only if ${\scriptsize\left|
 \begin{array}{ccc}
 i& -1 &0\\
 w_{11}+iw_{21}&w_{12}+iw_{22}&w_{13}+iw_{23}\\
 w_{31}&w_{32}&w_{33}
 \end{array}
 \right|}=0$, and $M_{ij}\in\mathcal{A}_2$ if and only if ${\scriptsize\left|
 \begin{array}{ccc}
 -i& -j &1\\
 w_{11}+iw_{31}&w_{12}+iw_{32}&w_{13}+iw_{33}\\
  w_{21}+jw_{31}&w_{22}+jw_{32}&w_{23}+jw_{33}
 \end{array}
 \right|}=0$.

\medskip
(4) It follows from similar arguments as (3). The details are
omitted.\qed

\begin{thm}\label{d=2-4}
$G$ has no abelian subgroup of index $p$, $G$ has at least two
distinct $\mathcal{A}_1$-subgroups of index $p$ and $d(G)=3$ if and
only if $G$ is isomorphic to one of the following pairwise
non-isomorphic groups:
\begin{enumerate}

\rr{Di} The type of $G/G'$ is $(p^2,p,p)$.
\begin{enumerate}

\rr{D1} $\langle a, b, c \di a^{p^{3}}=b^{p^2}=c^{p^2}=1,
[b,c]=a^{p^2},[c,a]=c^{-p},[a,b]=b^{p}c^{\nu p}\rangle$, where $p>2$
and $\nu=1$ or a fixed quadratic non-residue modulo $p$; moreover,
$|G|=p^7$, $\Phi(G)=Z(G)=\lg a^p,b^p,c^p\rg\cong C_{p^2}\times
C_{p}\times C_p$ and $G'=\lg a^{p^2},b^{p},c^p\rg \cong  C_p^3$.

\rr{D2} $\langle a, b, c \di a^{p^{3}}=b^{p^{2}}=c^{p^{2}}=1,
[b,c]=a^{p^2},[c,a]=b^{p},[a,b]=c^{\nu p}\rangle$,
where $p>2$, $\nu=1$ or a fixed quadratic non-residue modulo $p$; moreover, $|G|=p^7$,
$\Phi(G)=Z(G)=\lg a^p,b^p,c^p\rg\cong C_{p^2}\times C_{p}\times C_p$
and $G'=\lg a^{p^2},b^{p},c^p\rg\cong  C_p^3$.

\rr{D3} $\langle a, b, c \di a^{p^{3}}=b^{p^{2}}=c^{p^{2}}=1,
[b,c]=a^{p^2},[c,a]^{1+r}=b^{rp}c^{-p},[a,b]^{1+r}=b^pc^{p}\rangle$,
where $p>2$, $r=1,2,\dots,{p-2}$; moreover, $|G|=p^7$,
$\Phi(G)=Z(G)=\lg a^p,b^p,c^p\rg\cong C_{p^2}\times C_{p}\times C_p$
and $G'=\lg a^{p^2},b^{p},c^p\rg\cong  C_p^3$.

\rr{D4} $\langle a, b, c \di a^{8}=b^{4}=c^{4}=1,
[b,c]=a^{4},[c,a]=b^{2},[a,b]=c^{2}\rangle$; where $|G|=2^7$,
$\Phi(G)=Z(G)=\lg a^2,b^2,c^2\rg\cong C_{4}\times C_{2}\times C_2$
and $G'=\lg a^{4},b^{2},c^2\rg\cong  C_2^3$.

\rr{D5} $\langle a, b, c \di a^{8}=b^{4}=c^{4}=1,
[b,c]=a^{4},[c,a]=b^{2},[a,b]=b^{2}c^{2}\rangle$, where $|G|=2^7$,
$\Phi(G)=Z(G)=\lg a^2,b^2,c^2\rg\cong C_{4}\times C_{2}\times C_2$
and $G'=\lg a^{4},b^{2},c^2\rg\cong  C_2^3$.
\end{enumerate}

\rr{Dii} The type of $G/G'$ is $(p^2,p^2,p)$.

\begin{enumerate}
\rr{D6} $\langle a, b, c \di a^{p^3}=b^{p^3}=c^{p^2}=1,
[b,c]=a^{p^2},[c,a]=b^{\nu
p^{2}},[a,b]=c^{p}\rangle$, where $p>2$, $\nu=1$ or a fixed
quadratic non-residue modulo $p$ such that $-\nu\not\in (F_p)^2$; moreover, $|G|=p^8$,
$\Phi(G)=Z(G)=\lg a^p,b^p,c^p\rg\cong C_{p^2}\times C_{p^2}\times
C_p$ and $G'=\lg a^{p^2},b^{p^2},c^p\rg\cong  C_p^3$.

\rr{D7} $\langle a, b, c \di a^{p^3}=b^{p^3}=c^{p^2}=1,
[b,c]^{1+r}=a^{rp^2}b^{p^2},[c,a]^{1+r}=a^{-p^{2}}b^{p^{2}},[a,b]=c^{p}\rangle$,
where $p>2$, $r=1,2,\dots,p-2$ such that $-r\not\in (F_p)^2$; moreover, $|G|=p^8$,
$\Phi(G)=Z(G)=\lg a^p,b^p,c^p\rg\cong C_{p^2}\times C_{p^2}\times
C_p$ and $G'=\lg a^{p^2},b^{p^2},c^p\rg\cong  C_p^3$.

\rr{D8} $\langle a, b, c \di a^{8}=b^{8}=c^{4}=1,
[b,c]=a^{4}b^4,[c,a]=b^{4},[a,b]=c^{2}\rangle$; where $|G|=2^8$,
$\Phi(G)=Z(G)=\lg a^2,b^2,c^2\rg\cong C_{4}\times C_{4}\times C_2$
and  $G'=\lg a^{4},b^{4},c^2\rg \cong  C_2^3$.
\end{enumerate}

\rr{Diii} The type of $G/G'$ is $(p,p,p)$ where $p>2$.

\begin{enumerate}
\rr{D9} $\langle a, b, c \di a^{p^{2}}=b^{p^{2}}=c^{p^{2}}=1,
[b,c]=a^{p},[c,a]=b^{p},[a,b]=c^{p}\rangle$; where $|G|=p^6$,
$\Phi(G)=Z(G)=G'=\lg a^p,b^p,c^p\rg \cong C_p^3$.

\rr{D10} $\langle a, b, c \di a^{p^2}=b^{p^2}=c^{p^2}=1,
[b,c]=a^{p},[c,a]=c^{-p},[a,b]=b^{p}\rangle$;  where $|G|=p^6$,
$\Phi(G)=Z(G)=G'=\lg a^p,b^p,c^p\rg \cong C_p^3$.

\rr{D11} $\langle a, b, c \di a^{p^2}=b^{p^2}=c^{p^2}=1,
[b,c]=a^{p},[c,a]=c^{-p},[a,b]=b^{p}c^{\nu p}\rangle$, where $\nu=1$ or
is a fixed quadratic non-residue modulo $p$; moreover, $|G|=p^6$,
$\Phi(G)=Z(G)=G'=\lg a^p,b^p,c^p\rg \cong C_p^3$.

\rr{D12} $\langle a, b, c \di a^{p^2}=b^{p^2}=c^{p^2}=1,
[b,c]=a^{p},[c,a]^{1+r}=b^{rp}c^{-p},[a,b]^{1+r}=b^{p}c^{p}\rangle$, where
$r=1,2,\dots,p-2$. Moreover, $|G|=p^6$, $\Phi(G)=Z(G)=G'=\lg a^p,b^p,c^p\rg
\cong C_p^3$.

\rr{D13} $\langle a, b, c \di a^{p^2}=b^{p^2}=c^{p^2}=1,
[b,c]=a^{-p}b^pc^{p},[c,a]=a^{-p}b^{p},[a,b]=a^{p}\rangle$; where
$|G|=p^6$, $\Phi(G)=Z(G)=G'=\lg a^p,b^p,c^p\rg \cong C_p^3$.

\rr{D14} $\langle a, b, c; d \di a^{p^2}=b^{p^2}=c^{p}=d^p=1,
[b,c]=a^{p},[c,a]=b^{\nu
p},[a,b]=d,[d,a]=[d,b]=[d,c]=1\rangle$, where $p>2$, $\nu=1$ or a
fixed quadratic non-residue modulo $p$ such that $-\nu\not\in
(F_p)^2$; moreover, $|G|=p^6$,
$\Phi(G)=Z(G)=G'=\lg a^p,b^p,d\rg \cong C_p^3 $.

\rr{D15} $\langle a, b, c; d \di a^{p^2}=b^{p^2}=c^{p}=d^p=1,
[b,c]=a^{p},[c,a]=d,[a,b]=b^p,[d,a]=[d,b]=[d,c]=1\rangle$, where $p>2$; moreover, $|G|=p^6$,
$\Phi(G)=Z(G)=G'=\lg a^p,b^p,d\rg \cong C_p^3 $.

\rr{D16} $\langle a, b, c; d \di a^{p}=b^{p^2}=c^{p^2}=d^p=1,
[b,c]=d,[c,a]^{1+r}=b^{rp}c^{-p},[a,b]^{1+r}=b^pc^p,[d,a]=[d,b]=[d,c]=1\rangle$, where $p>2$, $r=1,2,\dots,p-2$ such that $-r\not\in
(F_p)^2$; moreover, $|G|=p^6$,
$\Phi(G)=Z(G)=G'=\lg b^p,c^p,d\rg \cong C_p^3 $.

\end{enumerate}

\rr{Div} The type of $G/G'$ is $(2,2,2)$.
\begin{enumerate}
\rr{D17} $\langle a, b,c; d \di a^{4}=b^{2}=c^{4}=d^2=1,
[b,c]=d,[c,a]=a^2,[a,b]=c^2,[d,a]=[d,b]=[d,c]=1\rangle$; where
$\Phi(G)=Z(G)=G'=\lg a^2,c^2,d\rg \cong C_2^3$.

\rr{D18} $\langle a,b,c; d \di a^{4}=b^{4}=c^{4}=d^2=1,
[b,c]=d,[c,a]=a^2,[a,b]=b^2=c^2,[d,a]=[d,b]=[d,c]=1\rangle$; where
$\Phi(G)=Z(G)=G'=\lg a^2,b^2,d\rg \cong C_2^3$.

\rr{D19} $\langle a,b,c; d \di a^{4}=b^{4}=c^{4}=d^2=1,
[b,c]=d,[c,a]=a^2b^2,[a,b]=a^2=c^2,[d,a]=[d,b]=[d,c]=1\rangle$;
$\Phi(G)=Z(G)=G'=\lg a^2,b^2,d\rg \cong C_2^3$.
\end{enumerate}
\end{enumerate}
Moreover, Table $\ref{table6}$ gives $(\mu_0,\mu_1,\mu_2)$  and
$\alpha_1(G)$ for the groups {\rm (D1)--(D19)}.
\begin{table}[h]
\centering
{
\tiny
\begin{tabular}[t]{c|c||c}
\hline
 $(\mu_0,\ \mu_1,\ \mu_2)$ &$\alpha_1(G)$     & types of  $\mathcal{A}$$_3$ groups \\

\hline $(0,\ p^2,\ p+1)$  &$p^3+2p^2$& (D1),(D4), (D6)--(D10), (D12)--(D14), (D16)  \\

\hline {$(0,\ p^2+p,\ 1)$} &$2p^2+p$ & (D2) where $-\nu\not\in (F_p^*)^2$;(D3) where $-r \not\in (F_p^*)^2$; (D5); (D11) where $-\nu \not\in (F_p^*)^2$;  \\

\hline $(0,\ p^2-p,\ 2p+1)$ &$2p^3+2p^2-p$ & (D2) where $-\nu \in (F_p^*)^2$; (D3) where $-\nu \in (F_p^*)^2$; (D11) where $-\nu \in (F_p^*)^2$; (D15) \\

\hline $(0,\ p^2-1,\ p+2)$ &$p^3+3p^2-1$ & (D17), (D18) \\

\hline $(0,\ p^2+1,\ p)$ &$p^3+p^2+1$ & (D19) \\

\hline

\end{tabular}}
\caption{The enumeration of (D1)--(D19)\label{table6}}
\end{table}
\end{thm}

\demo Let $A$ and $B$ be two distinct $\mathcal{A}_1$-subgroups of
index $p$. By Lemma \ref{G'}, $|G'|\le p|A'||B'|\le p^3$. It is easy
to see that $\Phi(A)=\Phi(B)=\Phi(G)$. By Lemma \ref{minimal
non-abelian equivalent conditions}, we have $\Phi(A)=Z(A)$ and
$\Phi(B)=Z(B)$. Since $[\Phi(G),A]=[\Phi(G),B]=1$ and $G=AB$, we
have $\Phi(G)\le Z(G)$. Moreover, $G'\le C_p^3$. If $|G'|\le p^2$,
then $G$ has an abelian subgroup of index $p$, a contradiction. Now
we have $d(G)=3, \Phi(G)\le Z(G)$ and $G'\cong C_p^3$. Thus $G$ is
one of the groups classified in \cite{QXA}. Let the type of $G/G'$
be $(p^{m_1},p^{m_2},p^{m_3})$, where $m_1\geq m_2\geq m_3$. By
\cite[Theorem 2.7]{QXA}, $m_3=1$. By \cite[Theorem 2.8]{QXA},
$m_1\le 2$.

If $m_1=2$ and $m_2=1$, then $G$ is either one of the groups
determined by \cite[Theorem 4.1]{QXA} for $p>2$ or one of the groups
determined by \cite[Theorem 7.1]{QXA} for $p=2$. By \cite[Theorem
4.3 \& 7.4]{QXA} , we get the groups (D1)--(D5).

If $m_1=m_2=2$, then $G$ is one of the groups determined by
\cite[Theorem 5.1]{QXA}. By \cite[Theorem 5.2]{QXA} , we get the
groups (D6)--(D8).

If $m_1=1$ and $p>2$, then $G$ is one of the groups determined by
\cite[Theorem 6.1]{QXA}. By \cite[Theorem 6.3]{QXA}, we get the
groups (D9)--(D16).

If $m_1=1$ and $p=2$, then $G$ is one of the groups listed in
\cite[Theorem 7.6]{QXA}. By \cite[Theorem 7.7]{QXA} , we get the
groups (D17)--(D19). Table \ref{table7} gives the correspondence.

\begin{table}[h]
  \centering
{
\tiny
\begin{tabular}{c|c||c|c}
\hline
\multirow{2}*{Groups} &Groups in & \multirow{2}*{Groups} & Groups in  \\

 &\cite[Theorem 4.1 \& 5.1 \& 6.1 \& 7.1 \& 7.6]{QXA}&  &  \cite[Theorem 4.1 \& 5.1 \& 6.1 \& 7.1 \& 7.6]{QXA} \\

\hline
(D1)&(D2) where $m_1=2$ and $m_2=1$ & (D11)&(J3) where $m_1=m_2=m_3=1$\\
(D2)&(D3) where $m_1=2$ and $m_2=1$ & (D12)&(J4) where $m_1=m_2=m_3=1$\\
(D3)&(D4) where $m_1=2$ and $m_2=1$ & (D13)&(J5) where $m_1=m_2=m_3=1$\\
(D4)&(M2) where $m_1=2$ & (D14)&(K1) where $m_1=m_2=m_3=1$ and $-\nu\not\in F_p^2$\\
(D5)&(M3) where $m_1=2$ & (D15)&(K2)\\
(D6)&(G3) where $-\nu\not\in F_p^2$, $m_1=m_2=2$ and $m_3=1$ & (D16)&(K6) where $m_1=m_2=m_3=1$ and $-r\not\in F_p^2$\\
(D7)&(G4) where $-r\not\in F_p^2$, $m_1=m_2=2$ and $m_3=1$ & (D17)&(S7)\\
(D8)&(G7) where $m_1=m_2=2$ and $m_3=1$ & (D18)&(S8)\\
(D9)&(J1) where $m_1=m_2=m_3=1$ &(D19) &(S9)\\
(D10)&(J2) where $m_1=m_2=m_3=1$ & &\\
\hline
\end{tabular}
\caption{{\small The correspondence from Theorem \ref{d=2-4} to
\cite[Theorem 4.1,5.1,6.1,7.1 \& 7.6]{QXA}}} \label{table7}}
\end{table}

\medskip
We calculate the $(\mu_0,\mu_1,\mu_2)$ and $\alpha_1(G)$ of those
groups in Theorem \ref{d=2-4} as follows.

By Lemma \ref{alpha of d=3}, $\alpha_1(G)=\mu_1+p^2\mu_2$. Hence we
only need to calculate $(\mu_0,\mu_1,\mu_2)$. we may assume that
$G=\lg a_1,a_2,a_3;x,y,z\rg$ with
$$x=[a_2,a_3],y=[a_3,a_1],z=[a_1,a_2],x^p=y^p=z^p=1,a_i^{p^{m_i}}=x^{w_{i1}}y^{w_{i2}}z^{w_{i3}}, i=1,2,3.$$
Let $M=\langle a_2, a_3, \Phi(G)\rangle$,
$M_i=\langle a_1a_2^i, a_3, \Phi(G)\rangle$ and
$M_{ij}=\langle a_1a_3^i, a_2a_3^j, \Phi(G)\rangle$, where $0\leq i,j\leq p-1$.
Then $M$, $M_i$ and $M_{ij}$ are all maximal subgroups of $G$. Let $w(G)=(w_{ij})$.

If $G$ is the group (D1), then $m_1=2$, $m_2=m_3=1$ and $w(G)={\scriptsize\left(
 \begin{array}{ccc}
1&0&0\\
0& \nu& 1\\
 0& -1& 0
 \end{array}
 \right)}$. By Lemma \ref{mu of d} (1), $M\in\mathcal{A}_2$, $M_i\in\mathcal{A}_2$
and $M_{ij}\in\mathcal{A}_1$. Hence  $(\mu_0, \mu_1, \mu_2)=(0,p^2,p+1)$ and $\alpha_1(G)=p^3+2p^2$.

If $G$ is the group (D2), then $m_1=2$, $m_2=m_3=1$ and $w(G)={\scriptsize\left(
 \begin{array}{ccc}
1&0&0\\
0& 1&0\\
 0& 0& \nu^{-1}
 \end{array}
 \right)}$. By Lemma \ref{mu of d} (1), $M\in\mathcal{A}_2$, $M_i\in\mathcal{A}_1$
and $M_{ij}\in\mathcal{A}_2$ if and only if $j^2=-\nu$. Hence  $(\mu_0, \mu_1, \mu_2)=(0,p^2+p,1)$ for
$-\nu\not\in (F_p^*)^2$ and $(\mu_0, \mu_1, \mu_2)=(0,p^2-p,2p+1)$ for  $-\nu\in (F_p^*)^2$.

If $G$ is the group (D3), then $m_1=2$, $m_2=m_3=1$ and $w(G)={\scriptsize\left(
 \begin{array}{ccc}
1&0&0\\
0& 1&1\\
 0& -1& r
 \end{array}
 \right)}$. By Lemma \ref{mu of d} (1), $M\in\mathcal{A}_2$, $M_i\in\mathcal{A}_1$
and $M_{ij}\in\mathcal{A}_2$ if and only if $j^2=-r^{-1}$. Hence  $(\mu_0, \mu_1, \mu_2)=(0,p^2+p,1)$ for
 $-r\not\in (F_p^*)^2$ and $(\mu_0, \mu_1, \mu_2)=(0,p^2-p,2p+1)$ for  $-r\in (F_p^*)^2$.

If $G$ is the group (D4), then $m_1=2$, $m_2=m_3=1$ and $w(G)={\scriptsize\left(
 \begin{array}{ccc}
1&0&0\\
0& 1& 0\\
 0& 0& 1
 \end{array}
 \right)}$. By Lemma \ref{mu of d} (1), $M\in\mathcal{A}_2$, $M_i\in\mathcal{A}_1$
and $M_{ij}\in\mathcal{A}_2$ if and only if $j=1$. Hence  $(\mu_0, \mu_1, \mu_2)=(0,p^2,p+1)$ and $\alpha_1(G)=p^3+2p^2$.

If $G$ is the group (D5), then $m_1=2$, $m_2=m_3=1$ and $w(G)={\scriptsize\left(
 \begin{array}{ccc}
1&0&0\\
0& 1& 0\\
 0& 1& 1
 \end{array}
 \right)}$. By Lemma \ref{mu of d} (1), $M\in\mathcal{A}_2$, $M_i\in\mathcal{A}_1$
and $M_{ij}\in\mathcal{A}_1$. Hence  $(\mu_0, \mu_1, \mu_2)=(0,p^2+p,1)$ and $\alpha_1(G)=2p^2+p$.

If $G$ is one of the groups (D6)--(D8), then $m_1=m_2=2$, $m_3=1$ and $w(G)={\scriptsize\left(
 \begin{array}{ccc}
1&0&0\\
0& \nu& 0\\
 0& 0& 1
 \end{array}
 \right)}$, ${\scriptsize\left(
 \begin{array}{ccc}
1&-1&0\\
1& r& 0\\
 0& 0& 1
 \end{array}
 \right)}$, or ${\scriptsize\left(
 \begin{array}{ccc}
1&1&0\\
0& 1& 0\\
 0& 0& 1
 \end{array}
 \right)}$. By Lemma \ref{mu of d} (2), $M\in\mathcal{A}_2$, $M_i\in\mathcal{A}_2$
and $M_{ij}\in\mathcal{A}_1$. Hence  $(\mu_0, \mu_1, \mu_2)=(0,p^2,p+1)$ and $\alpha_1(G)=p^3+2p^2$.

If $G$ is the group (D9), then $p>2$ and $m_1=m_2=m_3=1$ and $w(G)={\scriptsize\left(
 \begin{array}{ccc}
 1& 0 &0\\
 0&1&0\\
 0&0&1
 \end{array}
 \right)}$. By Theorem \ref{mu of d} (3), $M\in\mathcal{A}_1$, $M_i\in\mathcal{A}_2$ if and only if
 $i^2+1=0$ and $M_{ij}\in \mathcal{A}_2$ if and only if $i^2+j^2+1=0$.

 If $-1\in (F_p^*)^2$, then $i^2+1=0$ has two solutions, and by Lemma \ref{equivalent equation 1},
 $i^2+j^2+1=0$ has $p-1$ solutions. Hence $(\mu_0,\mu_1,\mu_2)=(0,p^2,p+1)$ .
 If $-1\not\in (F_p^*)^2$, then $i^2+1=0$ has no solution, and by Lemma \ref{equivalent equation 1},
 $i^2+j^2+1=0$ has $p+1$ solutions. Hence we also have $(\mu_0,\mu_1,\mu_2)=(0,p^2,p+1)$.

\medskip
If $G$ is the group (D10), then $p>2$ and $m_1=m_2=m_3=1$ and $w(G)={\scriptsize\left(
 \begin{array}{ccc}
 1& 0 &0\\
 0&0&1\\
 0&-1&0
 \end{array}
 \right)}$. By Theorem \ref{mu of d} (3), $M\in\mathcal{A}_1$, $M_i\in\mathcal{A}_2$
 if and only if $i=0$ and $M_{ij}\in \mathcal{A}_2$ if and only if $i=0$.
Hence $(\mu_0,\mu_1,\mu_2)=(0,p^2,p+1)$.

\medskip
If $G$ is the group (D11), then $p>2$ and $m_1=m_2=m_3=1$ and
$w(G)={\scriptsize\left|
 \begin{array}{ccc}
 1& 0 &0\\
 0&\nu&1\\
 0&-1&0
 \end{array}
 \right|}$. By Theorem \ref{mu of d} (3), $M\in\mathcal{A}_1$, $M_i\in\mathcal{A}_2$ if and only if
 $i=0$ and $M_{ij}\in \mathcal{A}_2$ if and only if $i^2=-\nu$.
Hence  $(\mu_0, \mu_1, \mu_2)=(0,p^2+p,1)$ for  $-\nu\not\in (F_p^*)^2$ and $(\mu_0, \mu_1, \mu_2)=(0,p^2-p,2p+1)$ for  $-\nu\in (F_p^*)^2$.

\medskip
If $G$ is the group (D12), then $p>2$ and $m_1=m_2=m_3=1$ and
$w(G)={\scriptsize\left(
 \begin{array}{ccc}
 1& 0 &0\\
 0&1&1\\
 0&-1&r
 \end{array}
 \right)}$. By Theorem \ref{mu of d} (3), $M\in\mathcal{A}_1$, $M_i\in\mathcal{A}_2$ if and only if
 $i^2(r+1)+r=0$ and $M_{ij}\in \mathcal{A}_2$ if and only if $i^2(r+1)+rj^2+1=0$.

 If $-(r+1)^{-1}r\in (F_p^*)^2$, then $i^2(r+1)+r=0$ has two solutions, and by Lemma \ref{equivalent equation 1},
 $i^2(r+1)+j^2r+1=0$ has $p-1$ solutions. Hence $(\mu_0,\mu_1,\mu_2)=(0,p^2,p+1)$.
 If $-(r+1)^{-1}r\not\in (F_p^*)^2$, then $i^2(r+1)+r=0$ has no solution, and by Lemma \ref{equivalent equation 1},
 $i^2(r+1)+j^2r+1=0$ has $p+1$ solutions. Hence we also have $(\mu_0,\mu_1,\mu_2)=(0,p^2,p+1)$.

\medskip
If $G$ is the group (D13), then $p>2$ and $m_1=m_2=m_3=1$ and
$w(G)={\scriptsize\left(
 \begin{array}{ccc}
 0& 0 &1\\
 0&1&1\\
 1&-1&0
 \end{array}
 \right)}$. By Theorem \ref{mu of d} (3), $M\in\mathcal{A}_1$, $M_i\in\mathcal{A}_2$ if and only if
 $i^2=1$ and $M_{ij}\in \mathcal{A}_2$ if and only if $(i+1)^2-j^2-1=0$. By Lemma \ref{equivalent equation 1},
 $(i+1)^2-j^2-1=0$ has $p-1$ solutions. Hence $(\mu_0,\mu_1,\mu_2)=(0,p^2,p+1)$.

\medskip
If $G$ is the group (D14), then $p>2$ and $m_1=m_2=m_3=1$ and
$w(G)={\scriptsize\left(
 \begin{array}{ccc}
 1& 0 &0\\
 0&\nu&0\\
 0&0&0
 \end{array}
 \right)}$. By Theorem \ref{mu of d} (3), $M\in\mathcal{A}_2$, $M_i\in\mathcal{A}_2$ and $M_{ij}\in \mathcal{A}_1$.
 Hence $(\mu_0,\mu_1,\mu_2)=(0,p^2,p+1)$.

\medskip
If $G$ is the group (D15), then $p>2$ and $m_1=m_2=m_3=1$ and
$w(G)={\scriptsize\left(
 \begin{array}{ccc}
 1& 0 &0\\
 0&0&1\\
 0&0&0
 \end{array}
 \right)}$. By Theorem \ref{mu of d} (3), $M\in\mathcal{A}_2$, $M_i\in\mathcal{A}_2$ and $M_{ij}\in \mathcal{A}_1$
 if and only if $j=0$. Hence $(\mu_0,\mu_1,\mu_2)=(0,p^2-p,2p+1)$.

\medskip
If $G$ is the group (D16), then $p>2$ and $m_1=m_2=m_3=1$ and
$w(G)={\scriptsize\left(
 \begin{array}{ccc}
 0& 0 &0\\
 0&1&1\\
 0&-1&r
 \end{array}
 \right)}$. By Theorem \ref{mu of d} (3), $M\in\mathcal{A}_1$, $M_i\in\mathcal{A}_2$ if and only if
 $i=0$ and $M_{ij}\in \mathcal{A}_1$ if and only if $i=0$. Hence $(\mu_0,\mu_1,\mu_2)=(0,p^2,p+1)$.

\medskip
If $G$ is the group (D17), then $p=2$ and $m_1=m_2=m_3=1$ and
$w(G)={\scriptsize\left(
 \begin{array}{ccc}
 0& 1 &0\\
 0&0&0\\
 0&0&1
 \end{array}
 \right)}$. By Theorem \ref{mu of d} (3), $M\in\mathcal{A}_2$, $M_i\in\mathcal{A}_2$ if and only if
 $i=0$ and $M_{ij}\in \mathcal{A}_1$ if and only if $j=0$. Hence $(\mu_0,\mu_1,\mu_2)=(0,3,4)$.

\medskip
If $G$ is the group (D18), then $p=2$ and $m_1=m_2=m_3=1$ and
$w(G)={\scriptsize\left(
 \begin{array}{ccc}
 0& 1 &0\\
 0&0&1\\
 0&0&1
 \end{array}
 \right)}$. By Theorem \ref{mu of d} (3), $M\in\mathcal{A}_2$, $M_i\in\mathcal{A}_2$ if and only if
 $i=0$ and $M_{ij}\in \mathcal{A}_1$ if and only if $j=0$. Hence $(\mu_0,\mu_1,\mu_2)=(0,3,4)$.

\medskip
If $G$ is the group (D19), then $p=2$ and $m_1=m_2=m_3=1$ and
$w(G)={\scriptsize\left(
 \begin{array}{ccc}
 0& 0 &1\\
 0&1&1\\
 0&0&1
 \end{array}
 \right)}$. By Theorem \ref{mu of d} (3), $M\in\mathcal{A}_1$, $M_i\in\mathcal{A}_2$ if and only if
 $i=0$ and $M_{ij}\in \mathcal{A}_1$ if and only if $i=j=0$. Hence $(\mu_0,\mu_1,\mu_2)=(0,5,2)$.
\qed

\begin{thm}\label{d=2-5}
$G$ has no abelian subgroup of index $p$ and $G$ has a unique
$\mathcal{A}_1$-subgroups of index $p$ if and only if $G$ is
isomorphic to one of the following pairwise non-isomorphic groups:
\begin{enumerate}
\rr{Ei} $p>2$. In this case $p=3$, $(\mu_0,\mu_1,\mu_2)=(0,1,3)$ and $\alpha_1(G)=10$ except {\rm (E2)}, in which $\alpha_1(G)=28$.
\begin{enumerate}
\rr{E1} $\lg {s}_1, s;{s}_2 \di
{{s}_1}^{9}={{s}_2}^{9}=1,s^3=s_2^{3\delta}, [{s}_1, s]={s}_2,
[{s}_2,s]={{s}_2}^{-3}{{s}_1}^{-3}, [{s}_2,{s}_1]=s_2^{3}\rg$ where
$\delta=0,1,2$; moreover, $|G|=3^{5}$, $c(G)=4$, $\Phi(G)=G'=\lg
s_1^3, s_2\rg\cong C_3\times C_9$ and $Z(G)=\lg s_2^3\rg\cong C_3$.

\rr{E2} $\lg a,b;c\di
a^{3^2}=b^{3^2}=c^3=1,[b,a]=c,[c,a]=a^3,[c,b]=b^{-3}\rg$; where
$|G|=3^{5}$, $c(G)=3$, $\Phi(G)=G'=\lg a^3, b^3,c\rg\cong C_3^3$ and
$Z(G)=\lg a^3, b^3\rg\cong C_3^2$.

\rr{E3}{\small  $\lg s_1, \beta;s_2,  x \di
{s_1}^{9}={s_2}^{9}=x^3=1,{\beta}^3=x, [s_1, \beta]=s_2$, $[s_2,
\beta]={s_2}^{-3}{s_1}^{-3},[s_1,s_2]=x, [x,s_1]=[x,\beta]=1\rg$;
where $|G|=3^{6}$, $c(G)=4$, $\Phi(G)=G'=\lg s_1^3, s_2,x\rg\cong
C_9\times C_3\times C_3$ and $Z(G)=\lg s_2^3,x\rg\cong C_3^2$.}

\rr{E4} {\small $\lg s_1, \beta;s_2,  x \di
{s_1}^{9}={s_2}^{9}=x^3=1,{\beta}^3={s_2}^{3}x, [s_1, \beta]=s_2,
[s_2, \beta]={s_2}^{-3}{s_1}^{-3}, [s_1,s_2]$.}

$=x, [x,s_1]=[x,\beta]=1\rg$; where $|G|=3^{6}$, $c(G)=4$,
$\Phi(G)=G'=\lg s_1^3, s_2,x\rg\cong C_9\times C_3\times C_3$ and
$Z(G)=\lg s_2^3,x\rg\cong C_3^2$.

\rr{E5} $\lg \alpha, \beta; s_1, s_2, x  \di
{s_1}^{9}={s_2}^{3}=x^3=1,{\beta}^3=x^2, {\alpha}^3={s_2}^{-1},
[\alpha, \beta]=s_1, [s_1,
\alpha]=x,[s_1,\beta]=s_2,[s_2,\beta]={s_1}^{-3},
[s_1,s_2]=[x,\alpha]=[x,\beta]=1\rg$; where  $|G|=3^{6}$, $c(G)=4$,
$\Phi(G)=G'=\lg s_1, s_2,x\rg\cong C_9\times C_3\times C_3$ and
$Z(G)=\lg s_1^3,x\rg\cong C_3^2$.

\rr{E6} $\lg \alpha, \beta;s_1, s_2, x  \di
{s_1}^{9}={s_2}^{3}=x^3=1,{\beta}^3=x, {\alpha}^3={s_2}^{-1}x,
[\alpha, \beta]=s_1, [s_1,
\alpha]=x,[s_1,\beta]=s_2,[s_2,\beta]={s_1}^{-3},
[s_1,s_2]=[x,\alpha]=[x,\beta]=1\rg$, where $|G|=3^{6}$, $c(G)=4$,
$\Phi(G)=G'=\lg s_1, s_2,x\rg\cong C_9\times C_3\times C_3$ and
$Z(G)=\lg s_1^3,x\rg\cong C_3^2$.
\end{enumerate}

\rr{Eii} $p=2$.

\begin{enumerate}

\rr{E7} $\lg a,b\di
a^{16}=1,b^{2^{s+t+2}}=a^{2^{s+t'+2}},[a,b]=a^2\rg$, where $s,t,t'$
are non-negative integers with $1\le s+t'\le 2$, $tt'=0$ and if
$t'\ne 0$, then $s+t'=2$; moreover, $|G|=2^{s+t+6}$, $c(G)=4$,
$\Phi(G)=\lg a^2, b^2\rg\cong C_8\times C_{2^{t-t'+3}}$, $G'=\lg
a^2\rg\cong C_8$, $Z(G)=\lg a^8,b^4\rg\cong C_2\times
C_{2^{t-t'+2}}$, $(\mu_0,\mu_1,\mu_2)=(0,1,2)$ and $\alpha_1(G)=5$;

\rr{E8} $\lg a,b,c\di
a^2=b^8=1,c^2=b^{4t},[a,b]=b^4,[a,c]=1,[c,b]=b^2\rg$ where $t=0,1$;
moreover, $|G|=2^5$, $c(G)=3$, $\Phi(G)=G'=\lg b^2\rg\cong C_4$,
$Z(G)=\lg b^4\rg\cong C_2$, $(\mu_0,\mu_1,\mu_2)=(0,1,6)$ and $\alpha_1(G)=9$;

\rr{E9} $\lg a,b,c;d\di
a^4=b^4=c^2=d^2=1,[a,b]=b^2,[c,a]=a^2b^2,[c,b]=d,[d,a]=[d,b]=[d,c]=1\rg$;
moreover $|G|=2^6$, $c(G)=2$, $\Phi(G)=Z(G)=G'=\lg a^2,b^2,d\rg\cong
C_2^3$, $(\mu_0,\mu_1,\mu_2)=(0,1,6)$ and $\alpha_1(G)=25$;

\rr{E10} $\lg a,b,c;d\di
a^4=b^4=c^2=d^2=1,[a,b]=a^2,[c,a]=a^2b^2,[c,b]=d,[d,a]=[d,b]=[d,c]=1\rg$,
moreover $|G|=2^6$, $c(G)=2$, $\Phi(G)=Z(G)=G'=\lg a^2,b^2,d\rg\cong
C_2^3$, $(\mu_0,\mu_1,\mu_2)=(0,1,6)$ and $\alpha_1(G)=25$.
\end{enumerate}
\end{enumerate}
\end{thm}

\demo Finite $p$-groups with a unique $\mathcal{A}_1$-subgroup were
classified in \cite{QYXA} for $p>2$ and \cite{QZGA} for $p=2$. If
$p>2$, then, by using the classification in \cite{QYXA} and
calculating the maximal index of $\mathcal{A}_1$-subgroups for these
groups, we get the groups (E1)--(E6). The details are omitted. Table
\ref{table 100} gives the correspondence.

\begin{table}[h]
  \centering
{
\tiny
\begin{tabular}{c|c||c|c||c|c}
\hline
Groups & Groups in \cite[Lemma 2.12]{QYXA} &Groups & Groups in \cite[Corollary 3.5]{QYXA}&Groups & Groups in \cite[Theorem 3.8]{QYXA} \\
\hline
(E1)&  (2) where $e=2$ & (E2)&  (2)  &(E3)&  (1) where $e=2$ and $k=1$\\
                                &        &  &  & (E4)&  (2)  where $e=2$ and $k=1$\\
                                 &       &     &  & (E5)& (3) where $e=2$ and $k=2$\\
                                  &        &    &  & (E6) &  (4) where $e=2$ and $k=1$\\

\hline
\end{tabular}
\caption{\footnotesize The correspondence from Theorem \ref{d=2-4}
to \cite[Lemma 2.12, Corollary 3.5, Theorem 3.8]{QYXA}} \label{table
100}}
\end{table}

If $p=2$, then, by lemma \ref{55}, $G$ is either metacyclic or
$|G|\le 2^8$. If $G$ is metacyclic, then, by Lemma \ref{metacyclic
An}, $|G'|=p^3$. By using the classification of metacyclic
$2$-groups in \cite{XZ} and checking the number of their
$\mathcal{A}_1$-subgroups of index $p$, we get the group (E7). The
details is omitted. If $G$ is not metacyclic, then $|G|\le 2^8$. By
using Magma to check the SmallGroup database, we get the groups
(E8)--(E10).

\medskip
We calculate the $(\mu_0,\mu_1,\mu_2)$ and $\alpha_1(G)$ of those
groups in Theorem \ref{d=2-5} as follows.

By hypothesis, $\mu_0=0$ and $\mu_1=1$. Hence $\mu_2=p$ for $d(G)=2$
and $\mu_2=p+p^2$ for $d(G)=3$. In the following, we calculate
$\alpha_1(G)$.

If $G$ is the group {\rm (E2)}, then $|G|=3^5$, and $|H'|=p$ for any maximal subgroup $H$.
By Lemma \ref{A_2-property} (7), $\alpha_1(H)=p^2$ for $H\in\mathcal{A}_2$. By Hall's enumeration principle,
$$\alpha_1(G)=\sum_{H\in
\Gamma_1} \alpha_1(H)=\mu_1+p^2\mu_2=p^3+1=28.$$

If $G$ is one of the group {\rm (E1)} and (E3)--(E6), then $c(G)=4$, and $c(H)=3$ for any $\mathcal{A}_2$-subgroup $H$.
Since $d(H)=2$, $\alpha_1(H)=p$. By Hall's enumeration principle,
$$\alpha_1(G)=\sum_{H\in
\Gamma_1} \alpha_1(H)=\mu_1+p\mu_2=p^2+1=10.$$

If $G$ is the group {\rm (E7)}, then $G$ is metacyclic. Hence $H$ is metacyclic for any $\mathcal{A}_2$-subgroup $H$.
Since $d(H)=2$, $\alpha_1(H)=p$. By Hall's enumeration principle,
$$\alpha_1(G)=\sum_{H\in
\Gamma_1} \alpha_1(H)=\mu_1+p\mu_2=p^2+1=5.$$

If $G$ is the group {\rm (E8)}, then $d(G)=3$ and $c(G)=3$.
All maximal subgroups of $G$ are
$M=\lg b, a\rg$, $M_0=\lg c,a,b^2\rg$, $M_1=\lg cb,a,b^2\rg$ and $M_{ij}=\lg ca^i, ba^j\rg$ where $i,j=0,1$.
Here $M$ is the unique $\mathcal{A}_1$-subgroup of index $2$, $|M_i'|=2$ and $c(M_{ij})=3$. By Lemma \ref{A_2-property},
$\alpha_1(M_i)=4$ and $\alpha_1(M_{ij})=2$. Let $H\in \Gamma_2(G)$. Then $H=\lg a^ib^jc^k,b^2\rg$ where $(i,j,k)\neq (0,0,0)$.
It is obvious that $H\in\mathcal{A}_1$ if and only if $k=1$.  By Hall's enumeration principle,
$$\alpha_1(G)=\sum_{H\in
\Gamma_1} \alpha_1(H)-2\sum_{H\in \Gamma_2} \alpha_1(H)=1+4\times
2+2\times 4-2\times 4=9.$$

If $G$ is the group {\rm (E9)} or (E10), then $d(G)=3$ and $\Phi(G)=Z(G)$.
By Theorem \ref{alpha of d=3}, $\alpha_1(G)=\mu_1+4\mu_2=25$.\qed

\section{$\mathcal{A}_3$-groups without $\mathcal{A}_1$-subgroup of index
$p$}

Assume $G$ is an $\mathcal{A}_3$-group without an
$\mathcal{A}_1$-subgroup of index $p$. We discuss according to $G$
has an abelian subgroup of index $p$ or not. In this case,
$\mathcal{A}_3$-groups have 151 non-isomorphic types. Theorem
\ref{d=3-1}, \ref{d=3-2}, \ref{d=3-3}, \ref{d=3-4} and \ref{d=3-5}
give the classification of $\mathcal{A}_3$-groups without an
$\mathcal{A}_1$-subgroup of index $p$ and with an abelian subgroup
of index $p$. They are the groups (F1)--(F8), (G1)--(G12),
(H1)--(H10), (I1)--(I11) and (J1)--(J9).  Theorem \ref{d=4-1},
\ref{d=4-2}, \ref{d=4-3}, \ref{d=4-4} and \ref{d=4-5} give the
classification of $\mathcal{A}_3$-groups without an
$\mathcal{A}_1$-subgroup of index $p$ and an abelian subgroup of
index $p$. They are the groups (K1)--(K5), (L1)--(L2), (M1)--(M62),
(N1)--(N26) and (O1)--(O6).

\subsection{$G$ has an abelian subgroup of index $p$}

In this section assume $G$ is an $\mathcal{A}_3$-group with an
abelian subgroup of index $p$, and $G$ without an
$\mathcal{A}_1$-subgroup of index $p$ in Theorem \ref{d=3-1},
\ref{d=3-2}, \ref{d=3-3}, \ref{d=3-4} and \ref{d=3-5}.

\begin{thm}\label{d=3-1}
$d(G)=2$ and $c(G)= 4$ if and only if  $G$ is isomorphic to one of
the following pairwise non-isomorphic groups:
\begin{enumerate}
\rr{F1} $\lg a,b\di a^{16}=b^{2^m}=1, [a,b]=a^{-2}\rg$; where
$|G|=2^{m+4}$, $\Phi(G)=\lg a^2, b^2\rg\cong C_8\times C_{2^{m-1}}$
if $m>1$, $\Phi(G)\cong C_8$ if $m=1$, $G'=\lg a^2\rg\cong C_8$,
$Z(G)=\lg a^8,b^2\rg\cong C_2\times C_{2^{m-1}}$ if $m>1$,
$Z(G)\cong C_2$ if $m=1$.

\rr{F2} $\lg a,b\di a^{16}=b^{2^m}=1, [a,b]=a^{6}\rg$; where
$|G|=2^{m+4}$, $\Phi(G)=\lg a^2, b^2\rg\cong C_8\times C_{2^{m-1}}$
if $m>1$, $\Phi(G)\cong C_8$ if $m=1$, $G'=\lg a^2\rg\cong C_8$,
$Z(G)=\lg a^8,b^2\rg\cong C_2\times C_{2^{m-1}}$ if $m>1$,
$Z(G)\cong C_2$ if $m=1$.

\rr{F3} $\lg a,b\di a^{16}=1, b^{2^m}=a^{8}, [a,b]=a^{-2}\rg$; where
$|G|=2^{m+4}$, $\Phi(G)=\lg a^2, b^2\rg\cong C_2^3$ if $m =1$,
$\Phi(G)\cong C_2^3\times C_2$ if $m =2$, $\Phi(G)\cong
C_{2^{m}}\times C_{2^2}$ if $m>2$, $G'=\lg a^2\rg\cong C_8$ and
$Z(G)=\lg b^2\rg \cong C_{2^{m}}$.

\rr{F4} $\lg a_1, b; a_2\di
a_1^{9}=a_2^{9}=b^{3^m}=1,[a_1,b]=a_2,[a_2,b]=a_1^{-3}a_2^{3t},[a_1,a_2]=1\rg$,
where $t=1,2$; moreover, $|G|=3^{m+4}$, $\Phi(G)=\lg a_2,a_1^3,
b^3\rg\cong C_3\times C_{9}\times C_{3^{m-1}}$ if $m>1$,
$\Phi(G)\cong C_3\times C_{9}$ if $m=1$, $G'=\lg a_2,a_1^3\rg\cong
C_3\times C_{9}$, $Z(G)=\lg a_2^3,b^3\rg\cong C_3\times C_{3^{m-1}}$
if $m>1$, $Z(G)\cong C_3$ if $m=1$.

\rr{F5} $\lg a_1, b; a_2\di a_1^{9}=a_2^{9}=1,
b^{3^m}=a_2^{-3},[a_1,b]=a_2,[a_2,b]=a_1^{-3}a_2^{-3},[a_1,a_2]=1\rg$;
where $|G|=3^{m+4}$, $\Phi(G)=\lg a_2,a_1^3, b^3\rg\cong C_9\times
C_3$ if $m=1$, $\Phi(G)\cong C_{3^m}\times C_3\times C_3$ if $m>1$,
$G'=\lg a_2,a_1^3\rg\cong C_9\times C_3$ and $Z(G)=\lg b^3\rg \cong
C_{3^m}$.

\rr{F6} $\lg a_1,b; a_2,a_3,a_4\di
a_i^p=b^{p^m}=1,[a_j,b]=a_{j+1},[a_4,b]=1,[a_i,a_j]=1\rg,$ where
$p\ge 5$, $1\leq i\leq 4$, $1\leq j\leq 3$; moreover, $|G|=p^{m+4}$,
$\Phi(G)=\lg a_2,a_3,a_4, b^p\rg \cong C_p^3\times C_{p^{m-1}}$ if
$m>1$, $\Phi(G)\cong C_p^3$ if $m=1$, $G'=\lg a_2,a_3,a_4\rg\cong
C_p^3$, $Z(G)=\lg a_4,b^p\rg\cong C_p\times C_{p^{m-1}}$ if $m>1$,
$Z(G)\cong C_p$ if $m=1$.

\rr{F7} $\lg a_1,b;a_2,a_3\di
a_i^p=b^{p^{m+1}}=1,[a_j,b]=a_{j+1},[a_3,b]=b^{p^m},[a_i,a_j]=1\rg,$
where $p\ge 5$, $1\leq i\leq 3$, $1\leq j\leq 2$; moreover,
$|G|=p^{m+4}$, $\Phi(G)=\lg a_2,a_3, b^p\rg\cong C_p^2\times
C_{p^m}$, $G'=\lg a_2,a_3,b^{p^m}\rg\cong C_p^3$ and  $Z(G)=\lg
b^p\rg\cong C_{p^m}$.

\rr{F8} $\lg a_1,b;a_2,a_3\di
a_1^{p^2}=a_i^p=b^{p^m}=1,[a_j,b]=a_{j+1},[a_3,b]=a_1^{tp},[a_i,a_j]=1\rg,$
where $2\leq i\leq 3,\ 1\leq j\leq 2$, and
$t=t_1,t_2,\dots,t_{(3,p-1)}$, where $p\ge 5$,
$t_1,t_2,\dots,t_{(3,p-1)}$ are the coset representatives of the
subgroup $(F_p^*)^3$ in $F_p^*$. Moreover, $|G|=p^{m+4}$,
$\Phi(G)=\lg a_2,a_3,a_1^p, b^p\rg\cong C_p^3\times C_{p^{m-1}}$ if
$m>1$, $\Phi(G)\cong C_p^3$ if $m=1$, $G'=\lg a_2,a_3,a_1^p\rg\cong
C_p^3$, $Z(G)=\lg a_1^p,b^p\rg\cong C_p\times C_{p^{m-1}}$ if $m>1$,
$Z(G)\cong C_p$ if $m=1$.
\end{enumerate}
Moreover, $(\mu_0,\mu_1,\mu_2)=(1,0,p)$ and $\alpha_1(G)=p^2$.
\end{thm}
\demo  Let $A$ be an abelian subgroup of index $p$ and $B$ be a
non-abelian subgroup of index $p$. By Lemma \ref{alj1}, $G_3=B'$,
$G_4=B_3$, $|G_3/G_4|=p$ and $|G'|=p^{c(G)-1}=p^3$. Since $c(G)= 4$,
$c(B)\ge 3$. Hence Lemma \ref{A_2-property} gives that $d(B)=2$ and
$c(B)=3$. By arbitrariness of $B$, all non-abelian subgroups of $G$
are generated by two elements. Such groups were classified by
\cite{Alj}. It follows from \cite[Main Theorem]{Alj} that $G$ is
either a $p$-groups of maximal class or one of the groups listed in
\cite[Theorem 3.12-3.13]{Alj}. Since $c(G)=4$, they are the groups
(F1)--(F8).

Since $d(G)=2$, $G$ has $1+p$ maximal subgroups. Hence
$(\mu_0,\mu_1,\mu_2)=(1,0,p)$. Since $d(H)=2$ for any non-abelian
maximal subgroup $H$, $\alpha_1(H)=p$. By Hall's enumeration
principle, $\alpha_1(G)=\sum_{H\in \Gamma_1}
\alpha_1(H)=p\mu_2=p^2.$\qed

\begin{thm}\label{d=3-2}
$d(G)=2$ and $c(G)=3$ if and only if  $G$ is isomorphic to one of
the following pairwise non-isomorphic groups:
\begin{enumerate}
\rr{Gi} $G'\cong C_{p^2}$.
\begin{enumerate}

\rr{G1}  $\langle a, b; c \di a^{8}=1, c^2=a^{4}=b^4,
[a,b]=c,[c,a]=1,[c,b]=c^2 \rangle$; where $|G|=2^{6}$, $\Phi(G)=\lg
a^2,
b^2,c\rg\cong C_2\times C_2\times C_4$, $G'=\lg c\rg$, $Z(G)=\lg ca^2, b^2\rg\cong C_2\times C_4$.

\rr{G2} $\langle a, b; c \di a^{8}=b^{4}=1, c^2=a^{4},
[a,b]=c,[c,a]=1,[c,b]=c^{2} \rangle$; where $|G|=2^{6}$,
$\Phi(G)=\lg a^2,
b^2,c\rg\cong C_2\times C_2\times C_4$, $G'=\lg c\rg$, $Z(G)=\lg ca^2, b^2\rg\cong C_2^2$.

\rr{G3} $\langle a, b; c \di a^{8}=b^{4}=1, c^2=a^{4},
[a,b]=c,[c,a]=c^2,[c,b]=1 \rangle$; where
$|G|=2^{6}$, $\Phi(G)=\lg a^2,
b^2,c\rg\cong C_{4}\times C_2^2$, $G'=\lg c\rg$ and $Z(G)=\lg a^2, cb^2\rg \cong C_{4}\times C_2$.

\rr{G4} $\langle a, b; c \di a^{2^{n+1}}=b^{4}=1, c^2=a^{2^n},
[a,b]=c,[c,a]=c^2,[c,b]=1 \rangle$, where $n>2$; moreover,
$|G|=2^{n+4}$, $\Phi(G)=\lg a^2,
b^2,c\rg\cong C_{2^n}\times C_2^2$, $G'=\lg c\rg$ and $Z(G)=\lg a^2, cb^2\rg \cong C_{2^n}\times C_2$.

\rr{G5} $\langle a, b; c \di a^{2^{n}}=b^{8}=1, c^2=b^{4},
[a,b]=c,[c,a]=c^{2},[c,b]=1 \rangle$, where $n>2$; moreover,
$|G|=2^{n+4}$, $\Phi(G)=\lg a^2,
b^2,c\rg\cong C_{2^{n-1}}\times C_2\times C_4$, $G'=\lg c\rg$ and $Z(G)=\lg a^2, cb^2\rg\cong C_{2^{n-1}}\times C_2$.

\rr{G6} $\langle a, b; c \di a^{2^{n}}=b^{4}=c^{4}=1,
[a,b]=c,[c,a]=c^2,[c,b]=1 \rangle$, where $n>2$. Moreover,
$|G|=2^{n+4}$, $\Phi(G)=\lg a^2,
b^2,c\rg\cong C_{2^{n-1}}\times C_2\times C_4$, $G'=\lg c\rg$ and $Z(G)=\lg a^2, cb^2\rg\cong C_{2^{n-1}}\times C_4$.

\end{enumerate}
\rr{Gii} $G'\cong C_p^2$.
\begin{enumerate}

\rr{G7}  $\langle a, b; c \di a^{p^{3}}=b^{p^2}=c^p=1,
[a,b]=c,[c,a]=1,[c,b]=a^{\nu p^2} \rangle$, where $p>2$ and $\nu=1$
 or a fixed quadratic non-residue modulo $p$; moreover, $|G|=p^{6}$, $\Phi(G)=\lg a^p,
b^p,c\rg\cong C_{p^2}\times C_p^{2}$, $G'=\lg a^{p^2},c\rg$ and $Z(G)=\lg a^p, b^p\rg\cong C_{p^2}\times C_{p}$.

\rr{G8}  $\langle a, b; c,d \di a^{p^{2}}=b^{p^{2}}=c^p=d^p=1,
[a,b]=c,[c,a]=d,[c,b]=1,[d,a]=[d,b]=1 \rangle$, where $p>2$. Moreover, $|G|=p^{6}$, $\Phi(G)=\lg a^p,
b^p,c,d\rg\cong C_{p}^4$, $G'=\lg c,d\rg$ and $Z(G)=\lg a^p, b^p,d\rg\cong C_{p}^3$.

\rr{G9}  $\langle a, b; c \di a^{p^{3}}=b^{p^2}=c^p=1,
[a,b]=c,[c,a]=a^{p^2},[c,b]=1 \rangle$, where $p>2$;
moreover, $|G|=p^{6}$, $\Phi(G)=\lg a^p,
b^p,c\rg\cong C_{p^2}\times C_p^{2}$, $G'=\lg a^{p^2},c\rg$ and $Z(G)=\lg a^p, b^p\rg\cong C_{p^2}\times C_{p}$.

\rr{G10}  $\langle a, b; c \di a^{p^{n+1}}=b^{p^2}=c^p=1,
[a,b]=c,[c,a]=a^{p^n},[c,b]=1 \rangle$, where $p>2$ and $n>2$;
moreover, $|G|=p^{n+4}$, $\Phi(G)=\lg a^p,
b^p,c\rg\cong C_{p^{n}}\times C_p^2$, $G'=\lg a^{p^n},c\rg$ and $Z(G)=\lg a^p, b^p\rg\cong C_{p^{n}}\times C_p$.

\rr{G11}  $\langle a, b; c \di a^{p^{n}}=b^{p^{3}}=c^p=1,
[a,b]=c,[c,a]=b^{\nu p^2},[c,b]=1 \rangle$, where $p>2$, $n>2$, and
$\nu=1$ or or a fixed quadratic non-residue modulo $p$; moreover,
$|G|=p^{n+4}$, $\Phi(G)=\lg a^p,
b^p,c\rg\cong C_{p^{n-1}}\times C_{p^2}\times C_p$, $G'=\lg b^{p^2},c\rg$ and $Z(G)=\lg a^p, b^p\rg\cong C_{p^{n-1}}\times C_{p^2}$.

\rr{G12}  $\langle a, b; c,d \di a^{p^{n}}=b^{p^{2}}=c^p=d^p=1,
[a,b]=c,[c,a]=d,[c,b]=1,[d,a]=[d,b]=1 \rangle$, where $p>2$ and
$n>2$. Moreover, $|G|=p^{n+4}$, $\Phi(G)=\lg a^p,
b^p,c,d\rg\cong C_{p^{n-1}}\times C_p^3$, $G'=\lg c,d\rg$ and $Z(G)=\lg a^p, b^p,d\rg\cong C_{p^{n-1}}\times C_p^2$.
\end{enumerate}
\end{enumerate}
Moreover, $(\mu_0,\mu_1,\mu_2)=(1,0,p)$ and $\alpha_1(G)=p^3$.
\end{thm}

\demo Let $A$ be an abelian subgroup of index $p$ and $B$ be a
non-abelian subgroup of index $p$. By Lemma \ref{alj1}, $G_3=B'$,
$G_4=B_3$, $|G_3/G_4|=p$ and $|G'|=p^{c(G)-1}$.

Since $c(G)=3$, $|B'|=|G_3|=p$ and $|G'|=p^2$. Let $\bar{G}=G/G_3$.
Then $|\bar{G}'|=p$ and $d(\bar{G})=2$. Hence $\bar{G}$ is an
$\mathcal{A}_1$-group by Lemma \ref{minimal non-abelian equivalent
conditions}. If $\bar{G}$ is metacyclic, then $G$ is also metacyclic
by Lemma \ref{metacyclic} and hence $G$ is an $\mathcal{A}_2$-group
by Lemma \ref{metacyclic An}. This contradicts that $G$ is an
$\mathcal{A}_3$-group. Thus $\bar{G}$ is a non-metacyclic
$\mathcal{A}_1$-group. If $G'\cong C_{p^2}$, then $G$ is one of the
groups listed in \cite[Theorem 3.5]{AHZ}. By \cite[Theorem 3.1 \&
3.6 ]{AHZ} we get groups (G1)--(G6). If $G'\cong C_p^2$, then $G$ is
one of the groups listed in \cite[Theorem 4.6]{AHZ}. By
\cite[Theorem 4.1 \& 4.7]{AHZ} we get groups (G7)--(G12).

\begin{table}[h]
  \centering
{
\begin{tabular}{c|c||c|c}
\hline
Groups & Groups in \cite[Theorem 3.5]{AHZ}&Groups & Groups in  \cite[Theorem 4.6]{AHZ} \\
\hline
(G1)&(C3) & (G7)&(G1) where $m=2$ and $p>2$\\
(G2)&(C4) & (G8)&(G2) where $m=2$ and $p>2$\\
(G3)&(C5) & (G9)&(G3) where $m=2$ and $p>2$\\
(G4)&(E2) where $p=m=2$ & (G10)&(J2) where $m=2$ and $p>2$\\
(G5)&(E5) where $p=m=2$ & (G11)&(J4) where $m=2$ and $p>2$\\
(G6)&(E8) where $p=m=2$ & (G12)&(J6) where $m=2$ and $p>2$\\
\hline
\end{tabular}
\caption{The correspondence from Theorem \ref{d=3-2} to
\cite[Theorem 3.5 \& 4.6]{AHZ}} \label{table8}}
\end{table}

Now we calculate the $(\mu_0,\mu_1,\mu_2)$ and $\alpha_1(G)$. Since
$d(G)=2$, $G$ has $1+p$ maximal subgroups. Hence
$(\mu_0,\mu_1,\mu_2)=(1,0,p)$. Since $|H'|=p$ for any non-abelian
maximal subgroup $H$, by Lemma \ref{A_2-property} (7),
$\alpha_1(H)=p^2$. By Hall's enumeration principle,
$\alpha_1(G)=\sum_{H\in \Gamma_1} \alpha_1(H)=p^2\mu_2=p^3.$\qed

\begin{thm}\label{d=3-3}
$d(G)=3$ and $\Phi(G)\le Z(G)$ if and only if $G$ is isomorphic to
one of the following pairwise non-isomorphic groups:
\begin{enumerate}
\rr{Hi} $\Phi(G)<Z(G)$, $G'\cong C_p$ and $c(G)=2$. Moreover, $(\mu_0,\mu_1,\mu_2)=(1+p,0,p^2)$ and $\alpha_1(G)=p^4$.
\begin{enumerate}
\rr{H1} $\langle a, b, c \di a^{p^{n+1}}=b^{p^m}=c^{p^2}=1,\
[a,b]=a^{p^{n}},\ [c,a]=[c,b]=1\rangle\cong M_p(n+1,m)\times
C_{p^2}$, where $\min\{n,m\}\ge 2$;  moreover, $|G|=p^{n+m+3}$,
$\Phi(G)=\lg a^p, b^p,c^p\rg\cong C_{p^{n}}\times C_{p^{m-1}}\times
C_{p}$, $G'=\lg a^{p^n}\rg$ and  $Z(G)=\lg a^p, b^p,c\rg\cong
C_{p^{n}}\times C_{p^{m-1}}\times C_{p^2}$.

\rr{H2} $\langle a, b, c; d \di a^{p^n}=b^{p^m}=c^{p^2}=d^p=1,\
[a,b]=d,\ [c,a]=[c,b]=1\rangle\cong M_p(n,m,1)\times C_{p^2}$, where
$n\ge m\ge 2$; moreover, $|G|=p^{n+m+3}$, $\Phi(G)=\lg a^p,
b^p,c^p,d\rg\cong C_{p^{n-1}}\times C_{p^{m-1}}\times C_{p}\times
C_{p}$, $G'=\lg d\rg$ and  $Z(G)=\lg a^p, b^p,c,d\rg\cong
C_{p^{n-1}}\times C_{p^{m-1}}\times C_{p^2}\times C_{p}$.

\rr{H3} $\langle a, b, c \di a^{p^n}=b^{p^m}=c^{p^{3}}=1,\
[a,b]=c^{p^2},\ [c,a]=[c,b]=1 \rangle\cong M_p(n,m,1)\ast C_{p^3}$,
where $n\ge m\ge 2$. Moreover,  $|G|=p^{n+m+3}$, $\Phi(G)=\lg a^p,
b^p,c^p\rg\cong C_{p^{n-1}}\times C_{p^{m-1}}\times C_{p^2}$,
$G'=\lg c^{p^2}\rg$ and $Z(G)=\lg a^p, b^p,c\rg\cong
C_{p^{n-1}}\times C_{p^{m-1}}\times C_{p^3}$.
\end{enumerate}

\rr{Hii} $\Phi(G)=Z(G)$, $G'\cong C_p^2$ and $c(G)=2$. Moreover, $(\mu_0,\mu_1,\mu_2)=(1,0,p+p^2)$ and $\alpha_1(G)=p^4+p^3$.

\begin{enumerate}

\rr{H4} $\langle a, b, c \di a^{p^{l}}=b^{p^{2}}=c^{p^{2}}=1,
[b,c]=1,[c,a]=c^{p},[a,b]=b^{-p}\rangle$, where $p>2$ and $l\ge 2$;
moreover, $|G|=p^{l+4}$, $\Phi(G)=Z(G)=\lg a^p, b^p,c^p\rg\cong
C_{p^{l-1}}\times C_{p}^2$, $G'=\lg b^p,c^{p}\rg$.

\rr{H5} $\langle a, b, c \di a^{p^{l}}=b^{p^{3}}=c^{p^3}=1,
[b,c]=1,[c,a]=b^{p^{2}}c^{tp^{2}},[a,b]=b^{-tp^2}c^{\nu
p^2}\rangle$, where $p>2$, $l\ge 2$, $\nu=1$ or a fixed quadratic
non-residue modulo $p$ such that $-\nu\not\in (F_p^*)^2$ and
$t\in\{0,1,\dots,\frac{p-1}{2}\}$; moreover, $|G|=p^{l+6}$,
$\Phi(G)=Z(G)=\lg a^p, b^p,c^p\rg\cong C_{p^{l-1}}\times
C_{p^2}\times C_{p^2}$, $G'=\lg b^{p^2},c^{p^2}\rg$.

\rr{H6}  $\langle a, b, c \di a^{2^{l}}=b^{4}=c^{4}=1,
[b,c]=1,[c,a]=c^{2},[a,b]=b^{2}\rangle$; moreover, $|G|=2^{l+4}$,
$\Phi(G)=Z(G)=\lg a^2, b^2,c^2\rg\cong C_{2^{l-1}}\times C_{2}\times
C_{2}$ if $l>1$, $G'=\lg b^{2},c^2\rg$ , $\Phi(G)=Z(G)\cong
C_{2}\times C_{2}$ if\  $l=1$.

\rr{H7}  $\langle a, b, c \di a^{2^{l}}=b^{8}=c^{8}=1,
[b,c]=1,[c,a]=b^{4},[a,b]=b^{4}c^{4}\rangle$, where $l\ge 2$;
moreover, $|G|=2^{l+6}$, $\Phi(G)=Z(G)=\lg a^2, b^2,c^2\rg\cong
C_{2^{l-1}}\times C_{4}\times C_{4}$, $G'=\lg b^{4},c^4\rg$.

\rr{H8} $\langle a, b, c \di
a^{p^{l+1}}=b^{p^{3}}=c^{p^{2}}=1,[b,c]=1,[c,a]=b^{p^{2}},
[a,b]=a^{p^l} \rangle$, where $l\ge 2$; moreover, $|G|=p^{l+6}$,
$\Phi(G)=Z(G)=\lg a^p, b^p,c^p\rg\cong C_{p^{l}}\times C_{p^2}\times
C_{p}$, $G'=\lg a^{p^l},b^{p^2}\rg$.

\rr{H9}  $\langle a, b, c; x\di a^{p^{l}}=b^{p^{2}}=c^{p^{2}}=x^p=1,
[a,b]=c^{p},[a,c]=x,[b,c]=[x,a]=[x,b]=[x,c]=1 \rangle$, where $l\ge
2$; moreover, $|G|=p^{l+5}$, $\Phi(G)=Z(G)=\lg a^p,
b^p,c^p,x\rg\cong C_{p^{l-1}}\times C_{p}^3$, $G'=\lg c^{p},x\rg$.

\rr{H10}  $\langle a, b, c; x,y \di a^{p^{l}}=b^{p}=c^{p}=x^p=y^p=1,
[a,b]=x,[a,c]=y,[b,c]=[x,a]=[x,b]=[x,c]=[y,a]=[y,b]=[y,c]=1
\rangle$, where $l\ge 2$ if $p=2$; moreover, $|G|=p^{l+4}$,
$\Phi(G)=Z(G)=\lg a^p, x,y\rg\cong C_{p^{l-1}}\times C_{p}^2$,
$G'=\lg x,y\rg$.

\end{enumerate}
\end{enumerate}
\end{thm}

\demo By Lemma \ref{Ber}, $|G'|\le p^2$. If $|G'|=p$, then $G$ is
one of the groups listed in \cite[Theorem 3.1]{ALQZ}. Since $G$ has
no $\mathcal{A}_1$-subgroup of index $p$, they are the groups
(H1)--(H3). If $|G'|=p^2$, then $G$ is one of the groups listed in
\cite[Theorem 4.8]{ALQZ}. By checking \cite[Table 4]{ALQZ} we get
groups (H4)--(H10).

\begin{table}[h]
  \centering
{\small
\begin{tabular}{c|c||c|c}
\hline
Groups&Groups in \cite[Theorem 4.8]{ALQZ}& Groups & Groups in \cite[Theorem 4.8]{ALQZ}\\
\hline
(H4)&(A1) where $l\ge 2$ and $m=1$ & (H8)&(A10) where $l\ge 2=m=n$\\
(H5)&(A3) where $l\ge m=2$ and $-\nu\not\in F_p^2$ & (H9)&(B2) where $m=2$ and $n=1$\\
(H6)&(A5) where $m=1$ & (H10)&(C) where $m=n=1$\\
(H7)&(A6) where $l\ge m=2$ and $-\nu\not\in (F_p^*)^2$& &\\
\hline
\end{tabular}
\caption{The correspondence from Theorem \ref{d=3-3} to
\cite[Theorem 4.8]{ALQZ}} \label{table9}}
\end{table}

Now we calculate the $(\mu_0,\mu_1,\mu_2)$ and $\alpha_1(G)$. If $G$
is one of the groups (H1)--(H3), then $\mu_0=1+p$. Since $\mu_1=0$,
$\mu_2(G)=p^2$. By Lemma \ref{alpha of d=3},
$\alpha_1(G)=p^2\mu_2=p^4.$

If $G$ is one of the groups {\rm(H4)--(H10)}, then $d(G)=3$ and $G'\cong C_p^2$. Since $\Phi(G)=Z(G)$,
$\mu_0=1$. Hence $(\mu_0,\mu_1,\mu_2)=(1,0,p+p^2)$. By Lemma \ref{alpha of d=3},
$\alpha_1(G)=p^2\mu_2=p^4+p^3.$\qed

\begin{thm}\label{d=3-4}
$d(G)=3$ and $\Phi(G)\not \le Z(G)$ if and only if  $G$ is
isomorphic to one of the following pairwise non-isomorphic groups:
\begin{enumerate}

\rr{Ii} $|G'|=p^2$ and $c(G)=3$. In this case, $Z(G)\not\le\Phi(G)$, $(\mu_0,\mu_1,\mu_2)=(1, 0, p^2+p)$ and $\alpha_1(G)=p^3$.

\begin{enumerate}

\rr{I1} $\lg a,b,x\di
a^8=b^{2^m}=x^2=1,[a,b]=a^{-2},[x,a]=[x,b]=1\rg=\lg a,b\rg\times\lg
x\rg$; where  $|G|=2^{m+4}$, $\Phi(G)=\lg a^2, b^2\rg\cong
C_{2^{m-1}}\times C_4$ if $m>1$, $\Phi(G)\cong C_4$ if $m=1$,
$G'=\lg a^{2}\rg\cong C_4$, $Z(G)=\lg a^4, b^2,x\rg\cong
C_{2^{m-1}}\times C_2^2$ if $m>1$, $Z(G)\cong C_2^2$ if $m=1$.

\rr{I2} $\lg a,b,x\di
a^8=b^{2^m}=1,x^2=a^4,[a,b]=a^{-2},[x,a]=[x,b]=1\rg=\lg
a,b\rg\ast\lg x\rg$; where $|G|=2^{m+4}$, $\Phi(G)=\lg a^2,
b^2\rg\cong C_{2^{m-1}}\times C_4$ if $m>1$, $\Phi(G)\cong C_4$ if
$m=1$, $G'=\lg a^{2}\rg\cong C_4$, $Z(G)=\lg  b^2,x\rg\cong
C_{2^{m-1}}\times C_4$ if $m>1$, $Z(G)\cong C_4$ if $m=1$.

\rr{I3} $\lg a,b,x\di
a^8=b^{2^m}=x^2=1,[a,b]=a^{2},[x,a]=[x,b]=1\rg=\lg a,b\rg\times\lg
x\rg$; where  $|G|=2^{m+4}$, $\Phi(G)=\lg a^2, b^2\rg\cong
C_{2^{m-1}}\times C_4$ if $m>1$, $\Phi(G)\cong C_4$ if $m=1$,
$G'=\lg a^{2}\rg\cong C_4$, $Z(G)=\lg a^4, b^2,x\rg\cong
C_{2^{m-1}}\times C_2^2$ if $m>1$, $Z(G)\cong C_2^2$ if $m=1$.

\rr{I4} $\lg a,b,x\di
a^8=x^2=1,b^{2^m}=a^4,[a,b]=a^{-2},[x,a]=[x,b]=1\rg=\lg
a,b\rg\times\lg x\rg$; where $|G|=2^{m+4}$, $\Phi(G)=\lg a^2, b^2\rg
\cong C_4$ if $m=1$, $\Phi(G)=\lg a^2, b^2\rg \cong C_{2^{m}}\times
C_2$ if $m>1$, $G'=\lg a^{2}\rg\cong C_4$ and $Z(G)=\lg
b^2,x\rg\cong C_{2^{m}}\times C_2$.

\rr{I5} $\lg a_1, b, x; a_2, a_3 \di
a_1^{p}=a_2^p=a_3^p=b^{p^m}=x^p=1,[a_1,b]=a_2,[a_2,b]=a_3,[a_2,a_1]=[a_3,a_1]=[a_3,a_2]=[x,a_1]=[x,b]=1\rg=\lg
a_1,b\rg\times\lg x\rg$, where $p>3$ for $m=1$ and $p>2$; moreover,
$|G|=p^{m+4}$, $\Phi(G)=\lg a_2,a_3, b^p\rg\cong C_{p^{m-1}}\times
C_p^2$ if $m>1$, $\Phi(G)\cong C_p^2$ if $m=1$, $G'=\lg
a_2,a_3\rg\cong C_p^2$, $Z(G)=\lg a_3,b^p,x\rg\cong
C_{p^{m-1}}\times C_p^2$ if $m>1$, $Z(G)\cong C_p^2$ if $m=1$.

\rr{I6} $\lg a_1,x, b; a_2 \di
a_1^{p}=a_2^p=b^{p^m}=x^{p^2}=1,[a_1,b]=a_2,[a_2,b]=x^p,[a_2,a_1]=[x,a_1]=[x,b]=1\rg=\lg
a_1,b\rg\ast\lg x\rg$, where $p>2$; moreover, $|G|=p^{m+4}$,
$\Phi(G)=\lg a_2, b^p,x^p\rg\cong C_{p^{m-1}}\times C_p^2$ if $m>1$,
$\Phi(G)\cong  C_p^2$ if $m=1$, $G'=\lg a_2,x^p\rg \cong C_p^2$,
$Z(G)=\lg x,b^p\rg\cong C_{p^{m-1}}\times C_{p^2}$ if $m>1$,
$Z(G)\cong C_{p^2}$ if $m=1$.

\rr{I7} $\lg a_1, x, b; a_2\di
a_1^{p}=a_2^p=b^{p^{m+1}}=x^p=1,[a_1,b]=a_2,[a_2,b]=b^{p^m},[a_2,a_1]=[x,a_1]=[x,b]=1\rg=\lg
a_1,b\rg\times\lg x\rg$, where $p>2$; moreover, $|G|=p^{m+4}$,
$\Phi(G)=\lg a_2, b^p\rg\cong C_{p^{m}}\times C_{p}$, $G'=\lg
a_2,b^{p^m}\rg\cong
 C_{p}^2$ and $Z(G)=\lg x,b^p\rg\cong C_{p^{m}}\times
C_{p}$.

\rr{I8} $\lg a_1, x, b; a_2 \di
a_1^{p^2}=a_2^p=b^{p^m}=x^p=1,[a_1,b]=a_2,[a_2,b]=a_1^{\nu
p},[a_2,a_1]=[x,a_1]=[x,b]=1\rg=\lg a_1,b\rg\times\lg x\rg$, where
$p>2$ and $\nu=1$ or a fixed quadratic non-residue modulo $p$;
Moreover, $|G|=p^{m+4}$, $\Phi(G)=\lg a_2,a_1^p, b^p\rg\cong
C_{p^{m-1}}\times C_p^2$ if $m>1$, $\Phi(G)\cong  C_p^2$ if $m=1$,
$G'=\lg a_2,a_1^p\rg \cong C_p^2$, $Z(G)=\lg a_1^p,b^p,x\rg\cong
C_{p^{m-1}}\times C_{p}^2$ if $m>1$, $Z(G)\cong C_p^2$ if $m=1$.

\rr{I9} $\lg a_1,x, b; a_2\di
a_1^9=a_2^3=x^3=1,b^3=a_1^3,[a_1,b]=a_2,[a_2,b]=a_1^{-3},[x,a_1]=[x,a_2]=1\rg=\lg
a_1,b\rg\times\lg x\rg$; where $|G|=3^{5}$, $\Phi(G)=G'=\lg
a_2,a_1^3\rg\cong C_3^2$ and $Z(G)=\lg a_1^3,x\rg \cong C_3^2$.
\end{enumerate}

\rr{Iii} $|G'|=p^3$  and $c(G)=3$. In this case, $Z(G)<\Phi(G)$, $(\mu_0,\mu_1,\mu_2)=(1, 0, p^2+p)$ and $\alpha_1(G)=p^3+p^2$.

\begin{enumerate}

\rr{I10} $\lg a,b,c\di
a^8=b^{2^{m+1}}=1,c^2=a^4b^{2^m},[a,b]=a^2,[c,a]=1,[c,b]=b^{2^m}\rg$;
where $|G|=2^{m+5}$, $\Phi(G)=\lg a^2, b^2\rg\cong C_{2^{m}}\times
C_{4}$, $G'=\lg a^2,b^{2^m}\rg\cong C_{4}\times C_{2}$ and $Z(G)=\lg
a^4,b^2\rg\cong C_{2^{m}}\times C_{2}$.

\rr{I11} $\lg a, b, c;d\di
a^{p^{n+1}}=b^{p}=c^{p^2}=d^p=1,[a,b]=d,[c,a]=a^{p^n},[d,a]=c^p,[c,b]=[d,b]=[d,c]=1\rg$,
where $p>2$. Moreover,  $|G|=p^{n+5}$, $\Phi(G)=\lg a^p,
c^p,d\rg\cong C_{p^{n}}\times C_{p}^2$, $G'=\lg a^{p^n},c^p,d\rg
\cong C_p^3$ and $Z(G)=\lg c^p,a^p\rg\cong C_{p^{n}}\times C_{p}$.
\end{enumerate}
\end{enumerate}
\end{thm}

\demo It is easy to see that $Z(G)\not\le \Phi(G)$ for groups
(I1)--(I9) and $Z(G)<\Phi(G)$ for (I10)--(I11). By checking maximal
subgroups of groups (I1)--(I11), we know that they are all
$\mathcal{A}_3$-groups and are pairwise non-isomorphic. Conversely,
in the following, we prove that $G$ is one of the groups
(I1)--(I11).

Let $A$ be an abelian subgroup  of index $p$ of $G$, and $B$ be a
non-abelian subgroup  of index $p$  of $G$. Then $A\cap B$ is an
abelian subgroup  of index $p$ of $B$. It follows that $Z(B)\le
A\cap B$ and hence $Z(B)\le Z(G)$.

We claim that there exists a non-abelian subgroup $B$  of index $p$
such that $d(B)=3$. If not, then all non-abelian proper subgroups of
$G$ are generated by two elements. By Lemma \ref{Dp3=>A_2},
$G\in\mathcal{A}_2$, a contradiction.

By Lemma \ref{A_2-property} (2), $B'\le C_p^2$ and $B'\le Z(B)\le
Z(G)$. If $G'=B'$, then we have $\Phi(G)\le Z(G)$, a contradiction.
Hence $G'>B'$. By Lemma \ref{G'}, $|G'|\le p|A'||B'|=p|B'|$. Thus
$|G'|=p|B'|$. Let $D$ be another non-abelian subgroup  of index $p$
of $G$ such that $d(D)=3$. Similarly we have $D'\le C_p^2$ and
$D'\le Z(D)\le Z(G)$. If $D'\not\le B'$, then $G'=B'D'\le Z(G)$ and
hence $\Phi(G)\le Z(G)$, a contradiction. Thus we have $D'\le B'$.
The same reason gives that $B'\le D'$ and hence $B'=D'$. By
arbitrariness of $D$, there exists a non-abelian subgroup $M$  of
index $p$ such that $d(M)=2$ and all non-abelian subgroups of $G/B'$
are generated by two elements.

\medskip

Case 1. $|G'|=p^2$

\medskip

By Lemma \ref{Ber}, $|G:Z(G)|=p^3$. Since $|G:\Phi(G)|=p^3$ and
$\Phi(G)\neq Z(G)$, there exists a non-abelian subgroup $M$  of
index $p$ such that $Z(G)\not\le M$. Let $x\in Z(G)\backslash M$.
Then $G=\lg M,x\rg$ and $G'=M'$. If $d(M)=3$, then, by Lemma
\ref{A_2-property} (2), $M'\le C_p^3$ and $M'\le Z(M)$. It follows
that $G'\le C_p^3$ and $G'\le Z(G)$. By calculation we get
$\Phi(G)\le Z(G)$, a contradiction. Hence $d(M)=2$. It follows that
$\Phi(M)=\Phi(G)$ and $M$ is one of the groups (1)--(7) in Lemma
\ref{A_2}.

\medskip

Subcase 1.1. $M$ is the group (1) in Lemma \ref{A_2}. That is,
$M=\lg a,b\di a^{8}=b^{2^m}=1, a^b=a^{-1}\rg.$

\medskip

Since $Z(M)=\lg a^4,b^2\rg$ and $x^2\in Z(G)\cap \Phi(G)=Z(M)$, we
may assume that $x^2=a^{4i}b^{2j}$. Thus $(xb^{-j})^2=a^{4i}$. If
$(2,j)=1$, then $|\lg a^2, xb^{-j}\rg|=8$. Since $G$ is an
$\mathcal{A}_3$-group, we have $|G|\le 2^5$ and $|M|\le 2^4$. Hence
$m=1$. It follows that $x^2=a^{4i}$. If $2\di j$, then, replacing
$x$ with $xb^{-j}$, we also have $x^2=a^{4i}$. Hence we get the
groups (I1) for $2\di i$ and (I2) for $(2,i)=1$, respectively.

\medskip

Subcase 1.2. $M$ is the group (2) in Lemma \ref{A_2}. That is,
$M=\lg a,b\di a^{8}=b^{2^m}=1, a^b=a^{3}\rg.$

\medskip

Since $Z(M)=\lg a^4,b^2\rg$ and $x^2\in Z(G)\cap \Phi(G)=Z(M)$, we
may assume that $x^2=a^{4i}b^{2j}$. Thus $(xb^{-j})^2=a^{4i}$. If
$(2,j)=1$, then $|\lg a^2, xb^{-j}\rg|=8$. Since $G$ is an
$\mathcal{A}_3$-group, we have $|G|\le 2^5$ and $|M|\le 2^4$. Hence
$m=1$. It follows that $x^2=a^{4i}$. If $2\di j$, then, replacing
$x$ with $xb^{-j}$, we also have $x^2=a^{4i}$. If $x^2=1$, then we
get the group (I3). If $x^2=a^4$, then, replacing $a$ with $ax$, we
get the group (I2).

\medskip

Subcase 1.3. $M$ is the group (3) in Lemma \ref{A_2}. That is,
$M=\lg a,b\di a^{8}=1, b^{2^m}=a^{4}, a^b=a^{-1}\rg.$

\medskip

Since $Z(M)=\lg a^4,b^2\rg$ and $x^2\in Z(G)\cap \Phi(G)=Z(M)$, we
may assume that $x^2=a^{4i}b^{2j}$. Thus $(xb^{-j})^2=a^{4i}$. If
$(2,j)=1$, then $|\lg a^2, xb^{-j}\rg|=8$. Since $G$ is an
$\mathcal{A}_3$-group, we have $|G|\le 2^5$ and $|M|\le 2^4$. Hence
$m=1$. It follows that $x^2=a^{4(i+1)}$. If $2\di j$, then,
replacing $x$ with $xb^{-j}$, we also have $x^2=a^{4i}$. If $x^2=1$,
then we get the group (I4). In the following, we may assume
that $x^2=a^4$. If $m=1$, then, replacing $b$ with $xb$, we get the
group (I2). If $m\ge 2$, then, replacing $x$ with $xb^{2^{m-1}}$, we
get the group (I4).

\medskip

Subcase 1.4. $M$ is the group (4) in Lemma \ref{A_2}. That is,
$M=\lg a_1,b; a_2,a_3\di
a_1^p=a_2^p=a_3^p=b^{p^m}=1,[a_1,b]=a_{2},[a_2,b]=a_3,[a_3,b]=1,[a_i,a_j]=1\rg,$
where $p>3$ for $m=1$, $p>2$ and $1\leq i,j\leq 3.$

\medskip

Since $Z(M)=\lg a_3,b^p\rg$ and $x^p\in Z(G)\cap \Phi(G)=Z(M)$, we
may assume that $x^p=a_3^{i}b^{jp}$. Thus $(xb^{-j})^p=a_3^{i}$. If
$(j,p)=1$, then $|\lg a_2, xb^{-j}\rg|=p^3$. Since $G$ is an
$\mathcal{A}_3$-group, we have $|G|\le p^5$ and $|M|\le p^4$. Hence
$m=1$. It follows that $x^p=a_3^{i}$. If $p\di j$, then, replacing
$x$ with $xb^{-j}$, we also have $x^p=a_3^{i}$. If $p\di i$, then we
get the group (I5). If $(i,p)=1$, then, replacing $x$ with
$x^{i^{-1}}$, we get the group (I6).

\medskip

Subcase 1.5. $M$ is the group (5) in Lemma \ref{A_2}. That is,
$M=\lg a_1,b;a_2\di
a_1^p=a_2^p=b^{p^{m+1}}=1,[a_1,b]=a_{2},[a_2,b]=b^{p^m},[a_1,a_2]=1\rg,\
{\rm where}\ p>2.$

\medskip

Since $Z(M)=\lg b^p\rg$ and $x^p\in Z(G)\cap \Phi(G)=Z(M)$, we may
assume that $x^p=b^{jp}$. Thus $(x^{-j^{-1}}b)^p=1$. If $p\di j$,
then, replacing $x$ with $x^{-j^{-1}}b$, we get the group (I7). If
$(j,p)=1$, then $|\lg a_2, x^{-j^{-1}}b\rg|=p^3$. Since $G$ is an
$\mathcal{A}_3$-group, we have $|G|\le p^5$ and $|M|\le p^4$. Hence
$m=1$. By replacing $b$ with $x^{-j^{-1}}b$ we get the group (I6).

\medskip

Subcase 1.6. $M$ is the group (6) in Lemma \ref{A_2}. That is,
$M=\lg a_1,b;a_2\di
a_1^{p^2}=a_2^p=b^{p^m}=1,[a_1,b]=a_{2},[a_2,b]=a_1^{\nu
p},[a_1,a_2]=1\rg,$ where $p>2$ and $\nu=1$ or a fixed quadratic
non-residue modulo $p$.

\medskip

Since $Z(M)=\lg a_1^p,b^p\rg$ and $x^p\in Z(G)\cap \Phi(G)=Z(M)$, we
may assume that $x^p=a_1^{ip}b^{jp}$. Thus
$(x^{-j^{-1}}b)^p=a_1^{ip}$. If $(j,p)=1$, then $|\lg a_2,
x^{-j^{-1}}b\rg|=p^3$. Since $G$ is an $\mathcal{A}_3$-group, we
have $|G|\le p^5$ and $|M|\le p^4$. Hence $m=1$. It follows that
$x^p=a_1^{ip}$. If $p\di j$, then, replacing $x$ with $xb^{-j}$, we
also have $x^p=a_1^{ip}$. If $p\di i$, then we get the group (I8).
If $(i,p)=1$, then, replacing $a_1$ and $x$ with $a_1x^{-i^{-1}}$
and $x^{i^{-1}\nu}$ respectively, we get the group (I6).

\medskip

Subcase 1.7. $M$ is the group (7) in Lemma \ref{A_2}. That is,
$M=\lg a_1,b;a_2\di
a_1^9=a_2^3=1,b^3=a_1^3,[a_1,b]=a_2,[a_2,b]=a_1^{-3}\rg.$

\medskip

If $x^3=1$, then $G$ is the group (I9). If $x^3\ne 1$, then we may
assume that $a_1^3=b^3=x^{-3}$. By replacing $a_1$ and $b$ with
$a_1x$ and $bx$ respectively, we get the group (I6).

\medskip

Case 2. $|G'|=p^3$

\medskip

In this case, $B'\cong C_p^2$ and $Z(B)=\Phi(B)$. By Lemma
\ref{Ber}, $|G:Z(G)|=p^4$. Since $|G:Z(B)|=p^4$ and $Z(B)\le Z(G)$,
we have $Z(G)=Z(B)=\Phi(B)<\Phi(G)$. By Lemma \ref{Dp3=>A_2},
$G/B'\in\mathcal{A}_2$. Hence $G/B'$ is one of the groups (8)--(12) in
Lemma \ref{A_2}.

\medskip

Subcase 2.1. $G/B'$ is the group (8) or (11) in Lemma
\ref{A_2}. That is, $G/B'=\langle \bar{a}, \bar{b}, \bar{c} \di
\bar{a}^4=\bar{c}^{2}=1, \bar{b}^2=\bar{a}^2=[\bar{a},\bar{b}],
[\bar{c},\bar{a}]=[\bar{c},\bar{b}]=1 \rangle \cong Q_8 \times
C_{2}$ or $G/B'=\langle \bar{a}, \bar{b}, \bar{c} \di \bar{a}^4=1,
\bar{b}^2=\bar{c}^{2}=\bar{a}^2=[\bar{a},\bar{b}],
[\bar{c},\bar{a}]=[\bar{c},\bar{b}]=1\rangle \cong Q_8 \ast C_{4}$.

\medskip

In this subcase, $\Phi(G)=G'=\lg a^2,B'\rg$. Since $a^2=b^2z$ where
$z\in B'\le Z(G)$, we have $[a^2,b]=1$. It follows that $\Phi(G)\le
Z(G)$, a contradiction.

\medskip

Subcase 2.2. $G/B'$ is the group (9) in Lemma \ref{A_2}. That
is, $G/B'=\langle \bar{a}, \bar{b}, \bar{c} \di
\bar{a}^{p^{n+1}}=\bar{b}^{p^m}=\bar{c}^{p}=1,
[\bar{a},\bar{b}]=\bar{a}^{p^{n}},
[\bar{c},\bar{a}]=[\bar{c},\bar{b}]=1\rangle\cong M(n+1,m)\times
C_{p}$.

\medskip

If $o(a)=p^{n+1}$, then, by calculation, we have $\Phi(G)\le Z(G)$,
a contradiction. Hence $o(a)=p^{n+2}$. If $p>2$, then $\lg a,b\rg$
has no abelian subgroup  of index $p$, a contradiction. Hence $p=2$.
If $b^{2^{m}}\in \lg a^{2^{n+1}}\rg$, then $\lg a^2,b\rg$ is
non-abelian and is of order $2^{n+m+1}$, which is contrary to
$G\in\mathcal{A}_3$. Hence $o(b)=2^{m+1}$ and $B'=\lg
a^{2^{n+1}},b^{2^m}\rg$. Since $\lg a,b\rg$ has an abelian subgroup
of index $p$, by calculation we have $n=1$ and $A\cap\lg a,b\rg=\lg
a,b^2\rg$, where $A$ is the abelian subgroup  of index $p$ of $G$.
Hence $A=\lg a,b^2,c\rg$. It follows that $[c,a]=1$. Assume that
$[c,b]=a^{4w}b^{u2^m}$. By replacing $c$ with $ca^{2w}$ we have
$[c,b]=b^{u2^m}$. Since $Z(G)<\Phi(G)$, we have $[c,b]=b^{2^m}$. If
$c^2\in \lg b^{2^m}\rg$, then $|\lg c,b\rg|=2^{m+2}$, which is
contrary to $G\in\mathcal{A}_3$. Thus we may assume that $$G=\langle
a, b, c\di a^{8}=b^{2^{m+1}}=1, c^{2}=a^{4}b^{i2^m},
[a,b]=a^{2}a^{4j}b^{k2^m}, [c,a]=1,[c,b]=b^{2^m}\rangle.$$ By
replacing $a$ with $ac^j$ we may assume that $[a,b]=a^2b^{k'2^m}$.
If $m\ge 2$, then, replacing $a$ and $c$ with $ab^{k'2^{m-1}}$ and
$cb^{(i+1)2^{m-1}}$ respectively, we have $[a,b]=a^2$ and
$c^2=a^4b^{2^m}$. Hence we get the group (I10). If $m=1$, then
$(ba)^2=b^2a^2[a,b]=a^4b^{2(k'+1)}$. Since $G\in\mathcal{A}_3$, we
have $|\lg a^2,ba\rg|=2^4$. It follows that $(ba)^2=a^4b^2$ and
$[a,b]=a^2$. Again since $G\in\mathcal{A}_3$, we have $|\lg
ca^2,ba\rg|=2^4$. It follows that $(ca)^2=b^2$ and $c^2=a^4b^2$.
Hence we also get the group (I10).

\medskip

Subcase 2.3. $G/B'$ is the group (10) in Lemma \ref{A_2}. That
is, $G/B'=\langle \bar{a}, \bar{b}, \bar{c};\bar{d} \di
\bar{a}^{p^n}=\bar{b}^{p^m}=\bar{c}^{p}=\bar{d}^p=1,
[\bar{a},\bar{b}]=\bar{d},
[\bar{c},\bar{a}]=[\bar{c},\bar{b}]=1\rangle\cong M(n,m,1)\times
C_{p}$, where $n\ge m$, and $n\ge 2$ if $p=2$.

\medskip

Since $\lg \bar{a},\bar{b}\rg$ is not metacyclic, $\lg a,b\rg$ is
not metacyclic. Since $\lg a,b\rg$ is a two-generator
$\mathcal{A}_2$-subgroup, by Lemma \ref{A_2-property} (6) and (7) we
have $p>2$ and $\exp(\lg a,b\rg')=p$. Hence $d^p=1$ and we may
assume $[a,b]=d$. If $m\ge 2$, then $$|G: \lg d,a\rg|\ge p^{m+1}\ge
p^3\ {\rm and}\ |G: \lg d,b\rg|\ge p^{n+1}\ge p^3.$$ Since
$G\in\mathcal{A}_3$,  $[d,a]=[d,b]=1$. It follows that $\lg
a,b\rg\in \mathcal{A}_1$, a contradiction. Hence we have $m=1$.
Since the index of $\lg b,c\rg$ is at least $p^{n+1}$, we have
$[b,c]=1$ for $n\ge 2$. If $n=1$, then, without loss of generality,
we may assume $[b,c]=1$. Since $Z(G)<\Phi(G)$, we have $[c,a]\ne 1$.
Hence $A=\lg b,c,\Phi(G)\rg$. It follows that $[d,b]=1$ and hence
$[d,a]\ne 1$. Since $G\in \mathcal{A}_3$ and $|G|=p^{n+5}$, we have
$|\lg a,d\rg|=p^{n+3}$. Hence $$o(a)=p^{n+1},\ [d,a]\not\in\lg
a^{p^n}\rg\ {\rm and}\ B'=\lg a^{p^n},[d,a]\rg.$$ Since
$Z(G)<\Phi(G)$, we have $[c,a]\not\in\lg [d,a]\rg$. Hence we may
assume $[c,a]=a^{ip^n}[d,a]^j$ where $(i,p)=1$. By replacing $c$
with $c^{i^{-1}}d^{-i^{-1}j}$ we have $[c,a]=a^{p^n}$. Hence $|\lg
c,a\rg|=p^{n+3}$ and $c^p\not\in\lg a^{p^n}\rg$. It follows that
$G'=\lg d,a^{p^n},c^p\rg$. Hence we may assume
$[d,a]=c^{up}a^{vp^n}$ where $(u,p)=1$. Since $|\lg
d,c^ua^{vp^{n-1}}\rg|=p^3\le p^{n+3}$,  $[d,c^ua^{vp^{n-1}}]=1$ and
hence $n\ge 2$ or $v=0$. By replacing $b$, $c$ and $d$ with
$b^{u^{-1}}$, $ca^{u^{-1}vp^{n-1}}$ and $d^{u^{-1}}$, respectively,
we get $[d,a]=c^p$. We may assume $b^p=a^{sp^n}(\mod \lg c^p\rg)$.
By calculation we have $(ba^{-sp^{n-1}})^p\in \lg c^p\rg$. Since
$|\lg c,ba^{-sp^{n-1}}\rg|\le p^{3}$, we have
$[c,ba^{-sp^{n-1}}]=1$. Assume $(ba^{-sp^{n-1}})^p=c^{tp}$. Since
$|\lg d,ba^{-sp^{n-1}}c^{-t}\rg|\le p^{3}$, we have
$[d,ba^{-sp^{n-1}}c^{-t}]=1$. By replacing $b$ with
$ba^{-sp^{n-1}}c^{-t}$ we have $b^p=1$. Hence $G$ is the group
(I11).

\medskip

Subcase 2.4. $G/B'$ is the group (12) in Lemma \ref{A_2}. That
is, $G/B'=\langle \bar{a}, \bar{b}, \bar{c} \di
\bar{a}^{p^n}=\bar{b}^{p^m}=\bar{c}^{p^{2}}=1,
[\bar{a},\bar{b}]=\bar{c}^{p}, [\bar{c},\bar{a}]=[\bar{c},\bar{b}]=1
\rangle\cong M(n,m,1)\ast C_{p^{2}}$, where $n\ge m$, and $n\ge 2$
if $p=2$.

\medskip

Since $\lg \bar{a},\bar{b}\rg$ is not metacyclic,  $\lg a,b\rg$ is
not metacyclic. Since $\lg a,b\rg$ is a two-generator
$\mathcal{A}_2$-group, by Lemma \ref{A_2-property} (6) and (7) we
have $p>2$ and $\exp(\lg a,b\rg')=p$. Hence $c^{p^2}=1$ and
$\exp(G')=p$. It follows that $[c^p,a]=[c^p,b]=1$. Thus $\lg
a,b\rg\in \mathcal{A}_1$, a contradiction.

\medskip
We calculate the $(\mu_0,\mu_1,\mu_2)$ and $\alpha_1(G)$ of those
groups in Theorem \ref{d=3-4} as follows.

It is easy to see that $(\mu_0,\mu_1,\mu_2)=(1, 0, p^2+p)$. We only
need to calculate $\alpha_1(G)$.

\medskip
{\bf Case 1.} $G$ is one of the groups (I1)--(I4).

 All maximal subgroups of $G$ are:

$N=\lg a,x,\Phi(G)\rg$;

$N_i=\lg ba^i, x, \Phi(G)\rg$, where $i=0,1$;

$N_{ij}=\lg ax^i, bx^j, \Phi(G)\rg$, where $i,j=0,1$.

It is easy to see that $N=\lg a,x,\Phi(G)\rg$ is the unique abelian maximal subgroup of $G$.
Since $|N_i'|=p$ and $d(N_i)=3$, by Lemma \ref{A_2-property} (7), $\alpha_1(N_i)=p^2=4$.
Since $N_{ij}=\lg ax^i,bx^j\rg$, $\alpha_1(N_{ij})=p=2$.

Let $H\in \Gamma_2(G)$. Then $H=\lg a^ib^jx^k,\Phi(G)\rg$ where $(i,j,k)\neq (0,0,0)$.
It is obvious that $H\in\mathcal{A}_1$ if and only if $j\neq 0$.  By Hall's enumeration principle,
$$\alpha_1(G)=\sum_{H\in
\Gamma_1} \alpha_1(H)-p\sum_{H\in \Gamma_2} \alpha_1(H)=p\times
p^2+p^2\times p-p\times p^2=p^3.$$

\medskip
{\bf Case 2.} $G$ is one of the groups (I5)--(I9).

 All maximal subgroups of $G$ are:

$N=\lg a_1,x,\Phi(G)\rg$;

$N_i=\lg ba_1^i, x, \Phi(G)\rg$, where $0\leq i\leq p-1$;

$N_{ij}=\lg a_1x^i, bx^j, \Phi(G)\rg$, where $0\leq i,j\leq p-1$.

It is easy to see that $N=\lg a_1,x,\Phi(G)\rg$ is the unique abelian maximal subgroup of $G$.
Since $|N_i'|=p$ and $d(N_i)=3$, by Lemma \ref{A_2-property}, $\alpha_1(N_i)=p^2$. Since $N_{ij}=\lg a_1x^i,bx^j\rg$, $\alpha_1(N_{ij})=p$.

Let $H\in \Gamma_2(G)$. Then $H=\lg a_1^ib^jx^k,\Phi(M)\rg$ where $(i,j,k)\neq (0,0,0)$.
It is obvious that $H\in\mathcal{A}_1$ if and only if $j\neq 0$.  By Hall's enumeration principle,
$$\alpha_1(G)=\sum_{H\in
\Gamma_1} \alpha_1(H)-p\sum_{H\in \Gamma_2} \alpha_1(H)=p\times
p^2+p^2\times p-p\times p^2=p^3.$$

\medskip
{\bf Case 3.} $G$ is the group (I10).

All maximal subgroups of $G$ are:

$N=\lg a,c,\Phi(G)\rg$;

$N_i=\lg ba^i, c, \Phi(G)\rg$, where $i=0,1$;

$N_{ij}=\lg ac^i, bc^j, \Phi(M)\rg$, where $i,j=0,1$.

It is easy to see that $N=\lg a,c,\Phi(G)\rg$ is the unique abelian maximal subgroup of $G$.
Since $|N_i'|=p^2=4$ and $d(N_i)=3$, by Lemma \ref{A_2-property}, $\alpha_1(N_i)=p^2+p$.
By calculation, $N_{ij}=\lg ac^i,bc^j\rg$. Hence $\alpha_1(N_{ij})=p$.

Let $H\in \Gamma_2(G)$. Then $H=\lg a^ib^jc^k,\Phi(G)\rg$ where $(i,j,k)\neq (0,0,0)$.
It is obvious that $H\in\mathcal{A}_1$ if and only if $j\neq 0$.  By Hall's enumeration principle,
$$\alpha_1(G)=\sum_{H\in
\Gamma_1} \alpha_1(H)-p\sum_{H\in \Gamma_2} \alpha_1(H)=p\times
(p^2+p)+p^2\times p-p\times p^2=p^3+p^2.$$

\medskip
{\bf Case 4.} $G$ is the group (I11).

 All maximal subgroups of $G$ are:

$N=\lg b,c,\Phi(G)\rg$;

$N_i=\lg ab^i, c, \Phi(G)\rg$, where $0\leq i\leq p-1$;

$N_{ij}=\lg ac^i, bc^j, \Phi(G)\rg$, where $0\leq i,j\leq p-1$.

It is easy to see that $N=\lg b,c,\Phi(G)\rg$ is the unique abelian maximal subgroup of $G$.
Since $|N_i'|=p^2$ and $d(N_i)=3$, by Lemma \ref{A_2-property}, $\alpha_1(N_i)=p^2+p$.
By calculation, $N_{ij}=\lg ac^i,bc^j\rg$. Hence $\alpha_1(N_{ij})=p$.

Let $H\in \Gamma_2$. Then $H=\lg a^ib^jc^k,\Phi(G)\rg$ where $(i,j,k)\neq (0,0,0)$.
It is obvious that $H\in\mathcal{A}_1$ if and only if $i\neq 0$.  By Hall's enumeration principle,
$$\alpha_1(G)=\sum_{H\in
\Gamma_1} \alpha_1(H)-p\sum_{H\in \Gamma_2} \alpha_1(H)=p\times
(p^2+p)+p^2\times p-p\times p^2=p^3+p^2.$$\qed

\begin{thm}\label{d=3-5}
$d(G)=4$ if and only if $G$ is isomorphic to one of the following
pairwise non-isomorphic groups:
\begin{enumerate}

\rr{Ji} $G'\cong C_p$  and $c(G)=2$. In this case, $(\mu_0,\mu_1,\mu_2)=(1+p,0,p^2+p^3)$ and $\alpha_1(G)=p^4$.

\begin{enumerate}

\rr{J1} $\langle {a},{b},{x},y \di {a}^4={x}^{2}=y^2=1,
{b}^2={a}^2=[{a},{b}], [{x},{a}]=[{x},{b}]=[y,a]=[y,b]=[x,y]=1
\rangle \cong Q_8 \times C_{2}\times C_2$; where $|G|=2^{5}$,
$\Phi(G)=G'=\lg a^{2}\rg$ and $Z(G)=\lg a^2,x,y\rg \cong C_2^3$.

\rr{J2}  $\langle {a}, {b},{x},y \di
{a}^{p^{n+1}}={b}^{p^m}={x}^{p}=y^p=1, [{a},{b}]={a}^{p^{n}},
[{x},{a}]=[{x},{b}]=[y,a]=[y,b]=[x,y]=1\rangle\cong M_p(n+1,m)\times
C_{p}\times C_p$; where $|G|=p^{n+m+3}$, $\Phi(G)=\lg a^p,
b^p\rg\cong C_{p^{n}}\times C_{p^{m-1}}$ if $m>1$, $\Phi(G)\cong
C_{p^{n}}$ if $m=1$, $G'=\lg a^{p^n}\rg$, $Z(G)=\lg a^p,b^p,x,y\rg
\cong C_{p^{n}}\times C_{p^{m-1}}\times C_{p}^2$ if $m>1$,
$Z(G)\cong C_{p^{n}}\times C_{p}^2$ if $m=1$.

\rr{J3}  $\langle {a}, {b},{x},y; {c} \di
{a}^{p^n}={b}^{p^m}={c}^{p}={x}^p=y^p=1, [{a},{b}]={c},
[{c},{a}]=[{c},{b}]=[{x},{a}]=[{x},{b}]=[y,a]=[y,b]=[x,y]=1\rangle\cong
M_p(n,m,1)\times C_{p}\times C_{p}$, where $n\ge m$, and $n\ge 2$ if
$p=2$; moreover, $|G|=p^{n+m+3}$, $\Phi(G)=\lg a^p, b^p,c\rg\cong
C_{p^{n-1}}\times C_{p^{m-1}}\times C_{p}$ if $m>1$, $\Phi(G)\cong
C_{p^{n-1}}\times C_{p}$ if $m=1$ and $n>1$,  $\Phi(G)\cong C_{p}$
if $m=n=1$, $G'=\lg c\rg$, $Z(G)=\lg a^p,b^p,c,x,y\rg \cong
C_{p^{n-1}}\times C_{p^{m-1}}\times C_{p}^3$ if\  $m>1$, $Z(G) \cong
C_{p^{n-1}}\times C_{p}^3$ if $m=1$ and $n>1$, $Z(G) \cong
 C_{p}^3$ if $m=n=1$.

\rr{J4}  $\langle {a}, {b}, {x},y \di {a}^4=y^2=1,
{b}^2={x}^{2}={a}^2=[{a},{b}],
[{x},{a}]=[{x},{b}]=[y,a]=[y,b]=[x,y]=1\rangle \cong Q_8 \ast
C_{4}\times C_{2}$; where $|G|=2^{5}$, $\Phi(G)=G'=\lg a^{2}\rg$ and
$Z(G)=\lg x,y\rg\cong C_4\times C_2$.

\rr{J5}  $\langle {a}, {b}, {x},y \di
{a}^{p^n}={b}^{p^m}={x}^{p^{2}}=y^p=1, [{a},{b}]={x}^{p},
[{x},{a}]=[{x},{b}]=[y,a]=[y,b]=[x,y]=1 \rangle\cong M_p(n,m,1)\ast
C_{p^{2}}\times C_{p}$, where $n\ge 2$ if $p=2$ and $n\ge m$.
Moreover, $|G|=p^{n+m+3}$, $\Phi(G)=\lg a^p, b^p,x^p\rg\cong
C_{p^{n-1}}\times C_{p^{m-1}}\times C_{p}$ if $m>1$, $\Phi(G)\cong
C_{p^{n-1}}\times C_{p}$ if $m=1$ and $n>1$,  $\Phi(G)\cong C_{p}$
if $m=n=1$, $G'=\lg x^p\rg$, $Z(G)=\lg a^p,b^p,x,y\rg \cong
C_{p^{n-1}}\times C_{p^{m-1}}\times C_{p}\times C_{p^2}$ if\  $m>1$,
$Z(G) \cong C_{p^{n-1}}\times C_{p}\times C_{p^2}$ if $m=1$ and
$n>1$, $Z(G) \cong
 C_{p}\times C_{p^2}$ if $m=n=1$.

\end{enumerate}

\rr{Jii} $G'\cong C_p^2$  and $c(G)=2$. In this ase, $(\mu_0,\mu_1,\mu_2)=(1,0,p+p^2+p^3)$ and $\alpha_1(G)=p^4+p^3$.
\begin{enumerate}
\rr{J6} $K\times C_2$, where $K=\lg {a},{b},{c}\di
{a}^{4}={b}^{4}=1,
{c}^2={a}^2{b}^2,[{a},{b}]={b}^2,[{c},{a}]={a}^2,[{c},{b}]=1\rg$;
and $|G|=2^{6}$, $\Phi(G)=G'=\lg a^2, b^2\rg\cong C_2^2$, $Z(G)\cong
C_2^3$.

\rr{J7} $K\times C_p$, where $K=\lg {a},{b},{d}\di
{a}^{p^{m}}={b}^{p^2}={d}^{p}=1,
[{a},{b}]={a}^{p^{m-1}},[{d},{a}]={b}^p,[{d},{b}]=1\rg$, where
$m\geq 3$ if $p=2$; moreover, $|G|=p^{m+4}$, $\Phi(G)=\lg a^p,
b^p\rg \cong C_{p^{m-1}}\times C_p$, $G'=\lg a^{p^{m-1}},b^p\rg$,
$Z(G)\cong C_{p^{m-1}}\times C_p^2$.

\rr{J8} $K\times C_p$, where $K=\lg {a},{b},{d}\di
{a}^{p^m}={b}^{p^2}={d}^{p^2}=1,
[{a},{b}]={d}^p,[{d},{a}]={b}^{jp},[{d},{b}]=1\rg$, where $(j,p)=1$,
$p>2$, $j$ is a fixed quadratic non-residue modulo $p$, and $-4j$ is
a quadratic non-residue modulo $p$; moreover, $|G|=p^{m+5}$,
$\Phi(G)=\lg a^p, b^p,d^p\rg\cong C_{p^{m-1}}\times C_p^2$ if $m>1$,
$\Phi(G) \cong C_p^2$ if $m=1$, $G'=\lg b^{p},d^p\rg$, $Z(G)\cong
C_{p^{m-1}}\times C_p^3$  if $m>1$, $Z(G)\cong C_p^3$  if $m=1$.

\rr{J9} $K\times C_p$, where $K=\lg {a},{b},{d}\di
 {a}^{p^m}={b}^{p^2}={d}^{p^2}=1,[{a},{b}]={d}^p,[{d},{a}]={b}^{jp}{d}^p,[{d},{b}]=1\rg$,
where if $p$ is odd, then $4j =1-\rho^{2r+1}$ with $1\le
r\le\frac{p-1}{2}$ and $\rho$ the smallest positive integer which is
a primitive root $(\mod p)$; if $p = 2$, then $j = 1$. Moreover,
$|G|=p^{m+5}$, $\Phi(G)=\lg a^p, b^p,d^p\rg\cong C_{p^{m-1}}\times
C_p^2$ if $m>1$, $\Phi(G) \cong C_p^2$ if $m=1$, $G'=\lg
b^{p},d^p\rg$, $Z(G)\cong C_{p^{m-1}}\times C_p^3$  if $m>1$,
$Z(G)\cong C_p^3$  if $m=1$.

\end{enumerate}
\end{enumerate}
\end{thm}

\demo By Lemma \ref{A_t+k}, groups (J1)--(J9) are
$\mathcal{A}_3$-groups. By checking maximal subgroups of groups
(J1)--(J9), we know that they are pairwise non-isomorphic.
Conversely, in the following, we prove that $G$ is one of the groups
(J1)--(J9).

By Lemma \ref{d=4}, we have $c(G)=2$, $\Phi(G)\le Z(G)$, all
$\mathcal{A}_1$-subgroups of $G$ contain $\Phi(G)$ and $G'\le
C_p^3$. Let $S$ be an abelian subgroup of index $p$ of $G$ and $t\in
G\setminus A$. By Lemma \ref{Ber}, $G'=\lg [s,t]\di s\in S\rg$.

We claim $G'\le C_p^2$. Otherwise, $G'\cong C_p^3$. Let $M=\lg
s,t\rg$ be an $\mathcal{A}_1$-subgroup of $G$. Since $G'\le M$, $M$
is not metacyclic. Hence we may assume that $$M=\lg s,t\di
s^{p^m}=t^{p^n}=r^p=1,[s,t]=r,[r,s]=[r,t]=1\rg$$ and $$G'=\lg
s^{p^{m-1}},t^{p^{n-1}},r\rg.$$ By Lemma \ref{Ber}, there exists
$a\in S$ such that $[a,t]=t^{p^{n-1}}$. Hence $\lg a,t\rg$ is
metacyclic, which is contrary to $G'\le \lg a,t\rg$. Hence $G'\le
C_p^2$.

By Lemma \ref{d(K) leq k+1 3}, there exists $K\le G$ such that
$K'=G'$ and $d(K)=3$. Since $K\in \mathcal{A}_2$,
 $K$ is maximal in $G$. By Lemma \ref{Ber},
$$|Z(G)|=\frac{|G|}{p|G'|}\ {\rm and}\ |Z(K)|=\frac{|K|}{p|K'|}.$$ It
follows that $Z(G)\not\le K$. Let $y\in Z(G)\setminus K$. Then
$G=\lg K,y\rg$. By Lemma \ref{A_2}, $K$ is one of the groups (8)--(16)
in Lemma \ref{A_2}.

\medskip

Case 1. $K$ is one of the groups (8)--(10) in Lemma \ref{A_2}. That
is, $K=\lg a,b\rg\times \lg x\rg$.

\medskip

Let $L=\lg a,b,y\rg$. Then $G=L\times \lg x\rg$ and $L$ is a group
of Type (8)--(12) in Lemma \ref{A_2}. Hence we get groups
(J1)--(J5).

\medskip

Case 2. $K$ is the groups (11) in Lemma \ref{A_2}. That is, $K=
\langle {a}, {b}, {x} \di {a}^4=1, {b}^2={x}^{2}={a}^2=[{a},{b}],
[{x},{a}]=[{x},{b}]=1\rangle \cong Q_8 \ast C_{4}$.

\medskip

Since $y^2\in\Phi(G)=\Phi(K)=\lg x^2\rg$, we may assume
$y^2=x^{2i}$. By replacing $y$ with $yx^i$ we have $y^2=1$. Hence we
get the group (J4).

\medskip

Case 3. $K$ is the group (12) in Lemma \ref{A_2}. That is, $K=
\langle {a}, {b}, {x} \di {a}^{p^n}={b}^{p^m}={x}^{p^{2}}=1,
[{a},{b}]={x}^{p}, [{x},{a}]=[{x},{b}]=1 \rangle\cong M_p(n,m,1)\ast
C_{p^{2}},$ where $n\ge 2$ if $p=2$ and $n\ge m$.

\medskip

Let $L=\lg a,b,y\rg$. Then $G=L\ast \lg x\rg$ and $L$ is one of the groups (8)--(12) in Lemma \ref{A_2}.
If $L$ is one of the groups
(8)--(11) in Lemma \ref{A_2}, then it is reduced to Case 1 or 2.
Hence we may assume that $L$ is also a group of Type (12) in Lemma
\ref{A_2}. Let $$L=\langle {a}, {b}, {y} \di
{a}^{p^n}={b}^{p^m}={y}^{p^{2}}=1, [{a},{b}]={y}^{p},
[{y},{a}]=[{y},{b}]=1 \rangle.$$ Since $x^p\in G'=K'$, we may assume
that $x^p=y^{ip}$. By replacing $x$ with $xy^{-i}$ we have $x^p=1$.
Hence we get the group (J5).

\medskip

Case 4. $K$ is the group (13) in Lemma \ref{A_2}. That is,
$K=\lg {a},{b},{c}\di {a}^{4}={b}^{4}=1,
{c}^2={a}^2{b}^2,[{a},{b}]={b}^2,[{c},{a}]={a}^2,[{c},{b}]=1\rg.$

\medskip

Since $y^2\in\Phi(G)=\Phi(K)=\lg a^2,b^2\rg$, we may assume
$y^2=a^{2i}b^{2j}$. Since $|\lg ya^i,b\rg|\le 8$, we have
$[ya^i,b]=1$ and hence $y^2=b^{2j}=c^{2j}a^{2j}$. Since $|\lg
yc^j,a\rg|\le 8$, we have $[yc^j,a]=1$ and hence $y^2=1$. Thus we
get the group (J6).

\medskip

Case 5. $K$ is the group (14) in Lemma \ref{A_2}. That is,
$K=\lg {a},{b},{d}\di {a}^{p^{m}}={b}^{p^2}={d}^{p}=1,
[{a},{b}]={a}^{p^{m-1}},[{d},{a}]={b}^p,[{d},{b}]=1\rg$, where
$m\geq 3$ if $p=2$.

\medskip

Since $y^p\in\Phi(G)=\Phi(K)=\lg a^p,b^p\rg$, we may assume
$y^p=a^{ip}b^{jp}$. Since $|\lg yb^{-j},a\rg|\le p^{m+1}$, we have
$[yb^{-j},a]=1$ and hence $y^p=a^{ip}$. Since $|\lg ya^{-i},d\rg|\le
p^3$, we have $[ya^{-i},d]=1$ and hence $a^{i}\in Z(K)$. By
replacing $y$ with $ya^{-i}$ we have $y^p=1$. Hence we get the group
(J7).

\medskip

Case 6. $K$ is the group (15) in Lemma \ref{A_2}. That is,
$K=\lg {a},{b},{d}\di {a}^{p^m}={b}^{p^2}={d}^{p^2}=1,
[{a},{b}]={d}^p,[{d},{a}]={b}^{jp},[{d},{b}]=1\rg$, where $(j,p)=1$,
$p>2$, $j$ is a fixed quadratic non-residue modulo $p$, and $-4j$ is
a quadratic non-residue modulo $p$.

\medskip

Since $y^p\in\Phi(G)=\Phi(K)=\lg b^p,d^p\rg$, we may assume that
$y^p=b^{rp}d^{sp}$. Since $|\lg yb^{-r},a\rg|\le p^{m+2}$, we have
$[yb^{-r},a]=1$ and hence $y^p=d^{sp}$. Since $|\lg yd^{-s},a\rg|\le
p^{m+2}$, we have $[yd^{-s},a]=1$ and hence $y^p=1$. Thus we get the
group (J8).

\medskip

Case 7. $K$ is the group (16) in Lemma \ref{A_2}. That is,
$K=\lg {a},{b},{d}\di
 {a}^{p^m}={b}^{p^2}={d}^{p^2}=1,[{a},{b}]={d}^p,[{d},{a}]={b}^{jp}{d}^p,[{d},{b}]=1\rg$,
where if $p$ is odd, then $4j =1-\rho^{2r+1}$ with $1\le
r\le\frac{p-1}{2}$ and $\rho$ the smallest positive integer which is
a primitive root $(\mod p)$; if $p = 2$, then $j = 1$.

\medskip

Since $y^p\in\Phi(G)=\Phi(K)=\lg b^p,d^p\rg$, we may assume
$y^p=b^{rp}d^{sp}$. Since $|\lg yb^{-r},a\rg|\le p^{m+2}$, we have
$[yb^{-r},a]=1$ and hence $y^p=d^{sp}$. Since $|\lg yd^{-s},a\rg|\le
p^{m+2}$, we have $[yd^{-s},a]=1$ and hence $y^p=1$. Thus we get the
group (J9).

\medskip
We calculate the $(\mu_0,\mu_1,\mu_2)$ and $\alpha_1(G)$ of those
groups in Theorem \ref{d=3-5} as follows.

Since $d(G)=4$ and $\mu_1=0$, $\mu_0+\mu_2=1+p+p^2+p^3$. Let $N$ be
an $\mathcal{A}_1$-subgroup of $G$. Since $G\in \mathcal{A}_3$,
$N\Phi(G)\in \Gamma_1$ or $\Gamma_2$. If $N\Phi(G)\in \Gamma_1$,
then $d(G)\le 3$, a contradiction. Hence $N\Phi(G)\in\Gamma_2$.
Since $G\in \mathcal{A}_3$, $N\Phi(G)\in\mathcal{A}_1$. That is,
$N=N\Phi(G)\in\Gamma_2$. Hence $\sum_{H\in \Gamma_2}
\alpha_1(H)=\alpha_1(G)$. By Hall's enumeration principle,
$\alpha_1(G)=\sum_{H\in \Gamma_1} \alpha_1(H)-p\sum_{H\in \Gamma_2}
\alpha_1(H)=\sum_{H\in \Gamma_1} \alpha_1(H)-p\alpha_1(G)$. Hence
$\alpha_1(G)=\frac{1}{1+p}\sum_{H\in \Gamma_1} \alpha_1(H).$

\medskip
{\bf Case 1.} $G$ is one of the groups (J1)--(J5).

Since $Z(G)$ is of index $p^2$, $\mu_0=1+p$ and hence
$\mu_2(G)=p^2+p^3$. Let $H\in\Gamma_1$. If $H$ is not abelian, then
$|H'|=p$. By Lemma \ref{A_2-property}, $\alpha_1(H)=p^2$. Hence
$$\alpha_1(G)=\frac{1}{1+p}\sum_{H\in \Gamma_1}
\alpha_1(H)=\frac{1}{1+p}\mu_2 p^2=p^4.$$

\medskip
{\bf Case 2.} $G$ is one of the groups (J6)--(J9).

Let $A$ be the unique abelian maximal subgroup of $K$. Then $A\times
\lg x\rg$ is the unique abelian maximal subgroup of $G$. Other
maximal subgroups of $G$ are:

$N_i=M_i\times \lg x\rg$, where $M_i$ are non-abelian maximal subgroups of $K$;

$N_{ijk}=\lg ax^i, bx^j, dx^k\rg$ (or $\lg ax^i, bx^j, dx^k\rg$ for (J6)), where $0\leq i,j,k\leq p-1$.

It is easy to see that $|N_i'|=p$. By Lemma \ref{A_2-property}, $\alpha_1(N_i)=p^2$. Since $N_{ijk}\cong K$,
$\alpha_1(N_{ijk})=p+p^2$. Hence
$$\alpha_1(G)=\frac{1}{1+p}\sum_{H\in
\Gamma_1}
\alpha_1(H)=\frac{1}{1+p}((p+p^2)p^2+p^3(p+p^2))=p^3+p^4.$$\qed

\subsection{$G$ has no abelian subgroup of index $p$}

In this section assume $G$ is an $\mathcal{A}_3$-group without an
abelian subgroup of index $p$ and an $\mathcal{A}_1$-subgroup of
index $p$ in  Theorem \ref{d=4-1}, \ref{d=4-2}, \ref{d=4-3},
\ref{d=4-4} and \ref{d=4-5}.

\begin{thm}\label{d=4-1}
Suppose that $G$ is a $p$-group all of whose maximal subgroups are
$\mathcal{A}_2$-groups, and every maximal subgroup of $G$ is
generated by two elements. Then $G$ is an $\mathcal{A}_3$-group if
and only if $G$ is isomorphic to one of the following pairwise
non-isomorphic groups:
\begin{enumerate}

\rr{Ki} $G$ is metacyclic. In this case, $(\mu_0,\mu_1,\mu_2)=(0,0,1+p)$ and $\alpha_1(G)=1+p+p^2$.

\begin{enumerate}
\rr{K1} $\lg a,b\di
a^{p^{r+3}}=1,b^{p^{r+s+t}}=a^{p^{r+s}},[a,b]=a^{p^r}\rg$, where
$p>2$ and $r,s,t$ are non-negative integers with $r\ge 1$ and
$r+s\ge 3$; $|G|=p^{2r+s+t+3}$, $c(G)=4$ for $r=1$; $c(G)=3$ for
$r=2$ and $c(G)=2$ for $r>2$, $\Phi(G)=\lg a^p, b^p\rg \cong
M_p(r+2, r+s+t-1)$ if $s\geq 3$, $\Phi(G)=\lg a^p, b^p\rg \cong
M_p(r+t+2, r+s-1)$ if $s< 3$; moreover, $G'=\lg a^{p^r}\rg\cong
C_{p^3}$, $Z(G)=\lg a^{p^3},b^{p^3}\rg \cong C_{p^r}\times
C_{p^{r+t+s-3}}$ if $s\geq 3$, $Z(G)=\lg a^{p^3},b^{p^3}\rg \cong
C_{p^{r+s-3}}\times C_{p^{r+t}}$ if $s<3$.

\rr{K2} $\lg a,b\di
a^{2^{r+3}}=1,b^{2^{r+s+t}}=a^{2^{r+s}},[a,b]=a^{2^r}\rg$, where
$r,s,t$ are non-negative integers with $r\ge 2$ and $r+s\ge 3$.
Moreover, $|G|=2^{2r+s+t+3}$,  $c(G)=3$ for $r=2$ and $c(G)=2$ for
$r>2$, $\Phi(G)=\lg a^2, b^2\rg \cong M_2(r+2, r+s+t-1)$ if $s\geq
3$, $\Phi(G)=\lg a^2, b^2\rg \cong M_2(r+t+2, r+s-1)$ if $s< 3$;
 $G'=\lg a^{2^r}\rg\cong C_{2^3}$,
$Z(G)=\lg a^{8},b^{8}\rg \cong C_{2^r}\times C_{2^{r+t+s-3}}$ if
$s\geq 3$, $Z(G)=\lg a^{2^3},b^{2^3}\rg \cong C_{2^{r+s-3}}\times
C_{2^{r+t}}$ if $s<3$.

\end{enumerate}

\rr{Kii} $G$ is not metacyclic. In this case, $p>2$, $|G|=p^6$, $(\mu_0,\mu_1,\mu_2)=(0,0,1+p)$ and $\alpha_1(G)=p+p^2$.

\begin{enumerate}
\rr{K3} $\lg a,b;c\di
a^{p^2}=b^{p^2}=c^{p^2}=1,[a,b]=c,[c,b]=a^{p}c^{mp},[c,a]=b^{\nu
p}c^{np}, [a,b^p]=[a^p,b]=c^p,[c,a^p]=[c,b^p]=[c^p,a]=[c^p,b]=1\rg$,
where $\nu$ is a fixed quadratic non-residue modulo $p$. The
parameters $m,n$ are the smallest positive integers satisfying
$(m-1)^2-\nu^{-1}(n+\nu)^2\equiv r (\mod p)$, for $r=0,1,\dots,p-1$;
moreover, $|G|=p^{6}$, $c(G)=4$, $\Phi(G)=G'=\lg a^p, b^p,c\rg\cong
C_p\times C_p \times C_{p^2}$ and $Z(G)=\lg c^p\rg\cong C_p$.

\rr{K4} $\lg a,b;c,d\di
a^9=b^9=c^3=d^3=1,[a,b]=c,[c,b]=a^3,[c,a]=b^{-3},[a^3,b]=[a,b^3]=d,[d,a]=[d,b]=1\rg$;
where $|G|=3^{6}$, $c(G)=4$, $\Phi(G)=G'=\lg a^3, b^3,c,d\rg\cong
C_3^4$ and $Z(G)=\lg d\rg\cong C_3$.

\rr{K5} $\lg a,b;c,d\di
a^9=b^9=c^3=d^3=1,[a,b]=c,[c,b]=a^3d,[c,a]=b^{-3}d,[a^3,b]=[a,b^3]=d,[d,a]=[d,b]=1\rg$,
where $|G|=3^{6}$, $c(G)=4$, $\Phi(G)=G'=\lg a^3, b^3,c,d\rg\cong
C_3^4$ and $Z(G)=\lg d\rg\cong C_3$.
\end{enumerate}
\end{enumerate}
\end{thm}

\demo Since all non-abelian proper subgroups of $G$ are generated by
two elements, $G$ is one of the groups classified in \cite{Alj}. It
follows from \cite[Main Theorem]{Alj}, $G$ is either a metacyclic
group or a $3$-group of maximal class or the group in \cite[Theorem
5.5 \& 5.6]{Alj}.

If $G$ is metacyclic, then $|G'|=p^3$ by Lemma \ref{metacyclic An}.
By using the classification of metacyclic $p$-groups in \cite{NX,XZ}
and checking their maximal subgroups,, we get the groups (K1)--(K2).
The details is omitted.

If $G$ is of maximal class, then, by Lemma \ref{jidalei3group}, $G$
does not satisfy the hypothesis.

If $G$ is the group in \cite[Theorem 5.5]{Alj}, then $G$ is the
group (K3). If $G$ is the group in \cite[Theorem 5.6]{Alj}, then $G$
is the group (K4)--(K5).

Since $d(G)=2$, $(\mu_0,\mu_1,\mu_2)=(0,0,1+p)$. If $G$ is one of
the groups (K1)--(K2), the $G$ is metacyclic and
$\Phi(G)\in\mathcal{A}_1$. Let $H\in\Gamma_1$. Then $H$ is also
metacyclic and $\alpha_1(H)=1+p$. By Hall's enumeration principle,
$$\alpha_1(G)=\sum_{H\in \Gamma_1} \alpha_1(H)-p\sum_{H\in \Gamma_2}
\alpha_1(H)=\mu_2(1+p)-p=p^2+p+1.$$

If $G$ is one of the groups {\rm(K3)--(K5)}, then $\Phi(G)$ is abelian.
Let $H\in\Gamma_1$. Then $\alpha_1(H)=p$.
By Hall's enumeration principle,
$$\alpha_1(G)=\sum_{H\in
\Gamma_1} \alpha_1(H)=\mu_2p=p^2+p.$$\qed

\begin{thm}\label{d=4-2}
Suppose that $|G|=p^n$ and $G$ is an $\mathcal{A}_3$-group. If $G$
has a three-generator maximal subgroup $M$ such that $d(M)=3$ and
$M'\not\le Z(G)$, then $p\ge 5$, $n=6$, all maximal subgroup of $G$
are $\mathcal{A}_2$-groups and $G$ is one of the following
non-isomorphic groups:
\begin{enumerate}
\rr{L1} $\lg x,m;a\di x^{p^2}=m^{p^2}=a^{p^2}=1,
[x,m]=a,[a,x]=x^p,[a,m]=m^{-p}\rg$; where $|G|=p^{6}$, $c(G)=4$,
$\Phi(G)=G'=\lg a, x^p, m^p\rg\cong C_p^2\times C_{p^2}$ and
$Z(G)=\lg a^p\rg\cong C_p$.

\rr{L2} $\lg x,m;a\di x^{p^2}=m^{p^2}=a^{p^2}=1,
[x,m]=a,[a,x]=x^pa^{p},[a,m]=m^{-p}a^{vp}\rg$, where $v\in F_p$.
Moreover, $|G|=p^{6}$, $c(G)=4$, $\Phi(G)=G'=\lg a, x^p, m^p\rg\cong
C_p^2\times C_{p^2}$ and $Z(G)=\lg a^p\rg\cong C_p$.
\end{enumerate}
Moreover, $(\mu_0,\mu_1,\mu_2)=(0,0,1+p)$ and $\alpha_1(G)=3p^2+p$.
\end{thm}

\demo By Lemma \ref{lem a0,b0,M' notin Z(G)}, we may assume $|M'|\le
p^2$. If $|M'|=p$, then $M'\le Z(G)$. Hence $|M'|=p^2$. It follows
from Lemma \ref{A_2-property}(4) that $c(M)=2$, $M'\cong C_p^2$,
$\Phi(M)=\mho_1(M)=Z(M)$, $M$ has a unique abelian subgroup $A$  of
index $p$ and $A/M'$ has type $(p^{n-6},p,p)$. Since $A$ is
characteristic in $M$ and $M\unlhd G$, we have $A\unlhd G$. We take
$x\in G\setminus M$ and let $H=\lg x,M'\rg$. Then $H$ is not abelian
since $M'\nleq Z(G)$. Moreover, $|H'|=p$, $H'\le M'$ and $H'\le
Z(G)$. By Lemma \ref{minimal non-abelian equivalent conditions} we
have $H$ is an $\mathcal{A}_1$-group. It follows that $|H|\ge
p^{n-2}$. Hence $o(xM')\ge p^{n-4}$. We will prove that

\medskip
(1) $\Phi(G)\le A$.
\medskip

Otherwise, $G/A$ is cyclic, $G=\lg x,A\rg$ and $M=\lg x^p,A\rg$.
Since $x^{p^2}\in \mho_1(M)=Z(M)$, we have $x^{p^2}\in Z(G)$. If
follows that $$|Z(G)/H'|=|Z(G)/Z(G)\cap M'|=|Z(G)M'/M'|\ge
o({x}^{p^2}M')\ge p^{n-6}.$$ Hence $|Z(G)|\ge p^{n-5}$. By Lemma
\ref{Ber}, $|G'|=|A/Z(G)|\le p^3$. Then $G_3\le M'$, $c(G)\le 4$ and
$G_4\le Z(G)$. Since $G'\ge M'\cong C_p^2$, $d(G')\ge 2$ and hence
$|\Phi(G')|\le p$. It follows that $\Phi(G')\le H'$.

We claim that $[x,w,x]^{p\choose 2}\in Z(G)$ for all $w\in G$.
Otherwise, we may assume that $p=2$ since $[x,w,x]^{p\choose 2}\in
\Phi(G')\le Z(G)$ for $p>2$. In this case, $[x,w,x,x]\neq 1$. Let
$L=\lg [x,w],x\rg$. Then $c(L)=3$. Since $G$ is an
$\mathcal{A}_3$-group, $L$ is an $\mathcal{A}_2$-group. By Lemma
\ref{A_2-property} (6), $L'\cong C_4$. On the other hand, $L'\le
G_3\le M'$, a contradiction.

Since $[x^p,w]\equiv [x,w]^p[x,w,x]^{p\choose 2}\ (\mod G_4)$ and
$\Phi(G')\le Z(G)$, $[x^p,w]\in Z(G)$. By Lemma \ref{Ber}, $M'=\{
[x^p,a]\di a\in A\}\le Z(G)$, a contradiction.

\medskip
(2) $|G|=p^6$, $o(xM')=p^2$ and $x^p\in A\setminus M'$.
\medskip

By (1), $x^p\in A$. Since $o(xM')\ge p^{n-4}$, we have $o(x^pM')\ge
p^{n-5}$ and hence $\exp(A/M')\ge p^{n-5}$. On the other hand,
$A/M'$ has type $(p^{n-6},p,p)$. If $n\ge 7$, then
$\exp(A/M')=p^{n-6}$, a contradiction. Hence $n=6$, $\exp(A/M')=p$,
$o(xM')=p^2$ and $x^p\in A\setminus M'$.

\medskip
(3) $G=\lg x,m;a\di x^{p^2}=m^{p^2}=a^{p^2}=1,
[x,m]=a,[a,x]=x^pa^{up},[a,m]=m^{-p}a^{vp}\rg$, where $u,v\in F_p$.
\medskip

Let $A=\lg x^p,a,M'\rg$ and $M=\lg m,A\rg$. Then $M'=\lg
[m,x^p]\rg\times \lg [m,a]\rg$. Since $[x^p,x,m]=1$ and
$[x,m,x^p]=1$,  we have $[m,x^p,x]=1$ by Witt's formula. Hence
$Z=\lg [m,x^p]\rg$. Since $M'\nleq Z(G)$, $[m,a,x]\neq 1$. Let
$[m,a,x]=[m,x^p]^i$, where $(i,p)=1$. By calculation we have
$$[[a,x]x^{-ip},m]=[a,x,m][x^{-ip},m]=[a,x,m][m,a,x]=1.$$ It follows
that $[a,x]x^{-ip}\in Z(M)=M'$.

We claim that $x^{p^2}=1$. Otherwise, $[a,x]^p\ne 1$. Hence $\lg
a,x\rg\in \mathcal{A}_2$. By Lemma \ref{A_2-property} (7), $\lg
a,x\rg$ is metacyclic. If $p>2$, then $[a,x^p]=[a,x]^p\ne 1$, a
contradiction. If $p=2$, then, by calculation,  we have
$1=[a,x^2]=[a,x]^2[a,x,x]$. Hence $[a,x,x]=x^4$. By calculation,
$[a^2,x]=[a,x]^2[a,x,a]=x^4$. Hence $a^2\not\in Z(G)$. It follows
that $M'=\lg a^2,x^4\rg$. Since $[a,x,x]\ne 1$, we have $[a,x]\equiv
x^2a^2 \ (\mod Z(G))$. By calculation we have $\lg a^2,ax\rg\cong
D_8$, which is contrary to that $G$ is an $\mathcal{A}_3$-group.
Hence $x^{p^2}=1$.

By Lemma \ref{A_2-property} (4), $\exp(A)>p$. It follows that
$a^p\ne 1$. Since $[a^p,x]=[a,x]^p=1$, $Z(G)=\lg a^p\rg$ and hence
$M'=\lg a^p,[a,m]\rg$.

We claim $p>2$. Otherwise, since $[a,x^2]=1$, we have $[a,x,x]=1$.
Hence we may assume $[a,x]=x^2a^{2u}$. By calculation we have $\lg
[a,m],ax\rg\cong D_8$, which is contrary to that $G$ is an
$\mathcal{A}_3$-group. Hence $p\ge 3$.

Since $[a,x,x]\in H'$, we may assume $[a,x,x]=a^{sp}$. Since $1\ne
[a,m,x]\in H'$, we may assume $[a,m,x]=a^{tp}$ where $(t,p)=1$. Let
$j=-2^{-1}t^{-1}s$. Then
$$[a,m^jx,m^jx]=[a,x,x][a,x,m]^j[a,m,x]^j=a^{(s+2jtp)}=1.$$ By
replacing $x$ with $xm^j$ we have $[a,x,x]=1$.

Since $[x^p,m]\ne 1$, we have $[x,m]^p\ne 1$ and $[x,m^p]\ne 1$.
Hence $M'=\lg a^p,m^p\rg$ and we may assume $[x,m]\equiv a^{v}\
(\mod \lg x^p,M'\rg)$, where $(v,p)=1$. By replacing $m$ with
$m^{v^{-1}}$ we have $[x,m]\equiv a\ (\mod \lg x^p,M'\rg)$. Without
loss of generality assume that $[x,m]=a$. Hence $G=\lg x,m\rg$
satisfy the following relations:
$$x^{p^2}=m^{p^2}=a^{p^2}=1,[x,m]=a,[a,x]\equiv x^{ip}\ (\mod \lg a^p\rg),[a,m]\equiv m^{jp}\ (\mod \lg
a^p\rg).$$ By replacing $m$ and $a$ with $m^{i^{-1}}$ and
$[x,m^{i^{-1}}]$ respectively, we have $[a,x]\equiv x^p\ (\mod \lg
a^p\rg)$. Since $G$ is metabelian, we have $[a,x,m]=[a,m,x]$. It
follows that $j=-1$. If $p=3$, then $\lg x^3, xm\rg$ is not abelian
and is of order $3^3$, which is contrary to that $G$ is an
$\mathcal{A}_3$-group. Hence $p\ge 5$ and we may assume that $$G=\lg
x,m;a\di x^{p^2}=m^{p^2}=a^{p^2}=1,
[x,m]=a,[a,x]=x^pa^{up},[a,m]=m^{-p}a^{vp}\rg,$$ where $u,v\in F_p$.

\medskip

(4) Two groups with parameters $(u,v)$ and $(u',v')$ are pairwise
isomorphic if and only if there exists $t\in F_p^*$ such that
$(u',v')=(u,v){\left[
                \begin{array}{c}
                  t \\
                  t^{-1} \\
                \end{array}
              \right]}
$ or $(u',v')=(v,u){\left[
                \begin{array}{c}
                  t \\
                  t^{-1} \\
                \end{array}
              \right]}
$.

\medskip
Suppose that $G=\lg x,m,a\rg$ and $\bar{G}=\lg
\bar{x},\bar{m},\bar{a}\rg$ with parameters $(u,v)$ and $(u',v')$
respectively. Let $\theta$ be an isomorphism from $\bar{G}$ to $G$.
Assume that $$\bar{x}^\theta=x^{s_{11}}m^{s_{12}}a^{s_{13}}r_1,\
\bar{m}^\theta=x^{s_{21}}m^{s_{22}}a^{s_{23}}r_2,$$ where $s_{ij}\in
F_p$ such that $s:=s_{11}s_{22}-s_{12}s_{21}\in F_p^*$ and
$r_1,r_2\in G_3$.

By calculation, $\bar{a}^\theta=[\bar{x},\bar{m}]^\theta\equiv
a^{s}\ (\mod G_3)$. Since $[\bar{a}^\theta,\bar{x}^\theta]\equiv
(\bar{x}^{p})^\theta\equiv (\bar{x}^\theta)^{p}\ (\mod G_4)$, we
have $[a^{s},x^{s_{11}}m^{s_{12}}]\equiv (x^{s_{11}}m^{s_{12}})^p \
(\mod G_4)$. Comparing indexes of $x^{p}$ and $m^{p}$ in two sides,
we have $ss_{11}=s_{11}$ and $ss_{12}=-s_{12}$. Since
$[\bar{a}^\theta,\bar{m}^\theta]=
(\bar{m}^{-p})^\theta=(\bar{m}^\theta)^{-p}$, we have
$[a^{s},x^{s_{21}}m^{s_{22}}]\equiv (x^{s_{21}}m^{s_{22}})^{-p}\
(\mod G_4)$. Comparing indexes of $x^{p}$ and $m^p$ in two sides, we
have $ss_{21}=-s_{21}$ and $ss_{22}=s_{22}$. It follows that
$$s_{11}s_{22}=1,s_{12}=s_{21}=0\ {\rm or}\
s_{11}=s_{22}=0,s_{12}s_{21}=1.$$

If $s_{11}s_{22}=1$ and $s_{12}=s_{21}=0$, then, by calculation,
$$\bar{a}^\theta=[\bar{x},\bar{m}]^\theta=[x^{s_{11}}a^{s_{13}},m^{s_{22}}a^{s_{23}}]\equiv
ax^{-s_{23}s_{11}p}m^{-s_{13}s_{22}p} \ (\mod G_4).$$ Since
$$[\bar{a}^\theta,\bar{x}^\theta]=(\bar{x}^{p}\bar{a}^{u'p})^\theta=
(\bar{x}^\theta)^{p}(\bar{a}^\theta)^{u'p},$$ we have
$$[ax^{-s_{23}s_{11}p}m^{-s_{13}s_{22}p},x^{s_{11}}]=
(x^{s_{11}}a^{s_{13}})^pa^{u'p}.$$ Comparing index of $a^{p}$ in two
sides, we have
\begin{equation}
s_{11}u=u'.
\end{equation}
 Since
$$[\bar{a}^\theta,\bar{m}^\theta]=(\bar{m}^{-p}\bar{a}^{v'p})^\theta=
(\bar{m}^\theta)^{-p}(\bar{a}^\theta)^{v'p},$$ we have
$$[ax^{-s_{23}s_{11}p},m^{s_{22}}]=
(m^{s_{22}}a^{s_{23}})^{-p}a^{v'p}.$$ Comparing index of $a^{p}$ in
two sides, we have
\begin{equation}
s_{22}v=v'.
\end{equation} Let $t=s_{11}$. Then $(u',v')=(u,v){\left[
                \begin{array}{c}
                  t \\
                  t^{-1} \\
                \end{array}
              \right]}
$.

On the other hand, if there exists $t\in F_p^*$ such that
$(u',v')=(u,v){\left[
                \begin{array}{c}
                  t \\
                  t^{-1} \\
                \end{array}
              \right]}
$, then $\theta: \bar{x}\mapsto x^t, \bar{m}\mapsto m^{t^{-1}}$ is
an isomorphism from $\bar{G}$ to $G$.

\medskip

If $s_{11}=s_{22}=0$ and $s_{12}s_{21}=1$, then, by calculation,
$$\bar{a}^\theta=[\bar{x},\bar{m}]^\theta=[m^{s_{12}}a^{s_{13}},x^{s_{21}}a^{s_{23}}]\equiv
a^{-1}x^{s_{13}s_{21}p}m^{s_{12}s_{23}p} \ (\mod G_4).$$ Since
$$[\bar{a}^\theta,\bar{x}^\theta]=(\bar{x}^{p}\bar{a}^{u'p})^\theta=
(\bar{x}^\theta)^{p}(\bar{a}^\theta)^{u'p},$$ we deduce that
$$[a^{-1}x^{s_{13}s_{21}p}m^{s_{12}s_{23}p},m^{s_{12}}]=
(m^{s_{12}}a^{s_{13}})^pa^{u'p}.$$ Comparing index of $a^{p}$ in two
sides, we have
\begin{equation}\label{i}
-s_{12}v=u'.
\end{equation}
Since
$$[\bar{a}^\theta,\bar{m}^\theta]=(\bar{m}^{-p}\bar{a}^{v'p})^\theta=
(\bar{m}^\theta)^{-p}(\bar{a}^\theta)^{v'p},$$ we deduce that
$$[a^{-1}x^{s_{13}s_{21}p}m^{s_{12}s_{23}p},x^{s_{21}}]=
(x^{s_{21}}a^{s_{23}})^{-p}a^{v'p}.$$ Comparing index of $a^{p}$ in
two sides, we have
\begin{equation}\label{ii}
-s_{21}u=v'.
\end{equation} Let $t=-s_{12}$. Then $(u',v')=(v,u){\left[
                \begin{array}{c}
                  t \\
                  t^{-1} \\
                \end{array}
              \right]}
$.

On the other hand, if there exists $t\in F_p^*$ such that
$(u',v')=(v,u){\left[
                \begin{array}{c}
                  t \\
                  t^{-1} \\
                \end{array}
              \right]}
$, then $\theta: \bar{x}\mapsto m^{-t}, \bar{m}\mapsto x^{-t^{-1}}$
is an isomorphism from $\bar{G}$ to $G$.

\medskip

(5) $G$ is one of the groups in the theorem. That is, we may assume
that $u=v=0$ or $u=1$.
\medskip

If $u=v=0$, then $G$ is the group (L1). If $u\ne 0$, then,
replacing $x$ and $m$ with $x^{u^{-1}}$ and $m^{u}$ respectively,
$G$ is the group (L2). If $v\ne 0$, then, replacing $x$ and
$m$ with $m^{v^{-1}}$ and $x^{v}$ respectively, $G$ is the group (L2).

By (4), groups in Theorem \ref{d=4-2} are pairwise non-isomorphic.

Now we calculate $(\mu_0,\mu_1,\mu_2)$ and $\alpha_1(G)$.

Since $d(G)=2$, $(\mu_0,\mu_1,\mu_2)=(0,0,1+p)$. All maximal
subgroups of $G$ are:

$M=\lg a,x,\Phi(G)\rg$;

$M_{i}=\lg a,mx^i, \Phi(G)\rg$, where $0\leq i\leq p-1$.

It is easy to see that $M=\lg a,x,m^p\rg$ such that $|M'|=p^2$ and $d(M)=3$. Hence $\alpha_1(M)=p^2+p$
and $\Phi(G)$ is the unique abelian maximal subgroup of $M$.
Similarly, $M_0=\lg a,m,x^p\rg$ and $\alpha(M_0)=p^2+p$. By calculation, $d(M_{i})=2$ for $i\neq 0$.
Hence $\alpha_1(M_{i})=p$ for $i\neq 0$.
By Hall's enumeration principle,
$$\alpha_1(G)=\sum_{H\in
\Gamma_1} \alpha_1(H)=2\times (p+p^2)+(p-1)\times p=3p^2+p.$$
 \qed

\begin{lem}\label{cong 0}
Suppose that $G(i, j)=\lg b,a_1;a_2,a_3\di
b^{p^2}=a_1^p=a_2^p=a_3^p=1,[a_1,b]=a_2,[a_2,b]=a_3,[a_3,b]=b^{ip},[a_2,a_1]=b^{jp},[a_3,a_1]=1\rg$,
where $i,j\in F_p^*$ and $p\ge 5$.
\begin{enumerate}
\rr1 If $G=G(i,j)$, then $|G|=p^{5}$, $G'=\lg a_2,a_3,b^p\rg$,
$G_3=\lg a_3,b^p\rg$, $G_4=\lg b^p \rg$, $G\in \mathcal{A}_3$, and
$M=\lg b^p,a_1,a_2,a_3\rg$ is the unique three-generator maximal
subgroup of $G$;

\rr2 $G(i,j)\cong G(i',j')$ if and only if there exist $r,s\in
F_p^*$ such that $j'=s^2j$ and $i'=r^2si$.
\end{enumerate}
\end{lem}

\demo (1) Let $$K=\lg a_2,a_1;d \di
a_2^{p}=a_1^{p}=d^p=1,[a_2,a_1]=d^j,[d,a_1]=[d,a_2]=1\rg\times \lg
a_3\rg.$$ We define an automorphism $\beta$ of $K$ as follows:
$$a_2^{\beta}=a_2a_3,\ a_1^\beta=a_1a_2\ {\rm and}\ a_3^\beta=a_3d^i.$$
Then $o(\beta)=p$. It follows from the cyclic extension theory that
$G=\lg K,b\rg$ is a cyclic extension of $K$ by $C_{p}$. Hence
$|G|=p^{5}$. It is easy to verify that $G'=\lg a_2,a_3,b^p\rg$,
$G_3=\lg a_3,b^p\rg$, $G_4=\lg b^p \rg$, $G\in \mathcal{A}_3$, and
$M=\lg b^p,a_1,a_2,a_3\rg$ is the unique three-generator maximal
subgroup of $G$.

(2) For convenience, let $G=\lg b,a_1,a_2,a_3\rg\cong G(i,j)$,
$\bar{G}=\lg \bar{b},\bar{a}_1,\bar{a}_2,\bar{a}_3\rg\cong G(i',j')$
and $\theta$ be an isomorphism from $\bar{G}$ to $G$. Since $M=\lg
b^p,a_1,a_2,a_3\rg$ and $\bar{M}=\lg
\bar{b}^p,\bar{a}_1,\bar{a}_2,\bar{a}_3\rg$ are the unique
three-generator maximal subgroups of $G$ and $\bar{G}$ respectively,
we have $\bar{M}^\theta=M$.  Hence we may assume that
$\bar{b}^\theta=b^{r}x$ and $\bar{a}_1^\theta=a_1^{s}y$, where
$r,s\in F_p^*$ and $x\in M$, $y\in \Phi(G)$.

By calculation we have
$\bar{a}_2^\theta=[\bar{a}_1,\bar{b}]^\theta\equiv a_2^{rs}\ (\mod
G_3)$, $\bar{a}_3^\theta=[\bar{a}_2,\bar{b}]^\theta\equiv
a_3^{r^2s}\ (\mod G_4)$.

Since $[\bar{a}_2^\theta,\bar{a}_1^\theta]=
(\bar{b}^{j'p})^\theta=(\bar{b}^\theta)^{j'p}$, we have
$[a_2^{rs},a_1^{s}]=b^{rj'p}$. Left side of the equation is
$b^{rs^2jp}$. By comparing index of $b^{p}$ in two sides we have
$j'=s^2j$.

Since
$[\bar{a}_3^\theta,\bar{b}^\theta]=(\bar{b}^{i'p})^\theta=(\bar{b}^\theta)^{i'p}$,
we have $[a_3^{r^2s},b^{r}]=b^{ri'p}$. Left side of the equation is
$b^{r^3si p}$. By comparing index of $b^{p}$ in two sides we have
$i'=r^2si$.

On the other hand, if there exist $r,s\in F_p^*$ such that $j'=s^2j$
and $i'=r^2si$, then, $\theta:$ $\bar{a_1}\rightarrow a_1^{s}$,
$\bar{b}\rightarrow b^{r}$ is an isomorphism from $\bar{G}$ to $G$.
\qed

\begin{lem}\label{cong 0a}
Suppose that $G(i, j)=\lg b,a_1;a_2,a_3\di
b^{p}=a_1^{p^2}=a_2^p=a_3^p=1,[a_1,b]=a_2,[a_2,b]=a_3,[a_3,b]=a_1^{ip},[a_2,a_1]=a_1^{jp},[a_3,a_1]=1\rg$,
where $i,j\in F_p^*$ and $p\ge 5$.
\begin{enumerate}
\rr1 If $G=G(i,j)$, then $|G|=p^{5}$, $G'=\lg a_2,a_3,a_1^p\rg$,
$G_3=\lg a_3,a_1^p\rg$, $G_4=\lg a_1^p \rg$, $G\in \mathcal{A}_3$,
and $M=\lg a_1,a_2,a_3\rg$ is the unique three-generator maximal
subgroup of $G$;

\rr2 $G(i,j)\cong G(i',j')$ if and only if there exist $r,s\in
F_p^*$ such that $j'=rsj$ and $i'=r^3i$.
\end{enumerate}
\end{lem}

\demo (1) Let $$K=\lg a_2,a_1 \di
a_2^{p}=a_1^{p^2}=1,[a_2,a_1]=a_1^{ip}\rg\times \lg a_3\rg.$$ We
define an automorphism $\beta$ of $K$ as follows:
$$a_2^{\beta}=a_2a_3,\ a_1^\beta=a_1a_2\ {\rm and}\
a_3^\beta=a_3a_1^{ip}.$$ Then $o(\beta)=p$. Hence $G=K\rtimes \lg
b\rg$ and $|G|=p^{5}$. It is easy to verify that $G'=\lg
a_2,a_3,a_1^p\rg$, $G_3=\lg a_3,a_1^p\rg$, $G_4=\lg a_1^p \rg$,
$G\in \mathcal{A}_3$, and $M=\lg a_1,a_2,a_3\rg$ is the unique
three-generator maximal subgroup of $G$.

\medskip
(2) For convenience, let $G=\lg b,a_1,a_2,a_3\rg\cong G(i,j)$,
$\bar{G}=\lg \bar{b},\bar{a}_1,\bar{a}_2,\bar{a}_3\rg\cong G(i',j')$
and $\theta$ be an isomorphism from $\bar{G}$ to $G$. Since $M=\lg
a_1,a_2,a_3\rg$ and $\bar{M}=\lg \bar{a}_1,\bar{a}_2,\bar{a}_3\rg$
are the unique three-generator maximal subgroups of $G$ and
$\bar{G}$ respectively, we have $\bar{M}^\theta=M$.  Hence we may
assume that $\bar{b}^\theta=b^{r}x$, $\bar{a}_1^\theta=a_1^{s}y$,
where $r,s\in F_p^*$ and $x\in \Omega_1(M)$, $y\in \Phi(G)$.

By calculation we have
$\bar{a}_2^\theta=[\bar{a}_1,\bar{b}]^\theta\equiv a_2^{rs}\ (\mod
G_3)$. $\bar{a}_3^\theta=[\bar{a}_2,\bar{b}]^\theta\equiv
a_3^{r^2s}\ (\mod G_4)$.

Since $[\bar{a}_2^\theta,\bar{a}_1^\theta]=
(\bar{a}_1^{j'p})^\theta=(\bar{a}_1^\theta)^{j'p}$, we have
$[a_2^{rs},a_1^{s}]=a_1^{sj'p}$. Left side of the equation is
$a_1^{rs^2jp}$. By comparing index of $a_1^{p}$ in two sides we have
$j'=rsj$.

Since
$[\bar{a}_3^\theta,\bar{b}^\theta]=(\bar{a}_1^{i'p})^\theta=(\bar{a}_1^\theta)^{i'p}$,
we have $[a_3^{r^2s},b^{r}]=a_1^{si'p}$. Left side of the equation
is $a_1^{r^3si p}$. By comparing index of $a_1^{p}$ in two sides we
have $i'=r^3i$.

On the other hand, if there exist $r,s\in F_p^*$ such that $j'=rsj$
and $i'=r^3i$, then $\theta:$ $\bar{a_1}\rightarrow a_1^{s}$,
$\bar{b}\rightarrow b^{r}$ is an isomorphism from $\bar{G}$ to $G$.
\qed

\begin{lem}\label{cong 1}
Suppose that $p\ge 3$ and $G=G(r,s)=\lg b,a_1;a_2\di
b^{p^{3}}=a_1^{p^2}=a_2^{p^2}=1,
[a_1,b]=a_2,[a_2,a_1]=b^{\nu_1p^{2}}a_2^{rp},[a_2,b]=a_1^{\nu_2p}a_2^{sp}\rangle$,
where $r,s\in F_p$ and $\nu_1,\nu_2=1$ or a fixed quadratic
non-residue modula $p$. Then
\begin{enumerate}
\rr1 $|G|=p^{7}$, $G'=\lg a_1^p,a_2,b^{p^2}\rg$, $G_3=\lg
a_1^p,a_2^p,b^{p^2}\rg$, $G_4=\lg a_2^p \rg$, $Z(G)\cap G'=\lg
a_2^p,b^{p^2}\rg$ and $M=\lg b^p,a_1,a_2\rg$ is the unique
three-generator maximal subgroup of $G$;

\rr2 $G$ is one of the following non-isomorphic groups:
\begin{enumerate}
\rr{i} $\lg b,a_1;a_2\di b^{p^{3}}=a_1^{p^2}=a_2^{p^2}=1,
[a_1,b]=a_2,[a_2,a_1]=b^{\nu_1p^{2}},[a_2,b]=a_1^{\nu_2p}a_2^{sp}\rangle$,
where $\nu_1,\nu_2=1$ or a fixed quadratic non-residue modula $p$
and $s=\nu_2,\nu_2+1,\dots,\nu_2+\frac{p-1}{2}$;

\rr{ii} $\lg b,a_1;a_2\di b^{p^{3}}=a_1^{p^2}=a_2^{p^2}=1,
[a_1,b]=a_2,[a_2,a_1]=b^{\nu_1p^{2}}a_2^{rp},[a_2,b]=a_1^{\nu_2p}\rangle$,
where $\nu_1,\nu_2=1$ or a fixed quadratic non-residue modula $p$
and $r=1,2,\dots,\frac{p-1}{2}$.
\end{enumerate}

\rr3 $[a_1^p,b]=a_2^p$;

\rr4 all maximal subgroups of $G$, except $M$, are
$\mathcal{A}_2$-groups;

\rr5 If $p=3$ and $\nu_2=-1$, then $|M'|=3$ and $G\in
\mathcal{A}_3$;

\rr6 If $p=3$ and $\nu_2=1$, then $|M'|=9$, and $G\in \mathcal{A}_3$
if and only if $r^2+4\nu_1$ is not a square;

\rr7 If $p\ge 5$, then $|M'|=p^2$, and $G\in A_3$ if and only if
$r^2-4\nu_1$ is not a square.

\end{enumerate}
\end{lem}

\demo (1) Let $$K=\lg a_2,a_1;d \di
a_2^{p^2}=a_1^{p^2}=d^p=1,[a_2,a_1]=d\rg.$$ We define an
automorphism $\beta$ of $K$ as follows:
$$a_2^{\beta}=a_2a_1^{\nu_2p}a_2^{sp}, \ a_1^\beta=a_1a_2.$$ Then
$$a_2^{\beta^p}=a_2, \ a_1^{\beta^{p^2}}=1,\ o(\beta)\le p^2.$$ It
follows from the cyclic extension theory that $G=\lg K,b\rg$ is a
cyclic extension of $K$ by $C_{p^2}$. Hence $|G|=p^{7}$. It is easy
to check that $G'=\lg a_1^p,a_2,b^{p^2}\rg$, $G_3=\lg
a_1^p,a_2^p,b^{p^2}\rg$, $G_4=\lg a_2^p \rg$, $Z(G)\cap G'=\lg
a_2^p,b^{p^2}\rg$ and $M=\lg b^p,a_1,a_2\rg$ is the unique
three-generator maximal subgroup of $G$.

\medskip
(2) For convenience, let $G=\lg b,a_1,a_2\rg\cong G(r,s)$,
$\bar{G}=\lg \bar{b},\bar{a}_1,\bar{a}_2\rg\cong G(r',s')$ and
$\theta$ be an isomorphism from $\bar{G}$ to $G$. Since $M=\lg
b^p,a_1,a_2\rg$ and $\bar{M}=\lg \bar{b}^p,\bar{a}_1,\bar{a}_2\rg$
are the unique three-generator maximal subgroups of $G$ and
$\bar{G}$ respectively, we have $\bar{M}^\theta=M$. Hence we may
assume that $\bar{b}^\theta=b^{l}a_1^mx$,
$\bar{a}_1^\theta=a_1^{i}a_2^{j}b^{kp}y$, where $l,i\in F_p^*$ and
$x\in \Phi(G)$, $y\in \mho_1(M)$.

By calculation, $\bar{a}_2^\theta=[\bar{a}_1,\bar{b}]^\theta\equiv
a_2^{il}\ (\mod G_3)$.

Since $[\bar{a}_2^\theta,\bar{a}_1^\theta]=
(\bar{b}^{\nu_1p^2}\bar{a}_2^{r'p})^\theta=(\bar{b}^\theta)^{\nu_1p^2}(\bar{a}_2^\theta)^{r'p}$,
we have $[a_2^{il},a_1^{i}]=b^{l\nu_1p^2}a_2^{ilr'p}$. Left side of
above equation is $b^{i^2l\nu_1p^{2}}a_2^{i^2lrp}$. Comparing
indexes of $b^{p^2}$ and $a_2^{p}$ in two sides, we have $i^2=1$ and
$r'=ir$.

Since $[\bar{a}_2^\theta,\bar{b}^\theta,\bar{b}^\theta]=
(\bar{a}_2^{\nu_2p})^\theta=(\bar{a}_2^\theta)^{\nu_2p}$, we have
$[a_2^{il},b^{l},b^{l}]=a_2^{il\nu_2p}$. Left side of above equation
is $a_2^{il^3\nu_2p}$. Comparing index of $a_2^{p}$ in two sides, we
have $l^2=1$.

By calculation we have
$$\bar{a}_2^\theta=[\bar{a}_1,\bar{b}]^\theta=[a_1^{i}a_2^{j}b^{kp}y,b^la_1^m]\equiv
[a_1^ia_2^j,b^l] \ (\mod Z(G)\cap G').$$

Moreover, by calculation we have
$$\bar{a}_2^\theta\equiv a_2^ia_1^{j\nu_2p} \ (\mod Z(G)\cap G')\ {\rm
for}\ l=1$$ and
$$\bar{a}_2^\theta\equiv a_2^{-i}a_1^{(i-j)\nu_2p} \ (\mod Z(G)\cap
G')\ {\rm for}\ l=-1.$$

We have $[\bar{a}_2^\theta,\bar{b}^\theta]=
(\bar{a}_1^{\nu_2p}\bar{a}_2^{s'p})^\theta=(\bar{a}_1^\theta)^{\nu_2p}(\bar{a}_2^\theta)^{s'p}$.

If $l=1$, then we have
$[a_2^{i}a_1^{j\nu_2p},ba_1^m]=(a_1^ia_2^jb^{kp})^{\nu_2p}a_2^{is'p}$.
Left side of above equation is
$b^{im\nu_1p^{2}}a_2^{imrp}a_1^{i\nu_2p}a_2^{isp}a_2^{jv_2p}$.
Comparing indexes of $b^{p^2}$, $a_1^p$ and $a_2^{p}$ in two sides,
we have $s'=s+ik\nu_2\nu_1^{-1}r$.

If $l=-1$, then we have
$[a_2^{-i}a_1^{(i-j)\nu_2p},b^{-1}a_1^m]=(a_1^ia_2^jb^{kp})^{\nu_2p}a_2^{-is'p}$.
Left side of above equation is
$b^{-im\nu_1p^{2}}a_2^{-imrp}a_1^{i\nu_2p}a_2^{isp}a_2^{-i\nu_2p}a_2^{(j-i)\nu_2p}$.
Comparing indexes of $b^{p^2}$, $a_1^p$ and $a_2^{p}$ in two sides,
we have $s'+s=2\nu_2+ik\nu_2\nu_1^{-1}r$.

On the other hand, if there exists $i,k$ such that $i^2=1$, $r'=ir$
and $s'+s=2\nu_2+ik\nu_2\nu_2^{-1}r$, then, $\theta:$
$\bar{a}_1\rightarrow a_1^{i}b^{kp}$, $\bar{b}\rightarrow
b^{-1}a_1^{-ik\nu_2\nu_1^{-1}}$ is an isomorphism from $\bar{G}$ to
$G$.

By above argument, if $r=0$, then $G$ is a group of Type (i), and if
$r\ne 0$, then $G$ is a group of Type (ii).

\medskip

(3) By calculation we have $[a_1^p,b]=[a_1,b]^p[a_1,b,a_1]^{p\choose
2}[a_1,b,a_1,a_1]^{p\choose 3}=a_2^p$.

\medskip

(4) Let $N$ be a maximal subgroup of $G$ such that $N\ne M$. Then we
may assume $N=\lg ba_1^i,\Phi(G)\rg$. By calculation we have
$[a_2,ba_1^i,ba_1^i]=[a_1^{\nu_2 p},ba_1^i]=a_2^{\nu_2 p}$. Hence
$N=\lg ba_1^i,a_2\rg$ is isomorphic to a group of Type (6) of Lemma
\ref{A_2}. Thus $N\in\A_2$.

\medskip

(5) By calculation,
$[a_1,b^3]=[a_1,b]^3[a_1,b,b]^3[a_1,b,b,b]=a_2^3[a_1^{-3},b]=1$.
Since $M'=\lg [a_1,b^3],[a_2,a_1]\rg$, $|M'|=3$. Moreover, $M=\lg
a_2,a_1\rg\ast\lg b^{3\nu_1 }a_2^r\rg$. By Lemma \ref{A_t+k} (2) we
get $M\in\A_2$. By (4) we get $G\in\A_3$.

\medskip

(6) By calculation,
$[a_1,b^3]=[a_1,b]^3[a_1,b,b]^3[a_1,b,b,b]=a_2^3[a_1^{3},b]=a_2^{-3}$.
Since $M'=\lg [a_1,b^3],[a_2,a_1]\rg$, $|M'|=9$. By (4) we get
$G\in\A_3$ if and only if $M\in\A_2$. Moreover, we claim that
$M\in\A_2$ if and only if $r^2+4\nu_1$ is not a square.

If $M\in\A_2$, then $\lg a_2a_1^x,b^3a_1^y,\Phi(M)\rg$ is either
abelian or minimal non-abelian, where $x,y\in F_3$. For $x\ne 0$ or
$y\ne 0$, since
$[a_2a_1^x,b^3a_1^y]=b^{9y\nu_1}a_2^{3yr}a_2^{-3x}\ne 1$, we have
$\lg a_2a_1^x,b^3a_1^y \rg$ is minimal non-abelian. Hence
$\Phi(M)\le \lg a_2a_1^x,b^3a_1^y \rg$. Moreover, $\Phi(M)=\lg
a_1^{3x}a_2^3, a_1^{3y}b^9, a_2^{3(yr-x)}b^{9y\nu_1}\rg$. It follows
that the equation $\left| {\begin{array}{*{20}{c}}
   x & 1 & 0  \\
   y & 0 & 1  \\
   0 & {yr - x} & {y{\nu _1}}  \\
\end{array}} \right| = 0$ only has solution $x=y=0$.
That is, the equation $x^2-xyr-y^2\nu_1=0$ only has solution $x= y=
0$. Hence $r^2+4\nu_1$ is not a square. Conversely, if $r^2+4\nu_1$
is not a square, then, by above argument, $\lg
a_2a_1^x,b^3a_1^y,\Phi(M)\rg$ is either abelian or minimal
non-abelian. It is easy to see that $\lg a_1,b^3,\Phi(M) \rg$ and
$\lg a_1,a_2b^{3z}, \Phi(M)\rg$ are minimal non-abelian. Hence all
maximal subgroups of $M$ are abelian or minimal non-abelian. That is
, $M\in \A_2$.
\medskip

(7) By calculation, $[a_1,b^p]=[a_1,b]^p=a_2^p$. Since $M'=\lg
[a_1,b^p],[a_2,a_1]\rg$, $|M'|=p^2$. By (4) we get $G\in\A_3$ if and
only if $M\in\A_2$. Moreover, we claim that $M\in\A_2$ if and only
if $r^2-4\nu_1$ is not a square.

If $M\in\A_2$, then $\lg a_2a_1^x,b^pa_1^y,\Phi(M)\rg$ is either
abelian or minimal non-abelian, where $x,y\in F_p$. For $x\ne 0$ or
$y\ne 0$, since $[a_2a_1^x,b^pa_1^y]=b^{y\nu_1 p^2}a_2^{yrp}a_2^{x
p}\ne 1$, we have $\lg a_2a_1^x,b^pa_1^y \rg$ is minimal
non-abelian. Hence $\Phi(M)\le \lg a_2a_1^x,b^pa_1^y \rg$. Moreover,
$\Phi(M)=\lg a_1^{xp}a_2^p, a_1^{yp}b^{p^2}, a_2^{(yr+x)p}b^{y\nu_1
p}\rg$. It follows that the equation $\left|
{\begin{array}{*{20}{c}}
   x & 1 & 0  \\
   y & 0 & 1  \\
   0 & {yr+x} & {y{\nu _1}}  \\
\end{array}} \right| =0$ only has solution $x= y=0$.
That is, the equation $x^2+xyr+y^2\nu_1=0$ only has solution $x= y=
0$. Hence $r^2-4\nu_1$ is not a square. Conversely, if $r^2-4\nu_1$
is not a square, then, by above argument, $\lg
a_2a_1^x,b^pa_1^y,\Phi(M)\rg$ is either abelian or minimal
non-abelian. It is easy to see that $\lg a_1,b^p,\Phi(M) \rg$ and
$\lg a_1,a_2b^{zp}, \Phi(M)\rg$ are minimal non-abelian. Hence all
maximal subgroups of $M$ are abelian or minimal non-abelian. That is
, $M\in \A_2$.
 \qed

\begin{lem}\label{cong 2}
Suppose that $G(\nu, t)=\lg b,a_1;a_2\di
a_1^{p^2}=a_2^{p^2}=1,b^{p^2}=a_2^{tp},
[a_1,b]=a_2,[a_2,a_1]=a_2^{p},[a_2,b]=a_1^{\nu p}\rangle$, where
$t\in F_p$, $\nu=1$ or a fixed quadratic non-residue modulo $p$.
\begin{enumerate}
\rr1 If $G=G(\nu,t)$, then $|G|=p^{6}$, $G'=\lg a_1^p,a_2\rg$,
$G_3=\lg a_1^p,a_2^p\rg$, $G_4=\lg a_2^p \rg$, $G\in \mathcal{A}_3$,
and $M=\lg b^p,a_1,a_2\rg$ is the unique three-generator maximal
subgroup of $G$;

\rr2 $G(\nu,t)\cong G(\nu',t')$ if and only if $(\nu',t')=(\nu,t)$.
\end{enumerate}
\end{lem}

\demo (1) Let $$K=\lg a_2,a_1 \di
a_2^{p^2}=a_1^{p^2}=1,[a_2,a_1]=a_2^p\rg.$$ We define an
automorphism $\beta$ of $K$ as follows: $$a_2^{\beta}=a_2a_1^{\nu
p},\ a_1^\beta=a_1a_2.$$ Then $a_2^{\beta^p}=a_2$,
$a_1^{\beta^{p^2}}=a_1$, $o(\beta)\le p^2$. It follows from the
cyclic extension theory that $G=\lg K,b\rg$ is a cyclic extension of
$K$ by $C_{p^2}$. Hence $|G|=p^{6}$. It is easy to check that
$G'=\lg a_1^p,a_2\rg$, $G_3=\lg a_1^p,a_2^p\rg$, $G_4=\lg a_2^p \rg$
and $M=\lg b^p,a_1,a_2\rg$ is the unique three-generator maximal
subgroup of $G$.

\medskip
(2) For convenience, let $G=\lg b,a_1,a_2\rg\cong G(\nu,t)$,
$\bar{G}=\lg \bar{b},\bar{a}_1,\bar{a}_2\rg\cong G(\nu',t')$ and
$\theta$ be an isomorphism from $\bar{G}$ to $G$. Since $M=\lg
b^p,a_1,a_2\rg$ and $\bar{M}=\lg \bar{b}^p,\bar{a}_1,\bar{a}_2\rg$
are the unique three-generator maximal subgroups of $G$ and
$\bar{G}$ respectively, we have $\bar{M}^\theta=M$.  Hence we may
assume that $\bar{b}^\theta=b^{l}x$ and $\bar{a}_1^\theta=a_1^{i}y$,
where $l,i\in F_p^*$ and $x\in M$, $y\in \Phi(G)$.

By calculation, $\bar{a}_2^\theta=[\bar{a}_1,\bar{b}]^\theta\equiv
a_2^{il}\ (\mod G_3)$.

Since $[\bar{a}_2^\theta,\bar{a}_1^\theta]=
(\bar{a}_2^{p})^\theta=(\bar{a}_2^\theta)^{p}$, we have
$[a_2^{il},a_1^{i}]=a_2^{ilp}$. Left side of above equation is
$a_2^{i^2lp}$. Comparing index of $a_2^{p}$ in two sides, we have
$i=1$.

Since
$(\bar{b}^\theta)^{p^2}=(\bar{b}^{p^2})^\theta=(\bar{a}_2^{t'p})^\theta=(\bar{a}_2^\theta)^{t'p}$,
we have $(b^l)^{p^2}=a_2^{t'lp}$. Left side of above equation is
$a_2^{ltp}$. Comparing index of $a_2^{3}$ in two sides, we have
$t'=t$.

Since $[\bar{a}_2^\theta,\bar{b}^\theta,\bar{b}^\theta]=
(\bar{a}_2^{\nu'p})^\theta=(\bar{a}_2^\theta)^{\nu'p}$, we have
$[a_2^{l},b^{l},b^{l}]=a_2^{\nu'lp}$. Left side of above equation is
$a_2^{l^3\nu p}$. Comparing index of $a_2^{p}$ in two sides, we have
$\nu'=\nu$. \qed

\begin{lem}\label{cong 3}
Suppose that $G(\nu,s,t)=\lg b,a_1;a_2\di
a_1^{p^2}=a_2^{p^2}=1,b^{p^m}=a_2^{tp},
[a_1,b]=a_2,[a_2,a_1]=1,[a_2,b]=a_1^{\nu p}a_2^{sp}\rangle$, where
$m\ge 2$, $t\in F_p$, $\nu=1$ or a fixed quadratic non-residue
modulo $p$.
\begin{enumerate}
\rr1 If $G=G(\nu,s,t)$, then $|G|=p^{m+4}$, $G'=\lg a_1^p,a_2\rg$,
$G_3=\lg a_1^p,a_2^p\rg$, $G_4=\lg a_2^p \rg$, $G\in\mathcal{A}_3$,
$M=\lg b^p,a_1,a_2\rg$ is the unique three-generator maximal
subgroup of $G$;

\rr2 $G$ is one of the following non-isomorphic groups:

\begin{enumerate}
\rr{i} $\lg b,a_1;a_2\di a_1^{p^2}=a_2^{p^2}=b^{p^m}=1,
[a_1,b]=a_2,[a_2,a_1]=1,[a_2,b]=a_1^{\nu p}a_2^{sp}\rangle$, where
$m\ge 2$, $\nu=1$ or a fixed quadratic non-residue modulo $p$,
$s=\nu,\nu+1,\dots,\nu+\frac{p-1}{2}$;

\rr{ii} $\lg b,a_1;a_2\di a_1^{p^2}=a_2^{p^2}=1,b^{p^m}=a_2^{p},
[a_1,b]=a_2,[a_2,a_1]=1,[a_2,b]=a_1^{\nu p}\rangle$, where $m\ge 2$,
$\nu=1$ or a fixed quadratic non-residue modulo $p$.
\end{enumerate}

\rr3 If $p=3$ and $\nu=-1$, then $M$ is abelian. If $\nu=1$ or $p\ge
5$, then $|M'|=p$.

\end{enumerate}
\end{lem}

\demo (1) Let $$K=\lg a_2,a_1 \di
a_2^{p^2}=a_1^{p^2}=1,[a_2,a_1]=a_2^p\rg.$$ We define an
automorphism $\beta$ of $K$ as follows: $$a_2^{\beta}=a_2a_1^{\nu
p}a_2^{sp},\ a_1^\beta=a_1a_2.$$ Then $a_2^{\beta^p}=a_2$,
$a_1^{\beta^{p^2}}=a_1$, $o(\beta)\le p^2$. It follows from the
cyclic extension theory that $G=\lg K,b\rg$ is a cyclic extension of
$K$ by $C_{p^m}$. Hence $|G|=p^{m+4}$. It is easy to verify that
$G'=\lg a_1^p,a_2\rg$, $G_3=\lg a_1^p,a_2^p\rg$, $G_4=\lg a_2^p
\rg$, $G\in\mathcal{A}_3$ and $M=\lg b^p,a_1,a_2\rg$ is the unique
three-generator maximal subgroup of $G$.

\medskip
(2) For convenience, let $G=\lg b,a_1,a_2\rg\cong G(\nu,s,t)$,
$\bar{G}=\lg \bar{b},\bar{a}_1,\bar{a}_2\rg\cong G(\nu',s',t')$ and
$\theta$ be an isomorphism from $\bar{G}$ to $G$. Since $M=\lg
b^p,a_1,a_2\rg$ and $\bar{M}=\lg \bar{b}^p,\bar{a}_1,\bar{a}_2\rg$
are the unique three-generator maximal subgroups of $G$ and
$\bar{G}$ respectively, we have $\bar{M}^\theta=M$.  Hence we may
assume that $\bar{b}^\theta=b^{l}x$ and
$\bar{a}_1^\theta=a_1^{i}a_2^jb^{kp^{m-1}}y$, where $l,i\in F_p^*$
and $x\in M$, $y\in G_3$.

By calculation we have
$\bar{a}_2^\theta=[\bar{a}_1,\bar{b}]^\theta\equiv a_2^{il}\ (\mod
G_3)$.

Since
$(\bar{b}^\theta)^{p^m}=(\bar{b}^{p^m})^\theta=(\bar{a}_2^{t'p})^\theta=(\bar{a}_2^\theta)^{t'p}$,
we have $(b^lx)^{p^m}=a_2^{t'ilp}$. Left side of the equation is
$a_2^{ltp}$. By comparing index of $a_2^{3}$ in two sides we have
$t'i=t$.

Since $[\bar{a}_2^\theta,\bar{b}^\theta,\bar{b}^\theta]=
(\bar{a}_2^{\nu'p})^\theta=(\bar{a}_2^\theta)^{\nu'p}$, we have
$[a_2^{l},b^{l},b^{l}]=a_2^{\nu'lp}$. Left side of the equation is
$a_2^{l^3\nu p}$. By comparing index of $a_2^{p}$ in two sides we
have $\nu'=\nu$ and $l^2=1$.

By calculation we have
$\bar{a}_2^\theta=[\bar{a}_1,\bar{b}]^\theta=[a_1^{i}a_2^{j}b^{kp^{m-1}}y,b^l]\equiv
[a_1^ia_2^j,b^l] \ (\mod G_4)$. Moreover, by calculation we have
$\bar{a}_2^\theta\equiv a_2^ia_1^{j\nu p} \ (\mod G_4)$ for $l=1$
and $\bar{a}_2^\theta\equiv a_2^{-i}a_1^{(i-j)\nu p} \ (\mod G_4)$
for $l=-1$.

 If $l=1$, then $$[a_2^{i}a_1^{j\nu
p},b]=[\bar{a}_2^\theta,\bar{b}^\theta]= (\bar{a}_1^{\nu
p}\bar{a}_2^{s'p})^\theta=(\bar{a}_1^\theta)^{\nu
p}(\bar{a}_2^\theta)^{s'p}=(a_1^ia_2^jb^{kp^{m-1}})^{\nu
p}a_2^{is'p}.$$ Left side of the equation is $a_1^{i\nu
p}a_2^{isp}a_2^{j\nu p}$. By comparing indexes of $a_1^p$ and
$a_2^{p}$ in two sides we have $s'=s-i^{-1}tk\nu$. On the other
hand, if there exists $i,k$ such that $t'i=t$ and
$s'=s-i^{-1}tk\nu$, then, $\theta:$ $\bar{a}_1\rightarrow
a^{i}b^{kp^{m-1}}$, $\bar{b}\rightarrow b$ is an isomorphism from
$\bar{G}$ to $G$.

If $l=-1$, then we have $$[a_2^{-i}a_1^{(i-j)\nu
p},b^{-1}]=[\bar{a}_2^\theta,\bar{b}^\theta]= (\bar{a}_1^{\nu
p}\bar{a}_2^{s'p})^\theta=(\bar{a}_1^\theta)^{\nu
p}(\bar{a}_2^\theta)^{s'p}=(a_1^ia_2^jb^{kp^{m-1}})^{\nu
p}a_2^{-is'p}.$$ Left side of the equation is $a_1^{i\nu
p}a_2^{isp}a_2^{-i\nu p}a_2^{(j-i)\nu p}$. By comparing indexes of
$a_1^p$ and $a_2^{p}$ in two sides we have $s'+s=2\nu +i^{-1}tk\nu$.
On the other hand, if there exists $i,k$ such that $t'i=t$ and
$s'+s=2\nu+i^{-1}tk\nu$, then, $\theta:$ $\bar{a}_1\rightarrow
a^{i}b^{kp^{m-1}}$, $\bar{b}\rightarrow b^{-1}$ is an isomorphism
from $\bar{G}$ to $G$.

By above argument we get that $G(\nu,s,0)\cong G(\nu',s',t')$ if and
only if $\nu=\nu'$, $t'=0$ and $s=s'$ or $s+s'=2\nu$. It follows
that $G$ is a group of Type (i)  for $t=0$. If $t\ne 0$, then, by
above argument, $G(\nu, s,t)\cong G(\nu,0,1)$. Hence $G$ is a group
of Type (ii).

\medskip
(3) By calculation, $[a_1^{\nu p},b]=a_2^{\nu p}$. Hence
$$[a_1,b^3]=[a_1,b]^3[a_1,b,b]^3[a_1,b,b,b]=a_2^3[a_1^{3\nu},b]=a_2^{3(1+\nu)}.$$
If $p\ge 5$, then $[a_1,b^p]=a_2^p$. Since $M'=\lg [a_1,b^p]\rg$,
the result is proved.
 \qed

\begin{thm}\label{d=4-3}
Suppose that $G$ is an $\mathcal{A}_3$-group all of whose maximal
subgroups are $\mathcal{A}_2$-groups and there exists a
three-generator maximal subgroup in $G$, and $M'\le Z(G)$ for every
three-generator maximal subgroup $M$. Then $d(G)=2$ if and only if
 $G$ is isomorphic to one of the following pairwise
non-isomorphic groups:
\begin{enumerate}

\rr{Mi} $G'\cong C_{p^2}$. In this case, $(\mu_0,\mu_1,\mu_2)=(0,0,1+p)$ and $\alpha_1(G)=p^3+p^2$.

\begin{enumerate}
\rr{M1}  $\langle a, b; c \di a^{8}=1, c^2=a^{4}=b^4,
[a,b]=c,[c,a]=1,[c,b]=1 \rangle$;  where  $|G|=2^{6}$, $c(G)=2$,
$\Phi(G)=\lg a^2, b^2,c\rg\cong C_4\times C_2^2$, $G'=\lg c\rg$ and
$Z(G)=\lg c\rg\cong C_4$.

\rr{M2} $\langle a, b; c \di a^{8}=b^{4}=1, c^2=a^{4},
[a,b]=c,[c,a]=1,[c,b]=1 \rangle$; moreover,
$|G|=2^{6}$, $c(G)=2$, $\Phi(G)=\lg a^2, b^2,c\rg\cong
C_{4}\times C_{2}^2$, $G'=Z(G)=\lg c\rg\cong C_{4}$.

\rr{M3} $\langle a, b; c \di a^{p^{3}}=b^{p^2}=1, c^p=a^{p^2},
[a,b]=c,[c,a]=1,[c,b]=c^{tp} \rangle$, where $p>2$ and $t\in F_p^*$;
moreover, $|G|=p^{6}$, $c(G)=3$, $\Phi(G)=\lg a^p, b^p,c\rg\cong
C_{p^2}\times C_p^2$, $G'=\lg c\rg$ and
$Z(G)=\lg c^p\rg \cong C_p$.

\rr{M4} $\langle a, b; c \di a^{p^{3}}=b^{p^2}=1, c^p=a^{p^2},
[a,b]=c,[c,a]=c^p,[c,b]=1 \rangle$, where $p>2$;
moreover, $|G|=p^{6}$, $c(G)=3$, $\Phi(G)=\lg a^p, b^p,c\rg\cong
C_{p^2}\times C_p^2$, $G'=\lg c\rg$ and
$Z(G)=\lg a^{p^2}\rg\cong C_{p}$.

\rr{M5} $\langle a, b; c \di a^{p^{n+1}}=b^{p^n}=1, c^p=a^{p^n},
[a,b]=c,[c,a]=1,[c,b]=1 \rangle$, where $n\ge 3$ for $p=2$ and $n\ge 2$; moreover,
$|G|=p^{2n+2}$, $c(G)=2$, $\Phi(G)=\lg a^p, b^p,c\rg\cong
C_{p^{n}}\times C_{p^{n-1}}\times C_p$, $G'=\lg c\rg$, $Z(G)=\lg
a^{p^2},b^{p^2},c\rg \cong C_{p^{n-1}}\times C_{p^{n-2}}\times C_p$
if $n>2$,
$Z(G)=\lg a^{p^2},b^{p^2},c\rg \cong C_p^2$ if $n=2$.

\rr{M6} $\langle a, b; c \di a^{p^{2}}=b^{p^{2}}=c^{p^2}=1,
[a,b]=c,[c,a]=c^p,[c,b]=1 \rangle$, where $p>2$;  moreover, $|G|=p^{6}$,
$c(G)=3$, $\Phi(G)=\lg a^p, b^p,c\rg\cong C_{p^2}\times C_p^2$,
$G'=\lg c\rg$,
$Z(G)=\lg c^{p}\rg \cong C_p$.

\rr{M7} $\langle a, b; c \di a^{p^{n}}=b^{p^{n}}=c^{p^2}=1,
[a,b]=c,[c,a]=1,[c,b]=1 \rangle$, where $n\ge 3$ for $p=2$ and $n\ge 2$; moreover, $|G|=p^{2n+2}$, $c(G)=2$, $\Phi(G)=\lg a^p,
b^p,c\rg\cong C_{p^{n-1}}^2\times C_{p^2}$, $G'=\lg
c\rg$, $Z(G)=\lg a^{p^2},b^{p^2},c\rg \cong C_{p^{n-2}}^2\times C_{p^2}$ if $n>2$, $Z(G) \cong
C_{p^2}$ if  $n=2$.

\rr{M8} $\langle a, b; c \di a^{p^{n+1}}=b^{p^2}=1, c^p=a^{p^n},
[a,b]=c,[c,a]=c^p,[c,b]=1 \rangle$, where $p>2$ and $n>2$;
moreover, $|G|=p^{n+4}$, $c(G)=3$, $\Phi(G)=\lg a^p, b^p,c\rg\cong
C_{p^n}\times C_p^2$, $G'=\lg c\rg$ and
$Z(G)=\lg a^{p^2}\rg\cong C_{p^{n-1}}$.

\rr{M9} $\langle a, b; c \di a^{p^{n}}=b^{p^{3}}=1, c^p=b^{p^2},
[a,b]=c,[c,a]=c^{tp},[c,b]=1 \rangle$, where $p>2$, $n>2$ and $t\in
F_p^*$; moreover, $|G|=p^{n+4}$, $c(G)=3$, $\Phi(G)=\lg a^p,
b^p,c\rg\cong C_{p^{n-1}}\times C_{p^2}\times C_{p}$, $G'=\lg c\rg$
and
$Z(G)=\lg a^{p^2},b^{p^2}\rg\cong C_{p^{n-2}}\times C_{p}$.

\rr{M10} $\langle a, b; c \di a^{p^{n}}=b^{p^{2}}=c^{p^2}=1,
[a,b]=c,[c,a]=c^p,[c,b]=1 \rangle$, where $p>2$ and $n>2$. Moreover,
$|G|=p^{n+4}$, $c(G)=3$, $\Phi(G)=\lg a^p, b^p,c\rg \cong
C_{p^{n-1}}\times C_{p^2}\times C_{p}$, $G'=\lg c\rg$ and
$Z(G)=\lg a^{p^2},c^{p}\rg \cong C_{p^{n-2}}\times C_{p}$.

\rr{M11} $\langle a, b; c \di a^{p^{n+1}}=b^{p^m}=1, c^p=a^{p^n},
[a,b]=c,[c,a]=1,[c,b]=1 \rangle$, where $n>m\ge 2$; moreover,
$|G|=p^{n+m+2}$, $c(G)=2$, $\Phi(G)=\lg a^p, b^p,c\rg\cong
C_{p^{n}}\times C_{p^{m-1}}\times C_p$, $G'=\lg c\rg$, $Z(G)=\lg
a^{p^2},b^{p^2},c\rg \cong C_{p^{n-1}}\times C_{p^{m-2}}\times C_p$
if $m>2$,
$Z(G)=\lg a^{p^2},b^{p^2},c\rg \cong C_{p^{n-1}}\times C_p$ if $m=2$.

\rr{M12} $\langle a, b; c \di a^{p^{n}}=b^{p^{m+1}}=1, c^p=b^{p^m},
[a,b]=c,[c,a]=1,[c,b]=1 \rangle$, where $n>m\ge 2$; $|G|=p^{n+m+2}$,
$c(G)=2$, $\Phi(G)=\lg a^p, b^p,c\rg \cong C_{p^{n-1}}\times
C_{p^{m}}\times C_p$, $G'=\lg c\rg$ and
$Z(G)=\lg a^{p^2},b^{p^2},c\rg \cong C_{p^{n-2}}\times C_{p^{m-1}}\times C_p$.

\rr{M13} $\langle a, b; c \di a^{p^{n}}=b^{p^{m}}=c^{p^2}=1,
[a,b]=c,[c,a]=1,[c,b]=1 \rangle$, where $n>m\ge 2$; moreover, $|G|=p^{n+m+2}$, $c(G)=2$, $\Phi(G)=\lg a^p,
b^p,c\rg\cong C_{p^{n-1}}\times C_{p^{m-1}}\times C_{p^2}$, $G'=\lg
c\rg$, $Z(G)=\lg a^{p^2},b^{p^2},c\rg \cong C_{p^{n-2}}\times
C_{p^{m-2}}\times C_{p^2}$ if $m>2$, $Z(G) \cong C_{p^{n-2}}\times
C_{p^2}$ if $m=2$.
\end{enumerate}

\rr{Mii} $c(G)=3$ and $G'\cong C_p^2$. In this case, $(\mu_0,\mu_1,\mu_2)=(0,0,1+p)$ and $\alpha_1(G)=p^3+p^2$.

\begin{enumerate}
\rr{M14}  $\langle a, b; c \di a^{8}=b^{4}=c^2=1,
[a,b]=c,[c,a]=1,[c,b]=a^{4} \rangle$;  where $|G|=2^{6}$,
$\Phi(G)=\lg a^2, b^2,c\rg\cong C_4\times C_2^2$, $G'=\lg
a^{4},c\rg$ and $Z(G)=\lg a^2\rg\cong C_4$.

\rr{M15} $\langle a, b; c,d \di a^{4}=b^{4}=c^2=d^2=1,
[a,b]=c,[c,a]=d,[c,b]=1,[d,a]=[d,b]=1 \rangle$;
moreover, $|G|=2^{6}$, $\Phi(G)=\lg a^2, b^2,c,d\rg\cong C_2^4$, $G'=\lg c,d\rg$ and  $Z(G)=\lg b^2,d\rg\cong C_2^2$.

\rr{M16} $\langle a, b; c \di a^{8}=b^{4}=c^2=1,
[a,b]=c,[c,a]=a^{4},[c,b]=1 \rangle$, where
$|G|=2^{6}$, $\Phi(G)=\lg a^2, b^2,c\rg\cong C_{4}\times
C_2^2$, $G'=\lg a^{4},c\rg$ and $Z(G)=\lg a^4, b^2\rg\cong
C_2^2$.

\rr{M17} $\langle a, b; c \di a^{2^{n+1}}=b^{4}=c^2=1,
[a,b]=c,[c,a]=a^{2^n},[c,b]=1 \rangle$, where $n>2$; moreover,
$|G|=2^{n+4}$, $\Phi(G)=\lg a^2, b^2,c\rg\cong C_{2^{n}}\times
C_2^2$, $G'=\lg a^{2^n},c\rg$ and $Z(G)=\lg a^4, b^2\rg\cong
C_{2^{n-1}}\times C_2$.

\rr{M18} $\langle a, b; c \di a^{2^{n}}=b^{8}=c^2=1,
[a,b]=c,[c,a]=b^{4},[c,b]=1 \rangle$, where $n>2$; moreover,
$|G|=2^{n+4}$, $\Phi(G)=\lg a^2, b^2,c\rg\cong C_{2^{n-1}}\times C_4\times
C_2$, $G'=\lg b^{4},c\rg$ and $Z(G)=\lg a^4, b^2\rg\cong
C_{2^{n-2}}\times C_4$.

\rr{M19} $\langle a, b; c,d \di a^{2^{n}}=b^{4}=c^2=d^2=1,
[a,b]=c,[c,a]=d,[c,b]=1,[d,a]=[d,b]=1 \rangle$, where $n>2$;
moreover, $|G|=2^{n+4}$, $\Phi(G)=\lg a^2, b^2,c,d\rg\cong
C_{2^{n-1}}\times C_2^3$, $G'=\lg c,d\rg$ and  $Z(G)=\lg a^4,
b^2,d\rg\cong C_{2^{n-2}}\times C_2^2$.

\end{enumerate}

\rr{Miii} $G_3\cong C_p$ and $\Phi(G')G_3\cong C_p^2$. In this case, $(\mu_0,\mu_1,\mu_2)=(0,0,1+p)$ and $\alpha_1(G)=p^3+2p^2$.

\begin{enumerate}
\rr{M20} $\langle a, b; c\di a^{8}=b^{8}=1,c^{2}=b^4,[a,b]=c,[c,a]=1,[c,b]=a^4b^4\rangle$; where $|G|=2^{7}$, $\Phi(G)=\lg
a^2, b^2,c\rg\cong C_{4}^2\times C_{2}$, $G'=\lg
a^{4},c\rg$ and  $Z(G)=\lg a^{4},
c^2\rg\cong C_2^2$.

\rr{M21} $\langle a, b; c\di a^{8}=c^{4}=1,b^{4}=c^2,
[a,b]=c,[c,a]=a^{4},[c,b]=1\rangle$; where $|G|=2^{7}$, $\Phi(G)=\lg
a^2, b^2,c\rg\cong C_{4}^2\times C_{2}$, $G'=\lg
a^{4},c\rg$ and  $Z(G)=\lg a^{4},
c^2\rg\cong C_2^2$.

\rr{M22} $\langle a, b; c\di a^{8}=c^{4}=1,b^{4}=c^2,
[a,b]=c,[c,a]=1,[c,b]=a^{4}\rangle$; where $|G|=2^{7}$, $\Phi(G)=\lg a^2, b^2,c\rg\cong C_4^2\times C_2$, $G'=\lg a^{4},c\rg$,
$Z(G)=\lg a^{4}, c^2\rg\cong C_2^2$.

\rr{M23} $\langle a, b; c\di a^{4}=b^{8}=c^{4}=1,
[a,b]=c,[c,a]=1,[c,b]=b^4\rangle$; where $|G|=2^{7}$, $\Phi(G)=\lg
a^2, b^2,c\rg\cong C_{4}^2\times C_{2}$, $G'=\lg
b^{4},c\rg$ and  $Z(G)=\lg b^{4},
c^2\rg\cong C_2^2$.

\rr{M24} $\langle a, b; c\di a^{p^{3}}=b^{p^3}=1,c^p=a^{p^2}b^{sp^2},
[a,b]=c,[c,a]=1,[b,c]=b^{p^2}\rangle$, where $p>2$, $s\in F_p$ and $1+4s\not\in F_p^2$; moreover, $|G|=p^{7}$, $\Phi(G)=\lg a^p, b^p,c\rg\cong
C_{p^{2}}\times C_{p^2}\times C_{p}$, $G'=\lg b^{p^2},c\rg$ and
$Z(G)=\lg a^{p^2},
b^{p^2}\rg\cong C_p^2$.

\rr{M25} $\langle a, b; c \di a^{p^{3}}=b^{p^{2}}=c^{p^2}=1,
[a,b]=c, [c,a]=1,[c,b]=c^{tp}a^{-tp^2}\rangle$, where $p>2$ and $t^2-4t\not\in F_p^2$; moreover, $|G|=p^{7}$, $\Phi(G)=\lg a^p, b^p,c\rg\cong
C_{p^{2}}\times C_{p^2}\times C_{p}$, $G'=\lg a^{p^2},c\rg$ and
$Z(G)=\lg a^{p^2},
c^{p}\rg\cong C_p^2$.

\rr{M26} $\langle a, b; c \di a^{p^{3}}=b^{p^{3}}=1,
[a,b]=c,c^p=b^{p^2}, [c,a]=1,[b,c]=a^{\nu p^2}\rangle$, where $p>2$, $\nu=1$ or a
fixed quadratic non-residue modular $p$, and $-\nu \not\in F_p^2$;
moreover, $|G|=p^{7}$, $\Phi(G)=\lg a^p, b^p,c\rg\cong
C_{p^{2}}\times C_{p^2}\times C_{p}$, $G'=\lg a^{p^2},c\rg$ and
$Z(G)=\lg a^{p^2},
b^{p^2}\rg\cong C_p^2$.

\rr{M27} $\langle a, b; c \di a^{p^{3}}=b^{p^{2}}=c^{p^2}=1,
[a,b]=c, [c,a]=1,[b,c]=a^{\nu p^2}\rangle$, where $p>2$, $\nu=1$ or a fixed quadratic non-residue modular $p$,
and $-\nu \not\in F_p^2$; moreover, $|G|=p^{7}$, $\Phi(G)=\lg a^p, b^p,c\rg\cong
C_{p^{2}}\times C_{p^2}\times C_{p}$, $G'=\lg a^{p^2},c\rg$ and
$Z(G)=\lg a^{p^2},
c^{p}\rg\cong C_p^2$.

\rr{M28} $\langle a, b; c \di a^{p^{2}}=b^{p^{3}}=c^{p^2}=1,
[a,b]=c, [c,a]=1,[c,b]=c^{p}b^{-p^2},[b^{p^2},a]=1\rangle$, where $p>2$; moreover, $|G|=p^{7}$, $\Phi(G)=\lg a^p, b^p,c\rg\cong
C_{p^{2}}\times C_{p^2}\times C_{p}$, $G'=\lg b^{p^2},c\rg$ and
$Z(G)=\lg c^{p},
b^{p^2}\rg\cong C_p^2$.

\rr{M29} $\langle a, b; c \di a^{p^{2}}=b^{p^{3}}=c^{p^2}=1,
[a,b]=c, [c,a]=1,[b,c]=b^{p^2}\rangle$, where $p>2$; moreover, $|G|=p^{7}$, $\Phi(G)=\lg a^p, b^p,c\rg\cong
C_{p^{2}}\times C_{p^2}\times C_{p}$, $G'=\lg b^{p^2},c\rg$ and
$Z(G)=\lg c^{p},
b^{p^2}\rg\cong C_p^2$.

\rr{M30} $\langle a,b;c\di a^{p^{n+1}}=b^{p^{3}}=1,
[a,b]=c,c^p=a^{p^n},[a,c]=b^{\nu p^2},[b,c]=1\rangle$, where $n>2$, $\nu=1$ or a fixed quadratic non-residue modular $p$; moreover, $|G|=p^{n+5}$,
$\Phi(G)=\lg a^p, b^p,c\rg\cong C_{p^n}\times C_{p^2}\times
C_{p}$, $G'=\lg b^{p^2},c\rg$ and  $Z(G)=\lg a^{p^2},b^{p^2}\rg\cong C_{p^{n-1}}\times C_p$.

\rr{M31} $\langle a,b;c\di a^{p^{n+1}}=b^{p^{3}}=1,
[a,b]=c,c^p=a^{p^n}b^{s\eta p^2},[a,c]=b^{\eta p^2},[b,c]=1\rangle$, where $p>2$, $n>2$, $s=1,2,\dots,\frac{p-1}{2}$,
and $\eta$ is a fixed quadratic non-residue modular $p$; moreover, $|G|=p^{n+5}$,
$\Phi(G)=\lg a^p, b^p,c\rg\cong C_{p^n}\times C_{p^2}\times
C_{p}$, $G'=\lg b^{p^2},c\rg$ and  $Z(G)=\lg a^{p^2},b^{p^2}\rg\cong C_{p^{n-1}}\times C_p$.

\rr{M32} $\langle a, b; c \di a^{p^{n+1}}=b^{p^{2}}=c^{p^2}=1,
[a,b]=c, [c,b]=1,[c,a]=c^{p}a^{-p^n},[a^{p^n},b]=1\rangle$, where $p>2$ and $n>2$; moreover, $|G|=p^{n+5}$,
$\Phi(G)=\lg a^p, b^p,c\rg\cong C_{p^n}\times C_{p^2}\times
C_{p}$, $G'=\lg a^{p^n},c\rg$ and  $Z(G)=\lg a^{p^2},c^p\rg\cong C_{p^{n-1}}\times C_p$.

\rr{M33} $\langle a, b; c \di a^{p^{n+1}}=b^{p^{3}}=1,
[a,b]=c,c^p=b^{p^2}, [b,c]=1,[a,c]=a^{p^n}\rangle$, where $n>2$; moreover, $|G|=p^{n+5}$,
$\Phi(G)=\lg a^p, b^p,c\rg\cong C_{p^n}\times C_{p^2}\times
C_{p}$, $G'=\lg a^{p^n},c\rg$ and  $Z(G)=\lg a^{p^2},c^p\rg\cong C_{p^{n-1}}\times C_p$.

\rr{M34} $\langle a, b; c \di a^{p^{n+1}}=b^{p^{2}}=c^{p^2}=1,
[a,b]=c, [b,c]=1,[a,c]=a^{p^n}\rangle$, where $n>2$; moreover, $|G|=p^{n+5}$,
$\Phi(G)=\lg a^p, b^p,c\rg\cong C_{p^n}\times C_{p^2}\times
C_{p}$, $G'=\lg a^{p^n},c\rg$ and  $Z(G)=\lg a^{p^2},c^p\rg\cong C_{p^{n-1}}\times C_p$.

\rr{M35} $\langle a, b; c \di a^{p^{n}}=b^{p^{3}}=c^{p^2}=1,
[a,b]=c, [c,b]=1,[c,a]=c^{tp}b^{-tp^2},[b^{p^2},a]=1\rangle$, where $p>2$, $n>2$ and $t^2+4t\not\in F_p^2$; moreover, $|G|=p^{n+5}$,
$\Phi(G)=\lg a^p, b^p,c\rg\cong C_{p^{n-1}}\times C_{p^2}\times
C_{p^2}$, $G'=\lg b^{p^2},c\rg$ and  $Z(G)=\lg a^{p^2},b^{p^2},c^p\rg\cong C_{p^{n-2}}\times C_p\times C_p$.

\rr{M36} $\langle a, b; c \di a^{p^{n}}=b^{p^{3}}=c^{p^2}=1,
[a,b]=c, [c,b]=1,[a,c]=b^{\eta p^2}\rangle$, where $p>2$, $n>2$, $\eta$ is a fixed quadratic non-residue modular $p$; moreover, $|G|=p^{n+5}$,
$\Phi(G)=\lg a^p, b^p,c\rg\cong C_{p^{n-1}}\times C_{p^2}\times
C_{p^2}$, $G'=\lg b^{p^2},c\rg$ and  $Z(G)=\lg a^{p^2},b^{p^2},c^p\rg\cong C_{p^{n-2}}\times C_p\times C_p$.

\end{enumerate}
\rr{Miv} $\Phi(G')\le G_3\cong C_p^2$. In this case, $(\mu_0,\mu_1,\mu_2)=(0,0,1+p)$ and $\alpha_1(G)=p^3+p^2$ except for (M37)--(M39).

\begin{enumerate}
\rr{M37} $\langle a, b; c \di a^{8}=b^{8}=c^{2}=1,
[a,b]=c,[c,a]=a^{4}b^{4},[c,b]=a^{4},[a^4,b]=1\rangle$; where $|G|=2^{7}$,
$\Phi(G)=\lg a^2, b^2,c\rg\cong C_4^2\times C_2$, $G'=\lg
a^{4}, b^{4},c\rg\cong C_2^3$, $Z(G)=\lg a^4,
b^4\rg\cong C_2^2$, and $\alpha_1(G)=18$;

\rr{M38} $\langle a, b; c \di a^{8}=b^{8}=1,
[a,b]=c,[c,a]=c^2=b^4,[c,b]=a^4\rangle$; where $|G|=2^{7}$,
$\Phi(G)=\lg a^2, b^2,c\rg\cong C_4^2\times C_2$, $G'=\lg
a^{4}, c\rg\cong C_4\times C_2$, $Z(G)=\lg a^4,
b^4\rg\cong C_2^2$, and $\alpha_1(G)=14$;

\rr{M39} $\langle a, b; c \di a^{8}=c^4=1,
[a,b]=c,[c,a]=c^2,[c,b]=a^4=b^4\rangle$; where $|G|=2^{7}$,
$\Phi(G)=\lg a^2, b^2,c\rg\cong C_4\times C_4\times C_2$, $G'=\lg
a^{4},c\rg\cong C_4\times C_2$, $Z(G)=\lg a^4,
c^2\rg\cong C_2^2$, and $\alpha_1(G)=18$;

\rr{M40} $\langle a, b; c \di a^{p^{2}}=b^{p^{2}}=c^{p}=1,
[a,b]=c,[c,a]=a^{p},[c,b]=b^{p}\rangle$, where $p>3$; moreover,
$|G|=p^{5}$, $\Phi(G)=G'=\lg a^p, b^p,c\rg\cong
C_p^3$ and $Z(G)=\lg a^p, b^p\rg\cong C_p^2$;

\rr{M41} $\langle a, b; c \di a^{p^{3}}=b^{p^{3}}=c^{p}=1,
[a,b]=c,[c,a]=b^{\nu p^{2}},[c,b]=a^{-p^2},[a^{p^2},b]=1\rangle$, where
$p>2$, $\nu=1$ or a
fixed quadratic non-residue modular $p$ such that $-\nu\not\in F_p^2$; moreover,
$|G|=p^{7}$, $\Phi(G)=\lg a^p, b^p,c\rg\cong C_{p^2}\times C_{p^2}\times C_{p}$, $G'=\lg a^{p^2},
b^{p^2},c\rg$ and $Z(G)=\lg a^p,
b^p\rg\cong C_{p^2}\times C_{p^2}$;

\rr{M42} {\small $\langle a, b; c \di a^{p^{3}}=b^{p^{3}}=c^{p}=1,
[a,b]=c,[c,a]^{1+r}=a^{p^2}b^{p^2},[c,b]^{1+r}=a^{-r p^2}b^{ p^2},\
[a^{p^2},b]=1\rangle$, where $p>3$ and  $-r\not\in (F_p)^2$;
moreover, $|G|=p^{7}$, $\Phi(G)=\lg a^p, b^p,c\rg\cong C_{p^2}\times
C_{p^2}\times C_{p}$, $G'=\lg a^{p^2}, b^{p^2},c\rg$ and $Z(G)=\lg
a^p,
b^p\rg\cong C_{p^2}\times C_{p^2}$;}

\rr{M43} $\langle a, b; c,d,e \di a^{p}=b^{p}=c^{p}=d^p=e^p=1,
[a,b]=c,[c,a]=d,[c,b]=e,[d,a]=[d,b]=[e,a]=[e,b]=1\rangle$, where
$p>3$; moreover,
$|G|=p^{5}$, $\Phi(G)=G'=\lg c,d,e\rg\cong
C_p^3$ and $Z(G)=\lg d, e\rg\cong C_p^2$;

\rr{M44} $\langle a, b; c \di a^{p^{3}}=b^{p^{3}}=1,
[a,b]=c,[c,a]=c^p=b^{sp^{2}},[c,b]=a^{-\nu p^2}b^{st\nu p^2}\rangle$, where $p>2$, $\nu=1$ or a
fixed quadratic non-residue modular $p$, $s\in F_p^*$, $t=0,1,\dots,\frac{p-1}{2}$
such that $(st\nu)^2-4\nu(s+1)\not\in F_p^2$; moreover,
$|G|=p^{7}$, $\Phi(G)=\lg a^p, b^p,c\rg\cong C_{p^2}\times C_{p^2}\times C_{p}$, $G'=\lg a^{p^2},c\rg\cong C_{p^2}\times C_p$ and $Z(G)=\lg a^{p^2},
b^{p^2}\rg\cong C_p^2$;

\rr{M45} $\langle a, b; c \di a^{p^{3}}=b^{p^{3}}=1,
[a,b]=c,[c,a]=c^p=b^{-p^{2}},[b,c]=a^{\nu p^2}\rangle$, where $p>2$, $\nu=1$ or a
fixed quadratic non-residue modular $p$; moreover,
$|G|=p^{7}$, $\Phi(G)=\lg a^p, b^p,c\rg\cong C_{p^2}\times C_{p^2}\times C_{p}$, $G'=\lg a^{p^2},c\rg\cong C_{p^2}\times C_p$ and $Z(G)=\lg a^{p^2},
b^{p^2}\rg\cong C_p^2$;

\rr{M46} $\langle a, b; c \di a^{p^{3}}=b^{p^{2}}=c^{p^2}=1,
[a,b]=c,[c,a]=c^{p},[c,b]=a^{-\nu p^2}c^{t\nu p}\rangle$, where $p>2$, $\nu=1$ or a
fixed quadratic non-residue modular $p$, $t=0,1,\dots,\frac{p-1}{2}$
such that $(t\nu)^2-4\nu\not\in F_p^2$; moreover,
$|G|=p^{7}$, $\Phi(G)=\lg a^p, b^p,c\rg\cong C_{p^2}\times C_{p^2}\times C_{p}$, $G'=\lg a^{p^2},c\rg\cong C_{p^2}\times C_p$ and $Z(G)=\lg a^{p^2},
c^{p}\rg\cong C_p^2$;

\rr{M47} $\langle a, b; c,d \di a^{p^{2}}=b^{p^{2}}=c^p=d^p=1,
[a,b]=c,[c,a]=b^{\nu p},[c,b]=d,[d,a]=[d,b]=1\rangle$, where $p>2$, $\nu=1$ or
a fixed quadratic non-residue modular $p$; moreover, $|G|=p^{6}$, $\Phi(G)=\lg a^p, b^p,c,d\rg\cong C_p^4$,
$G'=\lg b^{p},c,d\rg\cong C_{p}^3$ and $Z(G)=\lg a^{p},
b^{p},d\rg\cong C_p^3$.
\end{enumerate}

\rr{Mv} $c(G)=4$, $G$ has a unique three-generator maximal subgroup
$M$, and $G/M'$ is a group of Type $(4)$ in Lemma {\rm 2.5}. In this case,
$p\ge 5$, $|M'|=p$, $(\mu_0,\mu_1,\mu_2)=(0,0,1+p)$ and $\alpha_1(G)=2p^2$.

\begin{enumerate}
\rr{M48}  $\langle b,a_1; a_2, a_3 \di
b^{p^2}=a_1^p=a_2^p=a_3^p=1,[a_1,b]=a_2,[a_2,b]=a_3,[a_2,a_1]=b^p,[a_3,b]=b^p,[a_3,a_1]=1
\rangle$; where $|G|=p^{5}$, $\Phi(G)=G'=\lg a_2,a_3,b^p\rg\cong
C_p^3$ and $Z(G)=\lg b^p\rg\cong C_p$.

\rr{M49}  $\langle b,a_1; a_2, a_3 \di
b^{p^2}=a_1^p=a_2^p=a_3^p=1,[a_1,b]=a_2,[a_2,b]=a_3,[a_2,a_1]=b^{\eta
p},[a_3,b]=b^p,[a_3,a_1]=1 \rangle$, where $\eta$ is a fixed
quadratic non-residue modulo $p$;  moreover, $|G|=p^{5}$,
$\Phi(G)=G'=\lg a_2,a_3,b^p\rg\cong C_p^3$ and $Z(G)=\lg b^p\rg\cong
C_p$.

\rr{M50}  $\langle b,a_1; a_2, a_3 \di
b^{p^2}=a_1^p=a_2^p=a_3^p=1,[a_1,b]=a_2,[a_2,b]=a_3,[a_2,a_1]=b^p,[a_3,b]=b^{\eta
p},[a_3,a_1]=1 \rangle$, where $p\equiv 1(\mod 4)$ and $\eta$ is a
fixed quadratic non-residue modulo $p$;  moreover, $|G|=p^{5}$,
$\Phi(G)=G'=\lg a_2,a_3,b^p\rg\cong C_p^3$ and $Z(G)=\lg b^p\rg\cong
C_p$.

\rr{M51}  $\langle b,a_1; a_2, a_3 \di
b^{p^2}=a_1^p=a_2^p=a_3^p=1,[a_1,b]=a_2,[a_2,b]=a_3,[a_2,a_1]=b^{\eta
p},[a_3,b]=b^{\eta p},[a_3,a_1]=1 \rangle$, where $p\equiv 1 (\mod
4)$ and $\eta$ is a fixed quadratic non-residue modulo $p$;
 moreover, $|G|=p^{5}$, $\Phi(G)=G'=\lg a_2,a_3,b^p\rg\cong C_p^3$ and
$Z(G)=\lg b^p\rg\cong C_p$.

\rr{M52} $\langle b,a_1; a_2, a_3 \di
b^{p}=a_1^{p^2}=a_2^p=a_3^p=1,[a_1,b]=a_2,[a_2,b]=a_3,[a_2,a_1]=a_1^{p},[a_3,b]=a_1^{\nu
p},[a_3,a_1]=1 \rangle$, where $\nu=1$, $\eta_1$ or $\eta_2$,
$\{1,\eta_1,\eta_2\}$ is a transversal for $(F_p^*)^3$ in $F_p^*$;
 moreover, $|G|=p^{5}$, $\Phi(G)=G'=\lg
a_2,a_3,a_1^p\rg\cong C_p^3$ and $Z(G)=\lg a_1^p\rg\cong C_p$.

\rr{M53} $\langle b,a_1; a_2, a_3,a_4 \di
b^{p}=a_1^{p}=a_2^p=a_3^p=a_4^p=1,[a_1,b]=a_2,[a_2,b]=a_3,[a_2,a_1]=a_4,[a_3,b]=a_4,[a_3,a_1]=[a_4,a_1]=[a_4,b]=1
\rangle$, where $|G|=p^{5}$, $\Phi(G)=G'=\lg a_2,a_3,a_4\rg\cong
C_p^3$ and $Z(G)=\lg a_4\rg\cong C_p$.

\end{enumerate}

\rr{Mvi} $c(G)=4$, $G$ has a unique three-generator maximal subgroup
$M$, $|M'|=9$ and $G/M'$ is a group of Type $(6)$ in Lemma 2.5. In this case,
$(\mu_0,\mu_1,\mu_2)=(0,0,1+p)$ and $\alpha_1(G)=2p^2+p$.

\begin{enumerate}

\rr{M54} $\lg b,a_1; a_2\di b^{27}=a_1^{9}=a_2^{9}=1,
[a_1,b]=a_2,[a_2,a_1]=b^{-9},[a_2,b]=a_1^{3}a_2^{3s}\rangle$, where
$s=0,2$;  moreover, $|G|=3^{7}$, $\Phi(G)=\lg a_2,a_1^3,b^3\rg\cong
C_3\times C_9\times C_9$, $G'=\lg a_2,a_1^3,b^9\rg\cong C_3^2\times
C_9$ and $Z(G)=\lg a_2^3,b^9\rg\cong C_3^2$.

\rr{M55} $\lg b,a_1; a_2\di b^{27}=a_1^{9}=a_2^{9}=1,
[a_1,b]=a_2,[a_2,a_1]=b^{9}a_2^{3},[a_2,b]=a_1^{3}\rangle$, where
$|G|=3^{7}$, $\Phi(G)=\lg a_2,a_1^3,b^3\rg\cong C_3\times C_9\times
C_9$, $G'=\lg a_2,a_1^3,b^9\rg\cong C_3^2\times C_9$ and $Z(G)=\lg
a_2^3,b^9\rg\cong C_3^2$.

\end{enumerate}

\rr{Mvii} $c(G)=4$, $G$ has a unique three-generator maximal
subgroup $M$, $|M'|=p^2$ where $p\ge 5$ and $G/M'$ is a group of
Type $(6)$ in Lemma {\rm 2.5}. In this case, $(\mu_0,\mu_1,\mu_2)=(0,0,1+p)$ and $\alpha_1(G)=2p^2+p$.

\begin{enumerate}

\rr{M56} $\lg b,a_1; a_2\di b^{p^{3}}=a_1^{p^2}=a_2^{p^2}=1,
[a_1,b]=a_2,[a_2,a_1]=b^{\nu_1p^{2}},[a_2,b]=a_1^{\nu_2p}a_2^{sp}\rangle$,
where $\nu_1,\nu_2=1$ or a fixed quadratic non-residue modulo $p$
such that $-\nu_1$ is not a square, and
$s=2^{-1}\nu_2,2^{-1}\nu_2+1,\dots,2^{-1}\nu_2+\frac{p-1}{2}$;
 moreover, $|G|=p^{7}$, $\Phi(G)=\lg a_2,a_1^p,b^p\rg\cong C_{p^2}\times C_{p^2}\times C_{p}$, $G'=\lg
a_2,a_1^p,b^{p^2}\rg\cong C_p^2\times C_{p^2}$ and $Z(G)=\lg
a_2^p,b^{p^2}\rg\cong C_p^2$.

\rr{M57} $\lg b,a_1; a_2\di b^{p^{3}}=a_1^{p^2}=a_2^{p^2}=1,
[a_1,b]=a_2,[a_2,a_1]=b^{\nu_1p^{2}}a_2^{rp},[a_2,b]=a_1^{\nu_2p}\rangle$,
where $\nu_1,\nu_2=1$ or a fixed quadratic non-residue modulo $p$
and $r=1,2,\dots,\frac{p-1}{2}$ such that $r^2-4\nu_1$ is not a
square.  Moreover, $|G|=p^{7}$, $\Phi(G)=\lg a_2,a_1^p,b^p\rg\cong
C_{p^2}\times C_{p^2}\times C_{p}$, $G'=\lg
a_2,a_1^p,b^{p^2}\rg\cong C_p^2\times C_{p^2}$ and $Z(G)=\lg
a_2^p,b^{p^2}\rg\cong C_p^2$.
\end{enumerate}

\rr{Mviii} $c(G)=4$, $G$ has a unique three-generator maximal
subgroup $M$, $|M'|=p$ where $p\ge 3$ and $G/M'$ is a group of Type
$(6)$ in Lemma {\rm 2.5}. In this case, $(\mu_0,\mu_1,\mu_2)=(0,0,1+p)$ and $\alpha_1(G)=2p^2$.

\begin{enumerate}
\rr{M58} $\lg b,a_1; a_2\di a_1^{p^2}=a_2^{p^2}=1,b^{p^2}=a_2^{tp},
[a_1,b]=a_2,[a_2,a_1]=a_2^{p},[a_2,b]=a_1^{\nu p}\rangle$, where
$t\in F_p$, $\nu=1$ or a fixed quadratic non-residue modulo $p$;
moreover, $|G|=p^{6}$, $\Phi(G)=\lg a_2,a_1^p,b^p\rg\cong
C_{p^2}\times C_p^2$, $G'=\lg a_2,a_1^p\rg\cong C_{p^2}\times
 C_p$, $Z(G)=\lg a_2^p,b^{p^2}\rg\cong C_p$ except for
$p =3$ and $\nu=-1$, In case of $p =3$ and $\nu=-1$, we have
$Z(G)=\lg a_2^p,b^{p}\rg\cong C_p^2$ if $t\equiv0(\mod p)$ or
$Z(G)=\lg a_2^p,b^{p}\rg\cong C_{p^2}$ if $t\not\equiv0(\mod p)$.

\rr{M59} $\lg b,a_1; a_2\di a_1^{p^2}=a_2^{p^2}=b^{p^m}=1,
[a_1,b]=a_2,[a_2,a_1]=1,[a_2,b]=a_1^{\nu p}a_2^{sp}\rangle$, where
$m\ge 2$, $\nu=1$ for $p=3$, $\nu=1$ or a fixed quadratic
non-residue modulo $p$ for $p\ge 5$,
$s=\nu,\nu+1,\dots,\nu+\frac{p-1}{2}$;  moreover, $|G|=p^{m+4}$,
$\Phi(G)=\lg a_2,a_1^p,b^p\rg\cong C_p\times C_{p^{m-1}}\times
C_{p^{2}}$, $G'=\lg a_2,a_1^p\rg\cong C_p\times C_{p^{2}}$,
$Z(G)=\lg a_2^p,b^{p^2}\rg\cong C_p\times C_{p^{m-2}}$ for $m>2$,
$Z(G)=\lg a_2^p,b^{p^2}\rg\cong C_p$ for $m=2$.

\rr{M60} $\lg b,a_1; a_2\di a_1^{p^2}=a_2^{p^2}=1,b^{p^m}=a_2^{p},
[a_1,b]=a_2,[a_2,a_1]=1,[a_2,b]=a_1^{\nu p}\rangle$, where $m\ge 2$,
$\nu=1$ for $p=3$, $\nu=1$ or a fixed quadratic non-residue modulo
$p$ for $p\ge 5$.  Moreover, $|G|=p^{m+4}$, $\Phi(G)=\lg
a_2,a_1^p,b^p\rg\cong C_{p^{m}}\times C_p^2$, $G'=\lg
a_2,a_1^p\rg\cong C_p\times C_{p^{2}}$,  $Z(G)=\lg
a_2^p,b^{p^2}\rg\cong C_{p^{m-2}}$ for $m>2$,  $Z(G)=\lg
a_2^p,b^{p^2}\rg\cong C_p$ for $m=2$.

\end{enumerate}

\rr{Mix} $c(G)=4$, $G$ has a unique three-generator maximal subgroup
$M$, $|M'|=3$ and $G/M'\in\mathcal{A}_3$. In this case, $(\mu_0,\mu_1,\mu_2)=(0,0,1+p)$ and $\alpha_1(G)=2p^2$.

\begin{enumerate}

\rr{M61} $\lg b,a_1; a_2\di a_1^{9}=a_2^{9}=b^{27}=1,
[a_1,b]=a_2,[a_2,a_1]=b^{9s},[a_2,b]=a_1^{-3}a_2^{3t}\rangle$, where
$s,t=1,2$.  Moreover, $|G|=3^{7}$, $\Phi(G)=\lg
a_2,a_1^3,b^3\rg\cong C_3\times C_9\times C_9$, $G'=\lg
a_2,a_1^3,b^9\rg\cong C_9\times C_3\times C_3$ and  $Z(G)=\lg
a_2^3,b^{3}\rg\cong C_3\times C_9$.

\end{enumerate}

\rr{Mx} $c(G)=4$, $G$ has a unique three-generator maximal subgroup
$M$, $|M'|=p$ where $p\ge 5$ and $G/M'\in\mathcal{A}_3$. In this case, $(\mu_0,\mu_1,\mu_2)=(0,0,1+p)$ and $\alpha_1(G)=2p^2$.

\begin{enumerate}
\rr{M62} $\lg b,a_1; a_2\di b^{p^2}=a_1^{p^2}=a_2^{p}=a_3^{p}=1,
[a_1,b]=a_2,[a_2,a_1]=b^{\nu
p},[a_2,b]=a_3,[a_3,b]=a_1^{tp},[a_3,a_1]=1\rangle$, where $p\ge 5$,
$\nu=1$ or a fixed quadratic non-residue modulo $p$,
$t=t_1,t_2,\dots,t_{(3,p-1)}$, where $t_1,t_2,\dots,t_{(3,p-1)}$ are
the coset representatives of the subgroup $(F_p^*)^3$ in $F_p^*$.
 Moreover,$|G|=p^{6}$, $\Phi(G)=G'=\lg
a_2,a_3,a_1^p,b^p\rg\cong C_p^4$ and $Z(G)=\lg a_1^p,b^{p}\rg\cong
C_p^2$.
\end{enumerate}
\end{enumerate}
\end{thm}

\demo If $\Phi(G')G_3\le Z(G)$ and $\Phi(G')G_3\le C_p^2$, then $G$
is one of the groups listed in \cite[Theorem 3.5, 4.6, 5.5, 5.8 and
6.5]{AHZ}. By \cite[Theorem 3.6, 4.7, 5.1-5.2, 5.5-5.6 and
6.1]{AHZ}, we get the groups (M1)-(M47). For convenience, in Table
\ref{table10} we give the correspondence from those groups
(M1)-(M47) to \cite[Theorem 3.5, 4.6, 5.5, 5.8 and 6.5]{AHZ}. In the
following, we may assume that $\Phi(G')G_3\not\le Z(G)$ or
$\Phi(G')G_3\not\le C_p^2$.

\begin{table}[h]
  \centering
{
\tiny
\begin{tabular}{c|c||c|c}
\hline
Groups &Groups in \cite[Theorem 3.5, 4.6, 5.5, 5.8 and 6.5]{AHZ} & Groups&Groups in \cite[Theorem 3.5, 4.6, 5.5, 5.8 and 6.5]{AHZ} \\
\hline
(M1)&(C1) & (M25)&(L2) where $n=2$ and $t^2-4t\not\in F_p^2$ \\
(M2)&(C2) & (M26)&(L4) where $n=2$ and $-\nu\not\in F_p^2$\\
(M3)&(D1) where $p>2=n$ & (M27)&(L5) where $n=2$ and $-\nu\not\in F_p^2$\\
(M4)&(D2) where $p>2=n$ & (M28)&(L6) where $n=2$\\
(M5)&(D3) & (M29)& (L8) where $n=2$\\
(M6)&(I4) where $p>2=n$ & (M30)&(M1) where $m=2$ and $s=0$\\
(M7)&(D5) & (M31)&(M1) where $m=2$, $s=0$ and $\nu=\eta$\\
(M8)&(E2) where $p>2=m$ & (M32)&(M2) where $m=2$\\
(M9)&(E5) where $p>2=m$ & (M33)&(M4) where $m=2$\\
(M10)&(E8) where $p>2=m$ & (M34)&(M5) where $m=2$\\
(M11)&(E3) & (M35)&(M6) where $m=2$ and $t^2+4t\not\in F_p^2$\\
(M12)&(E6) & (M36)&(M8) where $m=2$ and $\nu=\eta$\\
(M13)&(E9) & (M37)&(N3)\\
(M14)&(G1) where $p=m=2$ & (M38)&(N8)\\
(M15)&(G2) where $p=m=2$ & (M39)&(N9)\\
(M16)&(G3) where $p=m=2$ & (M40)&(P1) where $n=1$ \\
(M17)&(J2) where $p=m=2$ & (M41)&(P3) where $n=2$ and $-\nu\not\in F_p^2$ \\
(M18)&(J4) where $p=m=2$ & (M42)&(p4) where $n=2$ and $-r\not\in F_p^2$\\
(M19)&(J6) where $p=m=2$ & (M43)&(P10) where $n=1$\\
(M20)&(K2) &(M44) &(Q1) where $n=2$ and $(st\nu)^2-4\nu(s+1)\not\in F_p^2$\\
(M21)&(K4) &(M45) &(Q1) where $n=2$, $s=-1$ and $t=0$\\
(M22)&(K8) &(M46) &(Q2) where $n=2$\\
(M23)&(K10) &(M47) &(S4) where $n=2$\\
(M24)&(L1) where $n=2$ and $1+4s\not\in F_p^2$ & &\\

\hline
\end{tabular}
\caption{{\small The correspondence from Theorem \ref{d=4-3} to
\cite[Theorem 3.5, 4.6, 5.5, 5.8 and 6.5]{AHZ} }} \label{table10}}
\end{table}

\medskip
If $G$ has two distinct three-generator maximal subgroup $M_1$ and
$M_2$, then, by Lemma \ref{A_2-property} (2), $M_1'\le C_p^3$ and
$M_2'\le C_p^3$. By hypothesis, $M_1'\le Z(G)$ and $M_2'\le Z(G)$.
Let $N=M_1'M_2'$. Then $\exp(N)=p$ and $N\le Z(G)$. Let
$\bar{G}=G/N$. Then $\bar{G}$ has two distinct abelian subgroups
$\bar{M}_1$ and $\bar{M}_2$  of index $p$. By Lemma \ref{minimal
non-abelian equivalent conditions} we have $\bar{G}\in
\mathcal{A}_1$. It follows that $\Phi(G')G_3=N$. By above
assumption, $|N|\ge p^3$. Hence $G$ is not metacyclic. By Lemma
\ref{metacyclic}, $\bar{G}$ is not metacyclic. Let $$\bar{G}=\lg
\bar{a},\bar{b};\bar{c}\di
\bar{a}^{p^n}=\bar{b}^{p^m}=\bar{c}^p=1,[\bar{a},\bar{b}]=\bar{c},[\bar{c},\bar{a}]=[\bar{c},\bar{b}]=1\rg,$$
where $n\ge m$. Then $N=\lg c^p,[c,a],[c,b]\rg$, $|N|=p^3$ and
$|G|=p^{n+m+4}$. Since $c^p\ne 1$, by calculation we have
$[a,b^p]\ne 1$. It follows that $m\ge 2$. Since $\lg c,b\rg\in
\mathcal{A}_1$, we have $|\lg c,b\rg|=p^{n+m+2}$. It follows that
$n\le 2$ and hence $n=m=2$. Moreover, $$N=\lg
c^p,a^{p^2},b^{p^2}\rg,\ {\rm where}\ [c,a]\not\in \lg
c^p,a^{p^2}\rg\ {\rm and}\ [c,b]\not\in \lg c^p,b^{p^2}\rg.$$ By
suitable replacement, we may assume that $[c,a]=b^{rp^2}c^{sp}$ and
$[c,b]=a^{u p^2}b^{vp^2}c^{wp}$, where $r,s,u,v,w\in F_p$, $r\ne 0$
and $u\ne 0$. If $p=2$, then we can prove that $G\in\mathcal{A}_4$.
The details are omitted. Hence $p\ge 3$. If $s^2-4r$ is a square,
then the equation $x^2-sx+r=0$ has a solution $x_1$. By calculation
we have $\lg a,cb^{x_1p}\rg$ is not abelian and of order $p^5$,
which contradicts that $G\in \mathcal{A}_3$. Hence $s^2-4r$ is not a
square. By Lemma \ref{equivalent equation}, the equation
$x^2+sxy+ry^2+wx+vy-u=0$ has a solution $(x_0,y_0)$. By calculation,
$\lg ca^{x_0p},ba^{y_0}\rg$ is not abelian and of order $p^5$, which
contradicts that $G\in \mathcal{A}_3$ again.

\medskip
By above argument, $G$ has a unique three-generator maximal subgroup
$M$. Let $\bar{G}=G/M'$. Then $\bar {G}$ has a three-generator
abelian subgroup $M/M'$  of index $p$, and every non-abelian
subgroup of $\bar{G}$ is generated by two elements. By Lemma
\ref{alj2} (2), $p\ge 3$. Since $\bar{G}\in\mathcal{A}_2$ or
$\mathcal{A}_3$, $\bar{G}$ is either one of the groups (4)--(7) in
Lemma \ref{A_2} or, by Theorem \ref{d=3-1}, one of the groups
(F4)--(F8).

We claim that $G'$ is abelian. Otherwise, $|G:G'|=p^2$ and $G'$ is
minimal non-abelian. By Lemma \ref{G' is abelian}, $G'$ is abelian,
a contradiction.

\medskip

Case 1: $\bar{G}$ is the group (4) in Lemma \ref{A_2}. That
is, $\bar{G}=\lg \bar{a}_1,\bar{b};\bar{a}_2,\bar{a}_3\di
\bar{a}_1^p=\bar{a}_2^p=\bar{a}_3^p=\bar{b}^{p^m}=1,[\bar{a}_1,\bar{b}]=\bar{a}_{2},[\bar{a}_2,\bar{b}]=\bar{a}_3,[\bar{a}_3,\bar{b}]=1,[\bar{a}_i,\bar{a}_j]=1\rg,$
 where $1\leq i,j\leq 3$.

\medskip

Since $d(\bar{M})=3$, we have $m=1$ and $M=\lg a_1,a_2,a_3\rg$.
Hence $p\ge 5$. Since $G'$ is abelian, $[a_2,a_3]=1$. By calculation
we have $$[a_3,a_1]=[a_2,b,a_1]=[a_2,a_1,b]=1.$$ It follows that
$$[a_1,a_2]\ne 1,\ |M'|=p\ {\rm and}\ |G|=p^5.$$ By calculation we
have $$a_3^p=[a_2,b]^p=[a_2^p,b]=1\ {\rm and}\ a_2^p=[a_1^p,b]=1.$$
Since $\lg a_2,b\rg\in\mathcal{A}_2$, $[a_3,b]\ne 1$ and hence $G$
is of maximal class.

If $\exp(G)=p^2$ and $\exp(M)=p$, then $a_1^p=1$ and $b^{p^2}=1$.
Hence $G=\lg b,a_1,a_2,a_3\rg$ has following relations:
$$b^{p^2}=a_1^p=a_2^p=a_3^p=1,[a_1,b]=a_2,[a_2,b]=a_3,[a_3,b]=b^{ip},[a_2,a_1]=b^{jp}$$
where $i,j\in F_p^*$. By Lemma \ref{cong 0} we get groups
(M48)--(M51).

If $\exp(G)=p^2$ and $\exp(M)=p^2$, then $a_1^{p^2}=1$. By suitable
replacement, we may assume that $b^{p}=1$. Hence $G=\lg
b,a_1,a_2,a_3\rg$ has following relations:
$$b^{p}=a_1^{p^2}=a_2^p=a_3^p=1,[a_1,b]=a_2,[a_2,b]=a_3,[a_3,b]=a_1^{ip},[a_2,a_1]=a_1^{jp}$$
where $i,j\in F_p^*$. By Lemma \ref{cong 0a} we get the group (M52).

If $\exp(G)=p$, then we get the group (M53).

\medskip

Case 2: $\bar{G}$ is the group (5) in Lemma \ref{A_2}. That
is, $\bar{G}=\lg \bar{a}_1,\bar{b};\bar{a}_2\di
\bar{a}_1^p=\bar{a}_2^p=\bar{b}^{p^{m+1}}=1,[\bar{a}_1,\bar{b}]=\bar{a}_{2},[\bar{a}_2,\bar{b}]=\bar{b}^{p^m},[\bar{a}_1,\bar{a}_2]=1\rg$.

\medskip

Since $d(\bar{M})=3$, $M=\lg a_1,a_2,b^p\rg$. Let $N=\lg
a_2,b,M'\rg$. Then $N$ is maximal in $G$ and $d(N)=2$. It follows
that $N\in \mathcal{A}_2$. Since $G'$ is abelian, we have
$[b^{p^m},a_2]=1$ and hence $c(N)=2$. By calculation we get
$b^{p^{m+1}}=[a_2,b]^p=[a_2^p,b]=1$ and hence $|N'|=p$. By Lemma
\ref{minimal non-abelian equivalent conditions} we have
$N\in\mathcal{A}_1$, a contradiction.

\medskip

Case 3: $\bar{G}$ is the group (6) in Lemma \ref{A_2}. That
is, $\bar{G}=\lg \bar{a}_1,\bar{b};\bar{a}_2\di
\bar{a}_1^{p^2}=\bar{a}_2^p=\bar{b}^{p^m}=1,[\bar{a}_1,\bar{b}]=\bar{a}_{2},[\bar{a}_2,\bar{b}]=\bar{a}_1^{\nu
p},[\bar{a}_1,\bar{a}_2]=1\rg,$ where $\nu=1$ or a fixed quadratic
non-residue modulo $p$.

\medskip

Since $d(\bar{M})=3$, we have $m\ge 2$ and $M=\lg a_1,a_2,b^p\rg$.
Let $N=\lg a_2,b,M'\rg$. Then $N$ is maximal in $G$ and $d(N)=2$. It
follows that $N\in \mathcal{A}_2$. By calculation, $a_1^{\nu
p^2}=[a_2,b]^p=[a_2^p,b]=1$. Hence $a_1^{p^2}=1$ and $\exp(N')=p$.
If $c(N)=2$, then $G\in\mathcal{A}_1$, a contradiction. Hence
$c(N)=3$. By calculation, $[a_2,b^p]=1$. Since $[a_1^p,a_2]=1$, we
have $[a_1^p,b]=a_2^p\ne 1$.

If $|M'|=p^2$, then $M'=\lg [a_1,a_2]\rg\times \lg [a_1,b^p]\rg$.
Since $G\in\mathcal{A}_3$, we have $|\lg a_1,a_2\rg|=|\lg
a_1,b^p\rg|=p^{m+3}$. It follows that $m=2$, $|G|=p^{7}$ and $M'=\lg
a_2^p,b^{p^2}\rg$. By suitable replacement, we may assume that
$$G=\lg b,a_1;a_2\di b^{p^3}=a_1^{p^2}=a_2^{p^2}=1,
[a_1,b]=a_2,[a_2,a_1]=b^{\nu_1p^2}a_2^{rp},[a_2,b]=a_1^{\nu_2p}a_2^{sp}\rg,$$
where $\nu_1,\nu_2=1$ or a fixed quadratic non-residue modulo $p$.
By Lemma \ref{cong 1} we get groups (M54)--(M57).

In the following, we may assume that $|M'|=p$. Hence $M'=\lg
a_2^2\rg$ and $|G|=p^{m+4}$.

If $[a_2,a_1]\ne 1$, then $|\lg a_1,a_2\rg|=p^{m+2}$ since
$G\in\mathcal{A}_3$. It follows that $m=2$ and $|G|=p^{6}$. By
suitable replacement, we may assume that $$G=\lg b,a_1;a_2\di
a_1^{p^2}=a_2^{p^2}=1, b^{p^2}=a_2^{tp},
[a_1,b]=a_2,[a_2,a_1]=a_2^{p},[a_2,b]=a_1^{\nu p}\rg,$$ where $t\in
F_p$, $\nu=1$ or a fixed quadratic non-residue modulo $p$. By Lemma
\ref{cong 2} we get the group (M58).

If $[a_2,a_1]=1$, then, by suitable replacement, we may assume that
$$G=\lg b,a_1;a_2\di a_1^{p^2}=a_2^{p^2}=1, b^{p^m}=a_2^{tp},
[a_1,b]=a_2,[a_2,a_1]=1,[a_2,b]=a_1^{\nu p}a_2^{sp}\rg,$$ where
$s,t\in F_p.$ By Lemma \ref{cong 3} we get groups (M59)--(M60).

\medskip

Case 4: $\bar{G}$ is the group (7) in Lemma \ref{A_2}. That
is, $\bar{G}=\lg \bar{a}_1,\bar{b};\bar{a}_2\di
\bar{a}_1^9=\bar{a}_2^3=1,\bar{b}^3=\bar{a}_1^3,[\bar{a}_1,\bar{b}]=\bar{a}_2,[\bar{a}_2,\bar{b}]=\bar{a}_1^{-3},[\bar{a}_2,\bar{a}_1]=1\rg$.

\medskip

In this case, we have $d(\bar{M})=d(\lg \bar{a}_2,\bar{a}_1\rg)=2$,
a contradiction.

\medskip

Case 5: $\bar{G}$ is the group (F4) in Main Theorem. That is,
$\bar{G}=\lg \bar{a}_1,\bar{b};\bar{a}_2\di
\bar{a}_1^{9}=\bar{a}_2^{9}=\bar{b}^{3^m}=1,[\bar{a}_1,\bar{b}]=\bar{a}_2,[\bar{a}_2,\bar{b}]=\bar{a}_1^{-3}\bar{a}_2^{3t},[\bar{a}_1,\bar{a}_2]=1\rg$,
where $t=1,2$.

\medskip

Since $d(\bar{M})=3$, we have $m\ge 2$ and $M=\lg a_1,a_2,b^3\rg$.
By Lemma \ref{A_2-property} (4), $|M'|=3$. By Lemma
\ref{A_2-property} (3), $m=2$ and hence $|G|=3^7$. Let $N=\lg
a_2,b,M'\rg$. Then $N$ is maximal in $G$ and $d(N)=2$. It follows
that $N\in \mathcal{A}_2$. By calculation,
$a_2^{9}=[a_1^3,b]^3=[a_1^9,b]=1$. Since $|\lg a_2^3,b\rg|\le 3^4$,
we have $[a_2^3,b]=1$. By calculation,
$a_1^{9}=[b,a_2]^3=[b,a_2^3]=1$. Hence $\exp(N')=3$. Let $L=\lg
a_1^3,b\rg$. Then $L\in\mathcal{A}_1$. It follows that $|L|=3^5$.
Hence $b^9\ne 1$. That is, $M'=\lg b^9\rg$. By calculation,
$[a_2,b^3]=1$ and $[a_1,b^3]=1$. Hence $[a_2,a_1]\ne 1$. Then we may
assume $$G=\lg b,a_1;a_2\di a_1^{9}=a_2^{9}=b^{27}=1,
[a_1,b]=a_2,[a_2,a_1]=b^{9s},[a_2,b]=a_1^{-3}a_2^{3t}b^{9r}\rg,$$
where $s,t=1,2$. By replacing $b$ with $a_1^{-rs}$ we get the group
(M61).

\medskip

Case 6: $\bar{G}$ is the group (F5) in Main Theorem. That is,
$\bar{G}=\lg \bar{a}_1,\bar{b};\bar{a}_2\di
\bar{a}_1^{9}=\bar{a}_2^{9}=1,
\bar{b}^{3^m}=\bar{a}_2^{-3},[\bar{a}_1,\bar{b}]=\bar{a}_2,[\bar{a}_2,\bar{b}]=\bar{a}_1^{-3}\bar{a}_2^{-3},[\bar{a}_1,\bar{a}_2]=1\rg$.

\medskip

Since $d(\bar{M})=3$, we have $m\ge 2$ and $M=\lg a_1,a_2,b^3\rg$.
By Lemma \ref{A_2-property} (4), $|M'|=3$. By Lemma
\ref{A_2-property} (3), $m=2$ and hence $|G|=3^7$. Let $N=\lg
a_2,b,M'\rg$. Then $N$ is maximal in $G$ and $d(N)=2$. It follows
that $N\in \mathcal{A}_2$. By calculation we have
$$a_2^{9}=[a_1^3,b]^3=[a_1^9,b]=1\ {\rm and}\
a_1^{9}=[b,a_2]^3=[b,a_2^3]=1.$$ Hence $\exp(N')=3$. Let $L=\lg
a_1^3,b\rg$. Then $L\in\mathcal{A}_1$. It follows that $|L|=3^5$.
Hence $b^9\ne a_2^{-3}$. That is, $M'=\lg a_2^3b^9\rg$. By
calculation, $[a_2,b^3]=1$ and $[a_1,b^3]=1$. Hence $[a_2,a_1]\ne
1$. Then we may assume $[a_2,a_1]=a_2^{3s}b^{9s}$, where $s=1,2$. By
calculation, $\lg a_2b^3,a_1\rg\in \mathcal{A}_1$ and is of order
$3^4$, which contradicts that $G\in\mathcal{A}_3$.

\medskip

Case 7: $\bar{G}$ is the group (F6) in Main Theorem. That is,
$\bar{G}=\lg \bar{a}_1,\bar{b};\bar{a}_2,\bar{a}_3,\bar{a}_4\di
\bar{a}_i^p=\bar{b}^{p^m}=1,[\bar{a}_j,\bar{b}]=\bar{a}_{j+1},[\bar{a}_4,\bar{b}]=1,[\bar{a}_i,\bar{a}_j]=1\rg,$
 where $p\ge 5$, $1\leq i\leq 4$, $1\leq j\leq 3$.

\medskip

In this case, we have $d(\bar{M})=d(\lg
\bar{a}_1,\bar{a}_2,\bar{a}_3,\bar{a}_4,\bar{b}^p\rg)=5$, a
contradiction.

\medskip

Case 8: $\bar{G}$ is the group (F7) in Main Theorem. That is,
$\bar{G}=\lg \bar{a}_1,\bar{b};\bar{a}_2,\bar{a}_3\di
\bar{a}_i^p=\bar{b}^{p^{m+1}}=1,[\bar{a}_j,\bar{b}]=\bar{a}_{j+1},[\bar{a}_3,\bar{b}]=b^{p^m},[\bar{a}_i,\bar{a}_j]=1\rg,$
 where $p\ge 5$, $1\leq i\leq 3$, $1\leq j\leq 2$.

\medskip

In this case, we  have $d(\bar{M})=d(\lg
\bar{a}_1,\bar{a}_2,\bar{a}_3,\bar{b}^p\rg)=4$, a contradiction.

\medskip

Case 9: $\bar{G}$ is the group (F8) in Main Theorem. That is,
$\bar{G}=\lg a_1,b;\bar{a}_2,\bar{a}_3\di
\bar{a}_1^{p^2}=\bar{a}_i^p=\bar{b}^{p^m}=1,[\bar{a}_j,b]=\bar{a}_{j+1},[\bar{a}_3,\bar{b}]=\bar{a}_1^{tp},[\bar{a}_i,\bar{a}_j]=1\rg,$
where $2\leq i\leq 3,\ 1\leq j\leq 2$, and
$t=t_1,t_2,\dots,t_{(3,p-1)}$, where $p\ge 5$,
$t_1,t_2,\dots,t_{(3,p-1)}$ are the coset representatives of the
subgroup $(F_p^*)^3$ in $F_p^*$.

\medskip

Since $d(\bar{M})=3$, we have $m=1$ and $M=\lg a_1,a_2,a_3\rg$.
Since $G'$ is abelian, $[a_2,a_3]=1$. By calculation we have
$$[a_3,a_1]=[a_2,b,a_1]=[a_2,a_1,b]=1.$$ It follows that $$[a_1,a_2]\ne
1,\ |M'|=p\ {\rm and}\ |G|=p^6.$$ By calculation we have
$$a_3^p=[a_2,b]^p=[a_2^p,b]=1,\ a_2^p=[a_1^p,b]=1\ {\rm and}\
a_1^{tp^2}=[a_3^p,b]=1.$$ Let $[a_2,a_1]=d$. Then we may assume

$G=\lg b,a_1;a_2,a_3,d\di a_1^{p^2}=a_2^{p}=a_3^p=1, b^{p}=d^{s},
[a_1,b]=a_2,[a_2,a_1]=d,[a_2,b]=\mbox{\hskip1.7in}a_3,[a_3,b]=a_1^{tp}d^r\rg.$

Since $[a_3,ba_1^x]=a_1^{tp}d^r\ne 1$, we have $|\lg
a_3,ba_1^x\rg|=p^4$ for any $x$. It follows that the equation
$st-rx\equiv 0\ (\mod p)$ about $x$ has no solution. Hence we have
$r=0$ and $s\ne 0$. By suitable replacement we get the group
(M62).

\medskip
We calculate the $(\mu_0,\mu_1,\mu_2)$ and $\alpha_1(G)$ of those
groups in Theorem \ref{d=4-3} as follows.

Since $d(G)=2$, $(\mu_0,\mu_1,\mu_2)=(0,0,1+p)$. In the following,
we calculate $\alpha_1(G)$.

\medskip
{\bf Case 1.} $G$ is one of the groups (M1)--(M19).

Since $|G'|=p^2$, $|M'|=p$ for any $M\in\Gamma_1$. By Lemma \ref{A_2-property} (7), $\alpha_1(M)=p^2$.
By Hall's enumeration principle,
$$\alpha_1(G)=\sum_{H\in
\Gamma_1} \alpha_1(H)=(1+p)p^2=p^2+p^3.$$

\medskip
{\bf Case 2.} $G$ is one of the groups (M20)--(M36).

In this case, $G=\lg a,b\rg$ such that $[a,b]=c$, $c^{p^2}=1$, $[c,a]=1$ and $[c,b]\not\in \lg c^p\rg$.
All maximal subgroups of $G$ are:

$M=\lg c,a,\Phi(G)\rg$;

$M_{i}=\lg c,ba^i, \Phi(G)\rg$, where $0\leq i\leq p-1$.

It is easy to see that $|M'|=p$. By Lemma \ref{A_2-property} (7), $\alpha_1(M)=p^2$.
By calculation, $d(M_{i})=3$ and $|M_i'|=p^2$. Hence $\alpha_1(M_i)=p^2+p$.
By Hall's enumeration principle,
$$\alpha_1(G)=\sum_{H\in
\Gamma_1} \alpha_1(H)=p^2+p\times (p^2+p)=p^3+2p^2.$$

\medskip
{\bf Case 3.} $G$ is one of the groups (M37).

All maximal subgroups of $G$ are:

$M=\lg c,a,\Phi(G)\rg$;

$M_{i}=\lg c,ba^i, \Phi(G)\rg$, where $i=0,1$.

By calculation, $d(M)=3$ and $|M'|=p^2$. Hence $\alpha_1(M)=p^2+p$.
Similarly, $\alpha_1(M_i)=p^2+p$.
By Hall's enumeration principle,
$$\alpha_1(G)=\sum_{H\in
\Gamma_1} \alpha_1(H)=(1+p)(p^2+p)=p^3+2p^2+p.$$

\medskip
{\bf Case 4.} $G$ is one of the groups {(M38)--(M39)}.

In this case, $G=\lg a,b\rg$ such that $[a,b]=c, c^4=1$, $[c,a]=c^2$ and $[c,b]\not\in\lg c^2\rg$.
All maximal subgroups of $G$ are:

$M=\lg c,a,\Phi(G)\rg$;

$M_{i}=\lg c,ba^i, \Phi(G)\rg$, where $i=0,1$.

By calculation, $M=\lg c,a,b^2\rg$ such that $d(M)=3$ and $|M'|=p^2=4$,
 $M_i=\lg c,ba^i,a^2\rg$ such that $d(M_i)=3$ and $|M_i'|=p=2$. By Lemma \ref{A_2-property}, $\alpha_1(M)=p^2+p$ and $\alpha_1(M_i)=p^2$.
By Hall's enumeration principle,
$$\alpha_1(G)=\sum_{H\in
\Gamma_1} \alpha_1(H)=(p^2+p)+p\times p^2=p^3+p^2+p.$$

\medskip
{\bf Case 5.} $G$ is one of the groups {(M40)--(M43)}.

Let $H\in \Gamma_1$. Then $|H'|=p$. By Lemma \ref{A_2-property},  $\alpha_1(H)=p^2$. By Hall's enumeration principle,
$$\alpha_1(G)=\sum_{H\in
\Gamma_1} \alpha_1(H)=(1+p)\times p^2=p^3+p^2.$$

\medskip
{\bf Case 6.} $G$ is one of the groups {(M44)--(M46)}.

 In this case, $G=\lg a,b\rg$ such that $[a,b]=c, c^{p^2}=1$, $[c,a]=c^p$ and $[c,b]\not\in\lg c^p\rg$ where $p>2$.
All maximal subgroups of $G$ are:

$M=\lg c,a,\Phi(G)\rg$;

$M_{i}=\lg c,ba^i, \Phi(G)\rg$, where $i=0,1$.

By calculation, $M=\lg c,a,b^2\rg$ such that $d(M)=3$ and $|M'|=p$,
 $M_i=\lg c,ba^i,a^2\rg$ such that $d(M_i)=3$ and $|M_i'|=p^2$. By Lemma \ref{A_2-property},
 $\alpha_1(M)=p^2$ and $\alpha_1(M_i)=p^2+p$.
By Hall's enumeration principle,
$$\alpha_1(G)=\sum_{H\in
\Gamma_1} \alpha_1(H)=p^2+p\times (p^2+p)=p^3+2p^2.$$

\medskip
{\bf Case 7.} $G$ is the group {(M47)}.

Let $H\in \Gamma_1$. Then $|H'|=p$. By Lemma \ref{A_2-property},  $\alpha_1(H)=p^2$.
By Hall's enumeration principle,
$$\alpha_1(G)=\sum_{H\in
\Gamma_1} \alpha_1(H)=(1+p)\times p^2=p^3+p^2.$$

\medskip
{\bf Case 8.} $G$ is one of the groups {(M48)--(M53)}.

In this case, $\Phi(G)$ is abelian and $|M'|=p$. Let $H\in \Gamma_1\setminus \{M\}$.
Then $d(H)=2$. Hence $\alpha_1(H)=p$ and $\alpha_1(M)=p^2$.
By Hall's enumeration principle,
$$\alpha_1(G)=\sum_{H\in
\Gamma_1} \alpha_1(H)=p^2+p\times p=2p^2.$$

\medskip
{\bf Case 9.} $G$ is one of the groups {(M54)--(M57)}.

In this case, $\Phi(G)$ is abelian and $|M'|=p^2$. Let $H\in \Gamma_1\setminus \{M\}$.
Then $d(H)=2$. Hence $\alpha_1(H)=p$ and $\alpha_1(M)=p^2+p$.
By Hall's enumeration principle,
$$\alpha_1(G)=\sum_{H\in
\Gamma_1} \alpha_1(H)=(p^2+p)+p\times p=2p^2+p.$$

\medskip
{\bf Case 10.} $G$ is one of the groups {(M58)--(M62)}.

In this case, $\Phi(G)$ is abelian and $|M'|=p$. Let $H\in \Gamma_1\setminus \{M\}$.
Then $d(H)=2$. Hence $\alpha_1(H)=p$ and $\alpha_1(M)=p^2$.
By Hall's enumeration principle,
$$\alpha_1(G)=\sum_{H\in
\Gamma_1} \alpha_1(H)=p^2+p\times p=2p^2.$$ \qed

\begin{lem}\label{cong 4}
Suppose that $p\ge 3$ and $G(\nu,k)=\lg a,b,x\di
a^{p^{3}}=b^{p^3}=x^{p}=1, [a,b]=a^p,[x,b]=a^{p^{2}},[x,a]=b^{\nu
p^2}a^{kp^2}\rangle$, where $k\in F_p$ and $\nu=1$ or a fixed
quadratic non-residue modula $p$.
\begin{enumerate}
\rr1 If $G=G(\nu,k)$, then $|G|=p^{7}$, $G'=\lg a^p,b^{p^2}\rg$,
$G_3=\lg a^{p^2}\rg$;

\rr2 $G(\nu',k')\cong G(\nu,k)$ if and only if $\nu'=\nu$ and
$k'=\pm k$.

\rr3 $G(\nu,k)\in A_3$ if and only if $k^2+4\nu$ is not a square.

\end{enumerate}
\end{lem}

\demo (1) Let $$K=\lg a,b \di
a^{p^3}=b^{p^3}=1,[a,b]=a^p\rg.$$ We define an
automorphism $\gamma$ of $K$ as follows:
$$a^{\gamma}=aa^{-kp^2}b^{-\nu p^2}, \ b^\gamma=ba^{-p^2}.$$ Then
$$a^{\gamma^p}=1, \ b^{\gamma^{p}}=1,\ o(\gamma)=p.$$ It
follows from the cyclic extension theory that $G=\lg K,x\rg$ is a
cyclic extension of $K$ by $C_{p}$. Hence $|G|=p^{7}$. It is easy
to check that $G'=\lg a^p,b^{p^2}\rg$, $G_3=\lg
a^{p^2}\rg$.

\medskip

(2) For convenience, let $G=\lg a,b,x\rg\cong G(\nu,k)$,
$\bar{G}=\lg \bar{a},\bar{b},\bar{x}\rg\cong G(\nu',k')$ and
$\theta$ be an isomorphism from $\bar{G}$ to $G$. Since
$M=\lg
x,b,a^p\rg$ is the unique three-generator maximal subgroups such that $M'\cong C_p^2$, $M\cha G$.
Similarly, $\bar{M}=\lg \bar{x},\bar{b},\bar{a}^p\rg\cha \bar{G}$.
Since $M,\Phi(G),\Omega_1(G)=\lg x,a^{p^2},b^{p^2}\rg$ and
$\bar{M},\Phi(\bar{G}),\Omega_1(\bar{G})=\lg \bar{x},\bar{a}^{p^2},\bar{b}^{p^2}\rg$
are characteristic in $G$ and $\bar{G}$ respectively, we may
assume that $\bar{a}^\theta=a^{l}x^nw$,
$\bar{b}^\theta=b^{i}x^jy$, $\bar{x}^\theta=x^rz$ where $l,i,r\in F_p^*$ and
$w,y\in \Phi(G)$, $z\in \Omega_1(M)$.

Since $[a^{l},b^{ip}]=[\bar{a}^\theta,(\bar{b}^p)^\theta]=[\bar{a},\bar{b}^p]^\theta=
(\bar{a}^{p^2})^\theta$,
we have $a^{lip^2}=a^{lp^2}$. Comparing
index of $a^{p^2}$ in two sides, we have $i=1$.

Since $[x^r,b]=[\bar{x}^\theta,\bar{b}^\theta]=
(\bar{a}^{p^2})^\theta=a^{lp^2}$,
we have $a^{rp^2}=a^{lp^2}$. Comparing
index of $a^{p^2}$ in two sides, we have $r=l$.

Since $[x^r,a^r]=[\bar{x}^\theta,\bar{a}^\theta]=
(\bar{b}^{\nu' p^2}\bar{a}^{k' p^2})^\theta=b^{\nu' p^2}a^{rk' p^2}$,
we have $b^{r^2\nu p^2}a^{r^2kp^2}=b^{\nu' p^2}a^{rk' p^2}$. Comparing
indexes of $a^{p^2}$ and $b^{p^2}$ in two sides, we have $\nu'=r^2\nu$ and $k'=rk$.

Since $\nu'=r^2\nu$, $\nu=\nu'=1$ or $\nu'=\nu$ is a fixed
quadratic non-residue modula $p$. Moreover, $r^2=1$ and hence $r=\pm 1$. It follows that $k'=\pm k$.

On the other hand, if $\nu'=\nu$ and $k'=-k$, then, $\theta:$
$\bar{a}\rightarrow a^{-1}$, $\bar{b}\rightarrow
b$, $\bar{x}\rightarrow x^{-1}$ is an isomorphism from $\bar{G}$ to
$G$.

\medskip

(3) All maximal subgroups of $G(\nu,k)$ are:

$N=\lg b,x,a^p\rg$;

$N_{i}=\lg ab^i,x,b^p\rg$ where $0\le i,j\le p-1$.

$N_{ij}=\lg ax^i,bx^j\rg$ where $0\le i,j\le p-1$.

It is easy to see that $N=\lg x,b\rg\ast\lg a^p\rg\in\mathcal{A}_2$.
Since $N_{ij}\cong \lg a,b\rg$, $N_{ij}\in\mathcal{A}_2$.
Hence $G(\nu,k)\in\mathcal{A}_3$ if and only if $N_i\in\mathcal{A}_2$ for any $i$.
Let $a_1=ab^i$ and $b_1=b^p$. Then
$N_i=\lg a_1,b_1,x\rg$ such that $$a_1^{p^3}=b_1^{p^2}=x^p=1, [a_1,b_1]=a_1^{p^2}b_1^{-ip},[x,a_1]=a_1^{(k+i)p^2}b_1^{(\nu-ik-i^2)p},[x,b_1]=1\rg.$$
By calculation, $[x,a_1]\neq 1$ for any $i$. Hence $\lg x,a_1\rg\in \mathcal{A}_1$.
Since $G(\nu,k)\in \mathcal{A}_3$, $|\lg x,a_1\rg|=p^5$.
It follows that $\nu-ik-i^2\neq 0$ for any $i$. That is,
equation $i^2+ik-\nu=0$ about $i$ has no solution. Hence $k^2+4\nu$ is not a square.

On the other hand, if $k^2+4\nu$ is not a square,
then it can be proved that $N_i\in\mathcal{A}_2$ for any $i$. Hence $G(\nu,k)\in A_3$ if and only if $k^2+4\nu$ is not a square.
 \qed

\begin{thm}\label{d=4-4}
Suppose that $G$ is an $\mathcal{A}_3$-group all of whose maximal
subgroups are $\mathcal{A}_2$-groups and there exists a
three-generator maximal subgroup in $G$, and $M'\le Z(G)$ for every
three-generator maximal subgroup $M$. Then $d(G)=3$ if and only if
 $G$ is isomorphic to one of the following pairwise
non-isomorphic groups:
\begin{enumerate}

\rr{Ni} $\Phi(G)\le Z(G)$ and $c(G)=2$. In this case, $(\mu_0,\mu_1,\mu_2)=(0,0,1+p+p^2)$ and
$\alpha_1(G)=p^4+p^3+p^2$.

\begin{enumerate}
\rr{N1}  $\langle a, b, c \di a^{p^{3}}=b^{p^{2}}=c^{p^{2}}=1,
[b,c]=a^{p^2},[c,a]=c^{p},[a,b]=b^{-p}\rangle$, where $p$ is odd;
 moreover, $|G|=p^{7}$,  $G'=\lg
a^{p^2},b^p,c^p\rg\cong C_p^3$ and $\Phi(G)=Z(G)=\lg a^p, b^p,c^p\rg
\cong C_{p^2}\times C_p\times C_p$.

\rr{N2} $\langle a, b, c; d \di a^{p^{2}}=b^{p^{2}}=c^{p^2}=d^p=1,
[b,c]=d, [c,a]=b^{p},[a,b]=c^{\nu
p},[d,a]=[d,b]=[d,c]=1\rangle$, where $p$ is odd, $\nu=1$ or a fixed
quadratic non-residue modulo $p$ such that $-\nu\not\in (F_p^*)^2$; moreover,
$|G|=p^{7}$, $\Phi(G)=Z(G)=\lg a^p, b^p,c^p,d\rg \cong C_p^4$
 and $G'=\lg
b^p,c^p,d\rg\cong C_p^3$ .

\rr{N3} $\langle a, b, c; d \di a^{p^{2}}=b^{p^{2}}=c^{p^2}=d^p=1,
[b,c]=d, [c,a]^{1+r}=b^{rp}c^{-p},[a,b]^{1+r}=b^{p}c^{p},[d,a]=[d,b]=[d,c]=1\rangle$,
where $p$ is odd, $-r\in F_p$ is not a square; moreover,
$|G|=p^{7}$, $\Phi(G)=Z(G)=\lg a^p, b^p,c^p,d\rg \cong C_p^4$
 and $G'=\lg
b^p,c^p,d\rg\cong C_p^3$ .

\rr{N4}  $\langle a, b, c \di a^{8}=b^{4}=c^{4}=1,
[b,c]=a^{4},[c,a]=c^{2},[a,b]=b^{2}\rangle$; where $|G|=2^{7}$,
$G'=\lg
a^{4},b^2,c^2\rg\cong C_2^3$ and  $\Phi(G)=Z(G)=\lg a^2, b^2,c^2\rg\cong C_4\times C_2\times C_2$.

\rr{N5} $\langle a, b, c; d\di a^{4}=b^{4}=c^{4}=d^2=1, [b,c]=d,
[c,a]=b^{2},[a,b]=b^{2}c^{2},[d,a]=[d,b]=[d,c]=1\rangle$, where
$|G|=2^{7}$, $\Phi(G)=Z(G)=\lg a^2, b^2,c^2,d\rg\cong C_2^4$ and
$G'=\lg
b^2,c^2,d\rg\cong C_2^3$.
\end{enumerate}

\rr{Nii} $Z(G)\nleq \Phi(G)$ and $c(G)=3$. In this case, $(\mu_0,\mu_1,\mu_2)=(0,0,1+p+p^2)$ and $\alpha_1(G)=p^3+p^2$.

\begin{enumerate}

\rr{N6} $\lg a, b,x\di
a^{p^{r+2}}=x^p=1,b^{p^{r+s+t}}=a^{p^{r+s}},[a,
b]=a^{p^r},[x,a]=[x,b]=1\rg=\lg a,b\rg\times \lg x\rg$, where $r\ge
2$ for $p=2$, $r\ge 1$ for $p\ge 3$, $t\ge 0$, $0\le s\le 2$ and
$r+s\ge 2$;  moreover, $|G|=p^{2r+s+t+3}$, $\Phi(G)=\lg a^p, b^p\rg
\cong C_{p^{r+t+1}}\times C_{p^{r+s-1}}$, $G'=\lg a^{p^r}\rg\cong
C_{p^2}$ and $Z(G)=\lg a^{p^2}, b^{p^2},x\rg \cong C_{p^{r+t}}\times
C_{p^{r+s-2}}\times C_p$.

\rr{N7} $\lg a, b,x\di a^{p^{3}}=b^{p^{t+3}}=1,x^p=a^{p^{2}},[a,
b]=a^{p},[x,a]=[x,b]=1\rg$, where $p\ge 3$ and $t\ge 0$;
 moreover, $|G|=p^{t+7}$, $\Phi(G)=\lg a^p, b^p\rg\cong C_{p^{t+2}}\times C_{p^2}$, $G'=\lg a^{p}\rg\cong C_{p^2}$,
$Z(G)=\lg a^{p^2}, b^{p^2},x\rg\cong C_{p^{t+1}}\times C_p^2$.

\rr{N8} $\lg a, b,x\di a^{p^{3}}=1, b^{p^{2}}=x^p=a^{p^{2}},[a,
b]=a^{p},[x,a]=[x,b]=1\rg$, where $p\ge 3$;  moreover, $|G|=p^{6}$,
$\Phi(G)=\lg a^p, b^p\rg\cong C_{p^2}\times C_p$, $G'=\lg
a^{p}\rg\cong C_{p^2}$ and $Z(G)=\lg x\rg\cong C_{p^2}$.

\rr{N9} $\lg a,b,x; c\di
a^{p^2}=b^{p^2}=c^p=x^p=1,[a,b]=c,[c,a]=b^{\nu p},
[c,b]=a^{p},[x,a]=[x,b]=1\rg=\lg a,b\rg\times \lg x\rg$, where $p\ge
5$, $\nu$ is a fixed quadratic non-residue modulo $p$;
 moreover, $|G|=p^{6}$, $\Phi(G)=G'=\lg a^p, b^p,c\rg\cong C_p^3$ and $Z(G)=\lg a^{p},
b^{p},x\rg\cong C_p^3$.

\rr{N10} $\lg a,b,x; c\di a^{p^2}=b^{p^2}=c^p=x^p=1,[a,
b]=c,[c,a]=a^{-p}b^{-lp},[c, b]=a^{-p},[x,a]=[x,b]=1\rg=\lg
a,b\rg\times \lg x\rg$, where $p\ge 5$, $4l=\rho^{2r+1}-1$, $r=1, 2,
\dots, \frac{1}{2}(p-1)$, $\rho$ is the smallest positive integer
which is a primitive root modulo $p$;  moreover, $|G|=p^{6}$,
$\Phi(G)=G'=\lg a^p, b^p,c\rg\cong C_p^3$ and $Z(G)=\lg a^{p},
b^{p},x\rg\cong C_p^3$.

\rr{N11} $\lg a, b,x; c\di a^9=b^9=c^3=x^3=1,[a,b]=c,[c,a]=b^{-3},
[c, b]=a^3,[a^3,b]=[x,a]=[x,b]=1\rg=\lg a,b\rg\times \lg x\rg$;
 moreover, $|G|=3^{6}$, $\Phi(G)=G'=\lg a^3, b^3,c\rg\cong C_3^3$ and $Z(G)=\lg
a^{3}, b^{3},x\rg\cong C_3^3$.

\rr{N12} $\lg a, b,x; c\di a^9=b^9=c^3=x^3=1, [a,b]=c, [c,a]=b^{-3},
[c, b]=a^{-3},[x,a]=[x,b]=1\rg=\lg a,b\rg\times \lg x\rg$; moreover,
$|G|=3^{6}$, $\Phi(G)=G'=\lg a^3, b^3,c\rg\cong C_3^3$ and $Z(G)=\lg
a^{3}, b^{3},x\rg\cong C_3^3$.
\end{enumerate}

\rr{Niii} $Z(G)<\Phi(G)$ and $G$ has at least two three-generator
maximal subgroups, in this case, $c(G)=2$ for ${\rm(N24)}$ and
$c(G)=3$ for else. Moreover, $(\mu_0,\mu_1,\mu_2)=(0,0,1+p+p^2)$
and $\alpha_1(G)=p^3+p^2$ except for {\rm (N14)}, {\rm (N18)} and {\rm (N22)--(N24)}.

\begin{enumerate}
\rr{N13} $\lg a,b,x\di
a^8=b^4=x^2=1,[a,b]=a^{-2},[x,a]=a^4,[x,b]=1\rg$; moreover, $|G|=2^{6}$,
$\Phi(G)=\lg a^2, b^2\rg\cong C_4\times C_2$, $G'=\lg a^{2}\rg\cong
C_4$, $Z(G)=\lg a^4, b^2\rg\cong C_2^2$, $(\mu_0,\mu_1,\mu_2)=(0,0,7)$ and $\alpha_1(G)=12$;

\rr{N14} $\lg a,b,x\di
a^8=b^8=x^2=1,[a,b]=a^{-2},[x,a]=b^4,[x,b]=a^4\rg$; moreover, $|G|=2^{7}$,
$\Phi(G)=\lg a^2, b^2\rg\cong C_4\times C_4$, $G'=\lg
a^{2},b^4\rg\cong C_4\times C_2$, $Z(G)=\lg a^4, b^2\rg\cong
C_4\times C_2$ and $\alpha_1(G)=14$;

\rr{N15} $\lg a,b,x\di
a^8=x^2=1,b^4=a^4,[a,b]=a^{-2},[x,a]=a^4,[x,b]=1\rg$;  where
$|G|=2^{6}$, $\Phi(G)=\lg a^2, b^2\rg\cong C_4\times C_2$, $G'=\lg
a^{2}\rg\cong C_4$ and $Z(G)=\lg b^2\rg\cong C_4$.

\rr{N16} $\lg a_1,b,x; a_2, a_3\di
a_1^p=a_2^p=a_3^p=b^{p}=x^p=1,[a_1,b]=a_{2},[a_2,b]=[x,a_1]=a_3,[a_3,b]=1,[x,b]=[a_i,a_j]=1\rg$,
where $p>3$ and $1\leq i,j\leq 3$;  where $|G|=p^{5}$,
$\Phi(G)=G'=\lg a_2,a_3\rg \cong C_p^2$ and $Z(G)=\lg a_3\rg\cong
C_p$.

\rr{N17} $\lg a_1,b,x; a_2\di
a_1^p=a_2^p=b^{p}=x^{p^2}=1,[a_1,b]=a_{2},[a_2,b]=[x,a_1]=x^p,[a_2,a_1]=[a_2,x]=[x,b]=1\rg$,
where $p>2$; moreover, $|G|=p^{5}$, $\Phi(G)=G'=\lg a_2,x^p\rg \cong
C_p^2$ and $Z(G)=\lg x^p\rg\cong C_{p}$.

\rr{N18} $\lg a_1,b,x; a_2\di
a_1^p=a_2^p=b^{p^2}=x^{p^2}=1,[a_1,b]=a_{2},[a_2,b]=x^p,[x,a_1]=b^p,[a_2,a_1]=[a_2,x]=[x,b]=1\rg$,
where $p>2$; moreover, $|G|=p^{6}$, $\Phi(G)=G'=\lg
a_2,x^p,b^p\rg\cong C_p^3$, $Z(G)=\lg b^p,x^p\rg\cong C_{p}\times
C_p$ and $\alpha_1(G)=p^2+2p^2-p$.

\rr{N19} $\lg a_1,b,x; a_2\di
a_1^p=a_2^p=b^{p^2}=x^{p}=1,[a_1,b]=a_{2},[a_2,b]=[x,a_1]=b^p,[a_2,a_1]=[a_2,x]=[x,b]=1\rg$,
where $p>2$; moreover, $|G|=p^{5}$, $\Phi(G)=G'=\lg a_2,b^p\rg \cong
C_p^2$ and  $Z(G)=\lg b^p\rg \cong C_p$.

\rr{N20} $\lg a_1,b,x; a_2\di
a_1^{p^2}=a_2^p=b^{p}=x^{p}=1,[a_1,b]=a_{2},[a_2,b]=a_1^{\nu p},
[x,a_1]=a_1^p,[a_2,a_1]=[a_2,x]=[x,b]=1\rg$, where $p>2$ and $\nu=1$
or a fixed quadratic non-residue modulo $p$; moreover, $|G|=p^{5}$,
$\Phi(G)=G'=\lg a_2,a_1^p\rg \cong C_p^2$ and  $Z(G)=\lg a_1^p\rg
\cong C_p$.

\rr{N21} $\lg a_1,b,x; a_2\di
a_1^{9}=a_2^3=x^3=1,b^{3}=a_1^{3},[a_1,b]=a_{2},[a_2,b]=a_1^{-3},
[x,a_1]=a_1^3,[a_2,a_1]=[a_2,x]=[x,b]=1\rg$;  moreover, $|G|=3^{5}$,
$\Phi(G)=G'=\lg a_2,a_1^3\rg \cong C_3^2$ and $Z(G)=\lg
a_1^3\rg\cong C_3$.

\rr{N22} $\lg a,b,x\di
a^{p^3}=b^{p^3}=x^p=1,[a,b]=a^p,[x,b]=a^{p^2},[x,a]=b^{\nu
p^2}a^{kp^2}\rg$, where $p>2$, $\nu=1$ or a fixed quadratic
non-residue modulo $p$, $0\le k\le \frac{p-1}{2}$ such that
$k^2+4\nu$ is not a square; moreover, $|G|=p^{7}$, $\Phi(G)=\lg a^p,
b^p\rg\cong C_{p^2}\times C_{p^2}$, $G'=\lg a^p,b^{p^2}\rg \cong
C_{p^2}\times C_p$, $Z(G)=\lg a^{p^2}, b^{p^2}\rg \cong C_p^2$, $(\mu_0,\mu_1,\mu_2)=(0,0,1+p+p^2)$ and $\alpha_1(G)=p^3+2p^2$.

\rr{N23} $\lg a,b,x\di
a^{p^3}=b^{p^{t+3}}=x^p=1,[a,b]=a^p,[x,b]=a^{p^2}b^{p^{t+2}},[x,a]=1\rg$,
where $p>2$ and $t\ge 1$; moreover, $|G|=p^{t+7}$, $\Phi(G)=\lg a^p,
b^p\rg \cong C_{p^2}\times C_{p^{t+2}}$, $G'=\lg
a^p,b^{p^{t+2}}\rg\cong C_{p^2}\times C_{p}$, $Z(G)=\lg a^{p^2},
b^{p^2}\rg \cong C_{p}\times C_{p^{t+1}}$, $(\mu_0,\mu_1,\mu_2)=(0,0,1+p+p^2)$ and $\alpha_1(G)=p^3+2p^2$.

\rr{N24} $\lg a,b,x\di
a^{p^{t+4}}=b^{p^3}=x^p=1,[a,b]=a^{p^{t+2}},[x,a]=b^{p^2},[x,b]=1\rg$,
where $t\ge 0$; moreover, $|G|=p^{t+8}$, $\Phi(G)=\lg a^p,
b^p\rg \cong C_{p^2}\times C_{p^{t+3}}$, $G'=\lg
a^{p^{t+2}},b^{p^2}\rg \cong C_{p^2}\times C_{p}$, $Z(G)=\lg a^{p^2},
b^{p^2}\rg \cong C_{p}\times C_{p^{t+2}}$, $(\mu_0,\mu_1,\mu_2)=(0,0,1+p+p^2)$ and $\alpha_1(G)=p^3+2p^2$.

\end{enumerate}

\rr{Niv} $Z(G)<\Phi(G)$, $c(G)=3$ and $G$ has a unique
three-generator maximal subgroup. In this case, $(\mu_0,\mu_1,\mu_2)=(0,0,1+p+p^2)$ and $\alpha_1(G)=2p^2$.

\begin{enumerate}

\rr{N25} $\lg a,b,d\di
a^{p^{m+1}}=b^{p^2}=d^p=1,[a,b]=a^{p^{m-1}},[d,a]=b^p,[d,b]=1\rg$,
where $p>2$ and $m\ge 2$; moreover, $|G|=p^{m+4}$, $\Phi(G)=\lg a^p,
b^p\rg\cong C_{p^m}\times C_{p}$, $G'=\lg a^{p^{m-1}},b^{p}\rg\cong
C_{p^2}\times C_{p}$ and $Z(G)=\lg a^{p^2}\rg\cong C_{p^{m-1}}$.

\rr{N26} $\lg a,b,d\di
a^{p^{3}}=b^{p^2}=d^p=1,[a,b]=a^{p},[d,a]=b^p,[d,b]=a^{\nu p^2}\rg$,
where $p>2$, $\nu=1$ or a fixed quadratic non-residue modulo $p$;
moreover, $|G|=p^{6}$, $\Phi(G)=G'=\lg a^p,b^{p}\rg\cong C_p\times
C_{p^2}$ and $Z(G)=\lg a^{p^2}\rg \cong C_p$.
\end{enumerate}
\end{enumerate}
\end{thm}

\demo If $\Phi(G)\le Z(G)$, then, since $G$ has no abelian maximal
subgroup, $\Phi(G)=Z(G)$. Hence $G$ is one of the groups listed in
\cite[Theorem 3.1,4.1,5.1,6.1,7.1-7.2,7.6]{QXA}. By hypothesis, the
minimal index of $\mathcal{A}_1$-subgroups and the maximal index of
$\mathcal{A}_1$-subgroups are $2$. By checking \cite[Theorem 3.3,
4.3,5.2,6.3,7.4-7.5,7.7]{QXA}, we get groups (N1)--(N5).

\begin{table}[h]
  \centering
{ 
\scriptsize
\begin{tabular}{c|c||c|c}
\hline
Groups & Groups in \cite[Theorem 4.1 \& 7.1]{QXA}& Groups &Groups in \cite[Theorem 4.1 \& 7.1]{QXA} \\
\hline
(N1)&(D1) where $m_1=2$ & (N4)&(M1) where $m_1=2$\\
(N2)&(E3) where $m_1=2$ and $-\nu\not\in F_p^2$ & (N5)&(N3) where $m_1=2$\\
(N3)&(E4) where $m_1=2$ and $-\nu\not\in F_p^2$ &&\\
\hline
\end{tabular}
\caption{ The correspondence from (N1)--(N5) to \cite[Theorem 4.1 \&
7.1]{QXA}} \label{table11}}
\end{table}

In the following, we may assume $\Phi(G)\not\le Z(G)$.

\medskip

{\bf Case 1:} $Z(G)\not\le \Phi(G)$.

\medskip

Then there exists a maximal subgroup $M$ such that $Z(G)\nleq M$.
Let $x\in Z(G)\setminus M$. Then $G=\lg M, x\rg$ and $x^p\in Z(M)$.
If $d(M)=3$, then $G'=M'\le Z(G)$ and $\exp(G')=p$. It follows that
$\Phi(G)\le Z(G)$, a contradiction. Hence $d(M)=2$. If $M$ has an
abelian subgroup  of index $p$, then $G$ also has an abelian
subgroup  of index $p$, a contradiction. Hence $\alpha_1(M)=1+p$.

\medskip

Since any maximal
subgroup of $M$ is an $\mathcal{A}_1$-group, $M$ is one of the
groups (17)--(21) in Lemma \ref{A_2}.

\medskip

Subcase 1.1: $M$ is the group of Type (17) in Lemma \ref{A_2}. That
is, $M=\lg a, b\di a^{p^{r+2}}=1,b^{p^{r+s+t}}=a^{p^{r+s}},[a,
b]=a^{p^r}\rg$, where $r\ge 2$ for $p=2$, $r\ge 1$ for $p\ge 3$,
$t\ge 0$, $0\le s\le 2$ and $r+s\ge 2$.

\medskip

If $s=2$, then $$Z(M)=\lg a^{p^2}\rg\times \lg b^{p^2}\rg,\
|M|=p^{2r+t+4}\ {\rm and}\ |G|=p^{2r+t+5}.$$ Since $x^p\in Z(G)$, we
may assume $x^p=a^{ip^2}b^{jp^2}$. Since $$|\lg
xb^{-jp},a\rg|=p^{r+3}<p^{2r+t+3}=\frac{|G|}{p^2},$$ we have
$[xb^{-jp},a]=1$ and hence $xb^{-jp}\in Z(G)$. By replacing $x$ with
$xb^{-jp}$ we get $x^p=a^{ip^2}$.

If $[xa^{-ip},b]=1$, then, replacing $x$ with $xa^{-ip}$, we get
$x^p=1$ and $G$ is a group of Type (N6).

If $[xa^{-ip},b]\ne 1$, then $(i,p)=1$. Without loss of generality,
we may assume $x^p=a^{p^2}$. Since $G\in \mathcal{A}_3$, we have
$$|\lg xa^{-ip},b\rg|=p^{r+t+4}=\frac{|G|}{p^2}=p^{2r+t+3}.$$ It
follows that $r=1$, $p\ge 3$ and $G$ is a group of Type (N7).

If $s=1$ and $r\ge 2$, then $$Z(M)=\lg a^{p^2}b^{-p^{t+2}}\rg\times
\lg b^{p^2}\rg,\ |M|=p^{2r+t+3}\ {\rm and}\ |G|=p^{2r+t+4}.$$ Since
$x^p\in Z(G)$, we may assume $x^p=(a^{p^2}b^{-p^{t+2}})^ib^{jp^2}$.
Since $$|\lg xa^{-ip},b\rg|=p^{r+t+3}<p^{2r+t+2}=\frac{|G|}{p^2},$$
we have $[xa^{-ip},b]=1$ and hence $xa^{-ip}\in Z(G)$. By replacing
$x$ with $xa^{-ip}$ we get $x^p=b^{j'p^2}$. Since $$|\lg
xb^{-j'p},a\rg|=p^{r+3}<p^{2r+t+2}=\frac{|G|}{p^2},$$ we have
$[xb^{-j'p},a]=1$ and hence $xb^{-j'p}\in Z(G)$. By replacing $x$
with $xb^{-j'p}$ we get $x^p=1$ and $G$ is a group of Type (N6).

If $s=1$ and $r=1$, then $$p\ge 3\ {\rm and}\ Z(M)=\lg b^{p^2}\rg,\
|M|=p^{t+5}\ {\rm and}\ |G|=p^{t+6}.$$ Since $x^p\in Z(G)$, we may
assume $x^p=b^{jp^2}$.

If $[xb^{-jp},a]=1$, then, replacing $x$ with $xb^{-jp}$, we get
$x^p=1$ and $G$ is a group of Type (N6).

If $[xb^{-jp},a]\ne 1$, then $(j,p)=1$. Without loss of generality,
we may assume $x^p=b^{p^2}$. Since $G\in \mathcal{A}_3$, we have
$$|\lg xb^{-p},a\rg|=p^{4}=\frac{|G|}{p^2}=p^{t+4}.$$ It follows that
$t=0$ and $G$ is a group of Type (N8).

If $s=0$ and $r\ge 3$, then $$Z(M)=\lg a^{p^2}b^{-p^{t+2}}\rg\times
\lg b^{p^2}\rg,\ |M|=p^{2r+t+2}\ a{\rm and}\ |G|=p^{2r+t+3}.$$ Since
$x^p\in Z(G)$, we may assume $x^p=(a^{p^2}b^{-p^{t+2}})^ib^{jp^2}$.
Since $$|\lg xa^{-ip},b\rg|=p^{r+t+3}<p^{2r+t+1}=\frac{|G|}{p^2},$$
we have $[xa^{-ip},b]=1$ and hence $xa^{-ip}\in Z(G)$. Replacing $x$
with $xa^{-ip}$, we get $x^p=b^{j'p^2}$. Since $$|\lg
xb^{-j'p},a\rg|=p^{r+3}<p^{2r+t+1}=\frac{|G|}{p^2},$$ we have
$[xb^{-j'p},a]=1$ and hence $xb^{-j'p}\in Z(G)$. Replacing $x$ with
$xb^{-j'p}$, we get $x^p=1$ and $G$ is a group of Type (N6).

If $s=1$ and $r\le 2$, then $$p\ge 3,\ r=2,\ Z(M)=\lg b^{p^2}\rg,\
|M|=p^{t+6}\ {\rm and}\ |G|=p^{t+7}.$$ Since $x^p\in Z(G)$, we may
assume $x^p=b^{jp^2}$. Since $$|\lg xb^{-jp},ab^{-p^t}\rg|\le
p^{4}<p^{t+5}=\frac{|G|}{p^2},$$ we have $[xb^{-jp},a]=1$ and hence
$xb^{-jp}\in Z(G)$. By replacing $x$ with $xb^{-jp}$ we get $x^p=1$
and $G$ is a group of Type (N6).

\medskip

Subcase 1.2: $M$ is the group of Type (18) in Lemma \ref{A_2}. That
is, $M=\lg a,b;c\di a^{p^2}=b^{p^2}=c^p=1,[a,b]=c,[c,a]=b^{\nu p},
[c,b]=a^{p}\rg$, where $p\ge 5$, $\nu$ is a fixed quadratic
non-residue modulo $p$.

\medskip

In this subcase, $Z(M)=\lg a^{p},b^p\rg$, $|M|=p^{5}$ and
$|G|=p^{6}$. Since $x^p\in Z(G)$, we may assume $x^p=a^{ip}b^{jp}$.
Since $|\lg c,xa^{-i}b^{-j}\rg|\le p^{3}$, we have
$[c,xa^{-i}b^{-jp}]=1$. It follows that $x^p=1$. Hence we get the
group (N9).

\medskip

Subcase 1.3: $M$ is one of the groups (19)--(21) in Lemma \ref{A_2}.

\medskip

By an argument similar to that of subcase 1.2, we get groups
(N10)--(N12) respectively. The details is omitted.

\medskip

{\bf Case 2:} $Z(G)<\Phi(G)$ and $G$ has at least two three-generator
maximal subgroups.

\medskip

Let $M_1,\cdots,M_s$ be all three-generator maximal subgroups of
$G$, $N=M_1'M_2'\dots M_s'$ and $\bar{G}=G/N$. Then $\bar{G}$ has at
least two abelian subgroups of index $p$ and every non-abelian
subgroup of $\bar{G}$ is generated by two elements. Since
$\Phi(G)\not\le Z(G)$, $\bar{G}$ is not abelian. By Lemma
\ref{Dp3=>A_2} we have $\bar{G}\in \mathcal{A}_2$. Hence $\bar{G}$
is a group of Type (8)--(12) in Lemma \ref{A_2}. Let $\bar{G}=\lg
\bar{a},\bar{b},\bar{x}\rg$ and $K=\lg a,b,\Phi(G)\rg$. Then
$d(K)=2$ since $\bar{K}$ is not abelian. Since $d(G)=3$, we have
$\Phi(K)=\Phi(G)$. It follows that $K=\lg a,b\rg$ is maximal in $G$
and $K\in \mathcal{A}_2$. Hence $K$ is a group of Type (1)--(7) or
(17)--(21) in Lemma \ref{A_2}.

\medskip

Subcase 2.1: $K$ is a group of Type (1) in Lemma \ref{A_2}. That is,
$K=\lg a,b\di a^{8}=b^{2^m}=1, a^b=a^{-1}\rg$.

\medskip

If $m=1$, then $\Phi(G)=\Phi(K)=\lg a^2\rg$ and $N=\lg a^4\rg$. By
calculation, a maximal subgroup $\lg a,c\rg$ of $G$ is abelian or
minimal non-abelian, a contradiction. Hence $$m\ge 2,\
\Phi(G)=\Phi(K)=\lg a^2,b^2\rg,\ \lg a^4\rg\le N\le \lg
a^4,b^{2^{m-1}}\rg\ {\rm and}\ |G|=2^{m+4}.$$ Let $M=\lg
x,a,b^2\rg$. Then $M\in\mathcal{A}_2$. Since $x^2\in Z(K)\cap
\Phi(M)$, we may assume that $x^2=a^{4i}b^{4j}$. By replacing $x$
with $xa^{2i}b^{-2j}$ we get $x^2=1$. Since $[x,b^2]=[a,b^2]=1$,
$[x,a]\ne 1$. Since $G\in \mathcal{A}_3$, $|\lg x,a\rg|=2^{m+2}$.
Since $|\lg x,a\rg|\le 2^5$, $m\le 3$.

If $m=2$, then $N=\lg a^4\rg$. We have $[x,a]=a^4$ and
$[x,b]=a^{4k}$. By replacing $b$ with $ba^{k}$ we get $[x,b]=1$.
Hence $G$ is the group (N13).

If $m=3$, then $|\lg x,a\rg|=2^5$. It follows that $N=\lg
a^4,b^4\rg$ and we may assume that $[x,a]=b^4a^{4i}$. If
$[x,a]=b^4a^4$, then $[x,b]=b^4$ or $[x,ba]=b^4$. It follows that
$|\lg x,b\rg|=2^4$ or $|\lg x,ba\rg|=2^4$, a contradiction. Hence
$[x,a]=b^4$. By suitable replacement we get $[x,b]=a^4$. Hence $G$
is the group (N14).

\medskip

Subcase 2.2: $K$ is a group of Type (2) in Lemma \ref{A_2}. That is,
$K=\lg a,b\di a^{8}=b^{2^m}=1, a^b=a^{3}\rg$.

\medskip

If $m=1$, then $\Phi(G)=\Phi(K)=\lg a^2\rg$ and $N=\lg a^4\rg$. By
calculation, a maximal subgroup $\lg a,c\rg$ of $G$ is abelian or
minimal non-abelian, a contradiction. Hence $$m\ge 2,\
\Phi(G)=\Phi(K)=\lg a^2,b^2\rg,\ \lg a^4\rg\le N\le \lg a^4,
b^{2^{m-1}}\rg\ {\rm and}\ |G|=2^{m+4}.$$ Let $M=\lg x,a,b^2\rg$.
Then $M\in\mathcal{A}_2$. Since $x^2\in Z(K)\cap \Phi(M)$, we may
assume $x^2=a^{4i}b^{4j}$. By replacing $x$ with $xa^{2i}b^{-2j}$ we
get $x^2=1$. Since $[x,b^2]=[a,b^2]=1$, $[x,a]\ne 1$. If
$[x,a]=a^4$, then, replacing $b$ with $bx$, it is reduced to Subcase
2.1. Hence $|\lg x,a\rg|=2^5$, $m=3$, $N=\lg a^4,b^4\rg$ and we may
assume $[x,a]=b^4a^{4i}$. If $[x,a]=b^4a^4$, then $[x,b]=b^4$ or
$[x,ba]=b^4$. It follows that $|\lg x,b\rg|=2^4$ or $|\lg
x,ba\rg|=2^4$, a contradiction. Hence $[x,a]=b^4$. By suitable
replacement we get $[x,b]=a^4$. By replacing $a$ and $b$ with $ax$
and $bx$, respectively, we get the group (N14) again.

\medskip

Subcase 2.3: $K$ is a group of Type (3) in Lemma \ref{A_2}. That is,
$K=\lg a,b\di a^{8}=1, b^{2^m}=a^{4}, a^b=a^{-1}\rg$.

\medskip

In this subcase, $\Phi(G)=\Phi(K)=\lg a^2,b^2\rg$ and $N=\lg
a^4\rg$. Let $M=\lg x,a,b^2\rg$. Then $M\in\mathcal{A}_2$. Since
$x^2\in Z(K)\cap \Phi(M)$, we may assume that $x^2=b^{4i}$. By
replacing $x$ with $xb^{-2j}$ we get $x^2=1$. Since
$[x,b^2]=[a,b^2]=1$, $[x,a]=a^4$. It follows that $m=2$. By suitable
replacement we get $[x,b]=1$. Hence $G$ is the group (N15).

\medskip

Subcase 2.4: $K$ is a group of Type (4) in Lemma \ref{A_2}. That is,
$K=\lg a_1,b;a_2,a_3\di
a_1^p=a_2^p=a_3^p=b^{p^m}=1,[a_1,b]=a_{2},[a_2,b]=a_3,[a_3,b]=1,[a_i,a_j]=1\rg,$
where $p\ge 5$ for $m=1$, $p\ge 3$ and $1\leq i,j\leq 3$.

\medskip

If $m=1$, then $p\ge 5$ and $N=\lg a_3\rg$. Let $M=\lg
x,a_1,a_2\rg$. Then $M\in\mathcal{A}_2$. Since
$[x,a_2]=[a_1,a_2]=1$,  $[x,a_1]\ne 1$. Without loss of generality
assume $[x,a_1]=a_3$. By suitable replacement we get $[x,b]=1$.

If $x^p=1$, then $G$ is the group (N16).

If $x^p\ne 1$, then we may assume that $a_3=x^{ip}$, where
$(i,p)=1$. By replacing $a_1$ and $a_2$ with $a_1^{i^{-1}}$ and
$[a_1^{i^{-1}},b]$, respectively, we get $[a_2,b]=[x,a_1]=x^p$.
Hence $G$ is the group (N17).

If $m\ge 2$, then $$\Phi(G)=\Phi(K)=\lg a_2,a_3,b^p\rg\ {\rm and}\
\lg a_3\rg\le N\le \lg a_3,b^{p^{m-1}}\rg.$$ Let $M=\lg
x,a_1,a_2,a_3,b^p\rg$. Then $M\in\mathcal{A}_2$. Since $x^p\in
Z(K)\cap \Phi(M)$, we may assume $x^p=a_3^{i}b^{jp}$. If $(j,p)=1$,
then $\lg a_1,xb^{-j}\rg\in \mathcal{A}_2$ and is of order $p^4$, a
contradiction. Hence $p\di j$. By replacing $x$ with $xb^{-j}$ we
get $x^p=a_3^i$. Since $[x,a_2]=[a_1,a_2]=1$, $[x,a_1]\ne 1$. Since
$|\lg x,a_1\rg|\le p^4$, we have $m=2$, $(i,p)=1$ and
$[x,a_1]\not\in\lg a_3\rg$. It follows that $N=\lg a_3,b^p\rg$.
Without loss of generality assume that $[a_2,b]=x^p$ and
$[x,a_1]=b^{jp}x^{kp}$, where $(j,p)=1$. Let $[x,b]=b^{sp}x^{tp}$.
By replacing $x$ and $b$ with $xa_2^{j^{-1}s-t}$ and
$ba_1^{-j^{-1}s}$, respectively, we get $[x,b]=1$. By replacing $b$
with $bx^{j^{-1}k}$ we get $[x,a_1]=b^{jp}$. By replacing $b$ and
$x$ with $b^{j^{-1}}$ and $x^{j^{-2}}$, respectively, we get
$[x,a_1]=b^p$. Hence $G$ is the group (N18).

\medskip

Subcase 2.5: $K$ is a group of Type (5) in Lemma \ref{A_2}. That is,
$K=\lg a_1,b;a_2\di
a_1^p=a_2^p=b^{p^{m+1}}=1,[a_1,b]=a_{2},[a_2,b]=b^{p^m},[a_1,a_2]=1\rg$,
where $p>2$.

\medskip

In this subcase, $N=\lg b^{p^m}\rg$. Let $M=\lg x,a_1,a_2,b^p\rg$.
Then $M\in\mathcal{A}_2$. Since $[x,a_2]=[a_1,a_2]=1$,  $[x,a_1]\ne
1$. Without loss of generality assume $[x,a_1]=b^{p^m}$. Since
$x^p\in Z(K)$, we may assume $x^p=b^{ip}$. If $(i,p)=1$, then $\lg
a_1,xb^{-i}\rg\in \mathcal{A}_2$ and is of order $p^4$. If $p\di i$,
then $\lg a_1,xb^{-i}\rg\in \mathcal{A}_1$ and is of order $p^3$.
Hence $m=1$ and $|G|=p^5$. By suitable replacement we get $[x,b]=1$.

If $x^p=1$, then $G$ is the group of Type (N19).

If $x^p=b^{ip}\ne 1$, then, replacing $a_1$ and $b$ with $a_1^{i}$
and $bx^{-i^{-1}}$, respectively, we get $b^p=1$ and
$[a_2,b]=[x,a_1]=x^p$. Hence $G$ is the group (N17).

\medskip

Subcase 2.6: $K$ is a group of Type (6) in Lemma \ref{A_2}. That is,
$K=\lg a_1,b\di
a_1^{p^2}=a_2^p=b^{p^m}=1,[a_1,b]=a_{2},[a_2,b]=a_1^{\nu
p},[a_1,a_2]=1\rg,$ where $p>2$ and $\nu=1$ or a fixed quadratic
non-residue modulo $p$.

\medskip

If $m=1$, then $N=\lg a_1^p\rg$. Let $M=\lg x,a_1,a_2\rg$. Then
$M\in\mathcal{A}_2$. Since $[x,a_2]=[a_1,a_2]=1$, $[x,a_1]\ne 1$.
Without loss of generality assume $[x,a_1]=a_1^p$. By suitable
replacement we get $[x,b]=1$.

If $x^p=1$, then $G$ is the group (N20).

If $x^p=a_1^{ip}\ne 1$, then, by suitable replacement, we get the
group (N17).

If $m\ge 2$, then $$\Phi(G)=\Phi(K)=\lg a_2,a_1^p,b^p\rg\ {\rm and}\
\lg a_1^p\rg\le N\le \lg a_3,b^{p^{m-1}}\rg.$$ Let $M=\lg
x,a_1,a_2,b^p\rg$. Then $M\in\mathcal{A}_2$. Since $x^p\in Z(K)\cap
\Phi(M)$, we may assume that $x^p=a_1^{ip}b^{jp}$. If $(j,p)=1$,
then $\lg a_1,xb^{-j}\rg\in \mathcal{A}_2$ and is of order $p^4$, a
contradiction. Hence $p\di j$. By replacing $x$ with $xb^{-j}$ we
get $x^p=a_1^{ip}$. If $x^p\ne 1$, then, replacing $a_1$ with
$a_1x^{-i^{-1}}$, it is reduced to Subcase 2.4. Hence we may assume
$x^p=1$. Since $[x,a_2]=[a_1,a_2]=1$, we have $[x,a_1]\ne 1$. Since
$|\lg x,a_1\rg|\le p^4$, we have $m=2$ and $[x,a_1]\not\in\lg
a_1^p\rg$. It follows that $N=\lg a_1^p,b^p\rg$ and
$[x,a_1]=b^{jp}a_1^{kp}$, where $(j,p)=1$. Let
$[x,b]=b^{sp}a_1^{tp}$. By replacing $x$ with $xa_2^{-\nu^{-1}t}$ we
get $[x,b]=b^{sp}$. Since $|\lg x,b\rg|=p^3$, $[x,b]=1$. Assume
$(ba_1)^p=b^pa_1^{rp}$. By calculation we have
$[x^{j^{-1}}a_2^{r-j^{-1}k},ba_1]=(ba_1)^p$. It follows that $\lg
x^{j^{-1}}a_2^{r-j^{-1}k},ba_1\rg\in \mathcal{A}_1$ and is of order
$p^3$, a contradiction.

\medskip

Subcase 2.7: $K$ is a group of Type (7) in Lemma \ref{A_2}. That is,
$K=\lg a_1,b;a_2\di
a_1^9=a_2^3=1,b^3=a_1^3,[a_1,b]=a_2,[a_2,b]=a_1^{-3}\rg$.

\medskip

In this subcase, $N=\lg a^3\rg=\lg b^3\rg$. If $x^p=b^{3i}\ne 1$,
then, replacing $b$ with $bx^{-i}$, it is reduced to Subcase 2.6.
Hence we may assume $x^3=1$. Let $M=\lg x,a_1,a_2\rg$. Then
$M\in\mathcal{A}_2$. Since $[x,a_2]=[a_1,a_2]=1$, $[x,a_1]\ne 1$.
Without loss of generality assume $[x,a_1]=a_1^3$. By suitable
replacement we get $[x,b]=1$. Hence $G$ is the group (N21).

\medskip

Subcase 2.8: $K$ is a group of Type (17) in Lemma \ref{A_2}. That
is, $K=\lg a, b\di a^{p^{r+2}}=1,b^{p^{r+s+t}}=a^{p^{r+s}},[a,
b]=a^{p^r}\rg$, where $r\ge 2$ for $p=2$, $r\ge 1$ for $p\ge 3$,
$t\ge 0$, $0\le s\le 2$ and $r+s\ge 2$.

\medskip

In this subcase, $$Z(K)=\lg a^{p^2}, b^{p^2}\rg,\ |K|=p^{2r+s+t+2}\
{\rm and}\ |G|=p^{2r+s+t+3}.$$ If $N=\lg a^{p^{r+1}}\rg$, then,
assuming that $[x,a]=a^{ip^{r+1}}$ and $[x,b]=a^{jp^{r+1}}$, we get
$xa^{-jp}b^{ip}\in Z(G)$, a contradiction. Hence $|N|=p^2$. Since
$x^p\in Z(G)$, we may assume that $x^p=a^{ip^2}b^{jp^2}$. By
replacing $x$ with $xa^{-ip}b^{-jp}$ we get $x^p=1$.

If $s=2$, then $$Z(K)=\lg a^{p^2}\rg\times \lg b^{p^2}\rg,\
|K|=p^{2r+t+4},\ |G|=p^{2r+t+5}\ {\rm and}\ N=\lg
a^{p^{r+1}},b^{p^{r+t+1}}\rg.$$ Without loss of generality assume
$[x,b]\ne 1$. Since $|\lg x,b\rg|\le p^{r+t+4}$, $r+t+4\ge 2r+t+3$.
It follows that $r=1$, $p\ge 3$ and $[x,b]\not\in\lg
b^{p^{t+2}}\rg$. We may assume $[x,b]=a^{p^2}b^{kp^{t+2}}$.

If $[x,a]\ne 1$, then $t=0$, $|G|=p^7$ and $[x,a]\not\in \lg
a^{p^2}\rg$. By suitable replacement we may assume that
$[x,b]=a^{p^2}$ and $[x,a]=b^{\nu p^2}a^{kp^2}$, where $\nu =1$ or a
fixed quadratic non-residue modulo $p$. By Lemma \ref{cong 4} we
have $G$ is the group (N22).

If $[x,a]=1$, then $(k,p)=1$ since $xa^{-p}\not\in Z(G)$. By
replacing $a$ and $x$ with $a^{k^{-1}}$ and $x^{k^{-1}}$,
respectively, we have $[x,b]=a^{p^2}b^{p^{t+2}}$. If $t=0$, then
$\lg x,ba\rg$ is not abelian and is of order $p^4$, a contradiction.
Hence $t\ge 1$ and $G$ is the group (N23).

If $s=1$, then $$|K|=p^{2r+t+3},\ |G|=p^{2r+t+4},\ r\ge 2\ {\rm
and}\ N=\lg a^{p^{r+1}},a^{-p^r}b^{p^{r+t}}\rg.$$ Since $|\lg
x,a^{-1}b^{p^t}|\le p^{r+3}<p^{2r+t+2}$, $[x,a^{-1}b^{p^t}]=1$.
Since $Z(G)<\Phi(G)$, $[x,b]\ne 1$. Since $|\lg
x,b\rg|=p^{r+t+4}=p^{2r+t+2}=\frac{|G|}{p^2}$, $r=2$. Without loss
of generality assume $[x,b]=a^{-p^2}b^{p^{t+2}}b^{ip^{t+3}}$. By
replacing $a$ with $axa^{-ip}$ we have $[a,b]=b^{p^{t+2}}$.  By
replacing $a$ and $b$ with $b$ and $a^{-1}b^{p^t}b^{ip^{t+1}}$,
respectively, we get $$G=\lg a,b,x\di
a^{p^{t+4}}=b^{p^{3}}=x^p=1,[a,b]=a^{p^{t+2}},[x,b]=1,[x,a]=b^{p^2}\rg.$$
Hence $G$ is a group of Type (N24).

If $s=0$, then $r\ge 2$. By suitable replacement we have $$K=\lg
a,b\di a^{p^{r+t+2}}=b^{p^r}=1,[a,b]=a^{p^{r+t}}\rg.$$ By
calculation we have $$Z(K)=\lg a^{p^2}\rg\times \lg b^{p^2}\rg,\
|K|=p^{2r+t+2},\ |G|=p^{2r+t+3},\ r\ge 3\ {\rm and}\ N=\lg
a^{p^{r+t+1}},b^{p^{r-1}}\rg.$$ Since $|\lg
x,a\rg|=p^{r+t+4}=\frac{|G|}{p^2}$, $r=3$. Without loss of
generality assume  $[x,a]=b^{p^{2}}a^{kp^{t+4}}$. By replacing $b$
with $ba^{kp^{t+2}}$ we get $[x,a]=b^{p^2}$. Hence $G$ is a group of
Type (N24).

\medskip

Subcase 2.9: $K$ is a group of Type (18) in Lemma \ref{A_2}. That
is, $K=\lg a,b\di a^{p^2}=b^{p^2}=c^p=1,[a,b]=c,[c,a]=b^{\nu p},
[c,b]=a^{p}\rg$, where $p\ge 5$, $\nu$ is a fixed quadratic
non-residue modulo $p$.

\medskip

In this subcase, $|G|=p^6$ and $N=Z(K)=\lg a^p,b^p\rg$. We claim
that $x^p=1$. Otherwise, since $x^p\in Z(K)$, we may assume that
$x^p=a^{ip}b^{jp}$, where $(i,j)\ne (0,0)$. If $i\ne 0$, then $\lg
xa^{-i},c\rg$ is not abelian and is of order $p^3$, a contradiction.
If $j\ne 0$, then $\lg xb^{-j},c\rg$ is not abelian and is of order
$p^3$, a contradiction again. Hence $x^p=1$. Let
$[x,a]=a^{sp}b^{tp}$. Since $|\lg xc^{-\nu^{-1}t},a\rg|=p^3$, we
have $[xc^{-\nu^{-1}t},a]=1$ and hence $[x,a]=b^{tp}$. By replacing
$x$ with $xc^{-\nu^{-1}t}$ we get $[x,a]=1$. The same reason gives
that $[x,b]=a^{up}$. Since $Z(G)<\Phi(G)$, $(u,p)=1$. Without loss
of generality assume $[x,b]=a^p$. Then $\lg cx^{\nu-1}, ab\rg$ is
not abelian and is of order $p^3$, a contradiction.

\medskip

Subcase 2.10: $K$ is a group of Type (19)--(21) in Lemma \ref{A_2}.
That is, $K=\lg a,b;c\di a^{p^2}=b^{p^2}=c^p=1,[a,b]=c,[c,a]=b^{\nu
p}, [c,b]=a^{p}\rg$, where $p\ge 5$, $\nu$ is a fixed quadratic
non-residue modulo $p$.

\medskip

By an argument similar to that of subcase 2.9, we get a
contradiction, respectively. The details is omitted.

\medskip

{\bf Case 3:} $Z(G)<\Phi(G)$ and $M$ is the unique three-generator maximal
subgroup of $G$.

\medskip

Then $G/M'$ has a unique abelian subgroup $M/M'$  of index $p$, and
all non-abelian subgroups of $G/M'$ are generated by two elements.
By Lemma \ref{Dp3=>A_2}, $G/M'\in \mathcal{A}_2$. Hence $G$ a group
of Type (13)--(16) in Lemma \ref{A_2}.

\medskip

Subcase 3.1: $G/M'$ is a group of Type (13) in Lemma \ref{A_2}. That
is, $G/M'=\lg \bar{a},\bar{b},\bar{c}\di \bar{a}^{4}=\bar{b}^{4}=1,
\bar{c}^2=\bar{a}^2\bar{b}^2,[\bar{a},\bar{b}]=\bar{b}^2,[\bar{c},\bar{a}]=\bar{a}^2,[\bar{c},\bar{b}]=1\rg$.

\medskip

Let $K=\lg a,b,\Phi(G)\rg$. Since $K/M'$ is not abelian, we have
$d(K)=2$ and $K\in\mathcal{A}_2$. Hence $\Phi(K)=\Phi(G)$ and $K=\lg
a,b\rg$. By calculation we get $$[c,a^2]=[c,a]^2[c,a,a]=a^4\ {\rm
and}\ [c^2,a]=[c,a]^2[c,a,c]=a^8=1.$$ Since $c^2\equiv a^2b^2\ (\mod
M')$, $[a,b^2]=1$. It follows that $b^4=[a,b^2]=1$ and hence $\lg
a,b\rg\in\mathcal{A}_1$, a contradiction.

\medskip

Subcase 3.2: $G/M'$ is a group of Type (14) in Lemma \ref{A_2}. That
is, $G/M'=\lg \bar{a},\bar{b},\bar{d}\di
\bar{a}^{p^{m}}=\bar{b}^{p^2}=\bar{d}^{p}=1,
[\bar{a},\bar{b}]=\bar{a}^{p^{m-1}},[\bar{d},\bar{a}]=\bar{b}^p,[\bar{d},\bar{b}]=1\rg$,
where $m\geq 3$ if $p=2$.

\medskip

Let $K=\lg a,b,\Phi(G)\rg$ and $L=\lg d,a,\Phi(G)\rg$. Since $K/M'$
and $L/M'$ is not abelian, we have $$d(K)=2,\ d(L)=2,\
K\in\mathcal{A}_2\ {\rm and}\ L\in\mathcal{A}_2.$$ Hence
$$\Phi(K)=\Phi(L)=\Phi(G),\ K=\lg a,b\rg\ {\rm and}\ L=\lg d,a\rg.$$ By
calculation, $b^{p^2}=[d^p,a]=1$. Since $L\in \mathcal{A}_2$,
$[a,b^p]\ne 1$. It follows that
 $a^{p^m}=[a,b^p]\ne 1$ and hence $o(a)=p^{m+1}$. Let $N=\lg
 a^p,b\rg$. Then $N\in \mathcal{A}_1$. Since $|N|=p^{m+2}$, we get
 $|G|=p^{m+4}$ and hence $M'=\lg a^{p^m}\rg$. We claim that $p>2$.
 Otherwise, $\lg d,a^2\rg$ is not abelian and is of order $2^{m+1}$,
 a contradiction. Hence $p>2$. Let $d^p=a^{ip^m}$.
 By replacing $d$ with $da^{-ip^{m-1}}$ we get $d^p=1$. By suitable
 replacement we get $[a,b]=a^{p^{m-1}}$ and
 $[d,a]=b^p$.

 If $[d,b]=1$ then $G$ is the group (N25).

If $[d,b]\ne 1$, then, since $|\lg d,b\rg|=p^4$, we have $m=2$.
Assume that $[d,b]=a^{jp^2}$, where $j=k^2\nu$, $\nu =1$ or a fixed
quadratic non-residue modulo $p$. By replacing $a$ and $d$ with
$a^k$ and $d^{k^{-1}}$, respectively, we get $[d,b]=a^{\nu p^2}$.
Hence $G$ is the group (N26).

\medskip

Subcase 3.3: $G/M'$ is a group of Type (15) in Lemma \ref{A_2}. That
is, $\lg \bar{a},\bar{b},\bar{d}\di
\bar{a}^{p^m}=\bar{b}^{p^2}=\bar{d}^{p^2}=1,
[\bar{a},\bar{b}]=\bar{d}^p,[\bar{d},\bar{a}]=\bar{b}^{jp},[\bar{d},\bar{b}]=1\rg$,
where $(j,p)=1$, $p>2$, $j$ is a fixed quadratic non-residue modulo
$p$, and $-4j$ is a quadratic non-residue modulo $p$.

\medskip

Let $K=\lg a,b,\Phi(G)\rg$ and $L=\lg d,a,\Phi(G)\rg$. Since $K/M'$
and $L/M'$ is not abelian, we have $$d(K)=2,\ d(L)=2,\
K\in\mathcal{A}_2\ {\rm and}\ L\in\mathcal{A}_2.$$ Hence
$$\Phi(K)=\Phi(L)=\Phi(G),\ K=\lg a,b\rg\ {\rm and}\ L=\lg d,a\rg.$$ By
Lemma \ref{A_2-property} (7), $\exp(K')=\exp(L')=p$. It follows that
$d^{p^2}=b^{p^2}=1$. By calculation we have $[d^p,b]=[d,b]^p=1$ and
$[d^p,a]=[d,a]^p=1$. Hence $K\in\mathcal{A}_1$, a contradiction.

\medskip

Subcase 3.4: $G/M'$ is a group of Type (16) in Lemma \ref{A_2}. That
is, $\lg \bar{a},\bar{b},\bar{d}\di
 \bar{a}^{p^m}=\bar{b}^{p^2}=\bar{d}^{p^2}=1,[\bar{a},\bar{b}]=\bar{d}^p,
 [\bar{d},\bar{a}]=\bar{b}^{jp}\bar{d}^p,[\bar{d},\bar{b}]=1\rg$,
where if $p$ is odd, then $4j =1-\rho^{2r+1}$ with $1\le
r\le\frac{p-1}{2}$ and $\rho$ the smallest positive integer which is
a primitive root $(\mod p)$; if $p = 2$, then $j = 1$.

\medskip
Let $K=\lg a,b,\Phi(G)\rg$ and $L=\lg d,a,\Phi(G)\rg$. Since $K/M'$
and $L/M'$ is not abelian, we have $$d(K)=2,\ d(L)=2,\
K\in\mathcal{A}_2\ {\rm and}\ L\in\mathcal{A}_2.$$ Hence
$$\Phi(K)=\Phi(L)=\Phi(G),\ K=\lg a,b\rg\ {\rm and}\ L=\lg d,a\rg.$$ By
Lemma \ref{A_2-property} (7), $\exp(K')=\exp(L')=p$. It follows that
$d^{p^2}=b^{p^2}=1$. By calculation we have $[d^p,b]=[d,b]^p=1$ and
$[d^p,a]=[d,a]^p=1$. Hence $K\in\mathcal{A}_1$, a contradiction.

\medskip
We calculate the $(\mu_0,\mu_1,\mu_2)$ and $\alpha_1(G)$ of those
groups in Theorem \ref{d=4-4} as follows.

Since $d(G)=3$, $(\mu_0,\mu_1,\mu_2)=(0,0,1+p+p^2)$. In the
following, we calculate $\alpha_1(G)$.

\medskip
{\bf Case 1.} $G$ is one of the groups (N1)--(N5).

Since $\Phi(G)\le Z(G)$, by Theorem \ref{alpha of d=3},
$\alpha_1(G)=\mu_1+\mu_2 p^2=p^4+p^3+p^2$.

\medskip
{\bf Case 2.} $G$ is one of the groups (N6)--(N12).

In this case, $G=M\ast \lg x\rg$ where $M=\lg a,b\rg$ such that $\alpha_1(M)=1+p$ and $x^p\in Z(M)$.
Other maximal subgroups of $G$ are:

$N_i=M_i \lg x\rg$ where $M_i\maxsgp M$;

$N_{ij}=\lg ax^i,bx^j,\Phi(M)\rg$ where $0\le i,j\le p-1$.

Since $d(N_i)=3$ and
$N_i'=M_i'$, by Lemma \ref{A_2-property} (7), $\alpha_1(N_i)=p^2$. Notice that $G=N_{ij}\ast \lg x\rg$.
If $d(N_{ij})=3$, then $G'=N_{ij}'\le Z(G)$ and $\exp(G')=p$. It follows that
$\Phi(G)\le Z(G)$, a contradiction. Hence $d(N_{ij})=2$. If $N_{ij}$ has an
abelian subgroup  of index $p$, then $G$ also has an abelian
subgroup  of index $p$, a contradiction. Hence $\alpha_1(N_{ij})=1+p$. Let $H\in\Gamma_2$.
Then $H=\lg xm,\Phi(G)\rg$ where $m\in M$ or $H=\lg m,\Phi(G)\rg$ where $m\in M\setminus\Phi(G)$.
It is obvious that $H'=1$ if and only if $H=\lg x,\Phi(G)\rg$. Thus $\sum_{H\in
\Gamma_2} \alpha_1(H)=p+p^2$.
 By Hall's enumeration principle,
$$\alpha_1(G)=\sum_{H\in
\Gamma_1} \alpha_1(H)-p\sum_{H\in \Gamma_2} \alpha_1(H)=(1+p)\times
p^2+p^2\times(1+p)-p\times (p+p^2)=p^3+p^2.$$

\medskip
{\bf Case 3.} $G$ is one of the groups (N13)--(N21) except for (N14)
and (N18).

In this case, $G=M\ast \lg x\rg$ where $\lg x\rg\cong C_p$ and $M=\lg a,b\rg$ such that
$d(M)=2$, $|M'|=p^2$, $[x,b]=1$, $[x,a]\in M_3$ and $\lg a,\Phi(M)\rg$
is the unique abelian maximal subgroup of $G$. Hence all maximal subgroups are:

$N_i=K_i\ast \lg x\rg$ where $K_i$ are maximal subgroups of $M$;

$N_{ij}=\lg ax^i,bx^j\rg$ where $0\le i,j\le p-1$.

It is easy to see that $K_i'=M_3$ and hence $|N_i'|=p$. By Lemma \ref{A_2-property},
$\alpha_1(N_i)=p^2$. Since $N_{ij}\cong M$, $\alpha_1(N_{ij})=p$.
Let $H\in\Gamma_2$. Then $H=\lg a^ib^jx^k,\Phi(G)\rg$.
It is obvious that $H\in \mathcal{A}_1$ if and only if $j\ne 0$. Hence $\sum_{H\in
\Gamma_2} \alpha_1(H)=p^2$.
 By Hall's enumeration principle,
$$\alpha_1(G)=\sum_{H\in
\Gamma_1} \alpha_1(H)-p\sum_{H\in \Gamma_2} \alpha_1(H)=(1+p)\times
p^2+p^2\times p-p\times p^2=p^3+p^2.$$

\medskip
{\bf Case 4.} $G$ is either the group (N14) or (N18).

In this case, $G=M\ast \lg x\rg$ where $\lg x\rg\cong C_p$ and $M=\lg a,b\rg$
such that $d(M)=2$, $|M'|=p^2$, $[x,b]\in M_3$, $[x,a]\not\in M'$ and $\lg a,\Phi(M)\rg$
 is the unique abelian maximal subgroup of $G$. Hence all maximal subgroups are:

$N=\lg b,x,\Phi(M)\rg$;

$N_{i}=\lg ab^i,x,\Phi(M)\rg$ where $0\le i,j\le p-1$.

$N_{ij}=\lg ax^i,bx^j\rg$ where $0\le i,j\le p-1$.

It is easy to see that $N'=M_3$. By Lemma \ref{A_2-property} (7),
$\alpha_1(N)=p^2$. By calculation, $|N_0'|=p$ and $|N_i'|=p^2$ for
$i\ne 0$. Hence $\alpha_1(N_0)=p^2$ and $\alpha_1(N_i)=p^2+p$ for
$i\ne 0$. Since $N_{ij}\cong M$, $\alpha_1(N_{ij})=p$. Let
$H\in\Gamma_2$. Then $H=\lg a^ib^jx^k,\Phi(G)\rg$. It is obvious
that $H\in \mathcal{A}_1$ if and only if $j\ne 0$. Hence $\sum_{H\in
\Gamma_2} \alpha_1(H)=p^2$. By Hall's enumeration principle,
$$\alpha_1(G)=\sum_{H\in
\Gamma_1} \alpha_1(H)-p\sum_{H\in \Gamma_2} \alpha_1(H)=2\times
p^2+(p-1)\times (p^2+p)+p^2\times p-p\times p^2=p^3+2p^2-p.$$

\medskip
{\bf Case 5.} $G$ is one of the groups (N22)--(N24).

In this case, $p>2$, $G=M\ast \lg x\rg$ where $\lg x\rg\cong C_p$ and $M=\lg a,b\rg$
such that $M$ is metacyclic, $[x,b]\in \Phi(M')$, $[x,a]\not\in M'$. Hence all maximal subgroups are:

$N=\lg b,x,\Phi(M)\rg$;

$N_{i}=\lg ab^i,x,\Phi(M)\rg$ where $0\le i,j\le p-1$.

$N_{ij}=\lg ax^i,bx^j\rg$ where $0\le i,j\le p-1$.

It is easy to see that $N'=\Phi(M')$. By Lemma \ref{A_2-property} (7), $\alpha_1(N)=p^2$.
By calculation, $|N_i'|=p^2$. Hence $\alpha_1(N_i)=p^2+p$.
Since $N_{ij}\cong M$, $\alpha_1(N_{ij})=1+p$. Let $H\in\Gamma_2$.
Then $H=\lg a^ib^jx^k,\Phi(G)\rg$. It is obvious that $H\not\in \mathcal{A}_1$ if and only if $i=j=0$. Hence $\sum_{H\in
\Gamma_2} \alpha_1(H)=p^2+p$.
 By Hall's enumeration principle,
$$\alpha_1(G)=\sum_{H\in
\Gamma_1} \alpha_1(H)-p\sum_{H\in \Gamma_2} \alpha_1(H)=p^2+p\times
(p^2+p)+p^2\times (1+p)-p\times(p^2+p)=p^3+2p^2.$$

\medskip
{\bf Case 6.} $G$ is one of the groups (N25)--(N26).

It is easy to verify that $M=\lg b,d,a^p\rg$ is the unique
three-generator maximal subgroup of $G$. Since $|M'|=p$, by Lemma
\ref{A_2-property}, $\alpha_1(N)=p^2$. Other maximal subgroups are:

$N_{i}=\lg ab^i,d\rg$ where $0\le i\le p-1$.

$N_{ij}=\lg ad^i,bd^j\rg$ where $0\le i,j\le p-1$.

It is easy to see that $N_i$ has a unique abelian maximal subgroup
$\lg d,a^p,b^p\rg$ and $N_{ij}$ has no abelian maximal subgroup.
Hence $\alpha_1(N_i)=p$ and $\alpha_1(N_{ij})=1+p$. Let
$H\in\Gamma_2$. Then $H=\lg a^ib^jd^k,\Phi(G)\rg$ where $(i,j,k)\ne
(0,0,0)$. It is obvious that $H\not\in \mathcal{A}_1$ if and only if
$i=j=0$, hence if and only if $H=\lg d,\Phi(G)\rg$. Thus $\sum_{H\in
\Gamma_2} \alpha_1(H)=p^2+p$. By Hall's enumeration principle,
$$\alpha_1(G)=\sum_{H\in
\Gamma_1} \alpha_1(H)-p\sum_{H\in \Gamma_2} \alpha_1(H)=p^2+p\times
p+p^2\times (1+p)-p\times(p^2+p)=2p^2.$$
 \qed

\begin{thm}\label{d=4-5}
Suppose that $G$ is an $\mathcal{A}_3$-group having no abelian
subgroup of index $p$. Then $d(G)=4$ if and only if  $G$ is
isomorphic to one of the following pairwise non-isomorphic groups:
\begin{enumerate}

\rr{Oi} $G'\cong C_p$ and $c(G)=2$.
In this case, $\Phi(G)=Z(G)=G'$, $(\mu_0,\mu_1,\mu_2)=(0,0,1+p+p^2+p^3)$ and $\alpha_1(G)=p^2+p^4$.

\begin{enumerate}
\rr{O1}  $D_8\ast Q_8$;

\rr{O2} $Q_8\ast Q_8$;

\rr{O3} $M_p(1,1,1)\ast M_p(2,1)$, where $p>2$;

\rr{O4} $M_p(1,1,1)\ast M_p(1,1,1)$, where $p>2$.
\end{enumerate}

\rr{Oii} $G'\cong C_p^2$  and $c(G)=2$.
\begin{enumerate}
\rr{O5} $G=\lg a,b,c,d\di
a^4=b^4=1,c^2=a^2,d^2=b^2,[a,b]=1,[a,c]=b^2,[b,c]=a^2,[a,d]=
a^2,[b,d]=a^2b^2,[c,d]=1\rg$;  here $|G|=2^6$, $\Phi(G)=Z(G)=G'=\lg
a^2,b^2\rg \cong C_2^2$, any two noncommutative elements of $G$ generate $M_2(2,2)$.
$(\mu_0,\mu_1,\mu_2)=(0,0,15)$ and $\alpha_1(G)=30$.
\end{enumerate}

\rr{Oiii} $G'\cong C_p^3$  and $c(G)=2$.
\begin{enumerate}
\rr{O6} $G=K\times \lg a_4\rg$ where $K=\lg a_1,a_2,a_3\di
a_1^4=a_2^4=a_3^4=1,[a_1,a_2]=a_3^2,[a_1,a_3]=a_2^2a_3^2,[a_2,a_3]=a_1^2a_2^2,[a_1^2,a_2]=[a_2^2,a_1]=1\rg$ and $\lg a_4\rg\cong C_2$;
here $|G|=2^7$, $\Phi(G)=G'=\lg a_1^2,a_2^2,a_3^2\rg$,
$Z(G)=\lg a_1^2,a_2^2,a_3^2,a_4\rg\cong C_2^4$, $(\mu_0,\mu_1,\mu_2)=(0,0,15)$ and $\alpha_1(G)=30$.
\end{enumerate}
\end{enumerate}
\end{thm}

\demo Let the type of $G/G'$ be $(p^{m_1},p^{m_2},p^{m_3},p^{m_4})$,
where $m_1\ge m_2\ge m_3\ge m_4$, and $G/G'=\lg a_1G'\rg\times\lg
a_2G'\rg\times\lg a_3G'\rg\times \lg a_4G'\rg$, where
$o(a_iG')=p^{m_i}$, $i=1,2,3,4$. Then $G=\lg a_1,a_2,a_3,a_4\rg$. By
Lemma \ref{d=4}, $c(G)=2$, $\Phi(G)\le Z(G)$,  $G'\le C_p^3$ and all
$\mathcal{A}_1$-subgroups of $G$ contain $\Phi(G)$. The last property gives that $\sum_{H\in
\Gamma_2} \alpha_1(H)=\alpha_1(G)$.

We claim that $m_1=1$. Otherwise, $m_1\ge 2$. Let $B=\lg
a_2,a_3,a_4\rg$. Since $|G:B|\ge p^2$, we deduce that $B'=1$. Hence
$A=\lg B,a_1^p\rg$ is an abelian subgroup  of index $p$ of $G$, a
contradiction. Hence $G/G'$ is elementary abelian.

\medskip

{\bf Case 1:} $G'\cong C_p$.

\medskip

In this case, $|G|=p^5$. We claim that $Z(G)=G'$. Otherwise, without
loss of generality assume $a_1\in Z(G)$. Let $B=\lg a_2,a_3,a_4\rg$.
Then $B'\cong C_p$. By Lemma \ref{A_2-property} (3), $|Z(B)|=p^2$.
Since $Z(G)\ge \lg a_1,Z(B)\rg$, $|Z(G)|\ge p^3$. Hence $G$ has an
abelian subgroup  of index $p$ of $G$, a contradiction.

By above argument, $G$ is an extraspecial $p$-group. Hence we get
groups (O1)--(O4).

\medskip

{\bf Case 2:} $G'\cong C_p^2$.

\medskip

In this case, $|G|=p^6$ and any two noncommutative elements generate
an $\mathcal{A}_1$-group of order $p^4$. Such groups were classified
by \cite{AP}. By checking those groups listed in \cite{AP} we get
the group (O5).

\medskip

{\bf Case 3:} $G'\cong C_p^3$.

\medskip

In this case, $|G|=p^7$ and any $\mathcal{A}_1$-subgroup is
$M_p(2,2,1)$. It follows that $\Omega_1(G)\le Z(G)$. Let $N\le G'$
such that $|N|=p$ and $\bar{G}=G/N$. Then $|\bar{G}'|=p^2$. By Lemma
\ref{d(K) leq k+1 3}, there exists $\bar{K}\le \bar{G}$ such that
$d(\bar{K})=3$ and $\bar{K}'=\bar{G}'$. Without loss of generality
assume $K=\lg a_1,a_2,a_3\rg$. It follows that $|K'|\ge p^2$ and
$G'=\lg a_1^p,a_2^p,a_3^p\rg$.

\medskip

Subcase 3.1: $K'\cong C_p^3$.

\medskip

By Lemma \ref{A_2}, $$K=\lg a_1,a_2,a_3\di a_1^4=a_2^4=a_3^4=1,
[a_1,a_2]=a_3^2, [a_1,a_3]=a_2^2a_3^2, [a_2,a_3]=a_1^2a_2^2\rg.$$
Since $a_4^2\in G'=K'$, there exists $c\in K$ such that $a_4^2=c^2$.
Since $|\lg a_4,c\rg|\le 16$, $[a_4,c]=1$. By replacing $a_4$ with
$a_4c$ we get $a_4^2=1$. Hence $a_4\in\Omega_1(G)\le Z(G)$ and we
get the group (O6).

\medskip

Subcase 3.2: $K'\cong C_p^2$.

\medskip

Without loss of generality assume $K'=\lg [a_1,a_2],[a_1,a_3]\rg$
and $[a_2,a_3]=1$. If $G'=\lg K',[a_1,a_4]\rg$, then there exists
$b\in\lg a_2,a_3,a_4\rg$ such that $[a_1,b]=a_1^p$, which
contradicts that $G'\le \lg a_1,b\rg$. Hence $[a_1,a_4]\in K'$. By
suitable replacement we get $[a_1,a_4]=1$. Without loss of
generality assume $G'=\lg K',[a_2,a_4]\rg$. By replacing $a_3$ with
$a_3a_4$ we have $K'=G'$. It is reduced to subcase 3.1.

\medskip
We calculate the $(\mu_0,\mu_1,\mu_2)$ and $\alpha_1(G)$ of those
groups in Theorem \ref{d=4-5} as follows.

Since $d(G)=4$, $(\mu_0,\mu_1,\mu_2)=(0,0,1+p+p^2+p^3)$. In the
following, we calculate $\alpha_1(G)$.

\medskip
{\bf Case 1.} $G$ is one of the groups (O1)--(O4).

Let $H\in\Gamma_1$. Then $|H'|=p$. By Lemma \ref{A_2-property} (7),
$\alpha_1(H)=p^2$. By Hall's enumeration principle,
$$\alpha_1(G)=\sum_{H\in \Gamma_1} \alpha_1(H)-p\sum_{H\in \Gamma_2}
\alpha_1(H)=(1+p+p^2+p^3)\times p^2-p\alpha_1(G).$$ Hence
$$\alpha_1(G)=\frac{(1+p+p^2+p^3)p^2}{1+p}=p^2+p^4.$$

\medskip
{\bf Case 2.} $G$ is the group (O5).

By calculation, $d(H)=3$ and $H'\cong C_2^2$ for any $H\in\Gamma_1$. Hence $\alpha_1(H)=6$.
By Hall's enumeration principle,
$$\alpha_1(G)=\sum_{H\in
\Gamma_1} \alpha_1(H)-p\sum_{H\in \Gamma_2} \alpha_1(H)=15\times
6-2\alpha_1(G).$$ Hence $\alpha_1(G)=30$.

{\bf Case 3.} $G$ is the group (O6).

Let $H\in\Gamma_1$. Then $H$ is one of the following types: (1)
$H=M\times \lg a_4\rg$ where $M\maxsgp K$; (2) $H=\lg
a_1a_4^i,a_2a_4^j,a_3a_4^k\rg$ where $i,j,k=0,1$. If $H$ is of Type
(1), then $|H'|=2$ and hence $\alpha_1(H)=4$ by Lemma
\ref{A_2-property}. If $H$ is of Type (2), then $H\cong K$ and hence
$\alpha_1(H)=7$. By Hall's enumeration principle,
$$\alpha_1(G)=\sum_{H\in \Gamma_1} \alpha_1(H)-p\sum_{H\in \Gamma_2}
\alpha_1(H)=7\times 4+8\times 7-2\alpha_1(G).$$ Hence
$\alpha_1(G)=28.$ \qed

\medskip

Now we list the triple $(\mu_0,\mu_1,\mu_2)$ and $\alpha_1(G)$ for
$\mathcal{A}_3$-groups in Table \ref{table 12} and \ref{table 13}
respectively. This solves Problem 893 in \cite{Ber2}, Problem 1595
in \cite{Ber3} and Problem 2829 in \cite{Ber4} respectively.

At the end of this paper, we list some enumeration properties of
$\mathcal{A}_3$-groups, which can easily be obtained from Table
\ref{table 13}.

\begin{thm}
Let $G$ be an $\mathcal{A}_3$-group. Then $p^2\le \alpha_1(G)\le
p^4+p^3+p^2+p$, and the following conclusions hold.

{\rm (1)} If $\alpha_1(G)=p^2$, then $d(G)=2$, $c(G)=4$ and $G$ has
an abelian maximal subgroup. In this case, non-abelian subgroups of
$G$ are generated by two elements;

{\rm (2)} If $\alpha_1(G)=p^4+p^3+p^2+p$, then $p=2$, $c(G)=2$,
$d(G)=4$ and $G=K\times C_2$ where $K$ is the smallest Suzuki
$2$-group;

{\rm (3)} If $p>2$, then $\alpha_1(G)\le p^4+p^3+p^2$;

{\rm (4)} If $d(G)=2$, then $\alpha_1(G)\le p^3+2p^2+p$;

{\rm (5)} If $d(G)=2$ and $\alpha_1(G)\le p^3+2p^2+p$, then $p=2$
and $G$ is the group {\rm (M37)} in Theorem {\rm \ref{d=4-3}};

{\rm (6)} If $d(G)\ge 3$, then $\alpha_1(G)\ge 2p^2-1$.
\end{thm}

{\bf Acknowledgments} During writing the paper, Professor Y.
Berkovich gave us many correspondences, we cordially express thanks
for  his encouragements and many valuable comments.  All these
helped us to improve the whole paper considerably.

\newpage
\begin{table}[h]
\centering { \scriptsize
\begin{tabular}[t]{c||c}

\hline
 $(\mu_0,\ \mu_1,\ \mu_2)$      & types of  $\mathcal{A}$$_3$ groups \\

\hline $(1,\ p-1,\ 1)$  &  (A1)--(A6)\\

\hline $(p+1,\ p^2-1,\ 1)$  &  (B1); (B2) where $m=n= 1$; (B3) where $n=1$; (B4); (B5)  where $n=1$\\

\hline $(p+1,\ p^2-p,\ p)$  &  (B2) where $m>1=n$ or $n>1=m$; (B3) where $n>1$; (B5)  where $n>1$ \\

\hline \multirow{2}*{$(1,\ p^2-1,\ p+1 )$}  &  (B6); (B9); (B11); (B13) where $p=2$ and $m=l=1$\\
 &(B17) where $l=1$; (B19) where $p=2$ and $m=l=1$\\

\hline \multirow{3}*{$(1,\ p^2,\ p )$ } &   (B7) where $l>1$; (B10); (B12) where $l>1$; (B13) where $m=1$ and $l>1$\\
& (B14) where $n=1$ or $m=1$; (B16); (B18) where $m=1$\\
                              &(B19) where $p=2$ and $l>1=m$; (B19) where $p>2$ and $m=1$\\

\hline $(1,\ p^2+p-1,\ 1)$  &  (B7) where $l=1$; (B12) where $l=1$; (B13) where $p>2$ and $m=l=1$; (B20)\\

\hline $(1,\ p^2+p-2,\ 2)$  &  (B8) where $l=1$\\

\hline \multirow{2}*{$(1,\ p^2-p,\ 2p)$}  &  (B8) where $l>1$; (B13) where $m=2$; (B14) where $n=m=2$\\
&(B15); (B17) where $l>1$; (B18) where $m=2$; (B19) where $m=2$\\

\hline $(0,\ p,\ 1)$  &  (C1)--(C6), (C8), (C9), (C11), (C15)--(C17)\\

\hline $(0,\ p-1,\ 2)$  &  (C7), (C10), (C12)--(C14)\\

\hline $(0,\ p^2,\ p+1)$  & (D1),(D4), (D6)--(D10), (D12)--(D14), (D16)  \\

\hline {$(0,\ p^2+p,\ 1)$}  & (D2) where $-\nu\not\in (F_p^*)^2$;(D3) where $-r \not\in (F_p^*)^2$; (D5); (D11) where $-\nu \not\in (F_p^*)^2$;  \\

\hline $(0,\ p^2-p,\ 2p+1)$  & (D2) where $-\nu \in (F_p^*)^2$; (D3) where $-\nu \in (F_p^*)^2$; (D11) where $-\nu \in (F_p^*)^2$; (D15) \\

\hline $(0,\ p^2-1,\ p+2)$  & (D17), (D18) \\

\hline $(0,\ p^2+1,\ p)$  & (D19) \\

\hline $(0,\ 1,\ p)$  & (E1)--(E7) \\

\hline $(0,\ 1,\ p^2+p)$  & (E8)--(E10) \\

\hline $(1,\ 0,\ p)$  & (F1)--(F8); (G1)--(G12) \\

\hline $(p+1,\ 0,\ p^2)$  & (H1)--(H3) \\

\hline $(1,\ 0,\ p^2+p)$  & (H4)--(H10), (I1)--(I11) \\

\hline $(p+1,\ 0,\ p^3+p^2)$  & (J1)--(J5) \\

\hline $(1,\ 0,\ p^3+p^2+p)$  & (J6)--(J9) \\

\hline $(0,\ 0,\ p+1)$  & (K1)--(K5); (L1)--(L2); (M1)--(M62)\\

\hline $(0,\ 0,\ p^2+p+1)$  & (N1)--(N26)\\

\hline $(0,\ 0,\ p^3+p^2+p+1)$  & (O1)--(O6)\\

\hline

\end{tabular}

\caption{The number of $\mathcal{A}$$_0$-, $\mathcal{A}$$_1$-,
$\mathcal{A}$$_2$-subgroups of index $p$ in
$\mathcal{A}$$_3$-groups}\label{table 12} }
\end{table}

\newpage
\begin{table}[h]
\centering { \scriptsize
\begin{tabular}[t]{c||c}
\hline $\alpha_1(G)$ &types of $\mathcal{A}_3$-groups $G$\\

\hline $p^2$& (F1)--(F8) \\

\hline $p^2+1$& (E1), (E3)--(E7)\\

\hline $p^2+p-1$& (A1)--(A6)\\

\hline $p^2+p$&  (C1)--(C6), (C8)--(C9), (C11), (C15), (C17), (K3)--(K5)\\

\hline $p^2+p+1$& (K1)--(K2)\\

\hline $p^2+2p$& (C16)\\

\hline $2p^2-1$& (B1), (B2) where $n=m=1$, (B3) where $n=1$, (B4), (B5) where $n=1$\\

\hline $2p^2$& (M48)--(M53), (M58)--(M62), (N25)--(N26)\\

\hline \multirow{2}*{$2p^2+p-1$}& (B7) where $l=1$, (B12) where
$l=1$,
(B13) where $p>2$ and $l=m=1$, (B20) \\
&(C7), (C10), (C12)--(C14)\\

\hline $2p^2+p$&   (D2) where $-\nu\not\in (F_p^*)^2$;(D3) where $-r \not\in (F_p^*)^2$;
(D5); (D11) where $-\nu \not\in (F_p^*)^2$; (M54)--(M57) \\

\hline $3p^2+1$& (D19)\\

\hline $3p^2+p-2$& (B8) where $l=1$ \\

\hline $3p^2+p$& (L1)--(L2)\\

\hline $p^3$&  (G1)--(G12), (I1)--(I9) \\
\hline $p^3+1$&     (E2), (E8) \\
\hline $p^3+p^2-p$&  (B2) where $m>1=n$ or $n>1=m$; (B3) where $n\ge 2$, (B5) where $n\ge 2$\\

\hline \multirow{4}*{$p^3+p^2$} & (B7) where $l\ge 2$, (B10), (B12)
where $l\ge 2$, (B13) where
$l\ge 2$ and $m=1$\\

& (B14) where $m=1$ or $n=1$, (B16), (B18) where $m=1$, (B19) where $m=1$ and $p>2$\\

& (B19) where $p=2$ and $m=1$ and $l>1$, (I10)--(I11), (M1)--(M19), (M40)--(M43), (M47)\\
&  (N6)--(N13), (N15)--(N17), (N19)--(N21) \\

\hline $p^3+p^2+p$& (M38)-(M39)\\

\hline $p^3+2p^2-p$& (N14), (N18)\\

\hline  \multirow{2}*{$p^3+2p^2-1$}& (B6), (B9), (B11), (B13) where $p=2$ and $l=m=1$\\
&(B17) where $l=1$, (B19) where $p=2$ and $l=m=1$\\

\hline $p^3+2p^2$& (D1),(D4), (D6)--(D10), (D12)--(D14), (D16), (M20)--(M36), (M44)--(M46), (N22)--(N24)\\

\hline   $p^3+2p^2+p$& (M37)\\

\hline   $p^3+3p^2-1$& (D17)--(D18)\\

\hline   \multirow{2}*{$2p^3+p^2-p$}&   (B8) where $l\ge 2$, (B13)
where $m=2$, (B14)
where $m=n=2$, (B15) \\
&(B17) where $l\ge 2$, (B18) where $m=2$, (B19) where $m=2$ \\

\hline   $2p^3+2p^2-p$& (D2) where $-\nu \in (F_p^*)^2$; (D3) where $-\nu \in (F_p^*)^2$; (D11) where $-\nu \in (F_p^*)^2$; (D15) \\

\hline   $p^4$&  (H1)--(H3), (J1)--(J5)             \\

\hline  $p^4+p^2$&(O1)--(O4) \\

\hline   $p^4+p^3+1$& (E9)--(E10) \\

\hline   $p^4+p^3$&(H4)--(H10), (J6)--(J9)\\

\hline   $p^4+p^3+p^2$& (N1)--(N5), (O6)\\

\hline   $p^4+p^3+p^2+p$&(O5)\\

\hline
\end{tabular}
\caption{The number of $\mathcal{A}_1$-subgroups in
$\mathcal{A}_3$-groups}\label{table 13} }

\end{table}

\medskip


\end{document}